\documentclass[11pt,draft]{article}
\usepackage{amsmath}
\usepackage{amssymb}
\usepackage{amscd}
\usepackage{xspace}
\usepackage{verbatim}
\newcommand{\mysection}[1]{\section{#1}\setcounter{equation}{0}}
\date{}
\begin{document}%
   
\title{\bf Elliptic Equations Involving Measures \footnote{{ To appear as a chapter in {Handbook for Differential Partial 
Differential Equations}, M. Chipot, P. Quittner Eds, Elsevier.}}}

\author{
{\bf Laurent V\'eron}\\
{\small Laboratoire de Math\'ematiques et Physique Th\'eorique}\\
{\small CNRS, UMR 6083}\\
{\small Facult\'e des Sciences, Universit\'e de Tours,}\\
{\small Parc de Grandmont, 37200 Tours, FRANCE}}

\maketitle

\newcommand{\txt}[1]{\;\text{ #1 }\;}
\newcommand{\tbf}{\textbf}
\newcommand{\tit}{\textit}
\newcommand{\tsc}{\textsc}
\newcommand{\trm}{\textrm}
\newcommand{\mbf}{\mathbf}
\newcommand{\mrm}{\mathrm}
\newcommand{\bsym}{\boldsymbol}
\newcommand{\scs}{\scriptstyle}
\newcommand{\sss}{\scriptscriptstyle}
\newcommand{\txts}{\textstyle}
\newcommand{\dsps}{\displaystyle}
\newcommand{\fnz}{\footnotesize}
\newcommand{\scz}{\scriptsize}
\newcommand{\be}{\begin{equation}}
\newcommand{\bel}[1]{\begin{equation}\label{#1}}
\newcommand{\ee}{\end{equation}}
\newcommand{\eqnl}[2]{\begin{equation}\label{#1}{#2}\end{equation}}
\newtheorem{subn}{\name}
\renewcommand{\thesubn}{}
\newcommand{\bsn}[1]{\def\name{#1}\begin{subn}}
\newcommand{\esn}{\end{subn}}
\newtheorem{sub}{\name}[section]
\newcommand{\dn}[1]{\def\name{#1}}   
\newcommand{\bs}{\begin{sub}}
\newcommand{\es}{\end{sub}}
\newcommand{\bsl}[1]{\begin{sub}\label{#1}}
\newcommand{\bth}[1]{\def\name{Theorem}\begin{sub}\label{t:#1}}
\newcommand{\blemma}[1]{\def\name{Lemma}\begin{sub}\label{l:#1}}
\newcommand{\bcor}[1]{\def\name{Corollary}\begin{sub}\label{c:#1}}
\newcommand{\bdef}[1]{\def\name{Definition}\begin{sub}\label{d:#1}}
\newcommand{\bprop}[1]{\def\name{Proposition}\begin{sub}\label{p:#1}}
\newcommand{\R}{\eqref}
\newcommand{\rth}[1]{Theorem~\ref{t:#1}}
\newcommand{\rlemma}[1]{Lemma~\ref{l:#1}}
\newcommand{\rcor}[1]{Corollary~\ref{c:#1}}
\newcommand{\rdef}[1]{Definition~\ref{d:#1}}
\newcommand{\rprop}[1]{Proposition~\ref{p:#1}}
\newcommand{\BA}{\begin{array}}
\newcommand{\EA}{\end{array}}
\newcommand{\BAN}{\renewcommand{\arraystretch}{1.2}
\setlength{\arraycolsep}{2pt}\begin{array}}
\newcommand{\BAV}[2]{\renewcommand{\arraystretch}{#1}
\setlength{\arraycolsep}{#2}\begin{array}}
\newcommand{\BSA}{\begin{subarray}}
\newcommand{\ESA}{\end{subarray}}
\newcommand{\BAL}{\begin{aligned}}
\newcommand{\EAL}{\end{aligned}}
\newcommand{\BALG}{\begin{alignat}}
\newcommand{\EALG}{\end{alignat}}
\newcommand{\BALGN}{\begin{alignat*}}
\newcommand{\EALGN}{\end{alignat*}}
\newcommand{\note}[1]{\textit{#1.}\hspace{2mm}}
\newcommand{\Proof}{\note{Proof}}
\newcommand{\qeda}{\hspace{10mm}\hfill $\square$}
\newcommand{\qed}{\\ ${}$ \hfill $\square$}
\newcommand{\Remark}{\note{Remark}}
\newcommand{\modin}{$\,$\\[-4mm] \indent}
\newcommand{\forevery}{\quad \forall}
\newcommand{\set}[1]{\{#1\}}
\newcommand{\setdef}[2]{\{\,#1:\,#2\,\}}
\newcommand{\setm}[2]{\{\,#1\mid #2\,\}}
\newcommand{\lra}{\longrightarrow}
\newcommand{\lla}{\longleftarrow}
\newcommand{\llra}{\longleftrightarrow}
\newcommand{\Lra}{\Longrightarrow}
\newcommand{\Lla}{\Longleftarrow}
\newcommand{\Llra}{\Longleftrightarrow}
\newcommand{\warrow}{\rightharpoonup}
\newcommand{\paran}[1]{\left (#1 \right )}
\newcommand{\sqbr}[1]{\left [#1 \right ]}
\newcommand{\curlybr}[1]{\left \{#1 \right \}}
\newcommand{\abs}[1]{\left |#1\right |}
\newcommand{\norm}[1]{\left \|#1\right \|}
\newcommand{\paranb}[1]{\big (#1 \big )}
\newcommand{\lsqbrb}[1]{\big [#1 \big ]}
\newcommand{\lcurlybrb}[1]{\big \{#1 \big \}}
\newcommand{\absb}[1]{\big |#1\big |}
\newcommand{\normb}[1]{\big \|#1\big \|}
\newcommand{\paranB}[1]{\Big (#1 \Big )}
\newcommand{\absB}[1]{\Big |#1\Big |}
\newcommand{\normB}[1]{\Big \|#1\Big \|}

\newcommand{\thkl}{\rule[-.5mm]{.3mm}{3mm}}
\newcommand{\thknorm}[1]{\thkl #1 \thkl\,}
\newcommand{\trinorm}[1]{|\!|\!| #1 |\!|\!|\,}
\newcommand{\bang}[1]{\langle #1 \rangle}
\def\angb<#1>{\langle #1 \rangle}
\newcommand{\vstrut}[1]{\rule{0mm}{#1}}
\newcommand{\rec}[1]{\frac{1}{#1}}
\newcommand{\opname}[1]{\mbox{\rm #1}\,}
\newcommand{\supp}{\opname{supp}}
\newcommand{\dist}{\opname{dist}}
\newcommand{\myfrac}[2]{{\displaystyle \frac{#1}{#2} }}
\newcommand{\myint}[2]{{\displaystyle \int_{#1}^{#2}}}
\newcommand{\q}{\quad}
\newcommand{\qq}{\qquad}
\newcommand{\hsp}[1]{\hspace{#1mm}}
\newcommand{\vsp}[1]{\vspace{#1mm}}
\newcommand{\ity}{\infty}
\newcommand{\prt}{\partial}
\newcommand{\sms}{\setminus}
\newcommand{\ems}{\emptyset}
\newcommand{\ti}{\times}
\newcommand{\pr}{^\prime}
\newcommand{\ppr}{^{\prime\prime}}
\newcommand{\tl}{\tilde}
\newcommand{\sbs}{\subset}
\newcommand{\sbeq}{\subseteq}
\newcommand{\nind}{\noindent}
\newcommand{\ind}{\indent}
\newcommand{\ovl}{\overline}
\newcommand{\unl}{\underline}
\newcommand{\nin}{\not\in}
\newcommand{\pfrac}[2]{\genfrac{(}{)}{}{}{#1}{#2}}

\def\ga{\alpha}     \def\gb{\beta}       \def\gg{\gamma}
\def\gc{\chi}       \def\gd{\delta}      \def\ge{\epsilon}
\def\gth{\theta}                         \def\vge{\varepsilon}
\def\vgf{\phi}       \def\vgf{\varphi}    \def\gh{\eta}
\def\gi{\iota}      \def\gk{\kappa}      \def\gl{\lambda}
\def\gm{\mu}        \def\gn{\nu}         \def\gp{\pi}
\def\vgp{\varpi}    \def\gr{\rho}        \def\vgr{\varrho}
\def\gs{\sigma}     \def\vgs{\varsigma}  \def\gt{\tau}
\def\gu{\upsilon}   \def\gv{\vartheta}   \def\gw{\omega}
\def\gx{\xi}        \def\gy{\psi}        \def\gz{\zeta}
\def\Gg{\Gamma}     \def\Gd{\Delta}      \def\vgf{\Phi}
\def\Gth{\Theta}
\def\Gl{\Lambda}    \def\Gs{\Sigma}      \def\Gp{\Pi}
\def\Gw{\Omega}     \def\Gx{\Xi}         \def\Gy{\Psi}

\def\CS{{\mathcal S}}   \def\CM{{\mathcal M}}   \def\CN{{\mathcal N}}
\def\CR{{\mathcal R}}   \def\CO{{\mathcal O}}   \def\CP{{\mathcal P}}
\def\CA{{\mathcal A}}   \def\CB{{\mathcal B}}   \def\CC{{\mathcal C}}
\def\CD{{\mathcal D}}   \def\CE{{\mathcal E}}   \def\CF{{\mathcal F}}
\def\CG{{\mathcal G}}   \def\CH{{\mathcal H}}   \def\CI{{\mathcal I}}
\def\CJ{{\mathcal J}}   \def\CK{{\mathcal K}}   \def\CL{{\mathcal L}}
\def\CT{{\mathcal T}}   \def\CU{{\mathcal U}}   \def\CV{{\mathcal V}}
\def\CZ{{\mathcal Z}}   \def\CX{{\mathcal X}}   \def\CY{{\mathcal Y}}
\def\CW{{\mathcal W}}
\def\BBA {\mathbb A}   \def\BBb {\mathbb B}    \def\BBC {\mathbb C}
\def\BBD {\mathbb D}   \def\BBE {\mathbb E}    \def\BBF {\mathbb F}
\def\BBG {\mathbb G}   \def\BBH {\mathbb H}    \def\BBI {\mathbb I}
\def\BBJ {\mathbb J}   \def\BBK {\mathbb K}    \def\BBL {\mathbb L}
\def\BBM {\mathbb M}   \def\BBN {\mathbb N}    \def\BBO {\mathbb O}
\def\BBP {\mathbb P}   \def\BBR {\mathbb R}    \def\BBS {\mathbb S}
\def\BBT {\mathbb T}   \def\BBU {\mathbb U}    \def\BBV {\mathbb V}
\def\BBW {\mathbb W}   \def\BBX {\mathbb X}    \def\BBY {\mathbb Y}
\def\BBZ {\mathbb Z}

\def\GTA {\mathfrak A}   \def\GTB {\mathfrak B}    \def\GTC {\mathfrak C}
\def\GTD {\mathfrak D}   \def\GTE {\mathfrak E}    \def\GTF {\mathfrak F}
\def\GTG {\mathfrak G}   \def\GTH {\mathfrak H}    \def\GTI {\mathfrak I}
\def\GTJ {\mathfrak J}   \def\GTK {\mathfrak K}    \def\GTL {\mathfrak L}
\def\GTM {\mathfrak M}   \def\GTN {\mathfrak N}    \def\GTO {\mathfrak O}
\def\GTP {\mathfrak P}   \def\GTR {\mathfrak R}    \def\GTS {\mathfrak S}
\def\GTT {\mathfrak T}   \def\GTU {\mathfrak U}    \def\GTV {\mathfrak V}
\def\GTW {\mathfrak W}   \def\GTX {\mathfrak X}    \def\GTY {\mathfrak Y}
\def\GTZ {\mathfrak Z}   \def\GTQ {\mathfrak Q}

\font\Sym= msam10 
\def\SYM#1{\hbox{\Sym #1}}
\newcommand{\bdw}{\prt\Gw\xspace}
\tableofcontents
\mysection {Introduction}
The role of measures in the study of nonlinear partial differential equations
has became more and more important in the last years, not only 
because it belongs to the mathematical spirit to try to extend the 
scope of a theory, but also because the extension from the 
function setting to the measure framework appeared to be the only 
way to bring into light nonlinear phenomena and to explain them. In a 
very similar process, the theory of linear equations shifted from the 
function setting to the distribution framework. 
The aim of this chapter is to bring into light several aspects of this 
interaction, in particular its connection with the singularity theory 
and the nonlinear trace theory. Our intention is not to present 
a truly self-contained text : clearly we shall assume that the reader is familiar with the 
standard second order linear elliptic equations regularity theory, as 
it is explained in Gilbarg and Trudinger's classical treatise \cite {GT}. 
Part of the results will be fully proven, and, for some of them, only 
the statements will be exposed. The starting point is the linear 
theory, in our case the study of 
\begin {equation}\label {Intr1}\left.\BA{rl}
Lu=\gl&\mbox { in }\Omega,\\[2mm]
u=\gm&\mbox { on }\prt\Omega, 
\EA\right.
\end {equation}
where $\Omega$ is a smooth bounded domain in $\mathbb R^n$, 
$L$ is a linear elliptic operator of second order, and $\gl$ and 
$\gm$ are Radon measures, respectively in $\Omega$ and $\prt\Omega$. 
Under some structural and regularity assumptions on $L$ (essentially 
that the maximum principle holds), it is
proven that (\ref {Intr1}) admits a unique solution. Moreover this 
solution admits a linear representation, i.e.
\begin {eqnarray}\label {Intr2}
u(x)=\int_{\Omega}G^\Omega_{L}(x,y)d\gl(y)+\int_{\prt\Omega}P^\Omega_{L}(x,y)d\gm(y),
\end {eqnarray}
for any $x\in\Omega$, where $G^\Omega_{L}$ and $P^\Omega_{L}$ are 
respectively the Green and the Poisson kernels associated to $L$ in 
$\Omega$. The presentation that we adopt is a combination of the 
classical regularity theory for linear elliptic equations and 
Stampacchia duality approach which provides the most powerful tool for the 
extension to semilinear equations.  In Section 3 we shall concentrate 
on semilinear equations 
with an absorption-reaction term of the following type
\begin {equation}\label {Intr3}\left.\BA{rl}
Lu+g(x,u)=\gl&\mbox { in }\Omega,\\[2mm]
u=0&\mbox { on }\prt\Omega, 
\EA\right.
\end {equation}
where $(x,r)\mapsto g(x,r)$ is a continuous function defined in 
$\Omega\times\mathbb R$, satisfying the absorption principle
\begin {eqnarray}\label {Intr4}
 {\rm sign}(r)g(x,r)\geq 0,\forevery (x,r)\in\Gw\times 
(-\infty,-r_{0}]\cup[r_{0},\infty),
\end {eqnarray}
for some $r_{0}\geq 0$. Under general assumptions on $g$, which are 
the natural generalisation of the Brezis-B\'enilan weak-singularity 
condition \cite{BB}, it is proven that for any Radon measure $\gl$ in $\Omega$
satisfying
\begin {eqnarray}\label {Intr5}
\int_{\Gw}\gr_{_{\prt\Gw}}^\ga d\abs\gl<\ity,
\end {eqnarray}
with $\gr_{_{\prt\Gw}}(x)=\dist (x,\prt\Gw)$ and $\ga\in [0,1]$, Problem (\ref {Intr3}) admits a solution. 
Notice that the assumption on $g$ depends both on $n$ and $\ga$. 
Furthermore, uniqueness holds if $r\mapsto g(x,r)$ is nondecreasing, 
for any $x\in\Gw$. However, the growth condition on $g$ is very 
restrictive. Thus the problem may not be solved for all the 
measures, but only for specific ones. A natural condition is to assume 
that the measure $\gl$ satisfies
\begin {eqnarray}\label {Intr6}
\int_{\Gw}g\left(x,\mathbb G_{L}^\Gw(\abs\gl)\right)\gr_{_{\prt\Gw}}dx<\ity,
\end {eqnarray}
where $\mathbb G_{L}^\Gw(\abs\gl)$, defined by
$$\mathbb G_{L}^\Gw(\abs\gl)(x)=\int_{\Omega}G^\Omega_{L}(x,y)d\abs 
{\gl}(y),\forevery x\in\Gw,
$$
is called the Green potential of $\abs\gl$. Under an additional 
condition on $g$, called the $\Gd_{2}$-condition, which excludes the exponential function, 
but not any 
positive power, it is shown that, in Condition (\ref{Intr6}), the measure $\gl$ can be replaced by its singular part
with respect to the $n$-dimensional Hausdorff measure in 
the Lebesgue decomposition, in order Problem (\ref {Intr3}) to be solvable. In the case where
$$g(x,r)=\abs r^{q-1}r,
$$
with $r>0$, Problem (\ref {Intr3}) can be solved for any bounded 
measure if $0<q<n/(n-2)$, but this is no longer the case if $q\geq 
n/(n-2)$. Baras and Pierre provide in \cite {BP1} a necessary and 
sufficient condition on the measure $\gl$ in terms of Bessel 
capacities. The solvability of nonlinear equations with measure is closely 
associated to removability question, the standard one being the 
following : assume $K$ is a compact subset of $\Gw$ and $u$ a solution 
of 
\begin {eqnarray}\label {Intr7}
Lu+g(x,u)=0\quad\mbox {in }\Gw\setminus K,
\end {eqnarray}
does it follows that $u$ can be extended, in a natural way, so that 
the equation is satisfied in all $\Gw$? The answer is positive if some
Bessel capacity, connected to the growth of $g$, of the set $K$ is zero. 
In Section 4 we give an overview of the semilinear problem 
with a source-reaction term of the following type
\begin {equation}\label {Intr8}\left.\BA{ll}
Lu=g(x,u)+\gl&\mbox { in }\Omega,\\[2mm]
\;u=0&\mbox { on }\prt\Omega, 
\EA\right.
\end {equation}
For this equation, not only the concentration of the measure is 
important, but also the total mass. The first approach, due to Lions 
\cite {Lio} is to construct a supersolution, the conditions are somehow 
restrictive. In the convex case, a rather complete presentation is 
provided by  Baras and Pierre \cite 
{BP2}, with the improvement of Adams and Pierre \cite {AP}. The idea 
is to write the solution $u$ of (\ref {Intr8}) under the form
\begin {eqnarray}\label {Intr8'}
u(x)=\int_{\Gw}G_{L}^\Gw g(y,u(y))dy+\mathbb G_{L}^\Gw(\gl)
\quad\mbox { in }\Omega.
\end {eqnarray}
The convexity of $r\mapsto g(x,r)$ gives a necessary condition 
expressed in term of the conjugate function $g^*(x,r)$. The difficulty 
is to prove that this condition is also sufficient and to link it to 
a functional analysis framework. An extension of this method is given 
by Kalton and Verbitsky \cite {KV} in connection with weighted 
inequalities in $L^q$ spaces. Finally, conditions for removability of 
singularities of positive solutions are treated by Baras and Pierre 
\cite {BP1}.  In Section 5 we consider the problem of 
solving boundary value problems with measures data for nonlinear equations 
with an absorption-reaction term, 
\begin {equation}\label {Intr9}\left.\BA{rl}
Lu+g(x,u)=0&\mbox { in }\Omega,\\[2mm]
u=\gm&\mbox { on }\prt\Omega, 
\EA\right.
\end {equation}
The first results in that direction are due to Gmira and V\'eron \cite 
{GmV} who prove that the B\'enilan-Brezis method can be adapted in a 
framework of weighted Marcinkiewicz spaces for obtaining existence of 
solutions in the so-called subcritical case : the case in which the problem is solvable with any boundary Radon measure. 
In a similar way as for Problem (\ref {Intr3}), it is shown that 
Problem (\ref {Intr9}) is solvable
 if the measure $\gm$ satisfies
\begin {eqnarray}\label {Intr10}
\int_{\Gw}g\left(x,\mathbb P_{L}^\Gw(\abs\gm)\right)\gr_{_{\prt\Gw}}dx<\ity,
\end {eqnarray}
where
$$\mathbb P_{L}^\Gw(\abs\gm)(x)=\int_{\prt\Omega}P^\Omega_{L}(x,y)d\abs 
{\gm}(y),\forevery x\in\Gw.
$$
It is also possible to extend the range of solvability if $\gm$ is replaced by its singular 
part with respect to the $(n-1)$-dimensional Hausdorff measure, 
for specific functions $g$ which verify a power like growth. In the 
last years the model 
case of equation
\begin {eqnarray}\label {Intr11}
L u+\abs u^{q-1}u=0,
\end {eqnarray}
acquired a central role because of its applications. The case 
$q=(n+2)/(n-2)$ is classical in Riemannian geometry and corresponds 
to conformal change of metric with prescribed constant negative 
scalar curvature \cite {LN}, \cite {RRV}. The case $1<q\leq 2$ is associated to 
superprocess in 
probability theory. It has been developed by Dynkin, \cite {Dy1}, \cite 
{Dy2} and Le Gall \cite {LG2} who introduced very powerful new tools for studying 
the properties of the positive solutions of this equation. The central 
idea is the discovery by Le Gall  \cite {LG1}, in the case $q=2=n$, and the extension 
by Marcus and V\'eron \cite {MV2}, in the general case $q>1$ and $n\geq 
2$, of the 
existence of a boundary trace of positive solutions of (\ref {Intr11}) 
in a smooth bounded domain $\Gw$. This boundary trace denoted by 
$Tr_{_{\prt\Gw}}(u)$ is no longer a Radon 
measure, but a $\sigma$-finite Borel measure which can takes infinite 
value on compact subsets of the boundary. The critical value for this 
equation, first observed by Gmira and V\'eron, is $q_{c}=(n+1)/(n-1)$. 
It is proven in \cite {LG1},  \cite {MV4} that for any positive 
$\sigma$-finite Borel measure $\gm$ on $\prt\Gw$ the problem
\begin {equation}\label {Intr12}\left.\BA{rl}
Lu+\abs u^{q-1}u=0&\mbox { in }\Omega,\\[2mm]
Tr_{_{\prt\Gw}}(u)=\gm&\mbox { on }\prt\Omega, 
\EA\right.
\end {equation}
admits a unique solution provided $1<q<q_{c}$. This is no longer the 
case when $q\geq q_{c}$. Although many results are now available for 
solving the super-critical case of Problem (\ref {Intr12}), the full 
theory is not yet completed. An important colateral problem deals with 
the question of boundary singularities, an example of which is the 
following : suppose $K$ is a compact subset of $\prt\Gw$, 
$u\in C^2(\Gw)\cap C(\overline\Gw\setminus K)$ is a solution of (\ref 
{Intr12}) in $\Gw$ which vanishes on $\prt\Gw\setminus K$ ; does it 
imply that $u$ is identically zero ? The answer to this question is 
complete, and expressed in terms of boundary Bessel capacities.
\mysection {Linear equations}
\subsection {Elliptic equations in divergence form}
We call $x=(x_{1},\ldots,x_{n}$ the variables in the space $\mathbb 
R^n$. Let $\Omega$ be a bounded domain in $\mathbb R^n$. The type of 
operators under consideration are linear 
second order differential operators in divergence form
\begin {eqnarray}\label {lin1}
Lu=-\sum_{i,j=1}^n \frac {\prt}{\prt x_{i}}\left(a_{ij}\frac {\prt u}{\prt x_{j}}\right)
+\sum_{i=1}^n b_{i}\frac {\prt u}{\prt x_{i}}-\sum_{i=1}^n \frac 
{\prt}{\prt x_{i}}\left(c_{i}u\right)+du
\end {eqnarray}
where the $a_{ij}$, $b_{i}$, $c_{i}$ and $d$ are at least bounded measurable functions satisfying
the uniform ellipticity condition in $\Omega$ :
\begin {eqnarray}\label {lin2}
\sum_{i,j=1}^n a_{ij}(x)\xi_{i}\xi_{j}
\geq \alpha \sum_{i=1}^n \xi_{i}^{2},\forevery 
\xi=(\xi_{1},\ldots,\xi_{n})\in \mathbb R^n,
\end {eqnarray}
for almost all $x\in\Omega$, where $\alpha>0$ is some fixed constant. 
It is classical to associate to $L$ the bilinear form $A_{L}$ 
\begin {eqnarray}\label {lin3}
A_{L}(u,v)=\int_{\Gw}a_{L}(u,v)dx,\forevery u,v\in 
W^{1,2}_{0}(\Omega),
\end {eqnarray}
where 
\begin {eqnarray}\label {lin4}
a_{L}(u,v)=\sum_{i,j=1}^n a_{ij}\frac {\prt u}{\prt x_{j}}\frac {\prt 
v}{\prt x_{i}}+\sum_{i=1}^n\left(b_{i}\frac {\prt u}{\prt x_{i}}v+
c_{i}\frac {\prt v}{\prt x_{i}}u\right)+duv.
\end {eqnarray}

An important uniqueness condition, symmetric in the $b_{i}$ and 
$c_{i}$, which also implies the maximum principle, 
is the following :
\begin {eqnarray}\label {uniq}
\int_{\Gw}\left(dv+\sum_{i=1}^n\frac {1}{2}(b_{i}+c_{i})
\frac {\prt v}{\prt x_{i}}\right)dx\geq 0,\forevery v\in C^1_{c}(\Gw),v\geq 0.
\end {eqnarray}

\blemma {LM} Let the coefficients of $L$ be bounded and measurable, 
and conditions (\ref {lin2}) and (\ref {uniq}) hold. Then for any 
$\phi\in W^{1,2}(\Gw)$ and $f_{i}\in L^2(\Gw)$ ($i=0,\ldots,n$) there 
exists a unique $u\in W^{1,2}(\Gw)$ solution of 
\begin {equation}\label {bvp}\left.\BA{ll}
Lu=f_{0}-\displaystyle{\sum_{i=1}^n}\myfrac {\prt f_{i}}{\prt x_{i}}&\mbox { in }\Omega,\\[2mm]
\;\;u=\phi&\mbox { on }\prt\Omega, 
\EA\right.
\end {equation}
\es 
\Proof By a solution, we mean $u-\phi\in W^{1,2}_{0}(\Gw)$ and 
\begin {eqnarray}\label {bvp1}
A_{L}(u,v)=\int_{\Gw}\left(f_{0}v+\sum_{i=1}^nf_{i}\frac {\prt v}{\prt 
x_{i}}\right)dx,\forevery v\in W^{1,2}_{0}(\Gw).
\end {eqnarray}
We put $\tilde u=u-\phi$. Then solving (\ref {bvp}) is equivalent to 
finding 
$\tilde u\in W^{1,2}_{0}(\Gw)$ such that 
\begin {eqnarray}\label {bvp1'}
A_{L}(\tilde u,v)=\int_{\Gw}\left(f_{0}v+\sum_{i=1}^nf_{i}\frac {\prt v}{\prt 
x_{i}}-a_{L}(\phi,v)\right)dx,\forevery v\in W^{1,2}_{0}(\Gw).
\end {eqnarray}
The bilinear form $A_{L}$ is clearly continuous on $W^{1,2}_{0}(\Gw)$ 
and
$$A_{L}(v,v)=\int_{\Gw}\left(\sum_{i,j=1}^n a_{ij}\frac {\prt v}{\prt x_{j}}\frac {\prt 
v}{\prt x_{i}}+dv^2+
\frac {1}{2}\sum_{i=1}^n(b_{i}+c_{i})\frac {\prt v^2}{\prt x_{i}}\right)dx.
$$
By (\ref {lin2}) and (\ref {uniq}),
$$A_{L}(v,v)\geq \ga \int_{\Gw}\abs{\nabla v}^2dx,\forevery v\in 
C^{1}_{0}(\Gw).
$$
By density  $A_{L}$ is coercive and thanks to Lax-Milgram's theorem, it defines an isomorphism 
between the 
Sobolev space $W^{1,2}_{0}(\Omega)$ and its dual space $W^{-1,2}(\Omega)$.\qeda \\

The celebrated De Giorgi-Nash-Moser regularity theorem asserts that, 
for $p>n$ and $f\in L_{loc}^p(\Omega)$, any $W^{1,2}_{loc}(\Omega)$ 
function $u$ which satisfies
\begin {eqnarray}\label {lin5}
\int_{\Gw}a_{L}(u,\phi) dx=\int_{\Gw}f\phi dx,\forevery \phi\in 
C_{0}^{\ity}(\Omega),
\end {eqnarray}
is locally H\" older continuous, up to a modification on a set of measure 
zero. Furthermore the weak maximum principle holds in the sense that if 
$u\in W^{1,2}(\Omega)$ satisfies
\begin {eqnarray}\label {subsol}
A_{L}(u,\phi)\leq 0,\forevery\phi\in C^{\infty}_{0}(\Omega),\,\phi\geq 
0,
\end {eqnarray}
such a $u$ is called a weak sub-solution, there holds
\begin {eqnarray}\label {subsol2}
\sup_{\Omega}u\leq\sup_{\prt\Gw}u.
\end {eqnarray}
In the above formula, 
$$ \sup _{\prt\Gw}v :=\inf\{k\in \mathbb R :\,(v-k)_{+}\in W^{1,2}_{0}(\Omega)\}.
$$
At end, the strong maximum principle holds : if for some 
ball $B\subset\bar B\subset\Gw$, 
\begin {eqnarray}\label {subsol3}
\sup_{B}u=\sup_{\Gw}u,
\end {eqnarray}
then $u$ is constant in the connected component of $\Omega$ 
containing $B$.\\

If the $a_{ij}$ and the $c_{i}$ are Lipschitz continuous, and the $b_{i}$ 
and $d$ are bounded measurable functions, the operator $L$ can be 
written in non-devergence form
\begin {eqnarray}\label {lin1'}
Lu=-\sum_{i,j=1}^n a_{ij}\frac {\prt^2 u}{\prt x_{i}\prt x_{j}}+
\sum_{j=1}^n b'_{j}\frac {\prt u}{\prt x_{j}}+d'u,
\end {eqnarray}
where
$$ b'_{j}=b_{j}-c_{j}-\sum_{i=1}^n\frac {\prt a_{ij}}{\prt x_{i}},
\quad d'=d-\sum_{i=1}^n\frac {\prt c_{i}}{\prt x_{i}}.
$$
Conversely, an operator $L$ in the non-divergence form (\ref {lin1'}) 
with Lipschitz continuous coefficients $a_{ij}$ and bounded and
measurable coefficients $b'_{i}$
can be written in divergence form
\begin {eqnarray}\label {lin1"}
Lu=-\sum_{i,j=1}^n \frac {\prt}{\prt x_{i}}\left(a_{ij}\frac {\prt u}{\prt x_{j}}
\right)+\sum_{j=1}^n\tilde b_{j}\frac {\prt u}{\prt x_{j}} +\tilde d,
\end {eqnarray}
with
$$\tilde b_{j}=b'_{j}+\sum_{i=1}^n\frac {\prt a_{ij}}{\prt x_{i}}.
$$
This duality between operators in divergence or in non-divergence form 
is very useful in the applications, in particular in the 
regularity theory of solutions of elliptic equations.
If $L$ is defined by (\ref {lin1}), the adjoint operator $L^{*}$ is defined by
\begin {eqnarray}\label {lin6}
L^{*}\phi=-\sum_{i,j=1}^n \frac {\prt}{\prt x_{j}}\left(a_{ij}\frac {\prt 
\phi}{\prt x_{i}}\right)
+\sum_{i=1}^nc_{i}\frac {\prt\phi}{\prt x_{i}}
-\sum_{i=1}^n\frac {\prt}{\prt x_{i}}\left(b_{i}\phi\right)
 +d\phi.
\end {eqnarray}
Under the mere assumptions that the coefficients $a_{ij}$, $b_{i}$,  
$c_{i}$ and $d$ are bounded and measurable in $\Gw$, the uniform ellipticity 
(\ref {lin2}), and the uniqueness condition (\ref {uniq}), the two operators 
$L$ and $L^*$ define 
an isomorphism between 
$W^{1,2}_{0}(\Gw)$ and $W^{-1,2}(\Gw)$. 
If the $a_{ij}$ and the $b_{i}$ are Lipschitz continuous, for any $u\in 
L^{1}_{loc}(\Gw)$, $Lu$ can be considered as a distribution in $\Omega$ 
if we define its action on test functions in the following way :
\begin {eqnarray}\label {lin7}
\langle Lu,\phi\rangle=\int_{\Gw}u L^{*}\phi dx,\forevery \phi\in 
C_{0}^{\ity}(\Omega).
\end {eqnarray}
\subsection {The $L^1$ framework}
Let $\Gw$ be a bounded domain with $C^2$ boundary and $L$ 
the operator given by (\ref {lin1}). 
\bdef {H}We say that the operator $L$ given by (\ref {lin1}) satisfies the 
condition (H), if the functions $a_{ij}$, $b_{i}$ and $c_{i}$ are Lipschitz continuous
in $\Gw$, $d$ is bounded and measurable, and if the uniform ellipticity 
condition (\ref {lin2}) and the uniqueness condition (\ref {uniq}) hold.
\es

Notice that this condition is symmetric in $L$ and $L^*$. We put
\begin {eqnarray}\label {distf}
\rho_{_{\prt\Gw}}(x)=\dist (x,\prt\Gw), \forevery x\in\overline\Gw.
\end {eqnarray}
We denote by $C_{c}^{1,L}(\overline \Omega)$ the space of 
$C^1(\overline \Gw)$ functions $\gz$, vanishing on $\prt\Gw$ and such that 
$L^*\gz\in 
L^{\ity}(\Gw)$, and by
\begin {eqnarray}\label {conorm}
\frac {\prt \gz}{\prt {\bf n}_{L^*}}=\sum_{i,j=1}^n
a_{ij}\frac {\prt \gz}{\prt x_{i}}{\bf n}_{j},
\end {eqnarray}
the co-normal derivative on the boundary following $L^*$ (here the ${\bf 
n}_{j}$ are the components of outward normal unit vector {\bf n} to $\prt\Gw$). 
\bdef {L1} Let $f\in L^1(\Gw;\rho_{_{\prt\Gw}}dx)$ and
$g\in L^1(\prt\Gw)$. We say that a function $u\in L^1(\Gw)$ is a very 
weak solution of the problem
\begin {equation}\label {L1-1}\left.\BA{rl}
Lu=f&\mbox { in }\Omega,\\[2mm]
\;u=g&\mbox { on }\prt\Omega, 
\EA\right.
\end {equation}
if, for any $\gz\in C_{c}^{1,L}(\overline \Omega)$, there holds
\begin {eqnarray}\label {L1-2}
\int_{\Gw}uL^*\gz dx=\int_{\Gw}f\gz dx-\int_{\prt\Gw}\frac 
{\prt\gz}{\prt{\bf n}_{L^*}}gdS.
\end {eqnarray}
\es
The next result is an adaptation of a construction,  essentially due to 
Brezis in the case of the Laplacian, although various forms 
of existence theorems were known for a long time.
\bth {L1-th} Let $L$ satisfy the condition (H). Then for any $f$ and $g$ 
as in \rdef {L1}, there exists 
one and only one very weak solution $u$ of Problem (\ref {L1-1}).
Furthermore, for any $\gz\in C_{c}^{1,L}(\overline \Omega)$, $\gz\geq 0$,
there holds
\begin {eqnarray}\label {L1-3}
\int_{\Gw}\abs u L^*\gz dx\leq \int_{\Gw}f{\rm sign} (u)\gz dx
-\int_{\prt\Gw}\frac 
{\prt\gz}{\prt{\bf n}_{L^*}}\abs gdS.
\end {eqnarray}
and
\begin {eqnarray}\label {L1-4}
\int_{\Gw} u_{+} L^*\gz dx\leq \int_{\Gw}f{\rm sign}_{+} (u)\gz dx
-\int_{\prt\Gw}\frac 
{\prt\gz}{\prt{\bf n}_{L^*}}g_{+}dS.
\end {eqnarray}
\es 
The following result shows the continuity of the process.
\blemma {L1-lem} There exists a positive constant $C=C(L,\Gw)$
such that if $f$ and $g$ are as in \rdef {L1} and $u$ is a very weak 
solution of (\ref {L1-1}),
\begin {eqnarray}\label {L1-5}
\norm u_{L^1(\Gw)}\leq C\left(\norm {\gr_{_{\prt\Gw}} f}_{L^1(\Gw)}+
\norm g_{L^1(\prt\Gw)}\right).
\end {eqnarray}
\es
\Proof We denote by $\eta_{u}$ the solution of
\begin {equation}\label {L1-1e}\left.\BA{ll}
L^*\eta_{u}={\rm sign}(u)&\mbox { in }\Omega,\\[2mm]
\phantom{L^*}\eta_{u}=0&\mbox { on }\prt\Omega, 
\EA\right.
\end {equation}
Notice that $\eta_{u}$ exists by \rlemma {LM}. Since the coefficients 
of $L$ are Lipschitz continuous, $\eta_{u}\in 
C_{c}^1(\overline\Gw)$ and $L^*\eta_{u}\in L^{\ity}(\Gw)$. Thus 
$\eta_{u}\in C_{c}^{1,L}(\overline \Omega)$. By the maximum principle
$$\abs {{\eta_{u}}}\leq \eta:=\eta_{1},
$$
thus
$$\abs {\frac {\prt\eta_{u}}{\prt{\bf n}_{L^*}}}
\leq -\frac {\prt\eta}{\prt{\bf n}_{L^*}}.
$$
Pluging this estimates into (\ref {L1-2}) one obtains
\begin {eqnarray}\label {L1-6}
\int_{\Gw}\abs udx\leq \int_{\Gw}\abs f\eta dx-\int_{\prt\Gw}\frac 
{\prt\eta}{\prt{\bf n}_{L^*}}\abs gdS,
\end {eqnarray}
from which (\ref {L1-5}) follows.\qeda \\

\noindent {\it Proof of \rth {L1-th}} -{\it Existence} Let 
$\{f_{n}\}$, $\{g_{n}\}$ be two sequences of $C^2$ functions defined 
respectively in $\Gw$ and $\prt\Gw$, $f_{n}$ with compact support, 
and such that
$$\norm {(f-f_{n})\gr_{_{\prt\Gw}}}_{L^1(\Gw)}+
\norm {g-g_{n}}_{L^1(\prt\Gw)}\to 0\quad\mbox {as } n\to\infty.
$$
Let $u_{n}$ be the classical solution (derived from \rlemma {LM} for 
example) of
\begin {equation}\label {L1-7}\left.\BA{ll}
Lu_{n}=f_{n}&\mbox { in }\Gw,\\[2mm]
\phantom{L}u_{n}=g_{n}&\mbox { on }\prt\Gw. 
\EA\right.
\end {equation}
Then $u_{n}\in W^{2,p}(\Omega)$ for any finite $p\geq 1$. By (\ref {L1-5}), 
$\{u_{n}\}$ is a Cauchy sequence in $L^1(\Gw)$. Because $u_{n}$ satisfies
\begin {eqnarray}\label {L1-8}
\int_{\Gw}u_{n}L^*\gz dx=\int_{\Gw}f_{n}\gz dx-\int_{\prt\Gw}\frac 
{\prt\gz}{\prt{\bf n}_{L^*}}g_{n}dS,
\end {eqnarray}
 for any 
$\gz\in C^{1,L}_{c}(\overline\Omega)$, letting $n\to\infty$ leads to (\ref {L1-2}).\\

\noindent {\it Estimates (\ref {L1-3}) and (\ref {L1-4})}. Let 
$\gamma$ be a smooth, odd and increasing function defined on $\mathbb R$ 
such that $-1\leq\gamma\leq 1$, and $\gz$ a nonnegative element of 
$C^{1,L}_{c}(\overline\Gw)$. Since 
\begin {eqnarray*}
\int_{\Gw}f_{n}\gamma (u_{n})\gz dx&=&\sum_{i,j=1}^n\int_{\Gw}a_{ij}
\frac {\prt u_{n}}{\prt x_{j}}\frac {\prt (\gamma (u_{n})\gz)}{\prt x_{i}}
dx\\
&+&\sum_{i}^n\int_{\Gw}\left(b_{i}\frac {\prt u_{n}}{\prt x_{i}}\gamma (u_{n})\gz
+c_{i}u_{n}\frac {\prt (\gamma (u_{n})\gz)}{\prt x_{i}}\right)dx
 +\int_{\Gw}du_{n}\gamma (u_{n})\gz dx\\
&\geq &\sum_{i,j=1}^n\int_{\Gw}a_{ij}
\frac {\prt u_{n}}{\prt x_{j}}\frac {\prt \gz}{\prt x_{i}}\gamma 
(u_{n})\gz dx\\
&+&\sum_{i}^n\int_{\Gw}\left(b_{i}\frac {\prt u_{n}}{\prt x_{i}}\gamma (u_{n})\gz
+c_{i}u_{n}\frac {\prt (\gamma (u_{n})\gz)}{\prt x_{i}}\right)dx
 +\int_{\Gw}du_{n}\gamma (u_{n})\gz dx.
\end {eqnarray*}
Put
$$j_{1}(r)=\int_{0}^r\gamma (s)ds,\;j_{2}(r)=r\gamma (r)\,\mbox { and }
j_{3}(r)=\int_{0}^rs\gamma' (s)ds.
$$
Then
\begin {eqnarray*}
\sum_{i,j=1}^n\int_{\Gw}a_{ij}
\frac {\prt u_{n}}{\prt x_{j}}\frac {\prt \gz}{\prt x_{i}}\gamma 
(u_{n})dx&=&\sum_{i,j=1}^n\int_{\Gw}a_{ij}
\frac {\prt j_{1}(u_{n})}{\prt x_{j}}\frac {\prt \gz}{\prt x_{i}}dx\\
&=&-\sum_{i,j=1}^n\int_{\Gw}j_{1}(u_{n})\frac {\prt}{\prt x_{j}}
\left(a_{ij}\frac {\prt \gz}{\prt x_{i}}\right)dx
+\int_{\prt\Gw}j_{1}(g_{n})\frac 
{\prt\gz}{\prt{\bf n}_{L^*}}dS,
\end {eqnarray*}
and
\begin {eqnarray*} 
\sum_{i=1}^n\int_{\Gw}\left(b_{i}\frac {\prt u_{n}}{\prt x_{i}}\gamma (u_{n})\gz
+c_{i}u_{n}\frac {\prt (\gamma (u_{n})\gz)}{\prt x_{i}}\right)dx
\qquad\qquad\qquad\qquad\qquad\qquad\qquad\qquad\qquad\\
=\sum_{i=1}^n\int_{\Gw}\left(b_{i}\frac {\prt j_{1}(u_{n})}{\prt 
x_{i}}\gz
+c_{i}\left(j_{2}(u_{n})\frac {\prt \gz}{\prt x_{i}}+
\gz\frac {\prt j_{3}(u_{n})}{\prt x_{i}}\right)\right)dx\qquad\quad \\
=\sum_{i=1}^n\int_{\Gw}\left(-j_{1}(u_{n})\frac {\prt}{\prt 
x_{i}}(b_{i}\gz)+c_{i}j_{2}(u_{n})\frac {\prt \gz}{\prt x_{i}}
-j_{3}(u_{n})\frac {\prt}{\prt x_{i}}(c_{i}\gz)\right)dx.
\end {eqnarray*}
Therefore
\begin {eqnarray*}
\int_{\Gw}f_{n}\gamma (u_{n})\gz dx-\int_{\prt\Gw}j_{1}(g_{n})\frac 
{\prt\gz}{\prt{\bf n}_{L^*}}dS\geq 
-\int_{\Gw}\left(\sum_{i,j=1}^n\frac {\prt}{\prt x_{j}}
\left(a_{ij}\frac {\prt\gz}{\prt x_{i}}\right)+\sum_{i=1}^n\frac 
{\prt}{\prt x_{i}}(b_{i}\gz)\right)j_{1}(u_{n})dx\\
+\int_{\Gw}\left(\sum_{i=1}^nc_{i}j_{2}(u_{n})\frac {\prt\gz}{\prt x_{i}}
-j_{3}(u_{n})\frac {\prt}{\prt 
x_{i}}(c_{i}\gz)+dj_{2}(u_{n})\gz\right)dx,
\end {eqnarray*}
and finally,
\begin {eqnarray*}
 \int_{\Gw} j(u_{n})L^*\gz dx\leq \int_{\Gw}f_{n}\gamma (u_{n})\gz dx-
 \int_{\prt\Gw} j(g_{n})\frac 
{\prt\gz}{\prt {\bf n}_{L^*}}dS.
\end {eqnarray*}
When $\gamma (r)\to {\rm sign}(r)$, $j_{1}(r)$ and $j_{2}(r)$ both converge to 
$\abs r$, and $j_{3}(r)$ converges to $0$ if, for example, we impose 
$0\leq \gamma'_{\epsilon}(r)\leq 
2\epsilon^{-1}\chi_{_{(-\epsilon,\epsilon)}}(r)$ and send $\epsilon$ to $0$. Letting successively $n\to\infty$ and 
$\gamma\to {\rm sign}$ yields to (\ref {L1-3}). We obtain (\ref {L1-4}) in the same way while 
 approximating sign$_{+}$ by $\gamma$.\qeda 
 
 \bcor {order} Under the assumptions of \rth {L1-th}, the mapping 
 $(f,g)\mapsto u$ defined by (\ref {L1-1}) is increasing.
 \es 

 For the regularity of solutions, the following result is due to 
 Brezis and Strauss \cite {BS} using  Stampacchia's duality method 
 \cite {St3}.
 \bth {L^1reg}Let $L$ satisfy the condition (H). Then for any  
 $1\leq q<n/(n-1)$, there exists $C=C(\Gw,q)>0$ such that for any $f\in 
 L^1(\Gw)$, the very weak solution $u$ of (\ref {L1-1}) with $g=0$ 
 satisfies
 \begin {equation}\label {reg1}
 {\norm u}_{W^{1,q}_{0}(\Gw)}\leq C\norm f_{L^1(\Gw)}.
 \end {equation}
 \es
 This theorem admits a local version.
  \bcor {L^1regloc} Let $L$ be the elliptic operator defined by (\ref 
  {lin1}), with Lipschitz continuous coefficients and satisfying 
  (\ref {lin2}). Let $u\in L^1_{loc}(\Gw)$ and $f\in L^1_{loc}(\Gw)$ 
  be such that 
   \begin {equation}\label {reg2}
\int_{\Gw}uL^*\gz dx=\int_{\Gw}f\gz dx,
 \end {equation}
 for any $\gz\in C^1_{c}(\Gw)$ such that $L^*\gz\in L^\ity(\Gw)$. Then for 
 any open subsets $G\subset\overline G\subset G'\subset \overline G'\subset 
 \Gw$, with $\overline G'$ compact and $1\leq q<n/(n-1)$, there exists 
 a constant $C=C(G,G',q,L)>0$ such that 
  \begin {equation}\label {reg3}
 {\norm u}_{W^{1,q}(G)}\leq C\left(\norm f_{L^1(G')}+\norm u_{L^1(G')}\right).
 \end {equation}
 \es 

\subsection {The measure framework}
We denote by $\mathfrak M(\Gw)$ and $\mathfrak M(\prt\Gw)$ the spaces of Radon 
measures on $\Omega$ and $\prt\Gw$ respectively, by $\mathfrak M_{+}(\Gw)$ 
and $\mathfrak M_{+}(\prt\Gw)$ their positive cones. For $0\leq\ga\leq 1$, we also denote by 
$\mathfrak M(\Gw;\gr^\ga_{_{\prt\Gw}})$ the 
subspace of the $\gm\in\mathfrak M(\Gw)$ satisfying
$$\int_{\Gw}\gr^\ga_{_{\prt\Gw}}d\abs\gm <\infty,
$$
and by $C(\overline\Omega;\gr^{-\ga}_{_{\prt\Gw}})$ the subspace of 
$C(\overline\Omega)$ of functions $\gz$ such that 
$$\sup_{\Gw}\abs{\gz}/\gr^\ga_{_{\prt\Gw}}<\infty.
$$
For the sake of clarity, we denote by
$$\mathfrak M(\Gw;\gr^0_{_{\prt\Gw}})=\mathfrak M^{b}(\Gw),$$
the space of bounded Radon measures in $\Gw$. 
Both $\mathfrak M(\Gw;\gr^\ga_{_{\prt\Gw}})$ and 
$C(\overline\Omega;\gr^{-\ga}_{_{\prt\Gw}})$ are endowed with the norm 
corresponding to their definition.
If $\gl\in \mathfrak M(\Gw;\gr_{_{\prt\Gw}})$ and 
$\gm\in \mathfrak M(\prt\Gw)$, the definition of a very weak solution to 
the measure data problem
\begin {equation}\label {M1}\left.\BA{ll}
Lu=\gl&\mbox { in }\Gw,\\[2mm]
\phantom{L}u=\gm&\mbox { on }\prt\Gw, 
\EA\right.
\end {equation}
is similar to \rdef {L1} : $u\in L^1(\Gw)$ and the equality
\begin {eqnarray}\label {M1-2}
\int_{\Gw}uL^*\gz dx=\int_{\Gw}\gz d\gl-\int_{\prt\Gw}\frac 
{\prt\gz}{\prt{\bf n}_{L^*}}d\gm,
\end {eqnarray}
holds for every $\gz\in C^{1,L}_{c}(\overline\Gw)$. 
\bth {M1-th} Let $L$ satisfy the condition (H). For every $\gl\in \mathfrak M(\Gw;\gr_{_{\prt\Gw}})$ and 
$\gm\in \mathfrak M(\prt\Gw)$ there exists a unique very weak solution 
$u$ to Problem (\ref {M1-2}). Furthermore the mapping $(\gl,\gm)\mapsto u$ 
is increasing.
\es
\Proof Uniqueness follows from \rlemma {L1-lem}. For existence, let 
$\{\gl_{n}\}$ be a sequence of smooth functions  
in $\overline\Omega$ such that 
$$\lim_{n\to\infty}\int_{\Gw}\gl_{n}\phi 
dx=\int_{\Gw}\phi d\gl,
$$
for every $\phi\in C(\overline\Gw;\gr^{-1}_{_{\prt\Gw}})$. Let $\{\gm_{n}\}$ be a sequence of $C^2$ 
functions on $\prt\Gw$ converging to $\gm$ in the weak sense of 
measures and $u_{n}$ denote the classical solution of
\begin {equation}\label {M1'}\left.\BA{ll}
Lu_{n}=\gl_{n}&\mbox { in }\Gw,\\[2mm]
\phantom{L}u_{n}=\gm_{n}&\mbox { on }\prt\Gw. 
\EA\right.
\end {equation}
Thus
\begin {eqnarray}\label {M1-n}
\int_{\Gw}u_{n}L^*\gz dx=\int_{\Gw}\gz \gl_{n} dx-\int_{\prt\Gw}\frac 
{\prt\gz}{\prt{\bf n}_{L^*}}\gm_{n}dS,
\end {eqnarray}
holds for every $\gz\in C^{1,L}_{c}(\overline\Gw)$. 
Since
$\norm{\gl_{n}\rho_{_{\prt\Gw}}}_{L^1(\Omega)}$ and 
$\norm{\gm_{n}}_{L^1(\prt\Omega)}$ are bounded
independently of $n$, it is the same with
$\norm{u_{n}}_{L^1(\Omega)}$ by \rlemma {L1-lem}. Let $\gw$ be a Borel subset of 
$\overline\Gw$, and $\gth_{\omega,n}$ the solution of
\begin {equation}\label {M1-3}\left.\BA{ll}
L^*\gth_{\omega,n}=\chi_{_{\omega}}{\rm sign}(u_{n})&\mbox { in }\Gw,\\[2mm]
\phantom{L^*}\gth_{\omega,n}=0&\mbox { on }\prt\Gw. 
\EA\right.
\end {equation}
Since $\gth_{\omega}$ is an admissible test function, 
$$
\int_{\omega}\abs {u_{n}} dx=\int_{\Gw}\gth_{\omega,n} \gl_{n}dx
-\int_{\prt\Gw}\frac 
{\prt\gth_{\omega,n}}{\prt{\bf n}_{L^*}}\gm_{n}dS.
$$
Moreover $-\gth_{\omega}\leq\gth_{\omega,n}\leq \gth_{\omega}$, where 
$\gth_{\omega}$ is the solution of 
\begin {equation}\label {M1-4}\left.\BA{ll}
L^*\gth_{\omega}=\chi_{_{\omega}}&\mbox { in }\Gw,\\[2mm]
\phantom{L^*}\gth_{\omega}=0&\mbox { on }\prt\Gw. 
\EA\right.
\end {equation}
Therefore
\begin {eqnarray}\label {M1-5}
\int_{\omega}\abs {u_{n}} dx
\leq\norm{\gl_{n}\rho_{_{\prt\Gw}}}_{L^1(\Omega)}{\norm 
{\gth_{\omega}/\gr_{_{\prt\Gw}}}}_{L^{\ity}(\Gw)}+
\norm{\gm_{n}}_{L^1(\prt\Omega)}{\norm 
{\prt\gth_{\omega}/\prt{\bf n}_{L^*}}}_{L^{\ity}(\prt\Gw)}.
\end {eqnarray}
By the $L^p$ regularity theory for elliptic equations and the
Sobolev-Morrey imbedding Theorem, for any $n<p<\infty$, there exists a constant $C=C(n,p)>0$ such that 
\begin {eqnarray}\label {M1-6}
\norm{\gth_{\omega}}_{C^1(\overline\Omega)}\leq 
C{\norm{\chi_{_{\omega}}}}_{L^p(\Gw)}=C{\abs\omega}^{1/p}.
\end {eqnarray}
This estimate, combined with (\ref {M1-5}), yields to
\begin {eqnarray}\label {M1-7}
\int_{\omega}\abs {u_{n}} dx\leq 
C(\norm{\gl_{n}\rho_{_{\prt\Gw}}}_{L^1(\Omega)}+
\norm{\gm_{n}}_{L^1(\prt\Omega)}){\abs\omega}^{1/p}\leq CM{\abs\omega}^{1/p},
\end {eqnarray}
for some $M$ independent of $n$. Therefore the 
sequence $\{u_{n}\}$ is uniformly integrable,
thus weakly compact in $L^1(\Omega)$ by the Dunford-Pettis Theorem, and there 
exist a subsequence $\{u_{n_{k}}\}$ and an integrable function $u$ 
such that $u_{n_{k}}\to u$, weakly in $L^1(\Gw)$. Passing to the limit in (\ref{M1-n})
leads to (\ref{M1-2}). Because of uniqueness the whole sequence 
$\{u_{n}\}$ converges weakly to $u$. The monotonicity assertion 
follows from uniqueness and \rcor {order}.\qeda
\medskip

\noindent \Remark Estimate (\ref {L1-4}) in the statement of \rth 
{L1-th} admits the following extension : Let the two 
measures $\gl$ and $\gm$ have Lebesgue decomposition
$$\gl=\gl_{r}+\gl_{s} \quad\mbox {and }\gm=\gm_{r}+\gm_{s},$$
$\gl_{r}$ and $\gm_{r}$ being the regular parts with respect to the $n$ 
and the $n-1$ dimensional Hausdorff measures and $\gl_{s}$ and $\gm_{s}$
the singular parts. If $\gl_{s}$ and $\gm_{s}$ are nonpositive, there 
holds
\begin {eqnarray}\label {M1-7'}
\int_{\Gw}u_{+}L^*\gz dx\leq 
\int_{\Gw}\gl_{r\,+}{sign}_{+}(u)\gz dx-\int_{\prt\Gw}\frac 
{\prt\gz}{\prt{\bf n}_{L^*}}\gm_{r\,+}dS,
\end {eqnarray}
for any $\gz\in C^{1,L}_{c}(\overline\Gw)$, $\gz\geq 0$.\medskip 

\noindent \Remark The above proof implies the following weak stability 
result. {\it If $\{\gl_{n}\}\subset\mathfrak M (\Gw;\gr_{_{\prt\Gw}})$ and 
$\{\gm_{n}\}\subset\mathfrak M (\prt\Gw)$
are sequences of 
measures wich converge respectively to $\gl$ in duality with 
$C(\overline\Gw;\gr^{-1}_{_{\prt\Gw}})$, and to $\gm$ in the weak 
sense of measures on $\prt\Gw$, the corresponding very weak solutions $u_{n}$ 
of (\ref{M1'}) converge weakly in $L^1(\Gw)$ 
to the very weak solution $u$ of (\ref{M1}).}
\subsection {Representation theorems and boundary trace}

If $\Gw$ is a bounded domain with a $C^2$ boundary, $L$ the elliptic 
operator defined by (\ref {lin1}), with Lispchitz countinuous 
coefficients and $u$ and $v$ two functions in $W^{2,p}(\Omega)$, 
with $p>n$, the Green formula implies
\begin {eqnarray}\label {R1}
\int_{\Omega}\left(vLu-uL^*v\right)dx=
\int_{\prt\Gw}\left(u\frac {\prt v}{\prt{\bf n}_{L^*}}-
v\frac {\prt u}{\prt{\bf n}_{L}}\right)dS,
\end {eqnarray}
where $L^*$ and $\prt v/\prt{\bf n}_{L^*}$ are respectively defined by 
(\ref {lin6}) and (\ref {conorm}), and

\begin {eqnarray}\label {norm}
\frac {\prt \gz}{\prt {\bf n}_{L}}=\sum_{i,j=1}^n
a_{ij}\frac {\prt \gz}{\prt x_{j}}{\bf n}_{i},
\end {eqnarray}
is the co-normal derivative following $L$. If we assume that 
condition (H) is fulfilled, and if $x\in\Gw$, we denote by $G^\Gw_{L}(x,.)$ the solution 
of 
\begin {equation}\label {Green}\left.\BA{ll}
L^*G^\Gw_{L}(x,.)=\gd_{x}&\quad\mbox {in }\Gw,\\[2mm]
\phantom{L^*}G^\Gw_{L}(x,.)=0&\mbox { on }\prt\Gw. 
\EA\right.
\end {equation}
The function $G^\Gw_{L}$ is the {\it Green function} of the operator $L$ in $\Gw$. 
Notice an ambiguity in terminology between $L$ and $L^*$, but it has 
no consequence because the condition (H) is invariant by duality and
the following symmetry result holds :
\begin {eqnarray}\label {Green sym}
G^\Gw_{L}(x,y)=G^\Gw_{L^*}(y,x),\forevery (x,y)\in 
\Gw\times\Gw,\;x\neq y.
\end {eqnarray}
The 
function $G^\Gw_{L}(x,.)$ is nonnegative by \rth {M1-th} and belongs to  
$W_{loc}^{2,p}(\Gw\setminus\{x\})$ for any $1<p<\infty$. Thus it is $C^1$ 
in $\overline\Gw\setminus\{x\}$. We denote 
\begin {eqnarray}\label {Poisson}
P^\Gw_{L}(x,y)=-\frac {\prt 
G^\Gw_{L}(x,y)}{\prt{\bf n}_{L^*}},\forevery (x,y)\in \Gw\ti\prt\Gw.
\end {eqnarray}
If $u\in C^2(\bar\Gw)$, the following Green representation formula derives from (\ref {R1})
\begin {eqnarray}\label {R2}
u(x)=\int_{\Omega}G^\Gw_{L}(x,y)Lu(y) 
dy+\int_{\prt\Gw}P^\Gw_{L}(x,y)u(y)dS(y),\forevery x\in\Gw.
\end {eqnarray}
By extension this representation formula holds almost everywhere if 
$(\gl,\gm)\in\mathfrak M(\Gw)\ti\mathfrak (\prt\Gw)$, and $u$ is
the very weak solution of (\ref {M1}), in the sense that
\begin {eqnarray}\label {R'2}
u(x)=\int_{\Gw}G_{L}^\Gw(x,y)d\gl(y)+
\int_{\Gw}G_{L}^\Gw(x,y)d\gm(y),\quad \mbox{a.e. in } \Gw.
\end {eqnarray}
Actually the representation formula is equivalent to the fact that $u$ 
is a very weak solution of Problem \ref {M1} (see \cite {BVV} for a 
proof). We set
\begin {eqnarray}\label {greenpot}
\mathbb G_{L}^\Gw(\gl)(x)=\int_{\Gw}G_{L}^\Gw(x,y)d\gl(y),
\end{eqnarray}
and call it the {\it Green potential} of $\gl$, and
\begin {eqnarray}\label {poispot}
\mathbb P_{L}^\Gw(\gl)(x)=\int_{\prt\Gw}P_{L}^\Gw(x,y)d\gl(y),
\forevery x\in\Gw,
\end{eqnarray}
the {\it Poisson potential }of $\gm$.
The Green kernel presents a singularity 
on the diagonal $D_{\Gw}=\{(x,x):x\in\Gw\}$, while the Poisson kernel 
becomes singular when the $x$ variable approaches the boundary point 
$y$. Many estimates on the singularities  have been obtained in the 
last thirty years \cite {Kr}, \cite {Mi}, \cite {GT}, \cite {Dy2}. We give below 
some useful estimates in which $\gr_{_{\prt\Gw}}$ is defined by (\ref 
{distf}).
\bth {estim} Assuming that $\Gw$ is bounded with a $C^2$ boundary and 
assumption (H) holds, then
\begin {eqnarray}\label {estG1}
G_{L}^\Gw(x,y)\leq C(L,\Omega)\frac {\min\left\{1,
\abs {x-y}\gr_{_{\prt\Gw}}(x)\right\}}
{\abs {x-y}^{n-2}}
,\forevery (x,y)\in(\Gw\times\Gw)\setminus D_{\Gw},
\end {eqnarray}
if $n\geq 3$, 
\begin {eqnarray}\label {estG2}
G_{L}^\Gw(x,y)\leq C(L,\Omega)\min\left\{1,
\abs {x-y}\gr_{_{\prt\Gw}}(x)\right\}
\ln_{+}\abs {x-y}
,\forevery (x,y)\in(\Gw\times\Gw)\setminus D_{\Gw},
\end {eqnarray}
if $n=2$. Moreover, for any $n\geq 2$,
\begin {eqnarray}\label {estP1}
K'(L,\Omega)\frac {\gr_{_{\prt\Gw}}(x)}{\abs{x-y}^{n}}\leq P_{L}^\Gw(x,y)\leq 
K(L,\Omega)\frac {\gr_{_{\prt\Gw}}(x)}{\abs{x-y}^{n}}
,\forevery (x,y)\in\Gw\times\prt\Gw.
\end {eqnarray}
\es 
\medskip 

 Another useful notion, from which some of the above estimates can be 
 derived is the notion of equivalence (see \cite{An},  \cite{Pi}). 
\bth {equivth} Assuming that $\Gw$ is bounded with a $C^2$ boundary and 
assumption (H) holds, there exists a positive constant $C$ such that 
\begin {eqnarray}\label {equivG}
C G_{-\Gd}^\Gw \leq G_{L}^\Gw\leq \frac {1}{C} G_{-\Gd}^\Gw\quad \mbox 
{in }\,(\Gw\times\Gw)\setminus 
D_{\Gw},
\end {eqnarray}
and
\begin {eqnarray}\label {equivP}
C P_{-\Gd}^\Gw \leq P_{L}^\Gw\leq \frac {1}{C} P_{-\Gd}^\Gw
\quad \mbox 
{in }\,\Gw\times\prt\Gw.
\end {eqnarray}
\es 
 
In order to study the boundary behaviour of harmonic functions, we 
introduce, for $\gb>0$,
\begin {eqnarray}\label {fol1}
\Gw_{\gb}=\{x\in\Gw:\gr_{_{\Gw}}(x)>\gb\},\quad 
G_{\gb}=\Gw\setminus\overline\Gw_{\gb},\quad
\Gs_{\gb}=\prt\Gw_{\gb}=\{x\in\Gw:\gr_{_{\Gw}}(x)=\gb\},
\end {eqnarray}
and $\Gs_{0}:=\Gs:=\prt\Gw$. Since $\Gw$ is $C^2$, there exists 
$\gb_{0}>0$ such that for every $0<\gb\leq\gb_{0}$ and $x\in G_{\gb}$ 
there exists a unique $\gs(x)\in\Gs$ such that $\abs 
{x-\gs(x)}=\gr_{_{\prt\Gw}}(x)$. We denote by $\Gp$ the mapping
 from $G_{\gb}$ to $(0,\gb)\times\Gs$ defined by
\begin {eqnarray}\label {fol12}
 \Gp(x)=(\gr_{_{\prt\Gw}}(x),\gs(x)).
\end {eqnarray}
The mapping $\Gp$ is a $C^2$ diffeomorphism with inverse given by
\begin {eqnarray}\label {fol13}
\Gp^{-1}(t,\gs)=\gs-t{\bf n},\forevery (t,\gs)\in (0,\gb)\times\Gs,
\end {eqnarray}
where ${\bf n}$ is the normal unit outward vector to $\prt\Gw$ at $x$ 
(see \cite {MV5} for details). If the distance coordinate is fixed 
in $(0,\gb_{0}]$, the mapping $\mathfrak H_{t}$ defined by 
$$\mathfrak H_{t}(x)=\gs(x)\forevery x\in\Gs_{t},$$ 
is the orthogonal 
projection from $\Gs_{t}$ to $\prt\Gw$. 
Thus $\mathfrak H_{t}^{-1}(.)=\Gp^{-1}(t,.)$
 is a $C^2$ diffeomorphism and the set 
$\{\Gs_{t}\}_{0<t\leq\gb_{0}}$ is a $C^2$ foliation of $G_{\gb_{0}}$.
For $0<t\leq\gb_{0}$ we can transfer naturally a Borel measure $\gm$, or a 
function $f$, on $\Gs_{t}$ into a Borel measure or a function on $\prt\Gw$ 
as follows :
\begin {equation}\label {fol14}\left.\BA{l}
\gm^t(E):=\gm (\mathfrak H_{t}^{-1}(E)), \quad\mbox {for every Borel subset 
E}\subset\prt\Gw,\\[2mm]
\,f^t(x):=f(\gs(x)-t{\bf n}(x)),\forevery x\in\prt\Gw. 
\EA\right.
\end {equation}
The Lebesgue classes on $\Gs_{t}$ and $\Gs$ are exchanged by this 
projection operator and actually
\begin {eqnarray}\label {fol15}
\gm\in \mathfrak M(\Gs_{t}),\,f\in L^1(\Gs_{t},\gm)\Longrightarrow
\left\{\BA{l}f^t\in L^1(\Gs,\gm^t),\\
\myint{\Gs_{t}}{}fd\gm=\myint{\Gs}{}f^td\gm^t.\EA\right.
\end {eqnarray}
\bdef {L-harm}{\rm Let $L$ be an elliptic operator defined by (\ref 
{lin1}) in $\Gw$, with bounded and measurable coefficients. We say that a 
function $u\in W^{1,2}_{loc}(\Gw)$ is weakly $L$-harmonic if}
\begin {eqnarray}\label {L-harm}
A_{L}(u,v)=0\forevery v\in C^1_{c}(\Gw).
\end {eqnarray}
\es 
\Remark If (\ref {lin2}) holds, any weakly $L$-harmonic 
function is H\" older continuous by the De Giorgi-Nash-Moser Theorem. 
If the coefficients of $L$ are Lispchitz continuous, the notion of 
weak $L$-harmonicity can be understood in the sense of distributions in 
$\Gw$, by assuming that 
$u$ is locally integrable in $\Gw$ and 
\begin {eqnarray}\label {L-harm2}
\int_{\Gw}u L^{*}\phi dx=0,\forevery \phi\in C_{0}^{\ity}(\Omega).
\end {eqnarray}
It can be verified that any locally integrable function $L$-harmonic in 
the sense of distributions in $\Gw$ is actually weakly $L$-harmonic. Therefore 
it belongs to $W^{2,p}_{loc}(\Gw)$, for any $1<p<\infty$, by the 
$L^p$-regularity theory of elliptic equations.

\bth {trace-th} Let $\Gw$ be a bounded domain of class $C^{2}$ and $L$ 
the elliptic operator defined by (\ref {lin1}) satisfying condition (H). 
Let $u$ be a nonnegative locally integrable $L$-harmonic function in $\Gw$.
Then there exists a unique nonnegative Radon measure 
$\gm$ on $\prt\Gw$ such that 
\begin {eqnarray}\label {L-harm3}
\lim_{t\to 0}\int_{\Gs_{t}}u(x)\gth(\gs(x)) dS
=\int_{\Gs}\gth d\gm,\forevery \gth\in C_{0}(\Gs).
\end {eqnarray}
Moreover $u$ is uniquely determined by $\gm$ and
\begin {eqnarray}\label {L-harm4}
u(x)=\int_{\prt\Gw}P^\Gw_{L}(x,y)d\gm(y),\forevery x\in\Gw.
\end {eqnarray}
\es
\Proof {\it Step 1 } The function $u$ is integrable. Let $0<\gb\leq\gb_{0}$. 
Since $u$ is continuous in 
$\overline\Gw_{\gb}$, its restriction to this set is the very weak solution of
\begin {equation}\label {L-harm5}\left.\BA{ll}
Lv=0&\mbox { in }\Gw_{\gb},\\[2mm]
\phantom{L}v=u_{\vline\Gs_{\gb}}&\mbox { on }\Gs_{\gb}.
\EA\right.
\end {equation}
Thus, if $\gz\in C_{c}^{1,L}(\overline\Gw_{\gb})$, there holds
\begin {eqnarray}\label {L-harm6}
\int_{\Gw_{\gb}}uL^{*}\gz dx=-\int_{\Gs_{\gb}} 
\frac{\prt\gz}{\prt{\bf n}_{L^{*}}}udS.
\end {eqnarray}
We fix $\gz=\eta_{1,\gb}$ as the solution of 
\begin {equation}\label {L-harm7}\left.\BA{ll}
L^*\eta_{1,\gb}=1&\mbox { in }\Gw_{\gb},\\[2mm]
\phantom{L^*}\eta_{1,\gb}=0&\mbox { on }\Gs_{\gb}. 
\EA\right.
\end {equation}
By Hopf's lemma, there exists $c>0$ such that 
$$
c\leq-\frac{\prt\eta_{1,\gb}}{\prt{\bf n}_{L^{*}}}\leq c^{-1}\quad\mbox {on }\Gs_{\gb}.
$$
Moreover, $c$ can be taken independent of $\gb\in(0,\gb_{0}]$. It 
follows 
\begin {eqnarray}\label {L-harm8}
\Psi(\gb)=\int_{\Gw_{\gb}}udx\geq 
c\int_{\Gs_{\gb}}udS=-c\Psi'(\gb),
\end {eqnarray}
from (\ref {L-harm6}), and 
$$\frac {d}{d\gb}\left(e^{\gb/c}\Psi(\gb)\right)\geq 0.
$$
Therefore
$$\lim_{\gb\to 0}\Psi(\gb)=\int_{\Gw}udx<\infty.
$$
Notice that (\ref{L-harm8}) implies that $\norm u_{L^{1}(\Gs_{\gb})}$ 
remains bounded independently of $\gb$.\smallskip 

\noindent {\it Step 2  } End of the proof. Let $\gth\in C^{2}(\prt\Gw)$,  
$w_{\gth}$ 
be the solution of 
\begin {equation}\label {L-harm9}\left.\BA{ll}
L^*w_{\gth}=0&\mbox { in }\Gw_{\gb},\\[2mm]
\phantom{L^*}w_{\gth}=\gth&\mbox { on }\Gs_{\gb}, 
\EA\right.
\end {equation}
and $h\in C(\Gs_{\gb})$ defined by
$$h=-\frac {\prt \eta_{1,\gb}}{\prt{\bf n}_{L^{*}}}.
$$
Then $\gz=\eta_{1,\gb}w_{\gth}h^{-1}$ belongs to $C^{1,L}_{c}(\overline\Gw_{\gb})$. 
Since $\prt\gz/\prt{\bf n}_{L^{*}}=\gth$ on $\Gs_{\gb}$,
$$\int_{\Gw_{\gb}}uL^{*}\gz dx=-\int_{\Gs_{\gb}}\frac 
{\prt\gz}{\prt{\bf n}_{L^{*}}}dS=\int_{\Gs_{\gb}}\gth u dS.
$$
It is easy to check that $L^{*}\gz $ is bounded in $L^{\ity}(\Gw_{\gb})$, 
independently of $\gb$. Therefore
$$\lim_{\gb\to 0}\int_{\Gw_{\gb}}uL^{*}\gz dx
$$
exists. The same holds true with
$$\lim_{\gb\to 0}\int_{\Gs_{\gb}}\gth udS,
$$
which defines a positive linear functional on $C^{2}(\prt\Gw)$. This 
characterizes the Radon measure $\gm$ in a unique way.\qeda 
\bdef {tr-def}{\rm The measure $\gm$ is called {\it the boundary trace} 
of $u$.}
\es
\Remark In the above theorem, many assumptions can be relaxed : the 
boundedness of $\Gw$ plays no role except that it allows a simpler statement of 
the result, and the integral represention (\ref {L-harm4}). The 
regularity of the boundary of the domain is not a key assumption, but 
in the case of a general domain, the boundary has to be replaced by 
the Martin boundary \cite {Ma}, and the Poisson kernel by the Martin 
kernel in order to have a representation formula valid for all the positive 
$L$-harmonic functions. \\

\noindent \Remark The Fatou Theorem asserts that for almost all $y\in\prt\Gw$ 
(for the $n-1$-dimensional Hausdorff measure $dH_{n-1}$) and for any cone $C_{y}$ 
interior to $\Gw$ the following limit exists,
\begin {equation}\label {fat}
\lim_{\tiny{\begin {array} {l}x\to y\\ x\in C_{y}\end {array}}} u(x)=\gm_{_{r}},
\end {equation}
where $\gm_{_{r}}$ is the regular part of the measure $\gm$ with respect 
to $dH_{n-1}$ in the Lebesgue decomposition. The proof of this result 
\cite {DaL}, \cite {Do},
is much more involved that the one of \rth {trace-th}. The trace in the sense of Radon measures 
is much more useful 
for our next considerations.\\

\bdef {L-supharm}{\rm A locally  integrable function $u$ defined in $\Gw$ is 
said super-$L$-harmonic if 
\begin {eqnarray}\label {L-supharm1}
\int_{\Gw}u L^{*}\phi dx\geq 0,\forevery \phi\in 
C_{0}^{\ity}(\Omega), \;\phi\geq 0.
\end {eqnarray}
}\es

\rth {trace-th} admits an extension to positive 
proven by Doob \cite {Lk}, \cite {Do}. 
\bth {trace-th2} Let $\Gw$ be a bounded domain of class $C^{2}$ and 
$L$ the elliptic operator defined by (\ref {lin1}). We assume that 
condition (H) holds. Let $u$ a nonnegative super-$L$-harmonic in $\Gw$. 
Then there exist two Radon measures $\gl\in\mathfrak M_{+}(\Gw)$ and 
$\gm\in \mathfrak M_{+}(\prt\Gw)$, such that
$$\int_{\Gw}\gr_{_{\prt\Gw}}d\gl<\ity,
$$
and $u$ is the unique
very weak solution to Problem (\ref {M1}). 
Furthermore (\ref {fat}) holds.
\es
\mysection {Semilinear equations with absorption}
In this section we consider the semilinear Dirichlet problem with 
right-hand side measure
\begin {equation}\label {SLDM}\left.\BA{rl}
Lu+g(x,u)=\gl&\mbox { in }\Omega,\\[2mm]
u=0&\mbox { on }\prt\Omega, 
\EA\right.
\end {equation}
in a bounded domain $\Gw$ of $\mathbb R^n$, where $g$ is a continuous 
function defined on $\mathbb R\times\Omega$, $\gl$ a Radon measure in 
$\Omega$ and $L$ the elliptic operator 
with Lipschitz continuous coefficients, defined by (\ref {lin1}).
\bdef {SLDMdef} Let $\gl\in \mathfrak M(\Gw;\gr_{_{\prt\Gw}})$. A 
function $u$ is a solution of (\ref {SLDM}), if $u\in L^1(\Gw)$, 
$g(.,u)\in L^1(\Gw;\gr_{_{\prt\Gw}}dx)$, and
if for any $\gz\in C_{c}^{1,L}(\overline \Omega)$, there holds
\begin {eqnarray}\label {SDLM1}
\int_{\Gw}\left( uL^*\gz +g(x,u)\gz\right)dx=\int_{\Gw}\gz d\gl.
\end {eqnarray}
\es 
The nonlinearity is understood as an absorption term, this means that 
$rg(x,r)$ is nonnegative for $\abs r\geq r_{0}$, uniformly with respect 
to $x\in\Gw$.
\bprop {genest}Let $L$ be the elliptic operator defined by (\ref {lin1}), 
satisfying the condition (H), and $\gl\in \mathfrak M(\Gw;\gr^\ga_{_{\prt\Gw}})$ for some $0\leq \ga\leq 
1$. If $g\in C(\Gw,\mathbb R)$ is an absorption nonlinearity 
satisfying 
$$rg(x,r)\geq 0,\forevery (x,r)\in\Gw\times \left((-\infty,- r_{0}]\cup 
[r_{0},\infty)\right),$$ 
and $g$ bounded on $\Gw\times 
(-r_{0},r_{0})$, any 
solution $u$ of (\ref {SLDM}) verifies $g(.,u)\in L^1(\Gw;\gr^\ga_{_{\prt\Gw}}dx)$.
\es
\Proof We set $h=g(.,u)$, then $u$ is the unique very weak solution of 
$$Lu=\gl-h\quad\mbox {in }\Gw,
$$
and, by assumption, $u\in L^1(\Gw)$, $h\in L^1(\Gw;\gr_{_{\prt\Gw}}dx)$. 
Let $\{\gl_{n}\}$ be a sequence of smooth functions in $\overline\Gw$ 
converging to $\gl$ in the weak sense of measures with duality with 
the space $C(\overline\Gw;\rho^{-\ga}_{_{\prt\Gw}})$ 
(thus ${\norm{\gl_{n}}}_{\mathfrak M(\Gw;\gr^\ga_{_{\prt\Gw}})}$ is 
bounded independently of $n$) and $\{u_{n}\}$ the 
corresponding sequence of solutions of 
$$Lu_{n}=\gl_{n}-h\quad\mbox {in }\Gw.
$$
By \rth {L1-th}, ${\norm{u_{n}}}_{L^1(\Gw)}$ is bounded independently 
of $n$, and for any $\gz\in 
C_{c}^{1,L}(\overline\Gw)$, $\gz\geq 0$, there holds
\begin {eqnarray*}\int_{\Gw}\left(\abs {u_{n}}L^*\gz +h\gz{\rm sign} (u_{n})\right)dx\leq 
\int_{\Gw}\gl_{n}\gz{\rm sign}  (u_{n}) dx.
\end {eqnarray*}
For test function $\gz$, we take 
$j_{\epsilon}(\eta_{1})=(\eta_{1}+\epsilon)^\ga-\epsilon^\ga$ 
where $\epsilon >0$. Then $0\leq j_{\epsilon}(\eta_{1})\leq 
\eta^\ga_{1}$, and, if we put $r_{1}=\sup_{\Gw}\eta_{1}$, one obtains

\begin {eqnarray*}
L^*(j_{\epsilon}(\eta_{1}))&=&-j'_{\epsilon}(\eta_{1})
\sum_{i,j=1}^n \frac {\prt}{\prt x_{j}}\left(a_{ij}\frac {\prt 
\eta_{1}}{\prt x_{i}}\right)
+j'_{\epsilon}(\eta_{1})\sum_{i=1}^nc_{i}\frac {\prt\eta_{1}}{\prt x_{i}}
-\sum_{i=1}^n\frac {\prt}{\prt x_{i}}\left(b_{i}j_{\epsilon}(\eta_{1})\right)
 +dj_{\epsilon}(\eta_{1})\\
 &\,&-j''_{\epsilon}(\eta_{1})\sum_{i,j=1}^n a_{ij}\frac 
 {\prt \eta_{1}}{\prt x_{j}}\frac{\prt \eta_{1}}{\prt x_{i}}\\
 &=&j'_{\epsilon}(\eta_{1})L^*\eta_{1}
 + \left (j_{\epsilon}(\eta_{1})-\eta_{1}j'_{\epsilon}(\eta_{1})\right)
 \left (d-\sum_{i=1}^n\frac {\prt b_{i}}{\prt x_{i}}\right)
 -j''_{\epsilon}(\eta_{1})\sum_{i,j=1}^n a_{ij}\frac 
 {\prt \eta_{1}}{\prt x_{j}}\frac{\prt \eta_{1}}{\prt x_{i}}\\
 &\geq &-\left (j_{\epsilon}(r_{1})-r_{1}j'_{\epsilon}(r_{1})\right)
 \abs {d-\sum_{i=1}^n\frac {\prt b_{i}}{\prt x_{i}}},
\end {eqnarray*}
since $L^*\eta_{1}=1$, $j_{\epsilon}$ is a concave 
and increasing function on $\mathbb R_{+}$, 
$r\mapsto j_{\epsilon}(r)-rj'_{\epsilon}(r)$ is positive and increasing, and 
ellipticity condition (\ref {lin2}) holds. Because 
$\left (j_{\epsilon}(r_{1})-r_{1}j'_{\epsilon}(r_{1})\right)$ is 
bounded when $0<\epsilon\leq 1$, and the coefficients $b_{i}$ and $d$ are 
respectively Lipschitz continuous and bounded in $\Gw$, there exists 
a constant $M>0$ independent of $\epsilon$ such that $L^*(j_{\epsilon}(\eta_{1}))\geq -M$ 
in $\Omega$. Therefore
$$-M\int_{\Gw}\abs{u_{n}}dx+\int_{\Gw}h j_{\epsilon}(\eta_{1}){\rm 
sign}(u_{n})dx\leq {\norm{\gl_{n}}}_{\mathfrak M(\Gw;\gr^\ga_{_{\prt\Gw}})}.
$$
Letting $n\to\infty$ yields to
\begin {eqnarray}\label {j-epsilon}
\int_{\Gw}g(x,u) j_{\epsilon}(\eta_{1}){\rm sign}
(u)dx\leq M\int_{\Gw}\abs{u}dx
+\sup_{n}{\norm{\gl_{n}}}_{\mathfrak M(\Gw;\gr^\ga_{_{\prt\Gw}})},
\end {eqnarray}
since $h\in L^1(\Gw;\rho_{_{\Gw}}dx)$. To be more precise, it is necessary to take a 
sequence of smooth approximations 
$\gg_{\gk}$ of the function ${\rm sign}$, then let $\gk\to 0$ and 
$\gg_{\gk}\to {\rm sign}$ as in the proof of \rth {L1-th}. Therefore there exists a 
positive constant $C$ such that 
\begin {eqnarray*}
\int_{\{x:\abs {u(x)}\geq r_{0}\}}g(x,u) j_{\epsilon}(\eta_{1}){\rm sign}(u)dx
\leq C+\int_{\{x:\abs {u(x)}< r_{0}\}}g(x,u) j_{\epsilon}(\eta_{1}){\rm sign}(u)dx.
\end {eqnarray*}
Using the fact that $g(x,r){\rm sign}(r)$ is positive for $\abs r\geq r_{0}$ 
and bounded for $\abs r<r_{0}$, we can let $\epsilon \to 0$ and 
conclude, thanks to Fatou's lemma, that
\begin {eqnarray}\label {j-epsilon3}
\int_{\Gw}\abs {g(x,u)}\eta_{1}^\ga dx
<C+\sup_{n}{\norm{\gl_{n}}}_{\mathfrak 
M(\Gw;\gr^\ga_{_{\prt\Gw}}dx)}<\infty,
\end {eqnarray}
which ends the proof.\qeda 

\subsection {The Marcinkiewicz spaces approach}
At first we recall some definitions and basic properties of the Marcinkiewicz spaces.
Let $G$ be an open subset of $\mathbb R^d$ and $\gl$ a positive Borel 
measure on $G$.
\bdef {marcin} 
{\rm For $p>1$, $p'=p/(p-1)$ and $u\in L^1_{loc}(G)$, we introduce}
\begin {eqnarray}\label {marsest}
\norm u_{M^p(G;d\gl)}=\inf \left\{c\in [0,\infty]:\int_{E}\abs ud\gl
\leq c\left(\int_{E}d\gl\right)^{1/p'},\;\forall E\subset G,\, 
EÊ\mbox { \rm Borel }\right\},
\end{eqnarray}
\es
and
\begin {eqnarray}\label {marsnorm}
M^p(G;d\gl)=\{u\in L^1(G;d\gl ):\,\norm u_{M^p(G;d\gl)}<\ity\}.
\end{eqnarray}
$M^p(G;d\gl)$ is called the Marcinkiewicz space of exponent $p$, or 
weak $L^p$-space. It is a Banach space and the following estimates 
can be found in \cite {BBC} and \cite {CC}.
\bprop {marcin2} Let $1\leq q<p<\ity$ and $u\in L^1_{loc}(G;d\gl)$. 
Then
\begin {eqnarray}\label {marsest1}
C(p)\norm u_{M^p(G;d\gl)}\leq \sup\left\{s>0:\,s^p\int_{\{x:\,\abs 
{u(x)}>s\}}d\gl\right\}\leq \norm u_{M^p(G;d\gl)}.
\end{eqnarray}
Furthermore
\begin {eqnarray}\label {marsest'}
\int_{E}\abs u^qd\gl\leq C(p,q)\norm 
u_{M^p(G;d\gl)}\left(\int_{E}d\gl\right)^{1-q/p},
\end{eqnarray}
for any Borel set $E\subset G$.
\es
The key role of Marcinkiewicz spaces is to give optimal estimates when 
solving elliptic equations in a measure framework. In particular, 
using (\ref {estG1}) and (\ref {estP1}) it is not difficult to prove 
the following result (see \cite {BVV} for a more general set of 
estimates in the case of the Laplacian operator).
\bth {estth}Let $\Gw\subset \mathbb R^n$, $n\geq 2$, be a $C^{2}$ 
bounded domain and $L$ an elliptic 
operator satisfying condition (H). Let $\ga\in 
[0,1]$, $\gl\in
\mathfrak M(\Gw;\rho^\ga_{_{\prt\Gw}})$, $\gm\in \mathfrak M (\prt\Gw)$. 
If  $n+\ga>2$, there holds
\begin {eqnarray}\label {marsest2}
{\norm{\mathbb G_{L}^\Gw(\gl)}}_{M^{(n+\ga)/(n+\ga-2)}(\Gw;\gr^\ga_{_{\prt\Gw}})}
\leq C{\norm \gl}_{\mathfrak M(\Gw;\rho^\ga_{_{\prt\Gw}})},
\end{eqnarray}
and
\begin {eqnarray}\label {marsest2'}
{\norm{\nabla{\mathbb G_{L}^\Gw(\gl)}}}_{M^{(n+\ga)/(n+\ga-1)}(\Gw;\gr^\ga_{_{\prt\Gw}})}
\leq C{\norm \gl}_{\mathfrak M(\Gw;\rho^\ga_{_{\prt\Gw}})}.
\end{eqnarray}
Furthermore, for any $\gg\in 
[0,1]$,
\begin {eqnarray}\label {marsest3}
{\norm{\mathbb P_{L}^\Gw(\gm)}}
_{M^{(n+\gg)/(n-1)}(\Gw;\gr^\gg_{_{\prt\Gw}})}
\leq C{\norm \gm}_{\mathfrak M(\prt\Gw)}.
\end{eqnarray}
\es 
The following definition is inspired by B\'enilan and Brezis classical 
work \cite {BB} (with $\ga=0$) and used by Gmira and V\'eron \cite 
{GmV} (with $\ga=1$).     

\bdef {BB} {\rm A real valued function $g\in C(\Gw\times \mathbb R)$ 
satisfies the $(n,\ga)$-{\it weak-singularity assumption}, $n\geq 2$, $\ga\in [0,1]$,
$n+\ga>2$, if there exists $r_{0}\geq 0$ such that 
\begin {eqnarray}\label {passive}
  rg(x,r)\geq 0, \forevery (x,r)\in\Gw\times 
(-\ity,-r_{0}]\cup[r_{0},\ity), 
\end{eqnarray}
and a nondecreasing function $\tilde g\in 
C([0,\ity))$ such that $\tilde g\geq 0$, 
\begin {eqnarray}\label {WWA}
\int_{0}^1\tilde g(r^{2-n-\ga})r^{n+\ga-1}dr<\ity,
\end{eqnarray}
and 
\begin {eqnarray}\label {WWA1}
\abs{g(x,r)} \leq \tilde g(\abs r),\forevery (x,r)\in 
\Gw\times \mathbb R.
\end{eqnarray}
}\es 
\bth {genMeas}Let $\Omega$ be a $C^2$ bounded domain in $\mathbb R^n$, $n\geq 
2$,  $L$ the elliptic operator defined by (\ref {lin1}) and 
$g\in C(\Gw\times \mathbb R)$ a real valued function. If $L$ 
satisfies assumptions (H) and $g$ the $(n,\ga)$-weak-singularity 
assumption (then $n\geq 3$ if $\ga=0$), for any 
$\gl\in \mathfrak M(\Gw;\gr^{\ga}_{_{\prt\Gw}})$ there 
exists a solution $u$ to Problem (\ref {SLDM}).
\es 
\Proof {\it Step 1 } Construction of approximate solutions. The 
technique developed below is adapted from Brezis and Strauss 
classical article
\cite {BS}. Let 
$\gl_{n}$ be a sequence of smooth functions, with compact support 
in $\Gw$, with uniformly bounded $L^1(\Gw;\gr_{_{\prt\Gw}}dx)$-norm, with the 
property
$$\lim_{n\to\ity}\int_{\Gw}\gl_{n}\gz dx\to\int_{\Gw}\gz d\gl,
$$
for any $\gz\in C(\overline\Gw)$ such that $\sup_{\Gw} 
(\gr^{-\ga}_{_{\prt\Gw}}\abs\gz)<\ity$. 
For $k>0$, we introduce the truncation $g_{k}(.,r)$ of $g(.,r)$ by
\begin {equation}\label {trunc}
g_{k}(x,r)=\left\{\BA {l}g(x,r)\phantom {k{\rm sign} ()}\quad \mbox{if 
}\abs {g(x,r)}\leq k,\\[2mm]
k\,{\rm sign} (g(x,r))\quad \mbox{if }\abs {g(x,r)}> k.\EA\right.
\end{equation}
By Lax-Milgram's theorem, for any $z\in L^2(\Gw)$, there exists a 
unique $w=\CT_{k}(z)$ such that
\begin {equation}\label {function3}
A_{L}(w,\phi)+\int_{\Gw}g_{k}(x,z)\phi dx
=\int_{\Gw}\gl_{n}\phi dx, \forevery \phi\in W^{1,2}_{0}(\Gw).
\end{equation}
Using (\ref {lin2}),
$$
\ga{\norm{\nabla w}}^2_{L^2(\Gw)}
\leq \left(k{\abs\Gw}^{1/2}+
{\norm{\gl_{n}}}_{L^2(\Gw)}\right){\norm{w}}_{L^2(\Gw)}.
$$
The mapping $\CT_{k}$ is continuous in $L^2(\Gw)$. By the above 
estimate and Rellich-Kondrachov's theorem, $\CT_{k}$ sends 
$L^2(\Gw)$ into a relatively compact subset of $L^2(\Gw)$. By
Schauder's theorem, it admits a fixed point, say $v=v_{k}$, and $v_{k}$ 
solves
\begin {equation}\label {approx1}
Lv_{k}+g_{k}(.,v_{k})=\gl_{n}\quad \mbox {in }\Gw.
\end{equation}
The functions $v_{k}$ belongs to $C^{1,L^*}_{c}(\overline\Gw)$, since 
$\gl_{n}$ and $g_{k}$ are bounded. Multiplying by $v_{k}$ and using (\ref {WWA1}) (one notices that the two
inequalities are uniform with respect to $k$), yields to
$$
\ga{\norm{\nabla v_{k}}}^2_{L^2(\Gw)}
\leq \left (\Gth{\abs\Gw}^{1/2}+
{\norm{\gl_{n}}}_{L^2(\Gw)}\right){\norm{v_{k}}}_{L^2(\Gw)},
$$
since $rg(x,r)\geq -\Gth\abs r$, for some  $\Gth$ verifying
\begin {equation}\label {theta}0\leq \Gth\leq \sup \{\abs 
{g(x,r)}:\,x\in\Gw,\,-r_{0}\leq r\leq r_{0}\}.
\end {equation}
Hence the set of functions $\{v_{k}\}$ remains bounded in $W^{1,2}_{0}(\Gw)$
independently of $k$. \smallskip

\noindent {\it Step 2 } Uniform estimates. In 
order to prove that there exists some $k$ such that $v_{k}$ satisfies 
\begin {equation}\label {approx2}
Lv_{k}+g(.,v_{k})=\gl_{n}\quad \mbox {in }\Gw.
\end{equation}
it is sufficient to prove that $v_{k}$ is uniformly bounded in $\Gw$. 
The technique used is due to Moser \cite {Mo}. For $\gth\geq 1$, 
$\abs{v_{k}}^{\gth-1}v_{k}$ belongs to $W^{1,2}_{0}(\Gw)$. For 
simplicity we denote it by $v^\gth_{k}$, thus
\begin {equation}\label {approx3}
A_{L}(v_{k},v^\gth_{k})+\int_{\Gw}g_{k}(x,v_{k})v^\gth_{k}dx
=\int_{\Gw}\gl_{n}v^\gth_{k}dx.
\end{equation}
But, using (\ref {lin2}) and (\ref {uniq}),
\begin {eqnarray*}
A_{L}(v_{k},v^\gth_{k})&\geq& \ga\gth\int_{\Gw}{\abs {\nabla 
v_{k}}}^2v^{\gth-1}_{k}dx +\sum_{i=1}^n\int_{\Gw}(b_{i}+\gth c_{i})v_{k}^{\gth}
\frac {\prt v_{k}}{\prt x_{i}}dx+\int_{\Gw}dv^{\gth+1}_{k}dx\\
&\geq&\frac {4\ga\gth}{(\gth+1)^2}\int_{\Gw}{\abs {\nabla \left(
{\abs {v_{k}}}^{(\gth+1)/2}\right)}}^2dx
+\frac {\gth-1}{2}\sum_{i=1}^n\int_{\Gw}(c_{i}-b_{i})\frac {\prt 
v_{k}}{\prt x_{i}}v_{k}^{\gth}dx\\
&\geq&\frac {4\ga\gth}{(\gth+1)^2}\int_{\Gw}{\abs {\nabla \left({\abs 
{v_{k}}}^{(\gth+1)/2}\right)}}^2dx-\frac {\gth-1}{2(\gth+1)}\int_{\Gw}
{\abs {v_{k}}}^{\gth+1}\rm {div} \CH dx
\end{eqnarray*}
where $ \CH_{i}=c_{i}-b_{i}$ and
$$\int_{\Gw}g_{k}(x,v_{k})v^\gth_{k}dx\geq -\Gth\int_{\Gw}{\abs {v_{k}}}^{\gth}dx.
$$
By using the previous estimates and Gagliardo-Nirenberg's inequality, 
it follows that, for some $\gs>0$ and $C_{i}\geq 0$ depending on 
$\gl_{n}$ but not on $k$, there holds
\begin{eqnarray*}
\frac {\gs\gth}{(\gth+1)^2}{\norm {v_{k}}}^{\gth+1}_{L^{(\gth+1)n/(n-2)}(\Gw)}
&\leq& C_{1}{\norm {v_{k}}}^{\gth}_{L^{\gth+1}(\Gw)}
+C_{2}{\norm {v_{k}}}^{\gth+1}_{L^{\gth+1}(\Gw)}\\
&\leq&C_{3}\max \{1,{\norm {v_{k}}}^{\gth+1}_{L^{\gth+1}(\Gw)}\}.
\end{eqnarray*} 
Putting $a=n/(n-2)$, $\gg=\gth+1$, 
$${\norm {v_{k}}}_{L^{a\gg}(\Gw)}\leq C^{1/\gg}_{4}\gg^{2/\gg}
\max \{1,{\norm {v_{k}}}_{L^{\gg}(\Gw)}\}.
$$
Iterating from $\gg=2$, we obtain
\begin{eqnarray*}{\norm {v_{k}}}_{L^{a^{m+1}\gg}(\Gw)}
    &\leq &C_{5}^{\Gs_{j=0}^ma^{-j}}
2^{\Gs_{j=0}^mja^{-j}}\max \{1,{\norm {v_{k}}}_{L^{2}(\Gw)}\}\\
&\leq &C_{6}\max \{1,{\norm {v_{k}}}_{L^{2}(\Gw)}\}.
\end{eqnarray*}
Consequently $\abs {v_{k}(x)}$ is uniformly bounded by some $k_{0}$. 
Taking $k>k_{0}$, $v_{k}$ is a solution of 
\begin {equation}\label {approx4}
Lv_{k}+g(.,v_{k})=\gl_{n}\quad \mbox {in }\Gw.
\end{equation}
In order to emphasize the fact that $v_{k}$ is actually independent of $k$, 
but not on $n$, we shall denote it by $u_{n}$.\smallskip 

\noindent {\it Step 3 } Uniform integrability. It follows from Step 2 
that $g(.,u_{n})u_{n}$ is integrable in $\Gw$ and 
the same is true with $g(.,u_{n})$, because of (\ref {WWA1}). 
The space $C_{c}^{1,L}(\overline\Gw)$ is a subspace of $W^{1,2}_{0}(\Gw)$, 
therefore (\ref {function3}) implies
\begin {eqnarray}\label {SDLM2}
\int_{\Gw}\left( u_{n}L^*\gz +g(x,u_{n})\gz\right)dx=\int_{\Gw}\gl_{n}\gz dx,
\end {eqnarray}
for every $\gz\in C_{c}^{1,L}(\overline\Gw)$. By \rth {L1-th}, for any 
 $\gz\in C_{c}^{1,L}(\overline\Gw)$, $\gz\geq 0$, one has
\begin {eqnarray}\label {SDLM2b}
\int_{\Gw}\left( \abs {u_{n}}L^*\gz +{\rm 
sign(u_{n})}g(x,u_{n})\gz\right)dx\leq\int_{\Gw}\abs{\gl_{n}}\gz dx.
\end {eqnarray} 
We take $\gz=\eta_{1}$ as \rlemma{L1-lem}, and derive from (\ref 
{passive}),
\begin {eqnarray}\label {SDLM3}
\norm {u_{n}}_{L^1(\Gw)}+\norm {\gr_{_{\prt\Gw}}g(.,u_{n})}_{L^1(\Gw)}
\leq \Gth\int_{\Gw}\gr_{_{\prt\Gw}} dx+C_{1}\norm {\gr_{_{\prt\Gw}}\gl_{n}}_{L^1(\Gw)}.
\end {eqnarray}
Consequently, by using (\ref {j-epsilon3}) in \rprop {genest} and (\ref {marsest2}) in \rth {estth},
\begin {eqnarray}\label {SDLM4}
{\norm{u_{n}}}_{M^{(n+\ga)/(n+\ga-2)}(\Gw;\gr^\ga_{_{\prt\Gw}})}
\leq C_{2}\norm {\gl_{n}-g(.,u_{n})}_{\mathfrak 
M(\Gw;\rho^\ga_{_{\prt\Gw}})}
\leq C_{3}\left(\Gth+\norm {\gr_{_{\prt\Gw}}\gl_{n}}_{L^1(\Gw)}\right).
\end {eqnarray}
By the local  regularity result  of \rcor {L^1regloc}, there exist a 
subsequence $\{u_{n_{k}}\}$ and a function $u\in W^{1,q}_{loc}(\Gw)$, 
for any $1\leq q<n/(n-1)$, such that $u_{n_{k}}\to u$ a.e. in $\Gw$ and 
weakly in  $W^{1,q}_{loc}(\Gw)$. Notice that $W^{1,q}_{loc}(\Gw)$ can 
be replaced by $W^{1,q}_{0}(\Gw)$ if $\ga=0$, by \rth {L^1reg}. Combining (\ref {SDLM3}) 
and  estimate (\ref {M1-7}) 
with $\mu_{n}=0$ and $\gl_{n}$ replaced by $\gl_{n}-g(.,u_{n})$, 
one obtains that, for any Borel subset $\gw\subset\Gw$, there holds
$$\int_{\gw}\abs {u_{n}}dx\leq 
\left(C'\abs \Gw+C'_{1}\norm 
{\gr_{_{\prt\Gw}}\gl_{n}}_{L^1(\Gw)}\right)\abs {\gw}^{1/p},
$$
if $p>n$. Thus, by the Vitali Theorem, it can also be assumed that 
$u_{n_{k}}\to u$ in $L^{1}(\Gw)$. Furthermore, for any $R\geq 0$,
\begin {eqnarray*}\int_{\gw}\abs {g(.,u_{n)}}\rho^{\ga}_{_{\prt\Gw}}dx
&\leq& \int_{\gw}\tilde g(\abs{u_{n}})\rho^\ga_{_{\prt\Gw}}dx\\
&\leq& \int_{\gw\cap\{\abs{u_{n}}\leq R\}}\tilde g(\abs{u_{n}})\rho^{\ga}_{_{\prt\Gw}}dx
+\int_{\gw\cap\{\abs{u_{n}}> R\}}
\tilde g(\abs{u_{n}})\rho^{\ga}_{_{\prt\Gw}}dx\\
&\leq& \tilde g(R)\int_{\gw}\rho^{\ga}_{_{\prt\Gw}}dx
-\int_{R}^\ity\tilde g(s)d\gth_{n}(s),
\end {eqnarray*}
where
\begin {eqnarray*}
\gth_{n}(s)=\int_{\{x\in\Gw:\abs {u_{n}}> s\}}\rho^{\ga}_{_{\prt\Gw}}(x)dx
&\leq& s^{-(n+\ga)/(n+\ga-2)}{\norm{u_{n}}}_
{M^{(n+\ga)/(n+\ga-2)}(\Gw;\gr^\ga_{_{\prt\Gw}})}\\
&\leq& Cs^{-(n+\ga)/(n+\ga-2)},
\end {eqnarray*}
by (\ref {marsest1}). Moreover
\begin {eqnarray*}-\int_{R}^\ity\tilde g(s)d\gth_{n}(s)&=&\tilde g(R)\gth_{n}(R)+
\int_{R}^\ity\gth_{n}(s)d\tilde g(s)\\
&\leq& \tilde g(R)\gth_{n}(R)+
C\int_{R}^\ity s^{-(n+\ga)/(n+\ga-2)}d\tilde g(s)\\
&\leq&\tilde g(R)\gth_{n}(R)-C\tilde g(R)R^{-(n+\ga)/(n+\ga-2)} \\
&\,&\qquad\qquad\qquad\qquad+\frac {C(n+\ga)}{n+\ga-2}\int_{R}^\ity\tilde 
g(s)s^{-2(n+\ga-1)/(n+\ga-2)}ds\\
&\leq&\frac {C(n+\ga)}{n+\ga-2}\int_{R}^\ity\tilde 
g(s)s^{-2(n+\ga-1)/(n+\ga-2)}ds.
\end {eqnarray*}
Since condition (\ref {marsest2}) is equivalent to
\begin {equation}\label {marsest2b}
\int_{1}^\ity\tilde 
g(s)s^{-2(n+\ga-1)/(n+\ga-2)}ds<\ity,
\end {equation}
given $\epsilon>0$, we first choose $R>0$ such that 
$$\frac {C(n+\ga)}{n+\ga-2}\int_{R}^\ity\tilde 
g(s)s^{-2(n+\ga-1)/(n+\ga-2)}ds\leq \epsilon/2.
$$
Then we put $\gd=\epsilon/(2(1+\tilde g(R))$ and derive
$$\int_{\gw}\rho_{_{\prt\Gw}}^\ga dx\leq \gd\Longrightarrow 
\int_{\gw}\abs {g(u_{n})}\rho_{_{\prt\Gw}}^\ga dx\leq \epsilon.
$$
Therefore $\{\rho_{_{\prt\Gw}}^\ga g(.,u_{n})\}$ 
is uniformly integrable, and we can assume that the previous 
sequence $\{n_{k}\}$ is such that
\begin {equation}\label {lim}
\lim_{n_{k}\to\ity}
\int_{\Gw} \abs{g_{n_{k}}(.,u_{n_{k}})-g(.,u)}\rho_{_{\prt\Gw}}^\ga dx=0
\Longrightarrow \int_{\Gw} 
\abs{g_{n_{k}}(.,u_{n_{k}})-g(.,u)}\rho_{_{\prt\Gw}}dx=0,
\end {equation}
since $\ga\in [0,1]$. Letting $n_{k}\to\ity$ in (\ref {SDLM2}), one obtains
\begin {eqnarray}\label {SDLM3b}
\int_{\Gw}\left( uL^*\gz +g(x,u)\gz\right)dx=\int_{\Gw}\gz d\gl.
\end {eqnarray}\qeda 
\medskip

Since the uniform integrability conditions depends only on the total 
variation norm of the measure $\rho_{_{\prt\Gw}}^\ga\gl$, the following stability result 
holds.
\bcor {genMeascor} Let $\Gw$ and $\ga$ be as in \rth {genMeas}, $g$ 
satisfy the $(n,\ga)$-weak-singularity assumption and $r\mapsto g(x,r)$
is nondecreasing, for any $x\in\Gw$. Then the solution $u$ is unique. 
If we assume that $\{\gl_{m}\}$ is a sequence of measures in  
$\mathfrak M(\Gw;\rho_{_{\prt\Gw}}^\ga)$ such that 
$$\lim_{m\to\infty}\int_{\Gw}\gz d\gl_{m}
=\lim_{m\to\infty}\int_{\Gw}\gz d\gl,
$$
for any $\gz\in C(\overline\Gw)$ satisfying 
$\sup_{\Gw}\gr_{_{\prt\Gw}}^{-\ga}\abs \gz<\infty$, then the 
corresponding solutions $u_{m}$ of problem 
\begin {equation}\label {SLDMm}\left.\BA{ll}
Lu_{m}+g(x,u_{m})=\gl_{m}&\mbox { in }\Omega,\\[2mm]
\phantom{Lu_{m}+g(x,)}u_{m}=0&\mbox { on }\prt\Omega, 
\EA\right.
\end {equation}
converge in $L^{1}(\Gw)$ to the solution $u$ of (\ref {SLDM}), when 
$m\to\ity$.
\es 
\Remark If $g(x,r)=\abs r^{q-1}r$,  the $(n,\ga)$-weak-singularity 
assumption is satisfied if and only if
\begin {eqnarray}\label {power q}
0<q<\frac {n+\ga}{n+\ga-2}.
\end{eqnarray}
\subsection {Admissible measures and the $\Delta_{2}$-condition}

\bdef {admdef}{\rm Let $\tilde g$ be a continuous real valued 
nondecreasing function defined in $\mathbb R_{+}$, $\tilde g\geq 0$. A measure $\gl$ 
in $\Gw$ is said $(\tilde g,k)$-admissible if
\begin {eqnarray}\label {g-adm}
\int_{\Gw}\tilde g(\mathbb G_{L}^\Gw(\abs\gl)+k)\rho_{_{\prt\Gw}}dx<\infty,
\end{eqnarray}
where $G_{L}^\Gw(\abs\gl)$ is the Green potential of $\gl$ and $k\geq 0$.}
\es
\bth {g-adm-th} Let $\Gw$ be a $C^2$ bounded domain in $\mathbb 
R^n$, $n\geq 2$, $L$ an elliptic operator defined by (\ref 
{lin1}), and $g\in C(\Gw\times\mathbb R)$. We assume that $L$ satisfies 
the condition (H), and $g$ (\ref {passive}) for some $r_{0} \geq 0$ and (\ref 
{WWA1}) for some function $\tilde g$ as in 
\rdef {admdef}. Then for any $(\tilde g,r_{0})$-admissible Radon measure 
$\gl\in\mathfrak M(\Gw;\rho_{_{\prt\Gw}})$, Problem (\ref {SLDM}) admits a 
solution.
\es
\Proof For $k>0$, we take the same truncation $g_{k}(.,r)$ of $g(.,r)$ 
defined by (\ref{trunc}). Since $g_{k}$ satisfies (\ref {WWA}) and (\ref {WWA1}), we 
denote by $u_{k}$ a solution of 
\begin {equation}\label {approx1b}\left.\BA{ll}
Lu_{k}+g_{k}(x,u_{k})=\gl&\mbox { in }\Omega,\\[2mm]
\phantom{Lu_{m}+g_{k}(x,)}u_{k}=0&\mbox { on }\prt\Omega, 
\EA\right.
\end {equation}
which exists by \rth {genMeas}. As in the proof of \rth {genMeas} the 
following estimates hold,
\begin {eqnarray}\label {SDLM3-bis}
\norm {u_{k}}_{L^1(\Gw)}+\norm {\gr_{_{\prt\Gw}}g_{k}(.,u_{k})}_{L^1(\Gw)}
\leq \Gth\int_{\Gw}\gr_{_{\prt\Gw}}dx+C_{1}
\norm {\gr_{_{\prt\Gw}}\gl}_{\mathfrak M(\Gw)},
\end {eqnarray}
where $\Gth$ is defined by (\ref {theta}), and
\begin {eqnarray}\label {SDLM4-bis}
{\norm{u_{k}}}_{M^{(n+1)/(n-1)}(\Gw;\gr_{_{\prt\Gw}})}
\leq C_{3}\left(\Gth+\norm {\gr_{_{\prt\Gw}}\gl}_{L^1(\Gw)}\right).
\end {eqnarray}
By \rcor {L^1regloc}, there exist a 
subsequence $\{u_{k_{j}}\}$ and a function $u\in W^{1,q}_{loc}(\Gw)$, 
for any $1\leq q<n/(n-1)$, such that $u_{k_{j}}\to u$ a.e. in $\Gw$ and 
weakly in  $W^{1,q}_{loc}(\Gw)$. Moreover $g_{k_{j}}(.,u_{k_{j}})\to 
g(.,u)$ almost everywhere in $\Gw$. Put 
$$w_{\gl_{+}}=\mathbb G_{L}^\Gw(\gl_{+})+r_{0}.
$$
Then
$$L(u_{k}-w_{\gl_{+}})+g_{k}(x,u_{k})=\gl-\gl_{+},
$$
and, for any $\gz\in C^{1,L}_{c}(\overline\Gw)$, $\gz\geq 0$,
\begin {equation}\label {SDLM5-bis}\BA {l}
\myint{\Gw}{}(u_{k}-w_{\gl_{+}})_{+}L^*\gz dx+
\myint{\Gw}{}g_{k}(x,u_{k})\rm {sign}_{+}(u_{k}-w_{\gl_{+}})\gz dx \\[2mm]
\phantom {\myint{\Gw}{}(u_{k}-w_{\gl_{+}})_{+}L^*\gz dx+
\myint{\Gw}{}g_{k}}
\leq-\myint{\prt\Gw}{}\myfrac 
{\prt\gz}{\prt{\bf n}_{L^*}}(u_{k}-w_{\gl_{+}})_{+}dS, 
\EA\end {equation}
by Inequality (\ref {M1-7'}). Since the boundary term in (\ref {SDLM5-bis}) vanishes, 
and $w_{\gl_{+}}\geq r_{0}$, there holds  
$g_{k}(x,u_{k})\rm {sign}_{+}(u_{k}-w_{\gl_{+}})\geq 0$, which implies
$$\int_{\Gw}(u_{k}-w_{\gl_{+}})_{+}L^*\gz dx\leq 0.
$$
Taking $\gz=\eta_{1}$ defined by (\ref {L1-1e}) (with $u=1$, hence $L^*\eta_{1}=1$), 
yields to $(u_{k}-w_{\gl_{+}})_{+}=0$ a.e. in $\Gw$. Thus
$$u_{k}\leq w_{\gl_{+}}=\mathbb G_{L}^\Gw(\gl_{+})+r_{0}.
$$
In the same way 
$$u_{k}\geq -\mathbb G_{L}^\Gw(\gl_{-})-r_{0}.
$$
Therefore 
\begin {eqnarray}\label {SDLM6}
\abs {u_{k}}\leq \mathbb G_{L}^\Gw(\abs \gl)+r_{0}
\Longrightarrow \abs {g_{k}(u_{k})}\leq \tilde g(\abs {u_{k}})\leq
\tilde g (\mathbb G_{L}^\Gw(\abs \gl)+r_{0}).
\end {eqnarray}
Because the right-hand side of (\ref {SDLM6}) belongs to 
$L^1(\Gw;\rho_{_{\prt\Gw}}dx)$, the sequence $\{g_{k}(.,u_{k})\}$ is 
uniformly integrable for the measure $\rho_{_{\prt\Gw}}dx$. As in the 
proof of \rth {genMeas}, we conclude by the Vitali Theorem that $u$ is a 
solution of (\ref {SLDM}).\qeda 
\medskip

The condition of ($g,r_{0}$)-admissibility on $\gl$ is too 
restrictive if the function $g$ has a strong power growth, in 
particular it leads to exclude some $\gl$ which are regular with respect the $n$-dimensional 
Hausdorff measure, even if we know, from the Brezis and Strauss Theorem 
(see \rth {genMeas}-Step 1), that Problem (\ref 
{SLDM}) is solvable for such measures. A natural extension is to impose only the 
($g,r_{0}$)-admissibility on the singular part $\gl_{s}$ of the 
measure. However, a generic power-like growth condition called the 
$\Delta_{2}$-condition is needed.

\bdef {delta2def}{\rm A real valued function $g\in C(\Gw\times \mathbb R)$ 
satisfies a uniform $\Delta_{2}$-condition if there exist two constants $\ell 
\geq 0$, $\gth>1$ such that}
\begin {eqnarray}\label {delta-2}
\abs {g(x,r+r')}\leq \gth (\abs {g(x,r)}+\abs {g(x,r')})+\ell,\forevery 
x\in\Gw,\,\forall (r,r')\in\mathbb R\times\mathbb R.
\end{eqnarray}
\es
\bth {g-adm-th2} Let $\Gw$ and $L$ be as in \rth {g-adm-th}. Assume
$g\in C(\Gw\times\mathbb R)$ satisfies the $\Delta_{2}$-condition, 
$r\mapsto g(x,r)$ is nondecreasing for any $x\in\Gw$ and
(\ref {WWA1}) holds for some function $\tilde g$ as in 
\rdef {admdef}. For any Radon measure 
$\gl\in\mathfrak M(\Gw;\rho_{_{\prt\Gw}})$, 
with $\gl=\tilde\gl+\gl^*$, where $\tilde\gl\in 
L^1(\Gw;\rho_{_{\prt\Gw}}dx)$, and $\gl^*$ is 
$(\tilde g,0)$-admissible and singular with respect to the 
$n$-dimensional Lebesgue measure, Problem (\ref {SLDM}) 
admits a unique solution.
\es
\Proof Uniqueness comes from the monotonicity of $r\mapsto g(x,r)$. \smallskip

\noindent {\it Step 1} If we write $g(x,r)=g(x,r)-g(x,0)+g(x,0)= \hat g(x,r)+g(x,0)$, then 
the equation is transformed into 
$$Lu+\hat g(x,u)=\gl -g(x,0)=\hat\gl ,
$$
where $r\mapsto \hat g(x,r)$ nondecreasing and $\hat g(x,0)=0$. 
Notice that $\abs {\hat g(x,0)}\leq \tilde g(0)$ by (\ref {WWA1}), and 
that $\gl^*$ is singular 
with respect to $\tilde\gl -g(x,0)$. Finally the new function $\hat g$ 
satisfies Inequality (\ref {delta-2}) with the same $\gth$ 
and $\ell$ replaced by $\hat\ell=\ell+(2\gth +1)\abs {\tilde g (0)}$, 
and (\ref {WWA1}) with $\tilde g$ replaced by $\tilde g +\abs {\tilde g 
(0)}$. From now we shall suppose that the function $g$ satisfies $g(x,0)=0$ for any 
$x\in\Gw$.
We introduce the truncation $g_{k}(.,r)$ by (\ref {trunc}). The truncated 
function $g_{k}$ satisfies also (\ref {delta-2}) (with $\gth$ replaced by 
$1+\gth$).\smallskip

\noindent {\it Step 2} We suppose that $\gl$ is nonnegative. Then 
$\tilde\gl$ and $\gl^*$ inherit the same property. Let $\{\tilde\gl_{i}\}$ 
be a sequence of smooth nonnegative functions with compact support in 
$\Gw$, converging to $\tilde\gl$ in the weak sense of 
$L^1(\Gw;\rho_{_{\prt\Gw}})$. Let $u_{i,k}$ be the solution of 
\begin {equation}\label {SLDM-i}\left.\BA{ll}
Lu_{i,k}+g_{k}(x,u_{i,k})=\tilde\gl_{i}+\gl^*&\mbox { in }\Omega,\\[2mm]
\phantom {Lu_{i,k}+g_{k}(x,)}
u_{i,k}=0&\mbox { on }\prt\Omega, 
\EA\right.
\end {equation}
and $v_{i,k}$ the one of 
\begin {equation}\label {SLDM-i2}\left.\BA{ll}
Lv_{i,k}+g_{k}(x,v_{i,k})=\tilde\gl_{i}&\mbox { in }\Omega,\\[2mm]
\phantom {Lv_{i,k}+g_{k}(x,)}
v_{i,k}=0&\mbox { on }\prt\Omega. 
\EA\right.
\end {equation}
Both solutions exist by \rth {g-adm-th}. By the maximum principle
\begin {eqnarray}\label {SLDM-k3}
0\leq u_{i,k}\leq v_{i,k}+\mathbb G_{L}^\Gw(\gl^*),
\end{eqnarray}
and by the monotonicity of $g_{k}$ and (\ref {delta-2}),
\begin {eqnarray}\label {SLDM-i3}
0\leq g_{k}(.,u_{i,k})\leq \gth \left(g_{k}(.,v_{i,k})+g_{k}(.,\mathbb 
G_{L}^\Gw(\gl^*))\right)+\ell 
\leq \gth \left(g_{k}(.,v_{i,k})+\tilde g(\mathbb 
G_{L}^\Gw(\gl^*))\right)+\ell.
\end{eqnarray}
By \rth {g-adm-th}), if $i$ is fixed and $k\to\ity$, the sequence $\{v_{i,k}\}$ converges weakly in 
$W^{1,q}_{loc}(\Gw)$ and a.e. in $\Gw$ to the solution $v_{i}$ of 
\begin {equation}\label {SLDM-i4}\left.\BA{ll}
Lv_{i}+g(x,v_{i})=\tilde\gl_{i}&\mbox { in }\Omega,\\[2mm]
\phantom {Lv_{i}+g(x,)}
v_{i}=0&\mbox { on }\prt\Omega. 
\EA\right.
\end {equation}
Since the $v_{i,k}$ are uniformly bounded with respect to $k$, the 
same property holds with the $g_{k}(v_{i,k})$, hence their convergence 
to $v_{i}$ and $g(.,v_{i})$ are uniform in $\overline\Gw$. Because of 
(\ref {SLDM-i3}) and the elliptic equations regularity theory, the 
sequence $\{u_{i,k}\}_{k\in \mathbb N_{*}}$ is relatively compact in 
the $W^{1,q}_{loc}(\Gw)$-topology. Thus there exist a subsequence
$\{u_{i,k_{j}}\}$ and a function $u_i$ such that $u_{i,k_{j}}\to 
u_{i}$ as $k_{j}\to\ity$ in this topology and a.e. in $\Gw$. 
By continuity, $g_{k_{j}}(.,u_{i,k_{j}})\to g(.,u_{i})$ a.e. in $\Gw$. 
Because of (\ref {SLDM-i3}) and the $(\tilde g,0)$-admissibility 
condition on $\gl^*$, Lebesgue's theorem implies that 
$$
\lim_{k_{j}\to\ity}g_{k_{j}}(.,u_{i,k_{j}})= g(.,u_{i})\quad\mbox {in }\, 
L^1(\Gw;\gr_{_{\prt\Gw}}dx).
$$
It follows from inequality (\ref {SLDM-k3}) that $u_{i,k_{j}}\to u_{i}$ 
in $L^1(\Gw)$ (we recall that $\mathbb G_{L}^\Gw(\gl^*)\in L^1(\Gw)$). 
Letting $k_{j}\to\ity$ in (\ref {SLDM-i}) we see that $u_{i}$ is the 
solution of
\begin {equation}\label {SLDM-i5}\left.\BA{ll}
Lu_{i}+g(x,u_{i})=\tilde\gl_{i}+\gl^*&\mbox { in }\Omega,\\[2mm]
\phantom {Lu_{i}+g(x,)}
u_{i}=0&\mbox { on }\prt\Omega. 
\EA\right.
\end {equation}
By uniqueness of $u_{i}$ the whole sequence $u_{i,k}$ converges to 
$u_{i}$ as $k\to\ity$. Moreover
\begin {eqnarray}\label {SLDM-6}\left.\BA {l}
(i)\, 0\leq u_{i}\leq v_{i}+\mathbb G_{L}^\Gw(\gl^*),\\[2mm]
(ii)\, 0\leq g(.,u_{i})\leq \gth \left(g(.,v_{i})+g(\mathbb 
G_{L}^\Gw(\gl^*))\right)+\ell 
\leq \gth \left(g(.,v_{i})+\tilde g(\mathbb G_{L}^\Gw(\gl^*))\right)+\ell.\EA\right.
\end{eqnarray}
By \rth {L1-th} with $\gz=G_{L}^{\Gw}(1)$,
\begin {eqnarray}\label {SLDM-7}
\norm {v_{i}-v_{j}}_{L^1(\Gw)}+
\norm {g(.,v_{i})-g(.,v_{j})}_{L^1(\Gw;\gr_{_{\prt\Gw}}dx)}\leq
C \norm {\tilde \gl_{i}-\tilde \gl_{j}}_{L^1(\Gw)}.
\end{eqnarray}
Therefore $v_{i}\to v$ in $L^1(\Gw)$ and $g(.,v_{i})\to g(.,v)$ in 
$L^1(\Gw;\gr_{_{\prt\Gw}}dx)$ where $v$ is the solution of 
\begin {equation}\label {SLDM-8}\left.\BA{ll}
Lv+g(x,v)=\tilde\gl&\mbox { in }\Omega,\\[2mm]
\phantom {Lv+g(x,)}
v=0&\mbox { on }\prt\Omega. 
\EA\right.
\end {equation}
By (\ref {SLDM-6})-(i) there exists a subsequence $\{u_{i_{j}}\}$ which 
converges in $L^1(\Gw)$ and a.e. in $\Gw$ to some function $u$. 
Because of (\ref {SLDM-6})-(ii), the admissibility condition on 
$\gl^*$ and the Vitali Theorem, the sequence $\{g(.,u_{i_{j}})\}$ converges to 
$g(.,u)$ in $L^1(\Gw;\gr_{_{\prt\Gw}}dx)$. Thus $u$ is the solution 
of (\ref {SLDM}).\smallskip

\noindent {\it Step 3 } In the general case we construct the solution 
$u_{i,k}$ of (\ref {SLDM-i}) and the
functions $U=\overline u_{i,k}$ and $U=\underline u_{i,k}$ solutions of
 \begin {equation}\label {SLDM-9}\left.\BA{ll}
LU+g_{k}(x,U)=\tilde\Gl&\mbox { in }\Omega,\\[2mm]
\phantom {LU+g_{k}(x,)}
U=0&\mbox { on }\prt\Omega. 
\EA\right.
\end {equation}
where $\Gl=\abs {\tilde\gl_{i}}+\abs {\gl^*}$ in the case of $\overline u_{i,k}$
and $\Gl=-\abs {\tilde\gl_{i}}-\abs {\gl^*}$ in the case of $\underline 
u_{i,k}$. We also construct the solutions 
$V=\overline v_{i,k}$ and $V=\underline v_{i,k}$ of the same equation 
with $\Gl=\abs {\tilde\gl_{i}}$ in the case of $\overline v_{i,k}$ and
 $\Gl=-\abs {\tilde\gl_{i}}$ in the case of $\underline 
v_{i,k}$. Since 
\begin {eqnarray}\label {SLDM-10}
\underline v_{i,k}-\mathbb G_{L}^\Gw(\abs{\gl^*})\leq u_{i,k}\leq \overline v_{i,k}+\mathbb G_{L}^\Gw(\abs{\gl^*}),
\end{eqnarray}
and 
\begin {eqnarray}\label {SLDM-11}
\gth \left(g_{k}(.,\underline v_{i,k})+g(.,\mathbb 
G_{L}^\Gw(-\abs{\gl^*}))\right)-\ell
\leq g_{k}(.,u_{i,k})\leq \gth \left(g_{k}(.,\overline v_{i,k})+g(.,\mathbb 
G_{L}^\Gw(\abs{\gl^*}))\right)+\ell, 
\end{eqnarray}
we conclude by using the Vitali Theorem and the convergence arguments 
of Step 2.\qeda 

\subsection {The duality method}

Let $\Gw$ be a domain in $\mathbb R^n$ and $L$ is an elliptic operator 
in $\Gw$. In this section we study the sharp solvability of 
Problem (\ref {SLDM}) when $g(x,r)=\abs u^{q-1}u$
with  $q>0$. For this type of nonlinearity, the 
$(n,0)$-weak-singularity assumption is satisfied if and only if 
$0<q<n/(n-2)$. Thus we shall concentrate on the case $n\geq 3$ and $q\geq n/(n-2)$
and for such a task the theory of Bessel capacities is needed. 
\medskip
\subsubsection {Bessel capacities}
Let $p>1$ be a real number and $p'=p/(p-1)$. If $m$ in an integer we endow the Sobolev 
space $W^{m,p}(\mathbb R^n)$ with the usual norm
$${\norm \phi}_{W^{m,p}(\mathbb R^n)}=\left(\sum_{\abs\gg\leq 
m}\int_{\Gw}{\abs{D^\gg\phi}}^pdx\right)^{1/p},
$$
and we introduce the associated capacity $C_{m,p}$ by
$$C_{m,p}(K)=\inf\left\{{\norm \phi}^p_{W^{m,p}(\mathbb R^n)}:\phi\in 
C^\infty_{c}(\mathbb R^n),\,\phi\geq 1\mbox { in a neighborhood of } 
K\,\right\},
$$
if $K$ is compact,
$$C_{m,p}(G)=\sup\left\{C_{m,p}(K):K\subset G,\,K\mbox { compact }\right\},
$$
if $G$ is open, and
$$C_{m,p}(E)=\inf\left\{C_{m,p}(G):E\subset G,\,G\mbox { open }\right\},
$$
for an arbitrary set $E$. The scale of Sobolev spaces is not accurate 
enough to describe the subets of $\mathbb R^n$ by means of their capacities. If $\ga$ is a real number, 
we introduce the Bessel kernel of order 
$\ga$ by 
\begin {eqnarray}\label {Bessel}
G_{\ga}=\mathcal F^{-1}\left((1+\abs\xi^2)^{-\ga/2}\right)
\end {eqnarray}
were $\mathcal F^{-1}$ is the inverse Fourier transform on the 
Schwartz space $\mathcal S'(\mathbb R^n)$. If 
$$\mathcal G_{\ga}=(I-\Gd)^{-\ga/2},
$$
there holds the Bessel potential representation 
\begin {eqnarray}\label {Bessel2}
f=\mathcal G_{\ga}g=G_{\ga}\ast g\Longleftrightarrow 
g=\mathcal G_{-\ga}g=G_{-\ga}\ast f \forevery f,\,g\in \mathcal 
S(\mathbb R^n).
\end {eqnarray}
\bdef{Bessel sp} {\rm Let $\ga$ and $p>1$ be two real numbers. The Bessel potential 
space of order $\ga$ and power $p$ is
$$L^{\ga,p}(\mathbb R^n)=\left\{f: f=G_{\ga}\ast g,\,g\in L^{p}(\mathbb 
R^n)\right\},
$$
with norm}
$${\norm f}_{L^{\ga,p}(\mathbb R^n)}=
{\norm g}_{L^{p}(\mathbb R^n)}=
{\norm {G_{-\ga}\ast f}}_{L^{p}(\mathbb R^n)}.
$$
\es 
As usual, $L_{0}^{\ga,p'}(\mathbb R^n)$ denotes the closure of 
$C^\ity_{c}(\mathbb R^n)$ in $L^{\ga,p'}(\mathbb R^n)$.
Thanks to a result due to Calderon, the 
functions in $W^{m,p}(\mathbb R^n)$ can be represented by mean of 
Bessel potentials. Actually for any $\ga\in\mathbb N_{*}$ and 
$1<p<\infty$, $W^{\ga,p}(\mathbb R^n)=L^{\ga,p}(\mathbb R^n)$ and their 
exists a positive constant $A$ such that 
\begin {eqnarray}\label {Bessel3}
A^{-1}{\norm f}_{L^{\ga,p}(\mathbb R^n)}\leq 
{\norm f}_{W^{\ga,p}(\mathbb R^n)}\leq A{\norm f}_{L^{\ga,p}(\mathbb R^n)}
,\forevery f\in W^{\ga,p}(\mathbb R^n).
\end {eqnarray}
\smallskip
By generalization (see \cite {Ch} for a general construction of capacities), the Bessel capacity of order $(\ga,p)$ ($\ga>0$, $p>1$) 
of a compact set $K$ is defined by
\begin {eqnarray}\label {Bessel4}
C_{\ga,p}(K)=\inf\left\{{\norm \phi}^p_{L^{\ga,p}(\mathbb R^n)}:\phi\in 
\mathcal S(\mathbb R^n),\,\phi\geq 1\mbox { in a neighborhood of } 
K\,\right\},
\end {eqnarray}
with the same extension to open sets and arbitrary sets as for 
Sobolev capacities. A dual definition involving measures is the 
following \cite {AH} :
\begin {eqnarray}\label {Bessel4'}
C_{\ga,p}(K)=\sup\left\{\left(\frac 
{\gm(K)}{{\norm{G_{\ga}\ast\gm}}_{L^{p'}(\mathbb R^n)}}\right)^p:\gm\in\mathfrak M_{+}(K)
\right\},
\end {eqnarray}
where $\mathfrak M_{+}(K)$ is the set of positive Radon measures with support in $K$.
An important result due to Maz'ya (see \cite {AH}) 
states that the following expression
\begin {eqnarray}\label {Bessel5}
\tilde C_{\ga,p}(K)=\inf\left\{{\norm \phi}^p_{L^{\ga,p}(\mathbb R^n)}:\phi\in 
\mathcal S(\mathbb R^n),\,\phi\equiv 1\mbox { in a neighborhood of } 
K\,\right\},
\end {eqnarray}
defines a new capacity which is equivalent to the $C_{\ga,p}$-capacity in 
the sense that there exists a positive constant $B$ such that 
\begin {eqnarray*}
B^{-1} C_{\ga,p}(K)\leq \tilde C_{\ga,p}(K)\leq B C_{\ga,p}(K),
\end {eqnarray*}
for any compact subset $K$. In the particular case of sets with 
zero capacity, the following useful result holds.
\bprop {zero cap} Let $\ga>0$, $1<p<\infty$, $K$ be a compact subset 
of $\mathbb R^n$ and $\CO$ an open subset containing $K$. If 
$C_{\ga,p}(K)=0$, there exists a sequence $\{\phi_{n}\}\subset 
C_{c}^{\infty}(\CO)$ such that $0\leq\phi_{n}\leq 1$, $\phi_{n}\equiv 
1$ in a neighborhood of $K$ and $\phi_{n}\to 0$ in $L^{\ga,p}(\mathbb 
R^n)$ as $n\to\infty$.
\es

By using smooth cut-off function with value in $[0,1]$, support in a 
neighborhood of $K$ and taking the value $1$ in a smaller neighborhood 
of $K$, the proof of this result is straightforward if $\ga$ is an 
integer, and more delicate if not (see \cite [Th. 3.3.3]{AH}).

\bdef {quasi a.e.} {\rm Let $\ga>0$ and $1<p<\infty$. \smallskip 

\noindent (i) A property is 
said to hold $C_{\ga,p}$-quasi everywhere if it holds everywhere but 
on a set of $C_{\ga,p}$-capacity zero.\smallskip 

\noindent (ii) A function $\phi$ defined in $\mathbb R^n$ 
is said to be $C_{\ga,p}$-quasicontinuous  if for any $\epsilon>0$, 
there is an open set $G\subset\mathbb R^n$ with $C_{\ga,p}(G)<\epsilon$
and $f\in C(G^c)$.
\smallskip 

\noindent (iii) Let $\CO$ be an open subset of $\mathbb R^n$ and 
$\gl\in\mathfrak M(\CO)$ . The measure $\gl$ is said not to charge 
 subsets of $\CO$ with $C_{\ga,p}$-capacity zero if
 $$\forall E\subset \CO, C_{\ga,p}(E)=0\Longrightarrow 
 \int_{E}d\abs\gl=0,$$
 where, $d\abs\gl$ denote in the same way the unique complete regular 
 Borel measure generated by the Radon measure $\abs\gl$.
}\es

It is proven in \cite {AH} that for any $\ga>0$, $1<p<\infty$ and 
$g\in L^p(\Gw)$, the function $G_{\ga}\ast g$ is 
$C_{\ga,p}$-quasicontinuous. Therefore, any element $\phi\in 
L^{\ga,p}(\mathbb R^n)$ admits a (unique) quasicontinuous 
representative, $\tilde \phi$. Furthermore, from any converging 
sequence $\{\phi_{n}\}\subset L^{\ga,p}(\mathbb R^n)$ it can be 
extracted a subsequence $\{\phi_{n}\}$ which converges 
$C_{\ga,p}$-quasi everywhere. The link between the measures which do not 
charge capacitary sets and elements of negative Bessel spaces is 
enlighted by three results. The first one is due essentially to 
Grun-Rehomme \cite {GhR} (see also \cite {AH}).
\bprop {BP1lem1} Let $\ga>0$ and $1<p<\infty$. If 
$\gl\in\mathfrak M(\Gw)\cap L^{-\ga,p}(\Gw)$, 
then $\gl$ does not charge sets with $C_{\ga,p'}$-capacity zero.
\es
\Proof By the Jordan 
decomposition Theorem of a measure, there exist two disjoint Borel 
subsets $A$ and $B$ such that 
$$A\cup B=\Gw,\quad \gl_{+}(B)=0,\quad \gl_{-}(A)=0.
$$ 
Let $E\subset \mathbb R^n$ with $C_{\ga,p'}(E)=0$. With no loss 
of generality $E$ can be assumed as being a Borel set. 
It is therefore sufficient that $\gl_{+}(A\cap E)=\gl_{-}(B\cap 
E)=0$. Because
$$\gl_{+}(A\cap E)=\sup\{\gl_{+}(K):K\mbox { compact },k\subset A\cup 
E\},
$$
it is sufficient to prove that for any compact subset $K\subset A\cap 
E$, $\gl_{+}(K)=0$. Let $\epsilon>0$, since $\gl_{-}(K)=0$, there 
exists an open subset $\gw$ of $\Gw$ containing $K$ such that 
$\gl_{-}(\gw)\leq \epsilon$. Let $\eta\in C^\infty_{c}(\gw)$, with 
value in $[0,1]$ and equal to $1$ on $K$. By \rprop{zero cap}, since $C_{\ga,p'}(K)=0$, 
there exists a sequence $\{\phi_{n}\}\subset C_{c}^\infty(\Gw)$, of 
functions with 
value in $[0,1]$, equal to $1$ in a neighborhood of $K$ and such that 
$\phi_{n}\to 0$ in $L^{\ga,p'}(\Gw)$ as $n\to\infty$. Then
\begin {eqnarray*}
\int_{K}d\gl_{+}\leq\int_{K}\phi_{n}\eta d\gl_{+}\leq \int_{\gg}\phi_{n}\eta d\gl_{+}
=\int_{\Gw}\phi_{n}\eta d\gl+\int_{\gw}\phi_{n}\eta d\gl_{-}.
\end {eqnarray*}
But
$$\int_{\gw}\phi_{n}\eta d\gl_{-}\leq \int_{\gw}d\gl_{-}\leq\epsilon,
$$
and
\begin {eqnarray*}\int_{\Gw}\phi_{n}\eta d\gl\leq \int_{\Gw}\phi_{n}d\gl
=\langle\gl,\phi_{n}\rangle_{[L^{-\ga,p},L^{\ga,p'}]}
\leq {\norm \gl}_{L^{-\ga,p}}{\norm \phi_{n}}_{L^{\ga,p'}},
\end {eqnarray*}
which goes to zero as $n\to\infty$. Therefore
$$\int_{K}d\gl_{+}\leq \epsilon.
$$
Since $\epsilon$ is arbitrary, $\gl_{+}(K)=0$. In the same way 
$\gl_{-}(B\cap E)=0$. Therefore $\abs\gl (E)=0$.\qeda 
\medskip

The second result is due to Feyel and de la Pradelle \cite {FdP}. It
shows the constructivity of certain measures which do not charge sets 
a given capacity of which vanishes.
\bprop {BP1lem2} Let $\ga>0$ and $1<p<\infty$. 
If $\gl\in\mathfrak M_{+}(\Gw)$ 
does not charge sets with $C_{\ga,p'}$-capacity zero, there exists an 
increasing sequence $\{\gl_{n}\}\subset \mathfrak M_{+}^{b}(\Gw)\cap L^{-\ga,p}(\Gw)$,  
$\gl_{n}$ with compact support in $\Gw$, which converges to $\gl$.
\es
\Proof We first assume that $\gl$ has compact support in $\Gw$. 
Let $\phi\in L_{0}^{\ga,p'}(\Gw)$ and $\tilde \phi$ its quasicontinuous 
representative. Since the function $\tilde \phi_{+}$ is 
quasicontinuous too, the following functional is well defined on 
$L_{0}^{\ga,p'}(\Gw)$, with values in $[0,\infty]$,
\begin {eqnarray}\label {FdP1}
P(\phi)=\int_{\Gw}\tilde \phi_{+}d\gl.
\end {eqnarray}
If $\{\phi_{n}\}$ converges to $\phi$ in $ L_{0}^{\ga,p'}(\Gw)$, 
there exists a subsequence $\{\phi_{n_{k}}\}$ which converges 
$C_{\ga,p'}$-quasi everywhere. Hence 
$$\int_{\Gw}\tilde \phi_{+}d\gl\leq \liminf_{n_{k}\to\infty}\int_{\Gw}\tilde 
\phi_{n\,+}d\gl,$$
by Fatou's lemma, and $\phi\mapsto P(\phi)$ is lower semicontinuous. 
Since $P$ is convex and positively homogeneous of order $1$, it is 
the upper hull of all the continuous linear functionals it dominates, by 
the Hahn-Banach Theorem. \smallskip

\noindent {\it Step 1 } Let $\epsilon>0$, and $\phi_{0}\in 
L^{\ga,p'}_{0}(\Gw)$. Then we claim that there exists a positive Radon 
measure $\gth$ belonging to $L^{-\ga,p}(\Gw)$ such that 
$0\leq\gth\leq\gl$, and
\begin {eqnarray}\label {FdP2}
\int_{\Gw}\phi_{0}d(\gn-\gth)<\epsilon.
\end {eqnarray}
Clearly 
$$(\phi_{0},P(\phi_{0})-\epsilon)\notin Epi (P)=\left\{(\phi,t)\in 
L^{\ga,p'}_{0}(\Gw)\times \mathbb R:\,t\geq P(\phi)\right\}.
$$
Since $Epi (P)$ is a closed convex subset of $L^{\ga,p'}_{0}(\Gw)\times 
\mathbb R$, it follows by the Hahn-Banach Theorem that there exist a 
continuous form $\ell$ on $L^{\ga,p'}_{0}(\Gw)$ and two constants $a$ 
and $b$ such that
\begin {eqnarray}\label {FdP3}
a+bt+\ell (\phi)\leq 0,\forevery (\phi,t)\in Epi(P),
\end {eqnarray}
and
\begin {eqnarray}\label {FdP4}
a+b(P(\phi_{0})-\epsilon)+\ell (\phi_{0})> 0.
\end {eqnarray}
But $(0,0)\in Epi(P)\Longrightarrow a\leq 0$. Thus (\ref {FdP4}) 
holds with $a=0$. If we apply (\ref {FdP3}) 
to $(\gt\phi,\gt t)$ with $\gt>0$ arbitrary (such a couple belongs to $Epi(P)$ 
since $P$ is positively homogeneous) and let $\gt\to\infty$, it follows 
\begin {eqnarray}\label {FdP5}
bt+\ell (\phi)\leq 0,\forevery (\phi,t)\in Epi(P).
\end {eqnarray}
In the particular case $\phi=0$ and $t>0$ (possible since $(0,t)\in Epi 
(P)$, $\forall t>0$), it gives $b\leq 0$. If $b$ were zero one would 
have $\ell (\phi)\leq 0$ for any $(\phi,t)\in Epi(P)$, and in 
particular $\ell (\phi_{0})\leq 0$, which would contradict (\ref {FdP4}) 
if we impose $b=0$. Since $b<0$, we define $\gth$ by
$$\gth (\phi)=-\frac {\ell (\phi)}{b}, \forevery\phi\in 
L^{\ga,p'}_{0}(\Gw),
$$
and derive 
\begin {eqnarray}\label {FdP6}
P(\phi)\geq \gth (\phi),
\end {eqnarray}
for any $\phi\in L^{\ga,p'}_{0}(\Gw)$, since $(P(\phi),\phi)\in Epi(P)$. 
In the particular case where $\phi\leq 0$, there holds $\gth (\phi)\leq 0$. This 
means that $\gth$ is a continuous positive linear functional on 
$L^{\ga,p'}_{0}(\Gw)$, dominated by $P$. It defines a unique Radon 
measure, still denoted by $\gth$, and (\ref {FdP2}) holds.\smallskip

\noindent {\it Step 2 } End of the proof. We assume now that $\gl$ has 
no longer a compact support in $\Gw$. There exists an exhaustive 
sequence of open subsets $\{\Gw_{k}\}$, compactely included in $\Gw$ 
such that
$$\Gw_{k}\subset\overline \Gw_{k}\subset \Gw_{k+1}\subset\overline 
\Gw_{k+1}\subset\ldots\Gw.
$$
We put $\gl_{k}=\gl\vline_{\Gw_k}$. We apply the result of step 1 to 
$\gl_{k}$, with $\epsilon=1/k$ and $\phi\equiv 1$ on $\Gw_{k}$ and 
derive the existence of a positive Radon measure $\gth_{k}\in L^{\ga,p'}(\Gw)$, 
with compact support in $\Gw$ satisfying $0\leq\gth_{k}\leq \gl$ and 
$$\int_{\Gw_{k}}d(\gl-\gth_{k})<1/k.
$$
The measure $\gl_{n}=\sup\{\gth_{1},\gth_{2},\ldots,\gth_{n}\}$ 
has compact support in $\Gw$, $\gl_{n}\leq \gl_{n+1}\leq\gl$ for any 
$n$, 
and 
$$\lim_{n\to\infty}\int_{\Gw}\gz d\gl_{m}=\int_{\Gw}\gz d\gl,\forevery 
\gz\in C_{c}(\Gw).
$$\qeda 

\bcor {decomposition}Let $\ga>0$ and $1<p<\infty$. 
If $\gl\in\mathfrak M^b (\Gw)$ 
does not charge sets with $C_{\ga,p'}$-capacity zero, there exist a 
 function $\gl^*\in L^1(\Gw)$ and a measure 
$\tilde\gl\in L^{-\ga,p}(\Gw)$ such that 
\begin {eqnarray}\label {sum alpha-p}
\gl=\tilde\gl+\gl^*.
\end {eqnarray}
\es
\Proof By assumption, both the positive and the negative parts of 
$\gl$ do not charge sets with $C_{\ga,p'}$-capacity zero. Therefore it 
is sufficient to prove (\ref {sum alpha-p}) with $\gl\in\mathfrak 
M_{b}^+(\Gw)$. Let $\{\gl_{n}\}\subset L^{-\ga,p}(\Gw)\cap\mathfrak 
M_{+}(\Gw)$ be the increasing sequence of measures with compact 
support in $\Gw$ which converges to $\gl$ weakly. We set
$$\gr_{j}=\gl_{j}-\gl_{j-1}, \quad\mbox {for }j\in\mathbb N_{*}
, \quad\mbox {and }\gr_{0}=\gl_{0}.
$$
Then 
$$\gl=\sum_{j=0}^\ity\gr_{j},
$$
and the series converges strongly in the space $\mathfrak M^b(\Gw)$. 
In particular
$$\sum_{j=0}^\ity{\norm{\gr_{j}}}_{\mathfrak M^b(\Gw)}<\infty.
$$
Let $\{\eta_{k}\}_{k\in\mathbb N_{*}}$ be a sequence of $C^\ity$ 
nonnegative functions in $\mathbb R^n$, with compact support in the 
open ball $B_{k^{-1}}(0)$, satisfying
$$\int_{\Gw}\eta_{k} dx=1.
$$
For any $j\in\mathbb N_{*}$ there exists $k^0_{j}\in \mathbb N_{*}$ 
such that for $k\geq k^0_{j}$, $\gr_{j,k}=\gr_{j}\ast\eta_{k}\in 
C_{c}^{\ity}(\Gw)$. Since $\gr_{j,k}\to \gr_{j}$ as $k\to\infty$, we 
fix $k_{j}\geq k^0_{j}$ such that 
$$\norm{\gr_{j,k_{j}}-\gr_{j}}_{L^{-\ga,p}(\Gw)}\leq 2^{-j}.
$$
 We set $\tilde\gr_{j,k_{j}}=\gr_{j}-\gr_{j,k_{j}}$. The series
 $\displaystyle {\sum_{j=0}^\infty\tilde\gr_{j,k_{j}}}$
is normaly convergent in $L^{-\ga,p}(\Gw)$ and, if $\tilde\gl$ denotes its sum, it 
belongs to $L^{-\ga,p}(\Gw)$. Moreover
$${\norm {\gr_{j,k_{j}}}}_{L^1(\Gw)}=
{\norm {\gr_{j}\ast\eta_{k_{j}}}}_{L^1(\Gw)}=
{\norm{\gr_{j}}}_{\mathfrak M^b(\Gw)}. 
$$
Thus the series $\displaystyle {\sum_{j=0}^\infty\gr_{j,k_{j}}}$ is 
normaly convergent in $L^1(\Gw)$ with sum $\gl^*$. The three series
$\displaystyle {\sum_{j=0}^\infty\gr_{j}}$, 
$\displaystyle {\sum_{j=0}^\infty\tilde \gr_{j,k_{j}}}$ and 
$\displaystyle {\sum_{j=0}^\infty\gr_{j,k_{j}}}$ converge in the sense 
of distributions in $\Gw$, therefore 
(\ref {sum alpha-p}) holds.
\qeda\medskip

\noindent \Remark If $\gl\geq 0$, it is the same with $\gl^{*}$. 
Unfortunately it is not clear that $\tilde \gl$ inherits the same 
property. Notice that $\gl^{*}$ and $\tilde \gl$ may not be mutually 
singular.\medskip

Another important and useful result concerning measures which do not 
charge sets with zero capacity is the following \cite {DaM}.
\bth {DAM} Let $\ga>0$ and $1<p<\infty$. If $\gl\in\mathfrak M_{+}(\Gw)$ 
does not charge sets with $C_{\ga,p'}$-capacity zero, there exist 
$\gn\in\mathfrak M_{+}(\Gw)\cap L^{-\ga,p}(\Gw)$ and a Borel function 
$f$ with value in $[0,\infty)$ such that 
\begin {equation}\label{dam1}
\gl (E)=\int_{E}fd\gn,\forevery E\subset \Gw,\; E\,\mbox { Borel}.
\end {equation}
\es 
\subsubsection {Sharp solvability}
The following theorem due to Baras and Pierre \cite {BP1} 
characterizes the bounded measures for which the problem
 \begin {equation}\label {power}\left.\BA{ll}
Lu+\abs u^{q-1}u=\gl&\mbox { in }\Omega,\\[2mm]
\phantom {Lu+\abs u^{q-1}\,}
u=0&\mbox { on }\prt\Omega, 
\EA\right.
\end {equation}
admits a solution.
\bth {BP1th} Let $\Gw$ be a $C^2$ bounded domain in $\mathbb R^n$, 
$n\geq 3$, $L$ the elliptic operator defined by (\ref {lin1}) 
satisfying the condition (H), $q\geq n/(n-2)$ and $\gl\in\mathfrak M^{b}(\Gw)$. 
Then Problem (\ref {power}) admits a solution if and only if $\gl$ 
does not charge sets with $C_{2,q'}$-capacity zero. The solution is 
unique and the 
mapping $\gl\mapsto u$ is nondecreasing.
\es
For proving this theorem we need the following regularity result.
\blemma {regularity} Let $\Gw$ and $L$ be as in \rth {BP1th}. Then for 
any $1<p<\ity$ and $\gl\in W^{-2,p}(\Gw)\cap \mathfrak M^b(\Gw)$, 
$\mathbb G_{L}^\Gw(\gl)\in L^p(\Gw)$. Moreover there exists 
$C=C(\Gw,L,p)>0$ such 
that
\begin {equation}\label {est-p}
{\norm{\mathbb G_{L}^\Gw(\gl)}}_{L^p(\Gw)}\leq C
{\norm{\gl}}_{W^{-2,p}(\Gw)}.
\end {equation}
\es 
\Proof Put $v=\mathbb G_{L}^\Gw(\gl)$, then
$$\int_{\Gw}vL^*\gz dx=\int_{\Gw}\gz d\gl,\forevery\gz\in 
C_{c}^{1,L}(\overline\Gw).
$$
Let $\phi\in C_{0}^\infty(\Gw)$, $\gz=\mathbb G_{L^*}^\Gw(\phi)$, then
$$\abs {\int_{\Gw}v\phi dx}\leq {\norm{\gl}}_{W^{-2,p}(\Gw)}
{\norm{\gz}}_{W^{2,p'}(\Gw)}\leq C{\norm{\gl}}_{W^{-2,p}(\Gw)}
{\norm{\phi}}_{L^{p'}(\Gw)},
$$
by the $L^p$-regularity theory of elliptic equations. Hence $v\in 
L^p(\Gw)$ and (\ref {est-p}) follows.\qeda 
\medskip

\noindent {\it Proof of \rth{BP1th}.} (i) Assume that $u$ is a solution of (\ref {power}). 
Since $\abs u^{q-1}u\in L^1(\Gw)$ by \rprop {genest}, it does not 
charge set with $C_{2,q'}$-capacity 
zero, which are negligible sets for the $n$-dimensional Hausdorff measure. 
Therefore $Lu\in \mathfrak M^b(\Gw)$, and
$$\abs{\langle Lu,\phi\rangle}=\abs{\int_{\Gw}u L^*\phi dx}\leq {\norm u}_{L^q(\Gw)}
{\norm {L^*\phi}}_{L^{q'}(\Gw)}\leq C{\norm u}_{L^q(\Gw)}
{\norm \phi}_{W^{2,q'}(\Gw)},
$$
for any $\phi\in C^{\infty}_{0}(\Gw)$. Therefore the measure $Lu$ 
defines a continuous linear functional on $W^{2,q'}_{0}(\Gw)$. 
Consequently $\gl$ is the sum of an integrable function and a measure 
in $W^{-2,q}(\Gw)$.\smallskip

(ii) Conversely, we first assume that $\gl$ is a positive measure. 
By \rprop {BP1lem2} there exists an increasing sequence of positive 
measures $\gl_{j}$ belonging to $W^{-2,q}$ converging to $\gl$ in the 
weak sense of measures. By \rth {g-adm-th} there exists $u_{j}$ 
solution to 
 \begin {equation}\label {power-j}\left.\BA{ll}
Lu_{j}+\abs {u_{j}}^{q-1}u_{j}=\gl_{j}&\mbox { in }\Omega,\\[2mm]
\phantom {Lu_{j}+\abs {u_{j}}^{q-1}\,}
u_{j}=0&\mbox { on }\prt\Omega. 
\EA\right.
\end {equation}
Moreover $u_{j}$ is nonnegative and $u_{j}\geq u_{j-1}$ for any 
$j\in\mathbb N_{*}$. For any $\gz\in C_{c}^{1,L}(\overline\Gw)$ there holds
\begin {eqnarray}\label {lpower-j-2}
\int_{\Gw}\left(u_{j}L^{*}\gz+u_{j}^q\gz\right)dx=\int_{\Gw}\gz 
d\gl_{j}.
\end {eqnarray}
Let $u=\lim_{j \to\infty}u_{j}$. If $\gz\geq 0$, we have, by the 
Beppo-Levi Theorem,
\begin {eqnarray}\label {power-j-2}
\int_{\Gw}\left(uL^{*}\gz+u^q\gz\right)dx=\int_{\Gw}\gz 
d\gl.
\end {eqnarray}
Hence $u\in L^{1}(\Gw)\cap L^q(\Gw;\rho_{_{\prt\Gw}}dx)$ and $u$ is 
the solution to Problem (\ref {power}). Because $\gl$ is bounded we 
have
$u\in L^{q}(\Gw)$ by \rprop {genest}. \smallskip

If $\gl$ is no longer 
positive, $\gl_{+}$ and $\gl_{-}$ do not charge Borel sets with 
 $C_{2,q'}$-capacity zero. Hence there exist two nondecreasing sequences of positive 
measures belonging to $W^{-2,q}(\Gw)$, $\{\gl_{j,+}\}$ and $\{\gl_{j,-}\}$, 
converging to $\gl_{+}$ and $\gl_{-}$ respectively. As in the proof 
of \rth {g-adm-th} we truncate the nonlinearity by putting $g_{k}(r)=$sign $ 
(r)\min\{k^q,\abs r^q\}$ for $k\in\mathbb N_{*}$, and we denote by 
$v_{k}$, (resp. $v_{k,+}$ and $v_{k,-}$) the solutions of 
 \begin {equation}\label {power-k}\left.\BA{ll}
Lv+g_{k}(v)=\gn&\mbox { in }\Omega,\\[2mm]
\phantom {Lv+g_{k}()\,}
v=0&\mbox { on }\prt\Omega, 
\EA\right.
\end {equation}
when $\gn=\gl_{j,+}-\gl_{j,-}$ (resp. $\gn=\gl_{j,+}$ and $\gn=\gl_{j,-}$). By \rth 
{genMeas}, $-v_{k,-}\leq v_{k}\leq v_{k,+}$, which implies $-g_{k}(v_{k,-})\leq 
g_{k}(v_{k})\leq g_{k}(v_{k,+})$. When $k\to \ity$, the sequences
$\{v_{k,+}\}$ and $\{v_{k,-}\}$ decrease and converge respectively to 
$u_{j,+}$ and $u_{j,-}$, the solutions of (\ref {power}) with respective 
right-hand side $\gl_{j,+}$ and $\gl_{j,-}$. Moreover
\begin {eqnarray}\label {power-l}
-\left(\mathbb G_{L}^\Gw(\gl_{j,-})\right)^q\leq -g_{k}(\mathbb G_{L}^\Gw(\gl_{j,-}))\leq g_{k}(v_{k})\leq 
g_{k}\left(\mathbb G_{L}^\Gw(\gl_{j,+})\right)\leq 
\left(\mathbb G_{L}^\Gw(\gl_{j,+})\right)^q.
\end {eqnarray}
Since the left and right-hand side terms are $L^{1}(\Gw)$-functions, the sequence 
$\{g_{k}(v_{k})\}$ is uniformly integrable. As in the proof of \rth 
{g-adm-th}, the sequence $\{v_{k}\}$ converges in $L^{q}(\Gw)$ to the 
solution $u_{j}$ of (\ref {power-j}) with right-hand side $\gl_{j,+}-\gl_{j,-}$. 
Furthermore 
$$-u_{j,-}\leq u_{j}\leq u_{j,+},\quad \mbox {and }\;-u^q_{j,-}\leq 
\abs{u_{j}}^{q-1}u_{j}\leq u^q_{j,+}.$$
Because $\{u_{j,+}\}$ and $\{u_{j,-}\}$ are monotone and converge 
in $L^q(\Gw)$, the sequence $\{u_{j,}\}$ is uniformly integrable in 
$L^q(\Gw)$ and converges a.e. in $\Gw$. Since $\gl_{j,+}-\gl_{j,-}$ 
converges weakly to $\gl$ in the sense of measures, there exists a 
function $u\in 
L^q(\Gw)$, solution of (\ref {power}).
\qeda 
\subsection {Removable singularities}
\subsubsection{Positive solutions}
In this section $\Gw$ is an arbitrary open set in $\mathbb R^n$. 
Let $L_{m}$ be a linear differential operator of order 
$m$ ($m\in \mathbb N_{*}$), defined by
\begin {equation}\label {m-operator}
L_{m}u=\sum_{0\leq \abs\ga\leq m}D^\ga (a_{\ga}u),
\end {equation}
where
\begin {equation}\label {m-operator2}
a_{\ga}\in L^\infty_{loc}(\Gw),\forevery\ga\in\mathbb N^n,\;\abs\ga\leq m.
\end {equation}
\bdef {weak-m} {\rm Let $G\subset\Gw$ be open, $u\in L^1_{loc}(G)$ and $T$ 
a distribution on $G$. We shall say that $u$ satisfies
\begin {equation}\label {m-operator3}
L_{m}u=T\quad \left(\mbox {resp. }L_{m}u\leq T\right )\quad \mbox {in }\CD'(G),
\end {equation}
or, equivalently, that $u$ is a distribution solution (resp. 
subsolution) of (\ref {m-operator3}), if
\begin {equation}\label {m-operator4}\left.\BA {l}
\myint{G}{}uL_{m}^*\gz dx=\langle T,\gz\rangle \quad (\mbox {resp. }
\myint{G}{}uL_{m}^*\gz dx\leq \langle T,\gz\rangle ),
\\[2mm]
\phantom{\myint{G}{}uL_{m}^*\gz dx=\langle T,\gz\rangle, 
\quad (\mbox {resp. }
}
\forevery\gz\in C_{c}^\infty(G)\quad (\mbox {resp.}\forevery\gz\in 
C_{c}^\infty(G)\,,\,\gz\geq 
0), 
\EA\right.\end {equation}
where $\langle.,.\rangle$ denote the duality pairing between 
$\CD'(G)$ and $\CD(G)$, and $L_{m}^*$ the formal adjoint of $L_{m}$ 
defined by
\begin {equation}\label {m-operator4b}
L^*_{m}\gz=\sum_{0\leq \abs\ga\leq m}(-1)^{\abs\ga}a_{\ga}D^\ga \gz.
\end {equation}
}\es

The following result is due to Baras and Pierre \cite {BP1}.
\bth{remov1}  Let $m\in\mathbb N_{*}$, $q>1$, $F$ be a relatively closed subset
of $G$, $\gl$ a Radon measure which does not charge sets with 
$C_{m,q'}$-capacity zero and $g$ a continuous real valued function 
which satisfies
\begin {equation}\label {m-operator5}
\liminf_{r\to\infty}g(r)/r^{q}>0.
\end {equation}
Let $u\in L^1_{loc}(\Gw\setminus F)$, such that $u\geq 0$ and
 $g(u)\in L^1_{loc}(\Gw\setminus F)$, be a solution of
 \begin {equation}\label {m-operator6}
 L_{m}u+g(u)\leq \gl\quad \mbox {in }\CD'(\Gw\setminus F).
 \end {equation}
 If $C_{m,q'}(F)=0$, then $u\in L^1_{loc}(\Gw)$, $g(u)\in L^1_{loc}(\Gw)$ and 
there holds
  \begin {equation}\label {m-operator7}
 L_{m}u+g(u)\leq \gl\quad \mbox {in }\CD'(\Gw).
 \end {equation}
\es 
\Proof {\it Step 1 } Let $\gz\in C_{c}^\infty(\Gw)$, and 
$K=$supp$(\gz)$. Since $K\cap F$ is a compact subset of $\Gw$ with 
$C_{m,q'}$-capacity zero, it follows by \rprop {zero cap} that 
there exists a sequence $\{\phi_{n}\}\subset 
C_{c}^\infty(\Gw)$ such that $0\leq\phi_{n}\leq 1$, $\phi_{n}\equiv 
1$ in a neighborhood of $K\cap F$ and $\phi_{n}\to 0$ as $n\to\infty$, in $W^{m,q'}(\Gw)$ 
 and $C_{m,q'}$-quasi everywhere. Therefore, 
$\gz_{n}=(1-\phi_{n})\gz$ satisfies :\smallskip

\noindent (i) $\gz_{n}\in C_{c}^\infty(\Gw\setminus F)$,\smallskip

\noindent (ii) $0\leq \gz_{n}\leq 1$,\smallskip

\noindent (iii) $\gz_{n}\to\gz$ in $W^{m,q'}(\Gw)$ 
 and $C_{m,q'}$-quasi everywhere as $n\to\infty$, and the 
 sequence $\{\gz_{n}\}$ is increasing.\smallskip
 
 \noindent {\it Step 2 } We claim that $g(u)\in L^1_{loc}(\Gw)$. We take $\gz\in C_{c}^\infty(\Gw)$, $\gz\geq 
0$ and $\{\gz_{n}\}$ be defined by the procedure in Step 1. Let 
$p\in\mathbb N$, $p\geq mq'$. Since $\gz_{n}^p\in C_{c}^\infty(\Gw\setminus 
F)$, (\ref {m-operator6}) implies
\begin {equation}\label {m-operator8}
\int_{\Gw}\left(u L_{m}^*(\gz_{n}^p)+g(u)\gz_{n}^p\right)dx\leq 
\int_{\Gw}\gz_{n}^pd\gl.
\end {equation}
Because $\gz_{n}^p\leq\gz$, there holds
\begin {equation}\label {m-operator9}
\int_{\Gw}g(u)\gz_{n}^pdx\leq 
\int_{\Gw}\gz d\abs {\gl}+\int_{\Gw}u \abs {L_{m}^*(\gz_{n}^p)}dx.
\end {equation}
Since the $a_{\ga}$ are locally bounded, 
$$\abs {L_{m}^*(\gz_{n}^p)}
\leq C\sum_{0\leq\abs\ga\leq m}{\abs{D^\ga(\gz_{n}^p)}}.
$$
The zero order term is estimated by
\begin {equation}\label {m-operator10}
\int_{\Gw}u\gz_{n}^pdx\leq \left(\int_{\Gw}u^q\gz_{n}^pdx\right)^{1/q}
\left(\int_{\Gw}\gz_{n}^pdx\right)^{1/p'}\leq 
\left(\int_{\Gw}u^q\gz_{n}^pdx\right)^{1/q'}
{\norm{\gz_{n}}}_{W^{m,q'}(\Gw)}.
\end {equation}
If $\abs\ga\geq 1$, 
\begin {equation*}
D^\ga(\gz_{n}^p)
=\sum_{j=1}^{ \abs\ga} c_{j}\gz_{n}^{p-j}\!\!\!\!\!\!\!\!\!\!\!\!
\sum_{\tiny{\BA{c}\abs{\gb_{i}}\geq 1\\
\gb_{1}+...+\gb_{j}=\ga\EA}}c_{\gb_{1},...,\gb_{j}}
D^{\gb_{1}}\gz_{n}...D^{\gb_{j}}\gz_{n},
\end {equation*}
where the $c_{j}$ and $c_{\gb_{1},...,\gb_{j}}$ are positive constants 
depending on the indices. Thus we are led to estimate a finite sum 
involving terms of the form
$$A=\int_{\Gw}u\gz_{n}^{p-j}
\abs{D^{\gb_{1}}\gz_{n}...D^{\gb_{j}}\gz_{n}} dx.
$$
By H\"older's inequality
$$A\leq \left(\int_{\Gw}u^q\gz_{n}^pdx\right)^{1/q}
\left(\int_{\Gw}\gz_{n}^{p-jq'}
\abs{D^{\gb_{1}}\gz_{n}...D^{\gb_{j}}\gz_{n}}^{q'}dx\right)^{1/q'}.
$$
Because $p\geq mq'\geq jq'$, it follows $0\leq \gz_{n}^{p-jq'}\leq 1$. By applying 
again H\"older's inequality, and using the fact that 
$\abs {\gb_{1}}+...+\abs {\gb_{j}}=\abs\ga$, it follows
$$A\leq \left(\int_{\Gw}u^q\gz_{n}^pdx\right)^{1/q}
\prod_{i=1}^j
\left(\int_{\Gw}\abs{D^{\gb_{i}}\gz_{n}}^{q'\abs{\ga}/\abs{\gb_{i}}}dx\right)
^{\abs{\gb_{i}}/\abs{\ga}q'}.
$$
By the Gagliardo-Nirenberg inequality, there holds
$$\abs {D^{\gb_{i}}\gz_{n}}^{q'\abs{\ga}/\abs{\gb_{i}}}
\leq C{\norm{\gz_{n}}}_{W^{\abs\ga,q'}(\Gw)}^{\abs{\gb_{i}}/\abs\ga}
\leq C{\norm{\gz_{n}}}_{W^{m,q'}(\Gw)}^{\abs{\gb_{i}}/\abs\ga}.
$$
Therefore
\begin {equation}\label {m-operator11}
A\leq C\left(\int_{\Gw}u^q\gz_{n}^pdx\right)^{1/q}
{\norm{\gz_{n}}}_{W^{m,q'}(\Gw)},
\end {equation}
from which derives
\begin {equation}\label {m-operator12}
\int_{\Gw}g(u)\gz_{n}^pdx\leq C_{1}+C_{2}\left(\int_{\Gw}u^q\gz_{n}^pdx\right)^{1/q}
{\norm{\gz_{n}}}_{W^{m,q'}(\Gw)}.
\end {equation}
By assumption, there exist two positive constants $a$ and $b$ such that
$$g(r)\geq ar^q-b,\forevery r\geq 0.
$$
Consequently, up to changing the constants $C_{i}$, 
\begin {equation}\label {m-operator13}
\int_{\Gw}(g(u)+b)\gz_{n}^pdx\leq C_{1}+C_{2}
\left(\int_{\Gw}(g(u)+b)\gz_{n}^pdx\right)^{1/q}
{\norm{\gz_{n}}}_{W^{m,q'}(\Gw)}.
\end {equation}
Finally, the left-hand side integral remains bounded independently 
of $n$ and we conclude by Fatou's lemma that $(g(u)+b)\gz^p\in 
L^1(\Gw)$.  Since $\gz$ is arbitrary, we find $g(u)\in L^1_{loc}(\Gw)$. The 
growth estimate on $g$ implies also $u\in L^q_{loc}(\Gw)$.\smallskip

\noindent {\it Step 3 } We claim that (\ref {m-operator7}) holds. 
Let $\gz\in C_{c}^\infty(\Gw)$, $\gz\geq 0$. By  constructing the same 
functions $\gz_{n}$ as above, we have
\begin {equation}\label {m-operator14}
\int_{\Gw}\left(u L_{m}^*\gz_{n}+g(u)\gz_{n}\right)dx\leq 
\int_{\Gw}\gz_{n}d\gl.
\end {equation}
Since $\abs\gl$ does not charge sets with $C_{m,q'}$-capacity zero and 
$\gz_{n}\to\gz$, $C_{m,q'}$-quasi everywhere in $\Gw$, this 
convergence holds also $\abs\gl$-a.e. in $\Gw$. By the Lebesgue Theorem
$$\lim_{n\to\infty}\int_{\Gw}\gz_{n}d\gl=\int_{\Gw}\gz d\gl.
$$
Because $g(u)$ is locally integrable in $\Gw$,
$$\lim_{n\to\infty}\int_{\Gw}g(u)\gz_{n}dx=\int_{\Gw}g(u)\gz dx,
$$
and finally, the convergence of $\{\gz_{n}\}$ to $\gz$ in 
$W^{m,q'}(\Gw)$ 
implies the convergence of $\{L_{m}^*\gz_{n}\}$ to $L_{m}^*\gz$ in 
$L^{q'}(\Gw)$. Passing to the limit in (\ref {m-operator14}) yields 
to (\ref {m-operator7}).\qeda \medskip

\noindent\Remark Contrary to the case of semilinear elliptic equations with 
an absorbing nonlinearity, which will be studied in next section, the removability of $F$ does not 
imply that the function $u$ is regular in whole $\Gw$ : the singularity 
is just not seen at the distributions level. 
\subsubsection {Semilinear elliptic equations with absorption}

The first result of unconditional removability of isolated sets 
for semilinear elliptic equations with absorption term is due to 
Brezis and V\'eron \cite {BV}. It deals with equation
\begin {equation}\label {BV1}
-\Gd u+g(u)=0,
\end {equation}
in $\Gw\setminus\{0\}$, where $\Gw$ is an open subset of $\mathbb R^n$ ($n\geq 3$) 
containing $0$ and $g$ a 
continuous function. They proved the following.
\bth {BVremth} Suppose $g$ satisfies
\begin {equation}\label {BV2}
\liminf_{r\to\infty}g(r)/r^{n/(n-2)}>0\quad\mbox{and }\;
\limsup_{r\to-\infty}g(r)/\abs r^{n/(n-2)}<0.
\end {equation}
If $u\in L^\infty_{loc}(\Gw\setminus\{0\})$ satisfies (\ref {BV1}) in 
the sense of distributions in $\Gw\setminus\{0\}$, there exists a 
 function $\tilde u \in C^1(\Gw)\cap W^{2,p}_{loc}(\Gw)$ for any 
 $1\leq p<\infty$, which coincides with $u$ a.e. in $\Gw$, and is a
solution of (\ref {BV1}) in whole $\Gw$.
\es

The proof of this result is settled upon a particular case of a 
general {\it a priori} estimate discovered  by Keller \cite {Ke} and Osserman \cite 
{Os} separately. In this particular case, and in assuming that 
$B_{R}(0)\subset\Gw$, it reads
\begin {equation}\label {BV3}
\abs {u(x)}\leq A\abs x^{2-n}+B,\forevery x\in 
B_{R/2}(0)\setminus\{0\},
\end {equation}
for some positive constants $A$ and $B$. From this estimate is derived 
the local integrability of $u$ in $\Gw$ and then of $g(u)$. Finally, it 
leads to the fact that Equation (\ref {BV1}) holds in the sense of distributions in $\Gw$. 
The conclusion follows by the maximum principle (which implies the 
boundedness of $u$ near $0$), and the elliptic equations regularity 
theory. Later on, this result was extended by V\'eron \cite {Ve7} as 
follows :
\bth {Vremth} Let $\Gs\subset\Gw$ be a complete and compact $d$-dimensional 
submanifold of class $C^2$ ($1\leq d< n-2$) and $g$ is a continuous 
real valued function such that
\begin {equation}\label {BV2b}
\liminf_{r\to\infty}g(r)/r^{(n-d)/(n-2-d)}>0\quad\mbox{and }\;
\limsup_{r\to-\infty}g(r)/\abs r^{(n-d/(n-2-d)}<0.
\end {equation}
If $u\in L^\infty_{loc}(\Gw\setminus\Gs)$ satisfies (\ref {BV1}) in 
the sense of distributions in $\Gw\setminus\Gs$, there exists a 
 function $\tilde u \in C^1(\Gw)\cap W^{2,p}_{loc}(\Gw)$ for any 
 $1\leq p<\infty$, which coincides with $u$ a.e. in $\Gw$ and is a
solution of (\ref {BV1}) in whole $\Gw$.
\es

Although more technical, the idea of the proof is similar to the one 
of \rth {BVremth}, except that the {\it a priori} estimate (\ref {BV3}) is replaced by
\begin {equation}\label {BV3b}
\abs {u(x)}\leq A\left(\dist(x,\Gs)\right)^{2-n-d}+B,\forevery x\in 
G\setminus\Gs,
\end {equation}
where $G$ is open and bounded and $\Gs\subset G\subset\overline G\subset 
\Gw$. The method developed by Baras and Pierre \cite {BP1} is settled 
upon integral identity, without using pointwise {\it a priori} estimates as 
the previous authors do.
\bth {remov2} Let $\Gw$ be a bounded open subset of $\mathbb R^n$, 
$n\geq 2$, with a $C^2$ boundary, $L$ an
elliptic operator defined by \ref {lin1} satisfying condition (H) and 
$q>1$. 
If $F$ is a compact subset of $\Gw$, any solution $u\in 
L^{q}_{loc}(\Gw\setminus K)$ of
\begin {equation}\label {2operator1}
 Lu+{\abs {u}}^{q-1}u=0,
\end {equation}
in ${\Gw\setminus K}$, belongs to $L^{q}_{loc}(\Gw)$ and 
satisfies (\ref {2operator1}) in whole $\Gw$, if and only if 
$C_{2,q'}(K)=0$. If this holds, $u\in W_{loc}^{2,p}(\Gw)$ for any 
$1\leq p<\infty$, and (\ref {2operator1}) is satisfied a.e. in $\Gw$.
\es 
\Proof (i) Let us assume that $C_{2,q'}(K)>0$. By (\ref {Bessel4'}),  there exists a 
positive Radon measure $\gl$ concentrated on $K$ such that 
$$\int_{\Gw}{\abs{G_{2}\ast\gm}}^qdx<\infty.
$$
This means that $\gl\in W^{-2,q}(\Gw)$. By \rth {BP1th}, Problem (\ref 
{power}) admits a solution in $\Gw$.\smallskip

\noindent (ii) Conversely we assume that $C_{2,q'}(K)=0$. By \rth 
{L1-th}, for any $\gz\in C^{1,L}_{c}(\Gw\setminus F)$, $\gz\geq 0$, 
there holds
$$\int_{\Gw}\left(\abs uL^*\gz+\abs u^{q}\gz\right)dx\leq 0.
$$
Therefore $v=\abs u$ is a subsolution of (\ref{2operator1}) in the sense 
of \rdef {weak-m}. Since $C_{2,q'}(K)=0$, we can extend $v$ 
as a solution of (\ref{2operator1}) in whole 
$\Gw$, and because $K$ has zero Lebesgue measure, $u\in L^q_{loc }(\Gw)$. 
Let $\gz_{n}=(1-\phi_{n})\gz$ be the functions defined in \rth{remov1} 
for an arbitrary $\gz\in C_{c}^\ity(\Gw)$ (we do not impose the 
positivity). Then $\gz_{n}\to\gz$ in $W^{2,q'}(\Gw)$ and 
$C_{2,q'}$-quasi everywhere. By assumption
$$\int_{\Gw}\left( uL^*\gz_{n}+\abs u^{q-1}u\gz_{n}\right)dx= 0.
$$
By Lebesgue's theorem, $\abs u^{q-1}u\gz_{n}\to \abs u^{q-1}u\gz$ in 
$L^1(\Gw)$. Moreover $L^*\gz_{n}\to L^*\gz$ in 
$L^{q'}(\Gw)$. Therefore, by letting $n\to \infty$, it is infered that
\begin {equation}\label {2operator2}
\int_{\Gw}\left( uL^*\gz+\abs u^{q-1}u\gz\right)dx= 0,
\end {equation}
which proves that (\ref {2operator1}) holds in $\Gw$. Let $G$ be any 
smooth open domain containing $K$ and such that $\overline G\subset\Gw$. 
For $\gb>0$ small enough we put $G_{\gb}=\{x\in G:\dist (x,\prt 
G>\gb\}$, and $\Gg_{\gb}=\{x\in G:\dist (x,\prt G)=\gb\}=\prt G_{\gb}$. 
There exists $\gb_{0}$ such that $\Gg_{\gb}$ is a smooth surface in 
$\mathbb R^n$. Because $u\in L^q(G\setminus\overline G_{\gb_{0}})$, it 
follows, by Fubini's theorem, that 
$u\vline_{\Gg_{\gb}}\in L^q(\Gg_{\gb})$ (endowed with the 
$(n-1)$-dimensional Hausdorff measure), for almost all $\gb\in [0,\gb_{0}]$. We fix a $\gb$ such that this 
property holds and denote by $V$ the Poisson potential of 
$u_{+}\vline_{\Gg_{\gb}}$ in $G_{\gb}$. By (\ref {L1-4}), for any 
$\gz\in C^{1,L}_{c}(\overline G_{\gb})$, $\gz\geq 0$, there holds
\begin {eqnarray}\label {2operator3}
\int_{G_{\gb}}\left( (u-V)_{+} L^*\gz+(u-V)_{+}\abs u^{q-1}u\gz \right)dx
\leq -\int_{\prt G_{\gb}}\frac{\prt\gz}{\prt{\bf n}_{L^*}}(u-u_{+})_{+}dS.
\end {eqnarray}
Taking $\gz=\mathbb G_{L}^{G_{\gb}}(1)$ implies $(u-V)_{+}\equiv 0$ 
in $G_{\gb}$. Thus $u\leq V$ in $G_{\gb}$. Since $V\in 
L_{loc}^{\infty}(G_{\gb})$, the same property holds with $u_{+}$. 
Since $G$ is arbitrary, $u_{+}\in L_{loc}^{\infty}(\Gw)$. In the same 
way $u_{-}\in L_{loc}^{\infty}(\Gw)$. We conclude with the elliptic 
equations regularity theory that $u\in W^{2,p}_{loc}(\Gw)$.\qeda 
\medskip

\noindent \Remark The following extension of \rth {remov2} is easy to 
establish : {\it
Let $g$ be a continuous real valued function 
which satisfies
\begin {eqnarray}\label {2operator4}
\liminf_{r\to\ity}g(r)/r^q>0\quad\mbox {and }\;
\limsup_{r\to-\ity}g(r)/\abs r^q<0,
\end {eqnarray}
for some $q>1$. Let $\gl\in \mathfrak M(\Gw)$ which does not charge sets with 
$C_{2,q'}$-capacity zero and $K$ a compact subset of 
$Gw$ with $C_{2,q'}$-capacity zero. Then  any function $u$, locally 
integrable in $\Gw\setminus K$ and such that $g(u)\in L_{loc}^1(\Gw\setminus 
K)$,
which verifies
\begin {eqnarray}\label {2operator5}
Lu+g(u)=\gl,
\end {eqnarray}
in $\CD'(\Gw\setminus K)$, can be extended as a solution of the same 
equation in $\CD'(\Gw)$. Furthermore $g(u)\in C(\Gw)$ and 
$u\in W^{2,p}_{loc}(\Gw)$, for any $1\leq p<\infty$.}
\subsection {Isolated singularities}

The description of the behaviour of solutions of semilinear elliptic 
equations near an isolated singularity deals with the following 
question : let $u$ be a solution of
\begin {eqnarray}\label {isolsing1}
Lu+g(u)=0\quad\mbox {in }\Gw\setminus \{0\},
\end {eqnarray}
where $\Gw$ is an open subset of $\mathbb R^n$ containing $0$, $L$ a elliptic 
operator under the form (\ref {lin2}) and $g$ a continuous 
real-valued function, can one describe the behaviour of 
$u(x)$ as $x\to 0$ ? When $L=-\Gd$ and $g=0$, it is known that $u$ admits an expansion 
in series of spherical harmonics. For the equation
\begin {eqnarray}\label {isolsing2}
-\Gd u+\abs u^{q-1}u=0\quad\mbox {in }\Gw\setminus \{0\},
\end {eqnarray}
($q>1$), much work on this subject has been done by V\'eron in \cite {Ve6}. Notice that if 
$q\geq n/(n-2)$ Brezis-V\'eron's result (see \rth {BVremth}) applies and the 
function $u$ is $C^2$ in whole $\Gw$. When $1<q<n/(n-2)$ this is no 
longer the case. For example there exists an explicit radial 
singular solution 
of (\ref {isolsing2}),
\begin {eqnarray}\label {isolsing3}
x\mapsto u_{s}(x)=\ell_{q,n}\abs x^{-2/(q-1)}
\end {eqnarray}
defined in $\mathbb R^n\setminus\{0\}$, where
\begin {eqnarray}\label {isolsing4}
\ell_{q,n}=\left(\left(\frac{2}{q-1}\right)\left(\frac{2q}{q-1}-n\right)\right)^{1/(q-1)}.
\end {eqnarray}
When $1<q<(n+1)/(n-1)$ there exist separable singular solutions. 
For expressing them, let $(r,\gs)$ be the spherical 
coordinates in $\mathbb R^n$ and $\Gd_{S^{n-1}}$ the Laplace-Beltrami 
operator on the unit sphere $S^{n-1}:=\{x\in \mathbb R^n:\abs x=1\}$. 
If $1<q<(n+1)/(n-1)$, one has $\ell_{q,n}> n-1=\gl_{1}(S^{n-1})$, the first 
nonzero eigenvalue of $\Gd_{S^{n-1}}$. Therefore, the classical 
variational analysis applies and there exist non-trivial solutions of 
\begin {eqnarray}\label {isolsing5}
-\Gd_{S^{n-1}}\gw-\ell_{q,n}\gw+\abs\gw^{q-1}\gw=0\quad\mbox {in }S^{n-1}.
\end {eqnarray}
Hence the function
\begin {eqnarray}\label {isolsing6}
x\mapsto u_{\gw}(x)=u_{\gw}(r,\gs)=r^{-2/(q-1)}\gw(\gs)
\end {eqnarray}
is a singular solution of (\ref {isolsing2}). Notice that $u_{s}$ is 
one of these solutions. Furthermore the constants $\ell_{q,n}$ and 
$-\ell_{q,n}$ are the only solutions of (\ref {isolsing5}) which have 
a constant sign. The following result is proven in \cite {Ve6}.
\bth {isolth1} Let $1<q<n/(n-2)$ ($q>1$ if $n=2$) and $u$ be positive solution of (\ref {isolsing2}) 
in some open set $\Gw$ containing $0$. Then, \smallskip

\noindent (i) either 
\begin {eqnarray}\label {isolsing7}
\lim_{x\to 0}\abs x^{2/(q-1)}u(x)=\ell_{q,n},
\end {eqnarray}

\noindent (ii) or there exists some $c\geq 0$ such that 
\begin {eqnarray}\label {isolsing8}
\lim_{x\to 0}\abs x^{n-2}u(x)=c,
\end {eqnarray}
if $n\geq 3$, and $\abs x^{n-2}$ replaced by $1/\ln (1/\abs x)$ in the 
above formula if $n=2$. Furthermore $u$ is a solution of 
\begin {eqnarray}\label {isolsing8c}
-\Gd u+u^q=C_{n}c\gd_{0} \quad\mbox {in }\CD'(\Gw),
\end {eqnarray}
for some positive constant  $C_{n}$ depending only on $n$.
\es

There are several proofs of this result, based either on a sharp use 
of the radial case and the Harnack inequality, or on a Lyapounov 
style analysis. If the function $u$ is no longer supposed 
to have constant sign, it is proven in \cite {Ve6} that the above 
dichotomy still holds provided $(n+1)/(n-1)\leq q<n/(n-2)$. However 
(i) has to be replaced by\smallskip

\noindent {\it \noindent (i') either 
\begin {eqnarray}\label {isolsing7'}
\lim_{x\to 0}\abs x^{2/(q-1)}u(x)=\ell\in \{\ell_{q,n},-\ell_{q,n}\},
\end {eqnarray}}
and (ii) by \smallskip

\noindent {\it \noindent (ii') or there exists some real number 
$c$ such that 
\begin {eqnarray}\label {isolsing8'}
\lim_{x\to 0}\abs x^{n-2}u(x)=c,
\end {eqnarray}
(if $n\geq 3$, with the classical modification if $n=2$). Moreover $u$ is a solution of 
\begin {eqnarray}\label {isolsing8c'}
-\Gd u+\abs u^{q-1}u=C_{n}c\gd_{0} \quad\mbox {in }\CD'(\Gw).
\end {eqnarray}
}
Actually, the Lyapounov analysis leads easily to a more general result \cite {CMV}.
\bth {isolth2} Let $1<q<n/(n-2)$ and $u$ be solution of (\ref {isolsing2}) 
in some open set $\Gw$ containing $0$. Then there exists a compact and 
connected subset $\CE$ of the set of solutions of (\ref {isolsing5}) 
such that 
\begin {eqnarray}\label {isolsing9}
\lim_{r\to 0} \dist_{C^2(S^{n-1})}(r^{2/(q-1)}u(r,.),\CE)=0,
\end {eqnarray}
where $\dist_{C^2(S^{n-1})}$ denotes the distance associated to the 
$C^2(S^{n-1})$-norm.
\es
This result leaves open two difficult questions :\smallskip

\noindent 1- Does it exist a particular element $\gw\in\CE$ such that 
\begin {eqnarray}\label {isolsing10}
\lim_{r\to 0} \norm {r^{2/(q-1)}u(r,.)-\gw}_{C^2(S^{n-1})}=0\,?
\end {eqnarray}
\smallskip

\noindent 2- What is the precise behaviour of $u$ when $\CE=\{0\}$ 
?\medskip

Besides the results above mentioned proven in \cite {Ve6}, the two 
questions have been thoroughly answered in \cite {CMV} in the 
two-dimensional case.
\bth {CMV1} Assume $n=2$, $q>1$ and $u$ is solution of (\ref {isolsing2})
in $\Gw\setminus\{0\}$. Then there exists a $2\gp$-periodic function 
$\gw$, solution of
\begin {eqnarray}\label {isolsing11}
-\frac {d^2\gw}{d\gs^2}-\left(\frac 
{2}{q-1}\right)^2\gw+\abs\gw^{q-1}\gw=0
\end {eqnarray}
such that (\ref {isolsing10}) holds on $S^1$.
\es

\bth {CMV2} Under the assumption of \rth {CMV1}, if $\gw=0$, let $k_{0}$ 
be the largest integer smaller than $2/(q-1)$. Then\smallskip

\noindent (i) either there exist an integer $k\in[1,k_{0}]$  and two
constants $A\neq 0$ and $\phi \in S^1$ such that 
\begin {eqnarray}\label {isolsing12}
\lim_{r\to 0} r^ku(r,\gs)=A\sin(k\gs+\phi),
\end {eqnarray}
in the $C^2(S^{n-1})$-topology,\smallskip

\noindent (ii) or there is a nonzero $c$ such that 
\begin {eqnarray}\label {isolsing13}
\lim_{r\to 0} u(r,\gs)/\ln(1/r)=c,
\end {eqnarray}
in the $C^2(S^{n-1})$-topology,\smallskip

\noindent (iii) or $u$ can be extended as a $C^2$
solution of (\ref {isolsing2}) in whole $\Gw$. \smallskip

\noindent In cases (ii) and (iii), $u$ is a solution of (\ref 
{isolsing8c'}) in $\CD'(\Gw)$.
\es

The proofs are extremely technical and use, in a fundamental manner, 
the Sturmian argument about the oscillations of solutions of $2$ 
dimensional elliptic 
equations jointly with the Jordan curve separation Theorem. \medskip

Many of the above results can be extended in a standard way to 
elliptic equations of the type
\begin {eqnarray}\label {isolsing14}
Lu+\abs u^{q-1}u=0,
\end {eqnarray}
where $L$ is the elliptic operator defined by (\ref {lin1}) subject 
to condition (H), and assuming $a_{ij}(x)=a_{ji}(x)$, an assumption 
which is not a real restriction. If we fix a linear change of 
variable in $\mathbb R^n$, $y=y(x)$, and write $u(x)=\tilde u(y)$, then
$$\frac {\prt^2u}{\prt x_{i}\prt x_{j}}=
\sum_{k,l} b_{li}b_{kj}\frac {\prt^2\tilde u}{\prt y_{l}\prt y_{k}},
$$
where
$$b_{\ga\gb}=\frac {\prt y_{\ga}}{\prt x_{\gb}}.
$$
Then
$$\sum_{i,j}a_{ij}(0)\frac {\prt^2u}{\prt x_{i}\prt x_{j}}=
\sum_{k,l}\frac {\prt^2\tilde u}{\prt y_{l}\prt y_{k}}
\sum_{i,j} a_{ij}(0)b_{li}b_{kj}.
$$
Since the matrix $\left(a_{ij}(0)\right)$ is symmetric, the $b_{\ga\gb}$ 
can be chosen such that 
$$\sum_{i,j} a_{ij}(0)b_{li}b_{kj}=\gd_{kl}.$$
With this transformation most of the above results can be restated 
with the variable $y$ replacing $x$. For example \rth {isolth1} transforms 
into
\bth {isolth1'} Let $1<q<n/(n-2)$ and $u$ be positive solution of  
(\ref{isolsing14}) 
in some open set $\Gw$ containing $0$. Then, \smallskip

\noindent (i) either 
\begin {eqnarray}\label {isolsing7b}
\lim_{y\to 0}\abs y^{2/(q-1)}\tilde u(y)=\ell_{q,n},
\end {eqnarray}

\noindent (ii) or there exists some $c\geq 0$ such that 
\begin {eqnarray}\label {isolsing8d}
\lim_{y\to 0}\abs y^{n-2}\tilde u(y)=c,
\end {eqnarray}
in which case $u$ is a solution of 
\begin {eqnarray}\label {isolsing8b}
L u+u^q=C_{n,L}c\gd_{0} \quad\mbox {in }\CD'(\Gw),
\end {eqnarray}
for some positive constant  $C_{n,L}$ depending only on $n$ and $L$.
\es
The description given by (\ref {isolsing7'}) of isolated singularities in the case 
of signed solutions of (\ref{isolsing14}) holds in the new unknown 
$\tilde u$ and variable 
$y$, provided $(n+1)/(n-1)<q<n/(n-2)$, and similarly the method which gives (\ref 
{isolsing9}) applies without restriction.
However the sharp analysis of the the limit case $q=(n+1)/(n-1)$ when the limit set is 
reduced to the zero function cannot be covered by this rough analysis.  
Moreover, the extension of the results given in \cite {CMV} (even in 
the non-critical cases where $2/(q-1)$ is not an integer) has not yet been done.
\subsection {The exponential and $2$-dimensional cases}
\subsubsection {Unconditional solvability}
As we have seen it above, the B\'enilan-Brezis weak-singularity 
assumption \cite {BB} is meaningless in the $2$-dimensional case for 
solving semilinear elliptic equations with bounded measures :
the $(n,0)$- weak-singularity assumption imposes $n\geq 3$ in 
\rdef{BB}. If $\Gw\subset\mathbb R^{2}$ is a smooth bounded domain, 
$L$ an elliptic operator, $g\in C(\Gw\times\mathbb R)$ is an absorbing nonlinearity and 
$\gl\in\mathfrak M^b(\Gw)$,  
a specific approach, developped by Vazquez \cite{Va1}, is 
needed, for solving
 \begin {equation}\label {2dim}\left.\BA{ll}
Lu+g(x,u)=\gl&\mbox { in }\Omega,\\[2mm]
\phantom {Lv+g(x,)\,}
u=0&\mbox { on }\prt\Omega. 
\EA\right.
\end {equation}

\bdef {EOG}{\rm Let $\tilde g\in C([0,\infty))$, $\tilde g\geq 0$. We denote by
\begin {eqnarray}\label {2dim1}
a_{+}(\tilde g):=\inf\left\{a\geq 0:\int_{0}^\ity\tilde 
g(s)e^{-as}ds<\infty\right\},
\end {eqnarray}
the {\it exponential order of growth of $\tilde g$ at infinity}. 
}\es 
 
If $ g^{*}\in C((-\infty,0])$, $g^{*}\leq 0$, 
the {\it exponential order of growth of $ g^{*}$ at minus infinity} 
is by definition the opposite of the exponential order of growth at infinity of the 
function $r\mapsto -g^{*}(-r)$, thus
\begin {eqnarray}\label {2dim1'}
a_{-}(g^{*}):=\sup\left\{a\leq 0:\int_{-\infty}^0
g^{*}(s)e^{as}ds>-\infty\right\}.
\end {eqnarray}
Those two quantities may be zero (for example if $\tilde g$ is a power), 
finite and nonzero ( if $\tilde g$ is an exponential) or infinite (if 
$\tilde g$ is a super-exponential). 
\bdef {VaV} {\rm A real valued function $g\in C(\Gw\times \mathbb R)$ 
satisfies the $2${\it-dimensional weak-singularity assumption}, if there 
exists $r_{0}\geq 0$ such that 
\begin {eqnarray}\label {passive'}
  rg(x,r)\geq 0, \forevery (x,r)\in\Gw\times 
(-\ity,-r_{0}]\cup[r_{0},\ity), 
\end{eqnarray}
and two nondecreasing functions $\tilde g_{1}\in 
C([0,\ity))$, $\tilde g_{1}\geq 0$, with zero exponential order of growth at 
infinity, and $\tilde 
g_{2}\in C((-\ity, 0])$ , $\tilde g_{2}\leq 0$, with zero exponential order of growth at 
minus infinity such that
\begin {eqnarray}\label {WWA1'}
 g(x,r) \leq \tilde g_{1}( r),\forevery (x,r)\in 
\Gw\times \mathbb R_{+},
\end{eqnarray}
and
\begin {eqnarray}\label {WWA1''}
\tilde g_{2}(r)\leq g(x,r) ,\forevery (x,r)\in 
\Gw\times \mathbb R_{-}.
\end{eqnarray}
}\es 

Notice that the zero exponential of growth assumptions can be written 
under the form
\begin {eqnarray}\label {2dim0}
\int_{0}^\infty\left(\tilde g_{1}( s)-\tilde g_{2}( -s)\right)e^{-as}ds<\infty,\forevery a>0.
\end{eqnarray}

\bth{VV1th} Let $\Gw\subset\mathbb R^{2}$ be a $C^2$ bounded domain and 
$g\in C(\Gw\times \mathbb R)$ satisfy the $2$-dimensional weak-singularity assumption. 
For any $\gl\in \mathfrak M_{b}(\Gw)$ Problem (\ref {2dim}) admits a solution. 
Furthermore, $\gd$ is invariant if we replace $g$ by $\ell g$, for 
any $\ell>0$.
\es

One of the tool of the proof is John-Nirenberg's theorem \cite [Th. 
7.21]{GT}.
\bth {JN}Let $G$ be a convex open domain in $\mathbb R^n$ and $v\in 
W^{1,1}(G)$. Assume that there exists $K>0$ such that
\begin {equation}\label {JN1}
\int_{G\cap B_{r}(a)}\abs {\nabla v}dx\leq Kr^{n-1},\forevery a\in G,\;\forall 
r>0.
\end {equation}
Then there exist two positive constants $C$ and $\gm_{0}$, depending 
only on $n$, such that
\begin {equation}\label {JN2}
\int_{G}\exp\left(\frac {\gm}{K}\abs {v-v_{G}}\right)dx
\leq C\left({\rm diam} (G)\right)^n,
\end {equation}
where $\gm=\gm_{0}\abs G\left({\rm diam} (G)\right)^{-n}$, and 
$v_{G}=\myfrac {1}{\abs G}\myint{G}{}v dx$. 
\es
Notice that for any bounded domain $G\subset \mathbb R^n$, 
${\rm diam} (G)={\rm diam} ({\rm conv\,} G)$. Then the following 
consequence of \rth {JN} is valid.

\bcor {JN2}Let $G$ be a bounded open domain in $\mathbb R^n$ and $v\in 
W_{0}^{1,1}(G)$. Assume that there exists $K>0$ such that (\ref {JN1}) 
holds. Then there exist two positive constants $C$ and $\gm_{0}$, depending 
only on $n$, such that (\ref {JN2}) holds with 
$\gm=\gm_{0}\abs {{\rm conv\;}G}\left({\rm diam} (G)\right)^{-n}$ and 
$v_{G}$ replaced by $v_{{\rm conv\;}G} =
\myfrac {1}{\abs {{\rm conv\;}G}}\myint{G}{}v dx$.
\es

\noindent {\it Proof of \rth {VV1th}. } {\it Step 1 } Approximation. 
First we multiply $\gl$ by the characteristic 
function $\chi_{_{\Gw_{n}}}$ of $\Gw_{n}=\{x\in\Gw: \rho_{_{\prt\Gw}}(x)>1/n\}$, and 
we regularize $\chi_{_{\Gw_{n}}}\gl$ by convolution with positive 
smooth functions with compact support and total mass $1$. By the property of 
convolution can replace $\gl_{+}$ and $\gl_{-}$ by 
$\gl_{n\,+}$ and $\gl_{n\,-}\in C_{c}^\ity(\Gw)$, and they satisfy, 
$${\norm {\gl_{n\,+}}}_{L^{1}(\Gw)}
\leq {\norm {\gl_{+}}}_{\mathfrak M^b(\Gw)},
$$
and 
$${\norm {\gl_{n\,-}}}_{L^{1}(\Gw)}
\leq {\norm {\gl_{-}}}_{\mathfrak M^b(\Gw)}.
$$
Let $u_{n}$ be the solution of 
 \begin {equation}\label {2dim4}\left.\BA{ll}
Lu_{n}+g(x,u_{n})=\gl_{n}&\mbox { in }\Omega,\\[2mm]
\phantom {Lu_{n}+g(x,)\,}
u_{n}=0&\mbox { on }\prt\Omega. 
\EA\right.
\end {equation}
Such a problem admits solutions (see Steps 1-3 of the 
proof of \rth {genMeas}). The following two estimates hold
\begin {eqnarray}\label {2dim5}
\norm {u_{n}}_{L^{1}(\Gw)}+\norm {\rho_{_{\prt\Gw}}g(.,u_{n})}_{L^{1}(\Gw)}
\leq \Gth\int_{\Gw}\rho_{_{\prt\Gw}}dx+
C_{1}\norm {\gl_{n}}_{L^{1}(\Gw)}\leq C_{2},
\end {eqnarray}
where $-\Gth\leq\min\{{\rm sign} (r)g(x,r): (x,r)\in \Gw\times \mathbb 
R\}$ is nonpositive, and
\begin {eqnarray}\label {2dim6}
\norm {\nabla u_{n}}_{M^{2}(\Gw)}\leq 
C_{4}(\Gth+\norm {\gl_{n}}_{L^{1}(\Gw)})\leq C_{5}.
\end {eqnarray}
Notice that (\ref {2dim6}), which replaces (\ref{SDLM4}), follows from 
(\ref {marsest2'}). As in the proof of \rth {genMeas} there exist a 
subsequence $\{u_{n_{k}}\}$ and a function $u\in 
W^{1,q}_{0}(\Gw)$, for any $1\leq q<2$, such that $u_{n_{k}}\to u$
in $L^{1}(\Gw)$ and a.e. in $\Gw$.\smallskip

\noindent{\it Step 2 } Convergence. Because (\ref {2dim6}) holds, 
\begin {eqnarray}\label {2dim7}
\int_{\Gw\cap B_{r}(a)}\abs{\nabla u_{n}}dx\leq C_{5}{\abs{\Gw\cap 
B_{r}(a)}}^{1/2}\leq C_{5}\sqrt\gp r,\forevery r>0,\;a\in\Gw,
\end {eqnarray}
and \rcor {JN2} implies 
\begin {eqnarray}\label {2dim8}
\int_{\Gw}\exp (\gm\abs {u_{n}}/C_{5}\sqrt\gp)dx\leq 
C_{6}\abs{\Gw}\exp (\gm\abs {u_{n\,{\rm 
conv\,}\Gw}}/C_{5}\sqrt\gp)\leq C_{7},
\end {eqnarray}
since ${\norm {u_{n}}}_{L^{1}(\Gw)}$ is uniformly bounded. If we set
$$\gth_{n}(s)=\int_{\{x\in\Gw:\abs {u_{n}(x)}>s\}}dx \quad \mbox 
{and }\;\gb=\frac {\gm}{C_{5}\sqrt \gp}, 
$$
then
\begin {eqnarray}\label {2dim9}
0\leq \gth_{n}(s)\leq C_{7}e^{-\gb s},\forevery s\geq 0.
\end {eqnarray}
Let $\gw$ be any Borel subset of $\Gw$. As in \rth {genMeas}-Step 3, for any $R>0$,
we have
\begin {eqnarray*}
\int_{\gw}\abs{g(x,u_{n}} dx &\leq&
\int_{\gw}\left(\tilde g_{1}(\abs{u_{n}})-\tilde g_{2}(-\abs{u_{n}})\right)dx,\\
&\leq&(\tilde g_{1}(R)-\tilde g_{2}(-R))\abs\gw-
\int_{R}^\ity(\tilde g_{1}(s)-\tilde g_{2}(-s))d\gth_{n}(s).
\end {eqnarray*}
Therefore, as in the proof of \rth {genMeas},
\begin {eqnarray*}
\int_{R}^\ity(\tilde g_{1}(s)-\tilde g_{2}(-s))d\gth_{n}(s)
&=&(\tilde g_{1}(R)-\tilde g_{2}(-R))\gth_{n}(R)
+\int_{R}^\ity\gth_{n}(s)d(\tilde g_{1}(s)-\tilde g_{2}(-s)),\\
&\leq&(\tilde g_{1}(R)-\tilde g_{2}(-R))\gth_{n}(R)+
C_{7}\int_{R}^\ity e^{-\gb s}d(\tilde g_{1}(s)-\tilde g_{2}(-s)),\\
&\leq &\frac {C_{7}}{\gb}\int_{R}^\ity(\tilde g_{1}(s)-\tilde 
g_{2}(-s))e^{-\gb s}ds.
\end {eqnarray*}
Let $\ge>0$ arbitrary. By (\ref {2dim0}) there exists $R>0$ such that
$$\frac {C_{7}}{\gb}\int_{R}^\ity(\tilde g_{1}(s)-\tilde 
g_{2}(-s))e^{-\gb s}ds\leq\epsilon/2.
$$
Now
$$\abs\gw\leq \epsilon /2(1+\tilde g_{1}(R)-\tilde g_{2}(-R))
\Longrightarrow \int_{\gw}\abs {g(x,u_{n})}dx\leq\epsilon.
$$
We conclude by the Vitali Theorem that $g(.,u_{n_{k}})\to g(.,u)$ in 
$L^{1}(\Gw)$, and we end the proof as for \rth{genMeas}.\qeda 
\medskip

If $g(x,r)=e^{ar}$ for some $a>0$, the previous 
result does not apply for any bounded measure $\gl$. However, if the 
constant $C_{5}$ is small enough, which means that $\Gth$ and 
${\norm \gl}_{\mathfrak M^b(\Gw)}$ are, accordingly, small, the 
uniform integrability may hold. The proof of the following variant 
 is parallel to the one of \rth {VV1th}.
 
\bth{VV1'th} Let $\Gw\subset\mathbb R^{2}$ be a $C^2$ bounded domain and 
$g\in C(\Gw\times \mathbb R)$ with finite exponential orders of growth 
at plus and minus infinity. Then there exists $\gd>0$ such that for 
any $\gl\in\mathfrak M^b(\Gw)$, if ${\norm \gl}_{\mathfrak 
M^b(\Gw)}\leq \gd$, Problem (\ref {2dim}) admits a solution.
\es

The monotonicity and uniform integrability arguments imply also the 
following stability result.

\bcor{VV1cor} Let $\Gw\subset\mathbb R^{2}$ be a $C^2$ bounded 
domain and
$g\in C(\Gw\times \mathbb R)$ satisfy the $2$-dimensional weak-singularity 
assumption. Assume also that $r\mapsto g(x,r)$ is nondecreasing for 
any $x\in\Gw$. Then, for any $\gl\in \mathfrak M_{b}(\Gw)$, 
the solution $u$ of Problem (\ref {2dim}) is unique and the mapping 
$\gl\mapsto u$ is nondecreasing. Furthermore, if $\{\gl_{m}\}$ is a 
sequence of bounded measures in $\Gw$ which converges in the sense of 
measures to $\gl$, the corresponding solutions $u_{m}$ to problem 
(\ref {2dim}) converge to $u$ in $L^{1}(\Gw)$.
\es
\subsubsection {Subcritical measures}
For simplicity we shall consider only nondecreasing absorption 
nonlinearities $g\in C(\mathbb R)$ in the problem
 \begin {equation}\label {2dim+1}\left.\BA{ll}
-\Gd u+g(u)=\gl&\mbox { in }\Omega,\\[2mm]
\phantom {-\Gd u+g()\,}
u=0&\mbox { on }\prt\Omega, 
\EA\right.
\end {equation}
where $\Gw$ is a smooth bounded domain of the plane, and 
$\gl\in\mathfrak M^b(\Gw)$.
\bdef {subcrit} {\rm Let $\gl$ be a bounded measure in $\Gw$, with 
Lebesgue decomposition $\gl=\gl^{*}+\gl_{s}+\sum_{j\in J}c_{j}\gd_{x_{j}}$
where $\gl^{*}$ is the absolutely continuous part with respect to 
the $2$-dimensional Hausdorff measure, $\gl_{s}$ the singular 
non-atomic part and $\{(c_{j},x_{j})\}_{j\in J}$ the set, at most 
countable, of atoms. Let
 $g$ be a continuous nondecreasing real valued function. We say that $\gl$ 
is {\it subcritical with respect to} $g$ if 
\begin {eqnarray}\label {2dim+2}
\frac {4\gp}{a_{-}(g)}\leq c_{j}\leq \frac {4\gp}{a_{+}(g)},\forevery 
j\in J.
\end {eqnarray}
}\es 
The following result is due to Vazquez \cite {Va1}.
\bth {VV2th} Let $\gl\in\mathfrak M^b(\Gw)$. Problem (\ref {2dim+1}) 
admits a solution if and only if $\gl$ is subcritical with respect 
to $g$.
\es
The local version of the necessary condition is the following.
\bprop {VV2prop} Assume $g$ has positive and finite exponential order 
of growth at infinity, $a_{+}(g)$. Let $R>0$ and 
$\gn\in\mathfrak M^b(B_{R}(0))$ with no atom. If $c>4\gp/a_{+}(g)$ 
there exists no function $u\in L^{1}(B_{R}(0))$ such that 
$g(u)\in L^{1}(B_{R}(0))$ and
\begin {eqnarray}\label {2dim+3}
\int_{B_{R}(0)}\left(-u\Gd\gz+g(u)\gz\right)dx=c\gz 
(0)+\int_{B_{R}(0)}\gz d\gn,\forevery\gz
\in C_{c}^\infty (B_{R}(0)).
\end {eqnarray}
\es
The next result is a particular case of a remarkable relaxation phenomenon which 
occurs above the critical level $4\gp/a_{+}(g)$. We denote by $B_{R}$ 
the ball of center $0$ and radius $R$ and by 
$B^{*}_{R}=B_{R\setminus\{0\}}$.

\blemma{VV3lem} Let $g$ be a continuous nondecreasing function with
positive and finite exponential order of growth at infinity $a_{+}(g)$ and, 
for $n\in\mathbb N_{*}$, $g_{n}(r)=\min \{g(r),g(n)\}$. Let $R>0$, 
$c>c_{+}(g)=4\gp/a_{+}(g)$ 
and $b$ 
be three constants, and $\gu_{n}$ the solution of 
 \begin {equation}\label {2dim+4}\left.\BA{ll}
-\Gd \gu_{n}+g_{n}(\gu_{n})=c\gd_{0}&\mbox { in }\CD'(B_{R}),\\[2mm]
\phantom {-\Gd \gu_{n}+g_{n}()\,}
\gu_{n}=b&\mbox { on }\prt B_{R}. 
\EA\right.
\end {equation}
When $n\to\ity$, $\{\gu_{n}\}$ decreases and converges, locally uniformly in 
$B^{*}_{R}$, to the solution $\gu_{c_{+}(g)}$ of
 \begin {equation}\label {2dim+5}\left.\BA{ll}
-\Gd \gu_{c_{+}(g)}+g(\gu_{c_{+}(g)})=c_{+}(g)\gd_{0}&\mbox { in }\CD'(B_{R}),\\[2mm]
\phantom {-\Gd \gu_{c_{+}(g)}+g()\,}
\gu_{c_{+}(g)}=b&\mbox { on }\prt B_{R}. 
\EA\right.
\end {equation}

\es 
\Proof  Since $a_{+}(g_{n})=0$, we know by \rth {VV1th}, that for any $c>0$, there exists a unique solution 
$\gu_{n}$ to (\ref {2dim+4}), which is therefore a radially symmetric 
function. Because $g_{n}$ is increasing, the sequence $\{\gu_{n}\}$ 
is nonincreasing. \smallskip

\noindent {\it Step 1}  Existence of a solution to problem ({\ref 
{2dim+5}}) in  the case $c<c_{+}(g)$. By comparing 
$\gu_{n}$ with the solution  $\Psi=\Psi_{c}$  of
 \begin {equation}\label {2dim+6}\left.\BA{ll}
-\Gd \Psi=c\gd_{0}+\abs{g(0)}&\mbox { in }\CD'(B_{R}),\\[2mm]
\phantom {-\Gd }
 \Psi=\abs b&\mbox { on }\prt B_{R}, 
\EA\right.
\end {equation}
there holds $\Psi\geq \max \{0,\gu_{n}\}$. But $\Psi$ has the explicit form
\begin {eqnarray}\label {2dim+7}
\Psi (x)=\frac {c}{2\gp}\ln(1/\abs x) +K.
\end {eqnarray}
for some constant $K$. The function $\gu_{n}$ is bounded from below by 
the solution $\Phi$ of 
\begin {equation}\label {2dim+8}\left.\BA{ll}
-\Gd \Phi+g(\Phi)=0&\mbox { in }\CD'(B_{R}),\\[2mm]
\phantom {-\Gd \Phi+g()}
 \Phi=b&\mbox { on }\prt B_{R}, 
\EA\right.
\end {equation}
and $\Phi$ is a bounded function. Therefore, for $n$ large enough,
$$ g(\Phi)\leq g_{n}(\gu_{n})\leq g(\gu_{n})\leq g(\Psi)
= g\left(\frac {c}{2\gp}\ln(1/\abs x) +K\right).
$$
But
$$\int_{B_{R}}g\left(\frac {c}{2\gp}\ln(1/\abs x) +K\right)dx
\leq \int_{B_{R}}g\left(\frac {c}{2\gp}\ln(k/\abs x) \right)dx
= \frac {2k\gp}{c}\int_{\gr}^\infty g(s)e^{-4\gp s/c}ds,
$$
for some $k>0$, $\gr>0$. This last integral is finite because 
$4\gp /c>a_{+}(g)$. We conclude with Lebesgue's theorem that $\gu_{n}$ 
converges to the solution $\gu_{c}$ to ({\ref {2dim+5}}).\smallskip

\noindent {\it Step 2}  Existence of a solution to problem 
(\ref {2dim+5}) in  the case $c=c_{+}(g)$. Let $\{c_{n}\}$ be a positive 
increasing sequence converging to $c_{+}(g)$. Then the sequence 
$\{\gu_{c_{n}}\}$ is increasing. Since $\Phi\leq 
\gu_{c_{n}}\leq\Psi_{c_{+}}$ (given by (\ref {2dim7}) and (\ref 
{2dim8})), the limit $\gu^\ast$ of the $\gu_{c_{n}}$ is attained in the 
$L^{1}(B_{R})$-norm, and
$$
\Phi\leq \gu^\ast\leq\Psi_{c_{+}}.
$$
The sequence $\{g(\gu_{c_{n}})\}$ is increasing and converges 
pointwise to $g(\gu^\ast)$. Let $\eta_{1}\in C^2_{c}(\overline {B_{R}})$ be the solution of
\begin {equation}\label {2dim+8'}\left.\BA{ll}
-\Gd \eta_{1}=1&\mbox { in }B_{R},\\[2mm]
\phantom {-\Gd }
 \eta_{1}=b&\mbox { on }\prt B_{R}. 
\EA\right.
\end {equation}
Hence $\eta_{1}\geq 0$ and
\begin {eqnarray}\label {2dim+9}
\myint{B_{R}}{}\left(-\gu_{c_{n}}\Gd\eta_{1}+g(\gu_{c_{n}})\eta_{1}\right)dx
=c_{n}\eta_{1}(0)-2\gp b\eta'_{1}(R).
\end {eqnarray}
Letting $n\to\infty$ and using the Beppo-Levi Theorem implies 
$$\lim_{n\to\infty}\norm 
{\left(g(\gu_{c_{n}})-g(\gu^\ast)\right)\eta_{1}}_{L^{1}(B_{R})}=0.$$
Thus $\gu^\ast$ is the solution of (\ref {2dim+5}) with 
$c=c_{+}$.\smallskip 

\noindent {\it Step 3 }  Nonexistence of a solution to problem 
(\ref {2dim+5}) in  the case $c>c_{+}(g)$. Suppose that such a 
solution $\gu_{c}$ exists. Because of uniqueness, it is a radial 
function, and $g(\gu_{c})\in L^{1}(B_{R})$. The function
$$r\mapsto w(r)-\frac {c}{2\gp}\ln(1/r),
$$
satisfies $(rw'(r))'=rg(\gu_{c})$ on $(0,R)$. Therefore $r\mapsto 
rw'(r)$ admits a limit when $r\to 0$. If the limit were not zero, 
say $\ga$, it would imply
$$w(r)=\ga\ln (1/r)(1+\circ (1))\quad {as }\;r\to 0,
$$
and 
$$\Gd w=rg(\gu_{c})-2\gp c\gd_{0},
$$
contradiction. Thus $rw'(r)\to 0$ as $r\to 0$, and by integration, 
\begin {eqnarray}\label {2dim+10}
\gu_{c}(r)=\frac {c}{2\gp}\ln (1/r)(1+\circ (1)).
\end {eqnarray}
Then, for any $0<\gg<c$, there exists $R_{\gg}\in (0,R]$ such that
$$\gu_{c}(r)\geq\frac {\gg}{2\gp}\ln (1/r),\quad \mbox{in 
}(0,R_{\gg}].
$$
Thus $g(\gu_{c})\geq g(\gg/(2\gp)\ln(1/r))$. Put $a=2\gp/\gg$. Since 
$g(\gu_{c})\in L^{1}(B)$, it implies
$$\int_{0}^\ity g(s)e^{-2as}ds<\infty\Longrightarrow 2a\geq a_{+}(g),
$$
and finally $c\leq c_{+}(g)$, a contradiction.\smallskip

\noindent {\it Step 4 } The relaxation phenomena when $c> c_{+}(g)$. For 
any $n$ and any $\ge>0$, the solution $\gu_{n}$ of (\ref {2dim+4}) is bounded from below 
by the solution $V_{n}$ of
\begin {equation}\label {2dim+11}\left.\BA{ll}
-\Gd V_{n}+g_{n}(V_{n})=(c_{+}(g)-\epsilon)\gd_{0}&\mbox { in }\CD'(B_{R}),\\[2mm]
\phantom {-\Gd V_{n}+g_{n}()}
V_{n}=b&\mbox { on }\prt B_{R}. 
\EA\right.
\end {equation}
Let $\tilde \gu$ be the limit of the $\gu_{n}$. 
Then $\tilde \gu$ is a solution of 
\begin {equation}\label {2dim+12}\left.\BA{ll}
-\Gd \tilde \gu+g(\tilde \gu)=0&\mbox { in }B^{*}_{R},\\[2mm]
\phantom {-\Gd \tilde \gu+g()}
\tilde \gu=b&\mbox { on }\prt B_{R}. 
\EA\right.
\end {equation}
Because $V_{n}$ converges to $\gu_{c_{+}(g)-\epsilon}$, there holds 
$\tilde \gu\geq \gu_{c_{+}(g)-\epsilon}$. Letting $\epsilon\to 0$ 
finally yields to $\tilde \gu\geq \gu_{c_{+}}(g)$. Taking the same 
test function $\eta_{1}$ defined by (\ref {2dim+8'}), one obtains
\begin {eqnarray}\label {2dim+12b}
\myint{B_{R}}{}\left(-\gu_n\Gd\eta_{1}+g_{n}(\gu_n)\eta_{1}\right)dx
=c\eta_{1}(0)-2\gp b\eta'_{1}(R).
\end {eqnarray}
Using the fact that $\gu_n\leq \Psi$ (see Step 1) and Fatou's lemma,
$$\myint{B_{R}}{}g(\tilde\gu)\eta_{1}dx\leq \liminf_{n\to\ity}
\myint{B_{R}}{}g_{n}(\gu_n)\eta_{1}dx<\ity.
$$
Thus $g(\tilde\gu)\in L^{1}(B_{R})$. Since $\tilde \gu\in 
L^{1}(B_{R})$, the distribution $T=-\Gd\tilde \gu+g(\tilde \gu)$ has the point $0$ for 
support, therefore there exist real numbers $c_{p}$, ($p\in\mathbb N^m$) such that
$$T=\sum_{\abs p\leq m}c_{p}D^p\gd_{0}.
$$
Let $\gz\in C^\ity_{c}(B)$ such that
$$(-1)^{\abs p}D^p\gz (0)=c_{p},\forevery p\in\mathbb N^m,\,\abs p\leq 
m,
$$
and for $\epsilon>0$, put $\gz_{\epsilon}(x)=\gz(x/\epsilon)$. Then
\begin {eqnarray}\label {2dim+13}
\int_{B}\left(-\tilde \gu\Gd\gz_{\epsilon}+g(\tilde \gu)\gz_{\epsilon}\right)dx
=\sum_{\abs p\leq m}\frac {c_{p}^{2}}{\epsilon^{\abs p}}.
\end {eqnarray}
But
\begin {eqnarray}\label {2dim+14}
\abs{\int_{B}\tilde \gu\Gd\gz_{\epsilon}dx}=
\frac {1}{\epsilon^{2}}\abs{\int_{B}\tilde \gu\Gd\gz (x/\epsilon)dx}\leq 
\frac {C}{\epsilon^{2}}\int_{0}^{R\epsilon}\ln(1/s)sds\leq 
C'\ln(1/\epsilon).
\end {eqnarray}
Comparing (\ref {2dim+13}) and (\ref {2dim+14}) 
implies $c_{p}=0$ for any $\abs p\geq 1$, from what is infered
\begin {eqnarray}\label {2dim+15}
-\Gd\tilde \gu+g(\tilde \gu)=c_{0}\gd_{0}\quad\mbox{in }\;\CD'(B).
\end {eqnarray}
By Step 3 and the inequality $\tilde \gu\geq \gu_{c_{+}}(g)$, one has
$c_{0}=c_{+}(g)$, which ends the proof.\qeda \medskip

\noindent {\it Proof of \rprop {VV2prop}}. Assume such a $u$ exists. 
By changing $R$, we can assume that $u\in L^{1}(\prt B_{R})$ and 
that $u$ is therefore the unique integrable function with $g(u)\in 
L^{1}(B_{R})$ which satisfies 
\begin {equation}\label {2dim+16}\left.\BA{ll}
-\Gd u+g(u)=&c\gd_{0}+\gn\quad\mbox {in }\CD'(B_{R}),\\[2mm]
\phantom {-\Gd u+g()}u
&\mbox {fixed on }\prt B_{R}. 
\EA\right.
\end {equation}
Put $g_{n}(r)=\min \{g(r),g(n)\}$, and let $\gu_{n}$ be the 
solution of
\begin {equation}\label {2dim+17}\left.\BA{ll}
-\Gd \gu_{n}+g_{n}(\gu_{n})=c\gd_{0}&\quad\mbox {in }\CD'(B_{R}),\\[2mm]
\phantom {-\Gd \gu_{n}+g_{n}()}
\gu_{n}=0&\quad\mbox {on }\prt B_{R}, 
\EA\right.
\end {equation}
and $v$ the one of
\begin {equation}\label {2dim+18}\left.\BA{ll}
-\Gd v=\gn_{+}&\quad\mbox {in }\CD'(B_{R}),\\[2mm]
\phantom {-\Gd }
v=u_{+}&\quad\mbox {on }\prt B_{R}. 
\EA\right.
\end {equation}
Since $g(\gu_{n}+v)\geq g_{n}(\gu_{n}+v)\geq g_{n}(\gu_{n})$, the function 
$U_{n}=\gu_{n}+v$ is a super-solution for Problem (\ref {2dim+16}). 
Therefore $u\leq \gu_{n}+v$. Letting $n\to\ity$  and using 
\rlemma {VV3lem} yields to
\begin {eqnarray}\label {2dim+19}
u\leq \gu_{c_{+}(g)}+v.
\end {eqnarray}
Writing again 
$$u(r,\gth)=u(x)=\frac {c}{2\gp}\ln(1/\abs x)+\gw(x),
$$
then
$$-\Gd \gw=\gn-g(u)\Longrightarrow -\Gd \bar\gw(r)=\overline{(\gn-g(u))(r)},
$$
where the overlining indicates the angular average. Because the 
measure $\gn$ has no atom and $g(u)\in L^{1}(B_{R})$, 
$$\int_{0}^r\overline{(\gn-g(u))(s)}ds\to 0,\quad\mbox {as }r\to 0.
$$
Thus
$$ \overline u(r)=\frac {c}{2\gp}\ln(1/r)(1+\circ (1)).
$$
In the same way
$$ \overline \gw(r)=\circ (\ln(1/r)),
$$
and, from \rlemma {VV3lem}-Step 2,
$$ \overline \gu_{c_{+}(g)}(r)=\gu_{c_{+}(g)}(r)=\frac {c_{+}(g)}{2\gp}\ln(1/r)(1+\circ (1)).
$$
Since $c>c_{+}(g)$, this contradicts (\ref {2dim+19}).\qeda 
\medskip                               %

\noindent {\it Proof of \rth {VV2th}. }  By replacing $\gl$ by 
$\gl-g(0)$, it is always possible to assume $g(0)=0$. The measure $\gl$ 
admits the decomposition
$$\gl=\sum_{j\in J}c_{j}\gd_{x_{j}}+\gn,
$$
where $\{x_{j}\}_{j\in J}$ is the set of atoms of $\gl$, and $\gn$ is the 
sum of a measure absolutely continuous with respect to the 
$2$-dimensional Hausdorff measure and a singular measure without atom.
\smallskip

\noindent {\it Step 1 } We assume 
that $\gl$ is positive with compact support in $\Gw$, and 
$c_{j}<c_{+}(g)$ for any $j\in J$. Let $\gd>0$ as in \rth {VV1'th}, 
 $J_{1}=\{j\in J:c_{j}\geq\gd/2\}$ (with $\#(J_{1})=K$), and $j_{2}=J\setminus 
 J'$. We denote 
 $$\gl_{\gd}=\gl-\sum_{j\in J_{1}}c_{j}\gd_{x_{j}}.$$
First, there exists a finite covering 
$\{\Gw_{i}\}_{i\in I}$ of $\Gw$ (with $\#(I)=N$) such that 
$\Gw_{i}\cap\Gw_{i'}=\emptyset$ if $i\neq i'$, and
\begin {equation}\label {2dim+20}\int_{\Gw_{i}}
d\gl_{\gd}<\gd.
\end {equation}
This covering can be chosen such that any $\overline\Gw_{i}$ contains at 
most one $x_{j}$ for $j\in J_{1}$, and actually $x_{j}\in 
\Gw_{i}$, we shall write $i=i(j)$ and this correspondence is one to 
one from $J_{1}$ into $I$. For such a $x_{j}$, there exists $\gs_{j}>0$ such that 
$\overline {B_{\gs_{j}}(x_{j})}\subset \Gw_{i(j)}$, and
\begin {equation}\label {2dim+20'}
\lim_{\gs\to 0}\int_{B_{\gs}(x_{j)}}d(\gl-c_{j}\gd_{x_{j}})=0.
\end {equation}
Let $R>0$ be such that 
$\overline\Gw\subset B_{R}(x_{j})$, $\forall j\in J_{1}$. For 
$0<\gs\leq \inf_{j\in J_{1}}\gs_{j}$ and $i=i(j)$ for some $j\in 
J_{1}$, we set
$$\Gw_{i(j)}=\overline {B_{\gs}(x_{j})}\cup \Gw'_{i(j),\gs}. 
$$
By \rlemma {VV3lem}-Step 1, each of the following equations admits a 
solution $u_{j}$,
\begin {eqnarray}\label {2dim+21}
-\Gd u_{j}+\frac {1}{2N}g(u_{j})&=&c_{j}\gd_{x_{j}}\quad\mbox {in 
}\;\CD'(B_{R}(x_{j})),\\
u_{j}&=&0\quad\mbox {on }\; \prt B_{R}(x_{j}),\notag
\end {eqnarray}
for $j\in J_{1}$. Let $\Gw_{i,\gs}=\{x\in\Gw_{i}:\dist (x,\Gw^c_{i})>\gs\}$. 
If $i\in I\setminus \{i(j):j\in J_{1}\}$, we set 
$\gl_{i,\gs}=\chi_{_{\Gw_{i,\gs}}}\gl_{\gd}$, and 
if $i=i(j)$ for some $j\in J_{1}$, we put $\gl_{i,\gs}=\chi_{_{\Gw'_{i,\gs}}}\gl_{\gd}$.
By \rth {VV1'th} there exist functions $\gu_{i,\gs}$ solutions of
\begin {equation}\label {2dim+21a}\left.\BA{ll}
-\Gd \gu_{i,\gs}+\myfrac {1}{2N}g(\gu_{i,\gs})=\gl_{i,\gs}
&\mbox {in}\CD'(\Gw),\\[2mm]
\phantom {-\Gd \gu_{i,\gs}+\frac {1}{2N}g()i}
\gu_{i,\gs}=0&\mbox {on }\; \prt \Gw,
\EA\right.
\end {equation}
for $i\in I$. Furthermore the $u_{j}$ and $\gu_{i,\gs}$ are 
respectively the limit of the $u_{j,n}$ and $\gu_{i,\gs,n}$ solutions of 
\begin {equation}\label {2dim+21b}\left.\BA{ll}
-\Gd u_{j,n}+\myfrac {1}{2N}g(u_{j,n})=c_{j}\gd_{x_{j}}\ast\gr_{n}
&\mbox {in }\;\CD'(B_{R}(x_{j})),\\[2mm]
\phantom {-\Gd u_{j,n}+\myfrac {1}{2N}g()}
u_{j,n}=0&\mbox {on }\; \prt B_{R}(x_{j}),
\EA\right.
\end {equation}
and
\begin {equation}\label {2dim+21'}\left.\BA{ll}
-\Gd \gu_{i,\gs,n}+\myfrac {1}{2N}g(\gu_{i,\gs,n})=\gl_{i,\gs}\ast\gr_{n}
&\mbox {in }\;\CD'(\Gw),\\[2mm]
\phantom {-\Gd \gu_{i,\gs,n}+\myfrac {1}{2N}g()}
\gu_{i,\gs,n}=0&\mbox {on }\; \prt \Gw,
\EA\right.
\end {equation}
where $\gr_{n}$ is a positive radial and smooth convolution kernel 
with shrinking compact support. Hence, for $n$ large enough and 
$\gs$ small enough, the support of the $c_{j}\gd_{x_{j}}\ast\gr_{n}$ 
and $\gl_{i,\gs}\ast\gr_{n}$ are all disjoint and included in 
$B_{\gs/2}(x_{j)}$ or in $\Gw_{i,\gs/2}$ (if $i\notin i(J_{1})$), or in $\Gw'_{i(j),\gs/2}$.
Finally, $g(u_{j,n})\to g(u_{j})$ 
in $L^{1}(B_{R}(x_{j}))$ (easy to check from \rlemma {VV3lem}-Step 1) 
and $g(\gu_{i,\gs,n})\to g(\gu_{i,\gs})$ in 
$L^{1}(\Gw)$, as $n\to\infty$ (by the proof of \rth {VV1th}).
Put
$$U_{n}=\sum_{j\in J_{1}}u_{j,n},\quad U=\sum_{j\in J_{1}}u_{j},
$$
both quantities defined in $\overline \Gw$, and
$$V_{n}=\sum_{i\in I}\gu_{i,\gs,n},\quad V_{\gs}=\sum_{i\in I}\gu_{i,\gs}.
$$
With the same convolution kernel $\gr_{n}$, we denote by $u_{n,\gs}$  
the solution to 
\begin {equation}\label {2dim+22}\left.\BA{ll}
-\Gd u_{\gs,n}+g(u_{\gs,n})=\gl_{\gs}\ast\gr_{n}
&\mbox {in }\;\CD'(\Gw),\\[2mm]
\phantom {-\Gd u_{\gs,n}+g()}
u_{\gs,n}=0&\mbox {on }\; \prt \Gw,
\EA\right.
\end {equation}
where
$$\gl_{\gs}=\sum_{j\in 
J_{1}}c_{j}\gd_{x_{j}}+\sum_{i\in I\setminus i(J_{1})}\chi_{_{\Gw_{i,\gs}}}\gl_{\gd}
+\sum_{i\in i(J_{1})}\chi_{_{\Gw'_{i,\gs}}}\gl_{\gd}.
$$
As in the proof of \rth {VV1th}, $u_{\gs,n}\to u_{\gs}$ in $L^{1}(\Gw)$ and a.e. 
in $\Gw$, $g(u_{\gs,n})$ is bounded in $L^{1}(\Gw)$, and $g(u_{\gs,n})\to 
g(u_{\gs})$ a.e. in $\Gw$. Because
\begin {eqnarray}\label {2dim+23}
-\Gd (U_{n}+V_{\gs,n})+g(U_{n}+V_{\gs,n})&=&-\sum_{j\in J_{1}}\Gd u_{j,n}
-\sum_{i\in I}\Gd \gu_{i,\gs,n}+g(U_{n}+V_{\gs,n})\notag \\
&\geq &\sum_{j\in J_{1}}\left(-\Gd u_{j,n}+
\frac {1}{2N}g(u_{j,n})\right)\\
&\,&\qquad\qquad+\sum_{i=1}^{N}\left(-\Gd \gu_{i,\gs,n}+\frac 
{1}{2N}g(\gu_{i,\gs,n})\right)\notag \\
&=&\gl_{\gs}\ast\gr_{n}\qquad\qquad \mbox {in }\; \CD'(\Gw),\notag 
\end {eqnarray}
and $U_{n}+V_{\gs,n}\geq 0$ on $\prt\Gw$, one obtains
$$0\leq u_{\gs,n}\leq U_{n}+V_{\gs,n}.
$$
The estimate of the uniform integrability of $\{g(U_{n}+V_{\gs,n})\}$ 
derives from the following argument : Let $\gw$ be a Borel subset of $\Gw$ 
and $\gw_{i}=\Gw_{i}\cap\gw$, $i\in I$. If $i\notin i(J_{1})$ we can write
$$U_{n}+V_{\gs,n}=\gu_{i,\gs,n}+ K(x),\forevery x\in\gw_{i},
$$
and, for $\gs$ fixed small enough, the function $x\mapsto K(x)$ is bounded
uniformly with respect to $n$ and $x\in\gw_{i}$, since the distance of 
the supports of 
the $\gl_{i',\gs}\ast\gr_{n}$ ($i'\neq i$), and the 
$c_{j}\gd_{x_{j}}\ast\gr_{n}$ ($j\in J_{1}$) to $\gw_{i}$ is larger or 
equal to $\gs/2$. As in the proof of \rth {VV1th}, we set
$$\gth_{n,i}(s)=\int_{\{x\in\gw_{i}:\abs {(U_{n}+V_{n,\gs})(x)}>s\}}dx,
$$
and 
$$\gth_{n,i}(s)\leq \int_{\{x\in\gw_{i}:\gu_{i,\gs,n}+K(x))>s\}}dx.
$$
The proof of \rth {VV1th} applies : for $\epsilon>0$ fixed, there exists 
$\gd>0$, 
such that
\begin {eqnarray}\label {2dim+24}
\abs{\gw_{i}}\leq \gd \Longrightarrow 
\int_{\gw_{i}}g(U_{n}+V_{n,\gs})dx<\epsilon/2N.
\end {eqnarray}
If $i=i(j)$ we put $\gw_{i}=\gw'_{i}\cup\gw''_{i}$, where 
$\gw'_{i}\subset\Gw'_{i(j),\gs}$ and $\gw''_{i}\subset B_{\gs}(x_{j})$.
On $\gw'_{i}$ we write 
$$U_{n}+V_{n,\gs}=\gu_{i(j),\gs,n}+K'(x),
$$
and $K'(x)$ is bounded independently of $n$, thus (\ref {2dim+24}) 
holds with $\gw'_{i}$ instead of $\gw_{i}$. On $\gw''_{i}$ there holds
$$U_{n}+V_{n,\gs}=u_{i(j),n}+K''(x),
$$
with $K''(x)$ bounded independently of $n$. Thus 
$$g(U_{n}+V_{n,\gs})\leq g(u_{i(j),n}+K''(x)).
$$
Because $g(u_{i(j),n})\to g(u_{i(j)})$ in $L^{1}(B_{R}(x_{i(j))}))$ 
as $n\to\ity$, $g(u_{i(j),n}+k)\to g(u_{i(j)}+k)$ for any $k>0$. Thus 
$\{g(u_{i(j),n}+k)\}$ is uniformly integrable. The same holds with
 $\{g(u_{i(j),n}+K"(x))\chi_{_{B_{\gs}(x_{i(j)})}}\}$, if we take 
 $k\geq K''$. Finally (\ref {2dim+24}) holds with $\gw''_{i}$ instead of 
 $\gw_{i}$. Consequently, 
 \begin {eqnarray}\label {2dim+25}
\forall\gw\subset\Gw,\;\gw\mbox { Borel },\abs{\gw}\leq \gd \Longrightarrow 
\int_{\gw}g(u_{n,\gs})dx\leq \int_{\gw}g(U_{n}+V_{n,\gs})dx<\epsilon.
\end {eqnarray}
We conclude by Vitali's theorem  that $g(u_{n,\gs})\to g(u_{\gs})$ in 
$L^{1}(\Gw)$, thus $u_{\gs}$ is the solution of 
\begin {equation}\label {2dim+26}\left.\BA{ll}
-\Gd u_{\gs}+g(u_{\gs})=\gl_{\gs}
&\mbox {in }\;\CD'(\Gw),\\[2mm]
\phantom {-\Gd u_{\gs}+g()}
u_{\gs}=0&\mbox {on }\; \prt \Gw.
\EA\right.
\end {equation}
In particular there holds
$$\int_{\Gw}(u_{\gs}+g(u_{\gs})\eta_{1}dx=\int_{\Gw}\eta_{1}d\gl_{\gs},
$$
if we take
\begin {eqnarray*}
-\Gd \eta_{1}&=&1
\quad\mbox {in }\;\Gw,\\
\eta_{1}&=&0\quad\mbox {on }\; \prt \Gw.\notag
\end {eqnarray*}
Letting $\gs\to 0$, $u_{\gs}$ increases to $u$ and
\begin {eqnarray}\label {2dim+27}
\int_{\Gw}(u+g(u)\eta_{1}dx=\int_{\Gw}\eta_{1}d\gl.
\end {eqnarray}
From this integrability property it follows that $u$ is the solution of (\ref {2dim+1}).\smallskip

\noindent {\it Step 2 } The case of a general positive bounded measure. 
We perform a double truncation, replacing $\gl$ by $\gl_{n}$ 
($n\in\mathbb N_{*}$), by putting
$$\gl_{n}=\sum_{j\in J_{c_{+}}}(c_{+}(g)-n^{-1})\gd_{x_{j}}+
\chi_{_{\Gw_{n}}}\left(\sum_{j\in J\setminus 
J_{c_{+}}}c_{x_{j}}\gd_{x_{j}}\gn\right),
$$
where $J_{c_{+}}\{j\in J:c_{j}=c_{+}(g)\}$, $\gn$ is the non-atomic 
part of $\gl$, and $\Gw_{n}=\{x\in\Gw:\dist (x,\prt\Gw)>1/n\}$. If 
$u_{n}$ is the solution corresponding to (\ref {2dim+1}), with $\gl$ 
replaced by $\gl_{n}$, the sequence $\{u_{n}\}$ is increasing and 
converges to some integrable function $u$. As in Step 1, we 
conclude, by Beppo-Levi's theorem and using Equality (\ref {2dim+27}) with $\gl_{n}$
 and $u_{n}$ instead of $\gl$ and $u$,
that $g(u_{n})$ converges to $g(u)$ a.e. and in 
$L^{1}(\Gw;\gr_{_{\prt\Gw}})$ and (\ref {2dim+27}) still holds at the 
limit. Furthermore $g(u)\in L^{1}(\Gw)$ by
\rprop {genest}.\smallskip

\noindent {\it Step 3 } The case of a general bounded measure. If 
$\gl=\gl_{+}-\gl_{_{}}$ is a bounded measure, subcritical with respect to $g$, 
we have
$$\gl_{+}=\sum_{j\in J_{+}}c_{j}\gd_{x_{j}}+\gn_{+},
$$
$$-\gl_{-}=\sum_{j\in J_{-}}c'_{j}\gd_{x'_{j}}-\gn_{-},
$$
where $\{(c_{j},x_{j})_{j\in J_{+}}\}$ (resp. $\{(c_{j},x'_{j})_{j\in 
J_{-}}\}$) 
is the set of positive atoms $c_{j}>0$ (resp.  $c'_{j}<0$). We 
trunctate the measures $\gl_{+}$ and $\gl_{-}$ as in Step 2, 
introduce the coverings $\{\Gw_{i}\}$ and $\{\tilde\Gw_{i}\}$ and the
separation parameter $\gs$ and 
construct the sets of solutions $u^{+}_{j}$, $\gu^{+}_{j,\gs}$, 
$u^{-}_{j}$ and $\gu^{-}_{j,\gs}$ such that 
\begin {eqnarray*}
-\Gd u^{+}_{j}+\frac {1}{2N}g(u^{+}_{j})&=&c_{j}\gd_{x_{j}}\quad\mbox {in 
}\;\CD'(B_{R}(x_{j})),\\
u^{+}_{j}&=&0\quad\quad\;\;\,\mbox {on }\; \prt B_{R}(x_{j}),\notag
\end {eqnarray*}
\begin {eqnarray*}
-\Gd \gu^{+}_{j,\gs}+\frac {1}{2N}g(\gu^{+}_{j,\gs})&=&\gl_{+\,i,\gs}\quad\mbox {in 
}\;\CD'(\Gw),\\
\gu^{+}_{j,\gs}&=&0\quad\quad\;\;\;\mbox {on }\; \prt \Gw,\notag
\end {eqnarray*}
   \begin {eqnarray*}
-\Gd u^{-}_{j}+\frac {1}{2N}g(u^{-}_{j})&=&c'_{j}\gd_{x'_{j}}\quad\mbox {in 
}\;\CD'(B_{R}(x'_{j})),\\
u^{-}_{j}&=&0\quad\quad\;\;\,\mbox {on }\; \prt B_{R}(x'_{j}),\notag
\end {eqnarray*}
and 
\begin {eqnarray*}
-\Gd \gu^{-}_{j,\gs}+\frac {1}{2N}g(\gu^{-}_{j,\gs})&=&\gl_{-\,i,\gs}\quad\mbox {in 
}\;\CD'(\Gw),\\
\gu^{-}_{j,\gs}&=&0\quad\quad\;\;\,\mbox {on }\; \prt \Gw,\notag
\end {eqnarray*}
and their approximations $u^{+}_{j,n}$, $\gu^{+}_{j,\gs,n}$, 
$u^{-}_{j,n}$ and $\gu^{-}_{j,\gs,n}$. We also construct $u_{n}$ 
solution of (\ref {2dim+22}). As in Step 1, we obtain
$$U_{-\,n}+V_{-\,\gs,n}\leq u_{\gs,n}\leq U_{+\,n}+V_{+\,\gs,n},
$$
where $U_{+\,n}$, $V_{+\,\gs,n}$, $U_{-\,n}$, $V_{-\,\gs,n}$ are 
defined as $U_{n}$ and $V_{\gs,n}$ as in Step 1, from the 
$u^{+}_{j,n}$, $\gu^{+}_{j,\gs,n}$, 
$u^{-}_{j,n}$ and $\gu^{-}_{j,\gs,n}$. Because
$$g(U_{-\,n}+V_{-\,\gs,n})\leq g(u_{n})\leq g(U_{+\,n}+V_{+\,\gs,n}),
$$
and the sets of functions $\{g(U_{-\,n}+V_{-\,\gs,n})\}$ and 
$\{g(U_{+\,n}+V_{+\,\gs,n})\}$ are uniformly integrable from Step 1, 
the same property is shared by the set $\{g(u_{n})\}$. We conclude by the Vitali Theorem as 
in Step 1, letting $n\to\ity$ and $\gs\to 0$. The other convergences, as 
in Step 2, follow by the same uniform integrability arguments and 
the monotonity. \qeda \medskip

The general approximation-relaxation result of \cite {Va1} is the 
following.
\bth{VV3th}Let $g$ be a continuous nondecreasing  function with 
finite exponential orders of growth at plus and minus infinity, and 
$\gl\in \mathfrak M^b(\Gw)$ with decomposition 
$$\gl=\gl^{*}+\gl_{s}+\sum_{j\in J}c_{j}\gd_{x_{j}},$$
$\gl^{*}$, $\gl_{s}$ being respectively the absolute continuous part 
and the singular non-atomic part of $\gl$. Let 
$$J^+=\{j\in J:c_{j}>c_{+}(g)\},\;\mbox {and }\;J^{-}=\{j\in J:c_{j}<c_{-}(g)\},$$
$\gr_{n}$ be a regularizing kernel and $u_{n}$ the solution of 
\begin {equation}\label {2dim+29}\left.\BA{ll}
-\Gd u_{n}+g(u_{n})=\gl\ast\gr_{n}&\mbox {in }\;\CD'(\Gw),\\[2mm]
\phantom {-\Gd u_{n}+g()}
u_{n}=0&\mbox {on }\;\prt\Gw.
\EA\right.
\end {equation}
Then $u_{n}\to u$ in $L^{1}(\Gw)$ where $u$ is the solution of 
\begin {equation}\label {2dim+30}\left.\BA{ll}
-\Gd u+g(u)=\gl^r&\mbox {in }\;\CD'(\Gw),\\[2mm]
\phantom {-\Gd u+g()}
u=0&\mbox {on }\;\prt\Gw,
\EA\right.
\end {equation}
and
$$\gl^r=\gl^{*}+\gl_{s}+\!\!\!\sum_{j\in J\setminus \{J^{+}\cup J^-\}}c_{j}\gd_{x_{j}}+
\sum_{j\in J^{+}}c_{+}(g)\gd_{x_{j}}+
\sum_{j\in J^{-}}c_{-}(g)\gd_{x_{j}}.$$
\es

The proof of this results follows by a combination of the arguments in 
\rprop {VV2prop} and \rth {VV2th}.

\mysection {Semilinear equations with source term}
\subsection {The basic approach}
The equation under consideration is written under the form
\begin {equation}\label {genspring}\left.\BA{ll}
Lu=g(x,u)+\gl&\mbox{in }\;\Gw,\\[2mm]
\phantom {L}
u=0&\mbox{on }\;\prt\Gw.
\EA\right.
\end {equation}
where $\Gw$ is a domain in $\mathbb R^n$, $L$ 
an elliptic operator defined in $\Gw$, $g$ a continuous function 
defined in $\mathbb R\times\Gw$ and $\gl$ a Radon measure in $\Gw$. 
The following general result plays an important role in proving existence
of solutions in presence of supersolutions and subsolutions (see
e.g. \cite {Ni}, \cite {RRV}).
\bth{subsup} Let $\Gw\subset\mathbb R^n$ be
any domain, $L$ a second order elliptic operator defined 
by the expression (\ref {lin1}) with
locally Lipschitz continuous coefficients. We assume that for any 
compact subset $K\subset\Gw$ there exists $\ga_{K}>0$ such that  
\begin {equation}\label {locellip}
\sum_{i,j=1}^na_{ij}(x)\xi_{i}\xi_{j}\geq 
\ga_{K}\sum_{i=1}^n\xi_{i}^{2},\forevery x\in K,\;\forall 
\xi=(\xi_{1},\ldots,\xi_{n})\in \mathbb R^n.
\end {equation}
Let $h^*,\,h^\dag\in C(\Gw\times \mathbb R)$ be such 
that $r\mapsto h^*(x,r)$ is nondecreasing for
every $x\in\Gw$, and
$(x,r)\mapsto h^\dag (x,r)$ is locally Lipschitz continuous with respect to the $r$
variable, uniformly when the $x$ variable stays in a compact subset
of $\Gw$, and put $h=h^*+h^\dag$. If there exist two
$C(\Gw)\cap W^{1,2}_{loc}(\Gw)$-functions $u_{*}$ and $u^*$ 
satisfying
\begin {equation}\label {sub/super}\left.\BA{ll}
(i)\quad L u_{*}+h(x,u_{*})\geq 0&\mbox {in }\Gw,\\[2mm]
(ii)\quad  L u^{*}+h(x,u^{*})\leq 0&\mbox {in }\Gw,\\[2mm]
(iii)\quad u_{*}\leq u^{*}&\mbox {in }\Gw,
\EA\right.
\end {equation}
where the equations are understood in the weak sense, then there is
a $C^{1}(\Gw)$-function $u$ which satisfies
\begin {equation}\label {sub=super}\left.\BA{ll}
(i)\quad L u+h(x,u)= 0&\mbox {in }\Gw,\\[2mm]
\phantom {}
(ii)\quad u_{*}\leq u\leq u^{*}&\mbox { in }\Gw.
\EA\right.
\end {equation}

\es

The following construction is at the origin of most of the methods 
for solving semilinear equations with reaction source term : if $\Gw$ is a bounded domain
in $\mathbb R^n$ with a $C^{2}$ boundary and $L$ the elliptic operator 
defined by (\ref {lin1}) satisfying condition (H), if $u$ is an 
integrable function solution of (\ref {genspring}) with 
$\gl\in \mathfrak M(\Gw;\gr_{_{\prt\Gw}})$ such that $g(.,u)\in 
L^{1}(\Gw;\rho_{_{\prt\Gw}}dx)$, there holds
\begin {eqnarray}\label {conv1}
u(x)=\int_{\Gw}G_{L}^\Gw(x,y)g(y,u(y)dy+
\int_{\Gw}G_{L}^\Gw(x,y)d\gl(y),\quad a.e.\mbox { in }\Gw.
\end {eqnarray}
\bth {supth} Assume $g(x,0)=0$, $r\mapsto g(x,r)$ is nondecreasing for any $x\in\Gw$ and 
$\gl\in \mathfrak M(\Gw;\gr_{_{\prt\Gw}})$ satisfies $\mathbb 
G_{L}^\Gw(\gl)\geq 0$. If there 
exists some $v\in L^{1}(\Gw)$, $v\geq 0$ such that $g(.,v)\in 
L^{1}(\Gw;\rho_{_{\prt\Gw}}dx)$ and
\begin {eqnarray}\label {conv3}
v\geq\mathbb G_{L}^\Gw(g(.,v)+\mathbb G_{L}^\Gw(\gl),
\end {eqnarray}
there exists a positive solution $u$ to Problem (\ref {genspring}).
\es 
\Proof 
The sequence 
$\{u_{n}\}_{n\in\mathbb N}$ defined by 
$u_{0}=0$ and
\begin {eqnarray}\label {conv2}
u_{n+1}=\mathbb G_{L}^\Gw(g(.,u_{n})+\mathbb G_{L}^\Gw(\gl),\forevery 
n\in\mathbb N,
\end {eqnarray}
is nondecreasing, as soon as $\mathbb G_{L}^\Gw(g(.,u_{n})$ exists, 
but the $u_{n}$ are well defined because it is easy to prove by induction 
that there holds 
\begin {eqnarray}\label {conv4}
0=u_{0}\leq u_{1}\leq u_{2}\leq \ldots\leq u_{n}\leq v.
\end {eqnarray}
Therefore there exists $u=\lim_{n\to\infty}u_{n}$ which satisfies 
$0\leq u\leq v$, $u\in L^{1}(\Gw)$, $g(.,u)\in 
L^{1}(\Gw;\rho_{_{\prt\Gw}}dx)$ and
\begin {eqnarray}\label {conv5}
u=\mathbb G_{L}^\Gw(g(.,u)+\mathbb G_{L}^\Gw(\gl).
\end {eqnarray}
This means that $u$ is a solution of (\ref {genspring}).\qeda 

\subsection {The convexity method}
The convexity method due to Baras and Pierre \cite {BP2} applies to 
a large variety of problems which contains Problem (\ref {genspring}).

\subsubsection {The general construction}
Let $(U,\gm)$ be a positive measured space with a $\gs$-finite 
measure $\gm$. We assume that $\{K_{n}\}_{n\in\mathbb N}$ is an increasing sequence of 
measurable subsets of $U$ such that
\begin {eqnarray}\label {conv6}
\gm(K_{n})<\infty,\forevery n\in\mathbb N,\quad \bigcup_{n\geq 0}K_{n}=U.
\end {eqnarray}
We denote by $L_{+}(U)$ (resp. $L_{+}(U\times U)$) the space of 
$\gm$-measurable (resp. $\gm\otimes\gm$-measurable) functions with 
value in $[0,\infty]$. We consider a kernel $N\in L_{+}(U\times U)$ 
and a fuction $j:U\times\mathbb R\mapsto [0,\infty]$, $\gm\otimes dx$-measurable such that 
\begin {equation}\label {conv7}\left.\BA{l}
(i)_{_{}}\, r\mapsto j(x,r)\mbox { is nondecreasing, convex and l.s.c., for 
almost all }x\in U,\cr\cr 
(ii)\, j(x,0)=0, \mbox { a.e. in }U. \EA\right.
\end {equation}
The conjugate function $j^\ast$, defined by
\begin {eqnarray}\label {conv8}
j^\ast (x,r)=\sup_{\ga\in\mathbb R}(r\ga-j(x,r))
\end {eqnarray}
satisfies (\ref {conv7}). If $u\in L_{+}(U)$, 
\begin {equation}\label {conv9}j(u)(x)=\left\{\BA{l} 
_{_{\,}}j(x,u(x))\quad\mbox {if } u(x)<\infty,\cr 
\displaystyle{\lim^{^{\,}}_{r\to\infty}j(x,r)\quad\mbox {if } u(x)=\infty}.
\EA\right.
\end {equation}
If $h\in L_{+}(U)$ we set
$$\mathbb N(h)(x)=\int_{U}N(x,y)h(y)d\gm(y),
$$
and
$$\mathbb N^{*}(h)(y)=\int_{U}N(x,y)h(x)d\gm(x).
$$
Notice that these two quantities are positive or infinite. All the 
$L^p(U)$-spaces ($1\leq p\leq \infty$) are relative to the measure 
$\gm$. We denote by $L_{+}^p(U)$ their positive cones,
\begin {equation}\label {conv10}
L_{c}^\infty(U)=\{h\in L^\ity(U):\,\exists n\in \mathbb N\,\mbox 
{ s.t. }\,h(x)=0,\,\mbox { a.e. in }U\setminus K_{n} \},
\end {equation}
and $L_{c_{+}}^\infty(U)=L_{c}^\infty(U)\cap L_{+}(U)$.
Being given $f\in L_{+}(U)$, the general problem lies in finding $u\in 
L_{+}(U)$ such that
\begin {equation}\label {gpb}
u=\mathbb N(j(u))+f.
\end {equation}
Multiplying (\ref {gpb}) by $h$ and integrating over $U$ implies
\begin {eqnarray}
\int_{U}fhd\gm=\int_{U}(u-\mathbb N^\ast(j(u)))hd\gm&=&\int_{U}(uh-j(u)
\mathbb N^\ast
(h))d\gm\notag \\
&=&\int_{U}\mathbb N^\ast (h)\left(u\frac {h}
{\mathbb N^\ast (h)}-j(u)\right)d\gm \label {conv11}\\
&\leq& \int_{U}j^\ast\left(\frac {h}{\mathbb N^\ast (h)}\right)
\mathbb N^\ast (h)d\gm,\notag 
\end {eqnarray}
provided $uh\in L^{1}(U)$. Therefore a necessary condition for 
existence of a solution to Equation (\ref {gpb}) is 
\begin {eqnarray}\label {conv12}
\int_{U}fhd\gm\leq \int_{U}j^\ast\left(\frac {h}{\mathbb N^\ast (h)}\right)
\mathbb N^\ast(h)d\gm,\forevery h\in L_{c_{+}}^\infty(U)\;
\mbox { such that } uh\in L^1(U).
\end {eqnarray}
Under a very mild additional assumption, this condition is also 
sufficient. Being given $C\geq 1$ and $h\in L_{c_{+}}^\infty(U)$, we 
denote
\begin {equation}F_{C}(h)=\left\{
\BA{l}\myint{U}{}j^\ast\left(\myfrac {h}{\mathbb N^\ast (h)}\right)\mathbb N^\ast 
(h)d\gm\mbox { if }\myfrac {h}{\mathbb N^\ast (h)}<\infty\;\mbox { a.e.}\\
\qquad\qquad\qquad\qquad \mbox {and  }j^\ast\left(\myfrac {h}{\mathbb N^\ast (h)}\right)
\mathbb N^\ast (h)\in L^1(U),\\
+\infty \qquad\qquad\mbox { if not.}
\EA\right.\label {conv13}\end {equation}
with the convention $h(x)/\mathbb N^\ast (h)(x)=0$ if $h(x)=\mathbb N^\ast (h)(x)=0$. 
If $C=1$, $F_{1}=F$. We put
$$X=\{h\in L_{c}^\ity(U):F(h)<\ity\},
$$
and
$$\hat X=\{h\in L_{c}^\ity(U):\exists C>1\mbox { s.t. }F_{C}(h)<\ity\}.
$$
In the sequel we adopt the convention $uh(x)=0$ if $h(x)=0$ and 
$u(x)=\infty$. The main existence result is as follows.
\bth {BPth} Let $f\in L_{+}(U)$. The following problem
\begin {equation}\left.
\BA{l} (i)\quad \,
u\in L_{+}(U),\quad u(x)=\mathbb N(j(u))(x)+f(x)\quad 
\gm\mbox{-a.e. in }U,\\[2mm]
(ii)\quad 
uh\in L^{1}(U),\forevery h\in\hat X,
\EA\right.\label {conv14}\end {equation}
admits a solution if and only if
\begin {eqnarray}\label {conv15}
\int_{U}fh d\gm\leq F(h),\forevery h\in\hat X.
\end {eqnarray}
\es
{\it Scheme of the proof. } For $\gg\in (0,1)$ we introduce the 
sequence $\{u_{n}\}$ defined by $u_{0}=\gg f$ and
\begin {eqnarray}\label {conv16}
u_{n+1}=\gg\left(\mathbb N(j(u_{n}))+f\right),\forevery n\in \mathbb N.
\end {eqnarray}
{\it Step 1 } We claim that
\begin {eqnarray}\label {conv17}
\int_{U}u_{n+1}h d\gm\leq \frac {\gg}{1-\gg}F(h),\forevery h\in\hat X.
\end {eqnarray}
For $1<C<1/\gg$, and $h\in\hat X$ such that $F_{C}(h)<\infty$, we suppose 
that there exists some $\psi\in L_{+}(U)$ such that
\begin {eqnarray}\label {conv18}
\psi(x)=\max\left\{\frac {1}{C}j'(u_{n})(x)\mathbb N^\ast(\psi)(x),h(x)\right\}.
\end {eqnarray}
It follows from (\ref {conv16}),
\begin {eqnarray}\label {conv19}
\int_{U}u_{n+1}\psi d\gm=\gg\int_{U}j(u_{n})\mathbb N^\ast(\psi)d\gm+\gg\int_{U}f\psi d\gm.
\end {eqnarray}
By assumption (\ref {conv15})
\begin {eqnarray*}
\int_{U}f\psi d\gm\leq F_{C}(\psi)&\leq &
\int_{U}j^{*}\left(\myfrac{\max\{j'(u_{n})\mathbb N^{*}(\psi),Ch\}}
{\mathbb N^{*}(\psi)}\right)\mathbb N^{*}(\psi)d\gm\\
&\leq &\int_{U}\max\left\{ 
j^*\left(j'(u_{n})\mathbb N^{*}(\psi)\right),
j^*\left(\myfrac {Ch}{\mathbb N^{*}(\psi)}\right)\mathbb N^{*}(\psi)\right\}d\gm.
\end {eqnarray*}
Since $\psi\geq h$, one has $\mathbb N^{*}(\psi)\geq \mathbb N^{*}(h)$. By convexity 
$j^{*}(\ga r)\leq \ga j^{*}( r)$, $\forall r\geq 0$, $\forall \ga\in 
[0,1]$, therefore
$$j^*\left(\myfrac {Ch}{\mathbb N^{*}(\psi)}\right)\mathbb N^{*}(\psi)\leq 
j^*\left(\myfrac {Ch}{\mathbb N^{*}(h)}\right)\mathbb N^{*}(h).
$$
By definition $j^*(j'(u_{n}))=u_{n}j'(u_{n})-j(u_{n})$. Thus, returning 
to (\ref {conv19}) implies
\begin {eqnarray*}
\int_{U}u_{n+1}\psi d\gm\leq \gg\int_{U}j(u_{n})\mathbb N^\ast(\psi)d\gm+
\gg\int_{U}\left(u_{n}j'(u_{n})-j(u_{n})\right)\mathbb N^\ast(\psi)d\gm
+\gg F_{C}(h).
\end {eqnarray*}
By combining this inequality with the definition of $\psi$, one derives
\begin {eqnarray*}
\int_{U}u_{n+1}\psi d\gm\leq \gg\int_{U}u_{n}\psi d\gm+\gg F_{C}(h).
\end {eqnarray*}
Because $u_{n+1}\geq u_{n} $ and $\psi\geq h$, we obtain
\begin {eqnarray*}
\int_{U}u_{n+1}h d\gm\leq\int_{U}u_{n+1}\psi d\gm\leq 
\frac {\gg}{1-\gg C}F_{C}(h).
\end {eqnarray*}
Letting $C\to 1$, (\ref {conv17}) follows. \medskip

\noindent {\it Step 2 } Convergence. Letting $n\to\infty$, $u_{n}$ 
increases and converges to some $u_{\gg}$ which satisfies
\begin {equation}\left.
\BA{l} (i)\quad\,
u_{\gg}\in L_{+}(U),\quad u_{\gg}=\gg\left(\mathbb N(j(u_{\gg}))+f\right)\quad 
\mbox{ in }U,\\[2mm]
(ii)\quad 
u_{\gg}h\in L^{1}(U),\forevery h\in\hat X,
\EA\right.\label {conv20}\end {equation}
This implies in particular
$$\int_{U}u_{\gg}h d\gm=\gg\int_{U}j(u_{\gg})\mathbb N^{*}(h)d\gm+
\gg\int_{U}fh d\gm,\forevery h\in \hat X.
$$
Let $C>1$ such that $F_{C}(h)<\infty$, then
$$\gg\int_{U}\left(u_{\gg}V\myfrac {h}
{\mathbb N^{*}(h)}-j(u_{\gg})\right)\mathbb N^{*}(h)d\gm
=(\gg C-1)\int_{U}u_{\gg}h d\gm+\gg\int_{U}fh d\gm,
$$
and consequently
\begin {eqnarray}\label {conv21}
\int_{U}u_{\gg}h d\gm\leq\frac {\gg}{\gg C-1}F_{C}(h).
\end {eqnarray}
Since the correspondence $\gg\mapsto u_{\gg}$ is increasing and, for almost 
all $x\in U$, $r\mapsto j(x,r)$ is continuous on the left, we can let 
$\gg\to 1$ in (\ref {conv21}) and (\ref {conv20})-(i) and deduce that 
the function $u=\lim_{\gg\to 1}u_{\gg}$ is a solution to problem 
(\ref {conv14}).\medskip

\noindent {\it Step 3 } Justification. The difficulties in the above 
proof are of two kinds :\smallskip

\noindent (1) It is not clear that $u_{n}<\infty$ on a set of 
positive measure. It is even not known if $u_{0}=\gg f$ satisfies 
$j(u_{0})<\infty$ a.e. in $U$. To go arround this difficulty we 
approximate $j(u_{n})$, formally equal to 
$u_{n}j'(u_{n})-j^{*}(j'(u_{n}))$, by $u_{n}\gb_{n}-j^{*}(\gb_{n})$ 
where the $\{\gb_{n}\}$ is an increasing sequence of regular enough 
fonctions converging to $j'(u_{n})$.\smallskip

\noindent (2) The existence of $\psi\in\hat X$ has to be proven.\smallskip

\noindent The full construction, which is extremely technical, is 
performed in \cite {BP2}.\qeda \medskip

In the presence of a subsolution $v$ to Problem (\ref{conv14}) it is 
possible to relax the assumption on the sign of $f$ and to produce a
signed solution $u$. More precisely, we assume that there exists a 
measurable function $v$ such that 
\begin {equation}\left.
\BA{l} (i)\quad\, 
v\in L^{1}(K_{n})\;\mbox { and}\quad N(.,.)j(v)(.)\in 
L^{1}(K_{n}\times U),\forevery n\in \mathbb N,\\[2mm]
(ii)\quad 
v(x)\leq \mathbb N(j(v))(x)+f(x)\quad \gm\mbox{-a.e. in }U,
\EA\right.\label {conv14'}\end {equation}
If $j:U\times \mathbb R\mapsto (-\ity,\ity]$ is a measurable function which 
satisfies (\ref {conv7}), we introduce $j^{*}_{v}$ and $\hat X_{v}$ :
$$j^{*}_{v}(x,r)=\sup_{\ga\geq v(x)}(r\ga-j(x,\ga),
$$
and
$$\hat X_{v}=\left\{h\in L_{c}^\ity(U):\exists C>1\mbox { s.t. }
j^{*}_{v}\left(\frac {Ch}{\mathbb N^{*}(h)}\right)\mathbb N^{*}(h)\in L^1(U)\right\}.
$$
\bcor {BPcor} There exists a measurable function $u:U\mapsto (-\ity,\ity]$ 
satisfying 
\begin {equation}\left.
\BA{l} (i)\quad\, 
u\geq v,\quad u(x)=\mathbb N(j(u))(x)+f(x)\quad 
\gm\mbox{-a.e. in }U,\\[2mm]
(ii)\quad 
uh\in L^{1}(U),\forevery h\in\hat X_{v},
\EA\right.\label {conv14''}\end {equation}
if and only if 
\begin {eqnarray}\label {conv15'}
\int_{U}fh d\gm\leq \int_{U}j^{*}_{v}\left(\frac 
{Ch}{\mathbb N^{*}(h)}\right)\mathbb N^{*}(h)d\gm,
\forevery h\in\hat X_{v}.
\end {eqnarray}
\es 
\Proof Put $w=u-v$ and define $\tilde j$ by \smallskip 

\indent $\qquad\tilde j(x,r)=0,\forevery (x,r)\in\Gw\times\mathbb 
R_{-}$,\smallskip 

\indent $\qquad\tilde j(x,r)=j(x,r+v(x))-j(x,v(x))\quad \mbox {if }j(x,v(x))<\infty\;\mbox {and } r>0$,
\smallskip

\indent $\qquad\tilde j(x,r)=\infty\quad \mbox {if }j(x,v(x))=\infty\;\mbox 
{and } r>0$.\smallskip

\noindent Thus $\tilde j$ takes nonnegative values and satisfies (\ref {conv7}). 
Moreover (\ref {conv14''}) is equivalent to
\begin {equation}\left.
\BA{l} (i)\quad\,
w\in L_{+}(U),\quad w=\mathbb N(\tilde j(w))+f+\mathbb N(j(v))-v\quad 
\gm\mbox{-a.e. in }U,\\[2mm]
(ii)\quad
wh\in L^{1}(U),\forevery h\in\hat X_{v}.
\EA\right.\label {conv15''}\end {equation}
Since

\indent $\qquad\tilde j^{*}(x,r)=j_{v}(x,r)+j(v(x))-rv(x)\quad$ if 
$\;j(x,v(x))<\infty$,\smallskip 

\indent $\qquad\tilde j^{*}(x,r)=0\quad$ if 
$\;j(v)=\infty$,\smallskip 

\noindent for any $h\in L_{c}^\infty(U)$, there holds
\begin {equation}\label {16'}
\tilde j^{*}\left(\frac {Ch}{\mathbb N^{*}(h)}\right)\mathbb N^{*}(h)=
j_{v}^{*}\left(\frac {Ch}{\mathbb N^{*}(h)}\right)\mathbb N^{*}(h)
+j(v)\mathbb N^{*}(h)-Chv,
\end {equation}
$\gm$-a.e. on $\{x\in U: j(v)(x)<\infty\}$. Therefore 
\begin {equation}\label {17'}
\tilde j^{*}\left(\frac {Ch}{\mathbb N^{*}(h)}\right)\mathbb N^{*}(h)\in 
L^{1}(U)\Longleftrightarrow j^{*}\left(\frac {Ch}{\mathbb N^{*}(h)}\right)
\mathbb N^{*}(h)\in 
L^{1}(U).
\end {equation}
The proof of \rcor {BPcor} follows from \rth {BPth} applied to 
Problem (\ref {conv15''}).\qeda 
\subsubsection {Application to elliptic semilinear equations}
Let $\Gw$ be a bounded domain in $\mathbb R^n$ with a $C^{2}$ 
boundary, $L$ an elliptic operator defined by (\ref {lin1}) satisfying 
(H) and $j:\Gw\times \mathbb R\mapsto [0,\infty]$ a 
measurable function (for the $(n+1)$-dimensional Hausdorff measure) such that
$j(x,r)= 0$, for almost all $x\in \Gw$ and every $r\leq 0$. The 
function $r\mapsto j(x,r)$ is also assumed to be convex, 
nondecreasing and l.s.c., thus it fulfills assumption (\ref {conv7}). 
If $\gl\in \mathfrak M_{+}(\Gw;\gr_{_{\prt\Gw}})$, 
$f=\mathbb G_{L}^{\Gw}(\gl)\in L^{1}(\Gw)$. We denote by
\begin {eqnarray}\label {conv21'}
Y(L)=\{\xi\in C_{c}^{1,L}(\overline\Gw)\}: L^*\xi\in L_{c}^\ity(\Gw)\cap 
L_{+}(\Gw),
\end {eqnarray}
the space $C^{1}$-functions $\xi$ vanishing on $\prt\Gw$ such that 
$L^{*}\xi $ has compact support and is essentially bounded. Notice 
that the elements of $Y(L)$ are nonnegative by the maximum principle.
\bth{BPth2} Assume there exist some $C>1$ and $\xi_{0}\in Y(L)$, 
$\xi\neq 
0$, such that 
\begin {eqnarray}\label {conv22}
j^*\left(C\frac {L^{*}\xi_{0}}{\xi_{0}}\right)\in L^{1}(\Gw).
\end {eqnarray}
If $\gl\in \mathfrak M_{+}(\Gw;\gr_{_{\prt\Gw}})$, there exists at 
least one $u\in L_{loc}^{1}(\Gw)$ such that 
$\mathbb G_{L}^\Gw(j(u))\in L_{loc}^{1}(\Gw)$ and
\begin {eqnarray}\label {conv23}
u=\mathbb G_{L}^\Gw(j(u))+\mathbb G_{L}^{\Gw}(\gl)\in 
L^{1}(\Gw),\mbox { a.e. in }\Gw,
\end {eqnarray}
if and only if
\begin {eqnarray}\label {conv24}
\int_{\Gw}\xi d\gl\leq \int_{\Gw}j^*\left(\frac 
{L^{*}\xi}{\xi}\right),\forevery \xi\in Y(L).
\end {eqnarray}
Moreover, if $\gm\geq 0$, there exists at least one positive solution.
\es
\Proof We put $\gm=dx$, the $n$-dimensional Hausdorff measure, and
$$N(x,y)=G_{L}^{\Gw}(x,y),\forevery (x,y)\in\Gw\times\Gw,\,x\neq y.
$$
Let $v$ be defined by
$$v(x)=\left\{\BA{l}0\quad\quad \mbox {if }\;f(x)\geq 0,\\[2mm]
f(x)\quad \mbox {if }\;f(x)\leq 0.
\EA\right.$$
Thus $v\in L^1(\Gw)$, $N^*(j(v))\equiv 0$ and (\ref {conv14'}) holds. 
Furthermore $j_{v}^*=j^*$ on $[0,\infty)$, $\hat X_{v}=\hat 
X\neq\{0\}$, because of (\ref {conv22}). If it exists, any solution $u$ of (\ref 
{conv23}) satisfies $u\geq v$, thus this problem is equivalent to 
$$\left\{\BA{l}u\geq v,\quad u=\mathbb N(j(u))+f,\\[2mm]
u\in L^1_{loc}(\Gw).
\EA\right.$$
If $\xi\in Y(L)$, we put $h=L^*\xi$, which means equivalently
 $$\xi=\mathbb G_{L^*}^\Gw(h)=\mathbb N^*(h).$$
By \rcor {BPcor} there exists a measurable function $u$ which 
satisfies $u=\mathbb N(j(u))+f$, $u\geq v$ and $uh\in L^1(\Gw)$, for 
every $h\in\hat X$. By (\ref {conv22}), $uL^*\xi_{0}\in L^1(\Gw)$, 
then $u(x_{0})$ is finite at least for one $x_{0}\in\Gw$, thus $N(x_{0},.)j(u)(.)\in L^1(\Gw)$, 
by the equation. For any compact 
$K\subset\Gw$ and any compact neighborhood 
$K_{0}$ of $K\cup\{x_{0}\}$, there exists a constant $C$ such that 
$$G_{L^*}^\Gw(x,y)\leq CG_{L^*}^\Gw(x_{0},y),\forevery (x,y)\in 
K\times (\Gw\setminus K_{0}). 
$$
Therefore
$$\int_{K}\int_{\Gw\setminus K_{0}}N(x,y)j(y,u(y))dydx
\leq C\abs K\int_{\Gw}N(x_{0},y)j(y,u(y))dy<\ity,
$$
from which it is infered that $\mathbb N(j(u))\in L_{loc}^1(\Gw)$, since $K$ is 
arbitrary. Furthermore $u\in L_{loc}^1(\Gw)$, from the equation.\qeda \medskip

When $j(x,r)=r_{+}^q$, for some $q>1$, the result is as follows.

\bcor {BPcor2} Let $q>1$, $\gl\in\mathfrak M(\Gw;\gr_{_{\prt\Gw}})$ and $\gs>0$. 
Then there exists a function 
$u\in L^1_{loc}(\Gw)$ such that 
$\mathbb G_{L}^\Gw(u_{+}^q)\in L^1_{loc}(\Gw)$ satisfying 
\begin {equation}\label {conv25}\left.\BA{ll}
Lu=u_{+}^q+\gs\gl&\mbox{in }\;\Gw,\\[2mm]
\phantom {L}
u=0&\mbox{on }\;\prt\Gw,
\EA\right.
\end {equation}
 if and only if
\begin {eqnarray}\label {conv26}
\gs\int_{\Gw}\xi d\gl\leq\myfrac {\gg-1}{\gg^{\gg'}}
\int_{\Gw}\frac {(L^*\xi)^{q'}}{\xi^{q'-1}},\forevery \xi\in 
Y(L),
\end {eqnarray}
where $q'=q/(q-1)$. Furthermore $u$ is nonnegative if 
$\mathbb G_{L}^\Gw(\gl)$ is so.
\es

Condition (\ref {conv26}) has two meanings : the first one is that the 
positive part of $\gl$ should not be too large, whatever is $q>1$, the 
second is that if $q$ is above some critical value, measure $\gl$ 
should not be too concentrated. This concentration is expressed in 
terms of Bessel capacities as for equations with absorption. If we 
assume for example that $\gl=\gl_{+}-\gl_{-}$ is a $L^p$-function, 
there holds,
\bcor {BPcor3} Let $q>1$, $\gl=\gl_{+}-\gl_{-}\in L^p(\Gw)$ Then there 
exists a function 
$u\in L^1_{loc}(\Gw)$ solution of Problem (\ref {conv25}) for $\gs>0$, small 
enough, if\smallskip

\noindent (i) $n=1,\,2\quad\mbox {and }1<q,\quad\mbox {or }\,n\geq 3 \quad\mbox {and }\, 
1<q<n/(n-2)$,\smallskip

\noindent or\smallskip

\noindent (ii) $n\geq 3,\quad q>n/(n-2)\quad\mbox {and }\,\gl_{+}\in 
L^p(\Gw)\quad\mbox {with }\,p\geq n(q-1)/2q$,\smallskip

\noindent or\smallskip

\noindent (ii) $n\geq 3,\quad q=n/(n-2)\quad\mbox {and }\,\gl_{+}\in 
L^p(\Gw)\quad\mbox {with }\,p>1$.
\es
\Proof Only condition (\ref {conv24}) is to be checked. If $\xi\in 
Y(L)$, 
we define $w$ by
\begin {equation}\label {conv27}
L^*\xi=w^{1/q'}\xi^{1/q}.
\end {equation}
If $\myfrac {1}{p}+\myfrac {1}{\gg}\leq 1$, there holds
\begin {equation}\label {conv28}\int_{\Gw}\xi d\gl\leq \int_{\Gw}\xi d\gl_{+}\leq C{\norm {\gl_{+}}}_{L^p}
{\norm {\xi}}_{L^\gg}.
\end {equation}
If we assume
\begin {equation}\label {conv29}\frac {1}{s}\leq \frac {1}{\gg}+\frac {2}{n},\quad {or }\;
\frac {1}{s}<\frac {2}{n},\quad {if }\; \gg=\ity,
\end {equation}
it follows, by (\ref {conv27}) and the Gagliardo and Sobolev inequalities,
$${\norm {\xi}}_{L^\gg}\leq C{\norm {\xi}}_{W^{2,s}}
\leq C{\norm {\Gd\xi}}_{L^{s}}\leq C\left(\int_{\Gw}w^{s/q'}\xi^{s/q}dx\right)^{1/s}.
$$
for any $1<s<\ity$. Furthermore, if 
\begin {equation}\label {conv30}s\leq q',
\end {equation}
one gets
$$\norm\xi_{L^\gg}\leq 
\left(\int_{\Gw}wdx\right)^{1/q'}\left(\int_{\Gw}\xi^{sq'/q(q'-s)}\right)^{(q'-s)/q's}.
$$
If 
\begin {equation}\label {conv31}\gg\geq sq'/q(q'-s),
\end {equation} 
we derive
$$\norm\xi_{L^\gg}\leq C\int_{\Gw}wdx.
$$
By combining this inequality with (\ref {conv28}), it is infered
$$\int_{\Gw}\xi d\gl \leq C\int_{\Gw}wdx.
$$
In order to get (\ref{conv29}), (\ref {conv30}), (\ref {conv31}), we choose $\gg=\infty$, $s<n/2$ if $n=1,\,2$ 
or $n\geq 3$. We take $\gg<\infty$ and $s$ such that equality holds 
in (\ref{conv29}), if $n\geq 3$, $q>n/(n-2)$, and $p\geq n(q-1)/2q$.\qeda 
\medskip

The next result expresses the condition of concentration which allows 
a measure to be admissible in Problem (\ref {conv25}).
\bprop {PBprop} Let $\gl_{\gs}=\gs\gl$ be a positive measure with compact 
support satisfying (\ref {conv26}). Then there exists 
$k=k(q,n,\gl_{\gs})$ such that
\begin {equation}\label {conv32}
\gl_{\gs}(K)\leq kC_{2,q'}(K),\forevery K\mbox {compact 
},\;K\subset\Gw.
\end {equation} 
\es 
\Proof We first notice that (\ref {conv26}) implies
\begin {equation}\label {conv33}
\int_{\Gw}v d\gl_{\gs}\leq \frac 
{q-1}{q^{q'}}\int_{\Gw}\frac{{\abs{L^* v}}^{q'}}{v^{q'-1}}dx,\forevery 
v\in C_{c}^{\infty}(\Gw),\; v\geq 0.
\end {equation} 
Indeed, if $v\geq 0$ belongs to $C_{c}^{\infty}(\Gw)$, we apply  (\ref 
{conv26}) to $\xi=\mathbb G_{L^*}^\Gw(\abs {L^* v})$ which is larger 
than $v$ by the maximum principle. We replace $v$ by $v^{2q'}$ in 
(\ref {conv33}). Since
\begin {eqnarray*}
L^* v^{2q'}=-2q'v^{2q'-1}\left[\sum_{i,j=1}^n \frac {\prt}{\prt x_{j}}\left(a_{ij}\frac {\prt 
v}{\prt x_{i}}\right)
+\sum_{i=1}^nc_{i}\frac {\prt v}{\prt x_{i}}
-\sum_{i=1}^n\frac {\prt}{\prt x_{i}}\left(b_{i}v\right)\right]\\
  -2q'(2q'-1)v^{2q'-2}\sum_{i,j=1}^na_{ij}
  \frac {\prt v}{\prt x_{j}}\frac {\prt v}{\prt x_{i}}
  +\left((2q'-1)v^{2q'}\frac {\prt b_{i}}{\prt x_{i}}+d\right)v^{2q'}.
\end {eqnarray*}
Then
$$\int_{\Gw}\frac {|L^* v^{2q'}|^{q'}}{v^{2q'(q'-1)}}dx
\leq C{\norm v}^{q'}_{L^\infty}{\norm v}^{q'}_{W^{2,q'}}+{\norm{\nabla 
v}}^{2q'}_{L^{2q'}},
$$
and finally  
\begin {eqnarray}\label {conv34}
\int_{\Gw}v^{2q'} d\gl_{\gs}\leq 
C{\norm v}^{q'}_{L^\infty}{\norm v}^{q'}_{W^{2,q'}},
\end {eqnarray}
by the Gagliardo-Nirenberg inequality. If $K\subset\Gw$ is compact, 
there exists a sequence $\{v_{k}\}\subset C_{c}^{\infty}(\Gw)$ such 
that $0\leq v_{k}\leq 1$, $v_{k}\equiv 1$ in a neighborhood of $K$ 
and ${\norm v_{k}}^{q'}_{W^{2,q'}}\to C_{2,q'}(K)$ when 
$k\to\infty$. Therefore (\ref {conv34}) implies (\ref 
{conv32}).\qeda \medskip

\noindent \Remark In the particular case where $K=B_{r}(x_{0})$ (for 
$0<r<\gr_{_{\prt\Gw}}(x_{0})$, the measure $\gl_{\gs}$ satisfies
\begin {equation}\label {conv36}
\gl_{\gs}(B_{r}(x_{0}))\leq C\left\{\BA {l}r^{n-2q'}\;\;\,\quad\quad\quad\quad\quad {if 
}\;q>n/(n-2),\\
\left(\ln (1+1/r)\right)^{1-q'}\quad {if 
}\;q=n/(n-2).\EA\right.
\end {equation}
Estimate (\ref {conv32}) can be understood in saying that the 
measure $\gl_{\gs}$ is Lipschitz continuous with respect to the 
capacity $C_{2,q'}$, although it must be noticed that a capacity is 
only an outer measure, not a regular one.\medskip

     Later on, Adams and Pierre \cite {AP} proved a series of 
remarquable equivalent properties linking estimates of type (\ref {conv32}) 
and Bessel capacities. 

\bth {APth} Let $n>2$, $p>1$ and $\gl$ be a nonnegative measure with 
compact support in $\Gw$. Then the following conditions are 
equivalent : \smallskip

\noindent (i) There exists $k_{1}>0$ such that for all compact 
subset $K\subset\Gw$, 
\begin {equation}\label {conv37}
\gl (K)\leq k_{1}C_{2,p}(K).
\end {equation}

\noindent (ii) There exists $k_{2}>0$ such that  
\begin {equation}\label {conv38}
\int_{\Gw}\xi^pd\gl\leq k_{2}\int_{\Gw}|\Gd \xi|^pdx,\forevery \xi\in 
Y(-\Gd).
\end {equation}

\noindent (iii) There exists $k_{3}>0$ such that  
\begin {equation}\label {conv39}
\int_{\Gw}\xi d\gl\leq k_{3}\int_{\Gw}|\Gd \xi|^p\xi^{1-p}dx,\forevery \xi\in 
Y(-\Gd).
\end {equation}

\noindent (iv) There exists $k_{4}>0$ such that  
\begin {equation}\label {conv40}
\int_{\Gw}\xi d\gl\leq k_{4}\int_{\Gw}|L^* \xi|^p\xi^{1-p}dx,\forevery \xi\in 
Y(L^*).
\end {equation}
\es
Their proof is performed with an elliptic operator with $C^1$ 
coefficients, but it can be adapted to an operator satisfying 
condition (H). It heavily relies on fine properties of real valued 
functions in connection with the Hardy-Littlewood maximal function 
and the Muckenhoupt weights. \medskip 

Usually a positive measure $\gl\in W^{-2,q}(\Gw)$ does not 
satisfies (\ref {conv37}), but only
\begin {equation} \label {conv37'}
\gl (G)\leq {\norm\gl}_{W^{-2,q}(\Gw)}C^{1/q}_{2,q'}(G),\forevery 
G\subset\Gw,\, G\mbox { compact}.
\end {equation}
However, the capacitary measure $\gl_{K}$ of a 
compact subset of $K\subset\Gw$ does 
verify it. This measure
 is the unique extremal for the dual definition of the capacity 
of $K$ given by (\ref {Bessel4'}). It is concentrated on $K$ and has 
the property that 
\begin {equation} \label {conv37''}
\gl_{K}(K)=C_{2,q'}(K),
\end {equation}
(see \cite [Th 2.2.7]{AH}). Moreover
\begin {equation} \label {conv37'''}
G_{1}\ast\gl_{K}\in L^{q}(\mathbb R^n)\,\mbox { and }
G_{1}\ast(G_{1}\ast\gl_{K})^{q-1}\in L^{\infty}(\mathbb R^n).
\end {equation}
where $G_{1}$ denotes the Bessel kernel of order $1$ defined by (\ref 
{Bessel}). The following result is proven in \cite {Pie}. 
\bprop {Pieth} Let $K\subset\Gw$ be compact subset with 
$C_{2,q'}(K)>0$ and $\gl_{K}$ the 
capacitary measure of $K$. Then there exists $k=k(n,q)$ such that
\begin {equation} \label {conv37''''}
\int_{\Gw}\xi d\gl_{K}\leq k
{\norm {G_{1}\ast(G_{1}\ast\gl_{K})^{q-1}}}_{L^{\infty}(\mathbb R^n)}
\int_{\Gw}\abs {\Gd\xi}^{q'}\xi^{1-q'}dx,\forevery\xi\in Y(-\Gd).
\end {equation}
\es
Hence, by \rcor {BPcor2}, Problem \ref {conv23} is solvable 
for any capacitary measure $\gl=\gl_{K}$, for $0<\gs\leq \gs_{0}$ 
for some $\gs_{0}>0$. 
Furthermore, since it is proven in \cite [Th. 3.1]{KM} that there 
exists a constant $k_{n,q}>0$ such that
$${\norm {G_{1}\ast(G_{1}\ast\gl_{K})^{q-1}}}_{L^{\infty}(\mathbb R^n)}
\leq k_{n,q}\forevery K\subset \Gw,\;K \mbox { compact},$$
it follows that $\gs_{0}=\gs_{0}(n,q)$.
\subsection {Semilinear equations with power source terms}
In this section we develop a direct methods for constructing explicit super 
solutions in order to apply \rth {supth}. We assume that $\Gw$ is a 
bounded open subset with a $C^2$ boundary and that $L$ 
defined by (\ref {lin1}) satisfies (H). 
\bth{KVth} Let $q>0$, $\gl\in \mathfrak M_{+}(\Gw;\gr_{_{\prt\Gw}})$. If 
there exists some $C_{0}>0$ such that 
\begin {eqnarray}\label {pst1}
\mathbb G_{L}^\Gw\left(\left(\mathbb G_{L}^\Gw(\gl)\right)^q\right)
\leq C_{0}\mathbb G_{L}^\Gw(\gl),\quad\mbox{ a.e. in } \Gw,
\end {eqnarray}
then problem
\begin {equation}\label {pst2}\left.\BA{ll}
Lu=\abs {u}^{q-1}u+\gs\gl&\mbox{in }\;\Gw,\\[2mm]
\phantom {L}
u=0&\mbox{on }\;\prt\Gw,
\EA\right.
\end {equation}
admits a positive solution $u\in L^{1}(\Gw)\cap 
L^q(\Gw;\gr_{_{\prt\Gw}}dx)$,\smallskip

\noindent (i) if  $0<\gs\leq\gs_{0}=\gs_{0}(q,C_{0})$, when $q>1$, \smallskip

\noindent (ii) for any $\gs>0$ when $0<q\leq 1$.
\es
\Proof Put $w=\gth \mathbb G_{L}^\Gw(\gs\gl)$, for some parameters 
$\gth,\,\gs>0>0$. Then, under condition (\ref {pst1}), 
$$\mathbb G_{L}^\Gw(w^q+\gs\gl)\leq (C_{0}\gth^q\gs^q+\gs)
\mathbb G_{L}^\Gw(\gl).
$$
Therefore
\begin {eqnarray}\label {pst3}
w\geq \mathbb G_{L}^\Gw(w^q)+ \mathbb G_{L}^\Gw(\gs\gl), 
\end {eqnarray}
as soon as
\begin {eqnarray}\label {pst4}
C_{0}\gth^q\gs^{q-1}+1\leq \gth.
\end {eqnarray}
If $q>1$ this is equivalent to
$$\gs\leq\max_{\gth>0}\left(\frac 
{\gth-1}{C_{0}\gth}\right)^{1/(q-1)}=\frac 
{1}{q(C_{0}q)^{1/(q-1)}},
$$
and we get (i) by \rth {supth}. If $0<q\leq 1$, for any $\gs>0$ one 
can find $\gth>0$ such that (\ref {pst4}) holds.\qeda\medskip

The next result due to \cite {KV} (\cite {BC} if $L=-\Gd$) points out 
how close to a necessary condition estimate (\ref {pst1}) is.
\bth{KVth2} Let $q>1$, $\gl\in \mathfrak 
M_{+}(\Gw;\gr_{_{\prt\Gw}})$, $\gs>0$. If 
there is a positive solution $u\in L^{1}(\Gw)$ to Problem (\ref 
{pst2}), there exists a constant $C_{1}>0$ such that
\begin {eqnarray}\label {pst5}
\mathbb G_{L}^\Gw\left(\left(\mathbb G_{L}^\Gw(\gs\gl)\right)^q\right)
\leq C_{1}\mathbb G_{L}^\Gw(\gs\gl),\quad\mbox{ a.e. in } \Gw.
\end {eqnarray}
If $L=-\Gd$, $C_{1}=1/(q-1)$.
\es
\blemma{BClem} Let $h\in L^{1}(\Gw;\rho_{_{\prt\Gw}}dx)$, $h\geq 0$, 
and $\gm,\,\eta\in \mathfrak M_{+}(\Gw;\gr_{_{\prt\Gw}})$, 
$\gm\neq 0$, such that $\gm-\eta\geq h$.  If $\phi\in C^2([0,\ity))$ is a concave 
nondecreasing function such that $\phi(1)\geq 0$, there holds
\begin {eqnarray}\label {pst6}
h\phi'\left(\frac {\mathbb G_{-\Gd}^\Gw(\gm)}{\mathbb G_{-\Gd}^\Gw(\eta)}\right)
\in L^{1}(\Gw;\rho_{_{\prt\Gw}}dx),
\end {eqnarray}
and
\begin {eqnarray}\label {pst7}
-\Gd\left(\phi\left(\frac {\mathbb G_{-\Gd}^\Gw(\gm)}
{\mathbb G_{-\Gd}^\Gw(\eta)}\right)\mathbb G_{-\Gd}^\Gw(\eta)\right)\geq
h\phi'\left(\frac {\mathbb G_{-\Gd}^\Gw(\gm)}{\mathbb G_{-\Gd}^\Gw(\eta)}\right).
\end {eqnarray}
\es
\Proof Put $z=\mathbb G_{-\Gd}^\Gw(\gm)$ and $w=\mathbb 
G_{-\Gd}^\Gw(\eta)$. We write $\eta=h+\gm+\gs$ where $\gs$ is a 
positive Radon measure. Let $h_{n}$, $\gm_{n}$ and $\gs_{n}$ be elements 
of $C_{c}^\ity(\Gw)$ such that $h_{n}\to h$ in 
$L^{1}(\Gw;\rho_{_{\prt\Gw}}dx$, and $\gm_{n}\to\gm$ and 
$\gs_{n}\to\gs$, in the weak sense of $\mathfrak 
M_{+}(\Gw;\gr_{_{\prt\Gw}})$. Put $z_{n}=\mathbb G_{-\Gd}^\Gw(\gm_{n})$ and 
$w_{n}=\mathbb G_{-\Gd}^\Gw(h_{n}+\gm_{n}+\gs_{n})$, then
$z_{n}\to z$ and $w_{n}\to w$ in $L^1(\Gw)$ as $n\to\infty$, and 
a.e. (after extraction of a subsequence). Thus $z_{n}>0$ in $\Gw$, 
for $n$ large enough. Because of the concavity, 
$\phi(1)\geq 0$ and $\phi'\geq 0$, there  holds
$$-\Gd\left(z_{n}\phi\left(\frac {w_{n}}{z_{n}}\right)\right)
\geq \phi'\left(\frac {w_{n}}{z_{n}}\right)(h_{n}+\gs_{n})
\geq \phi'\left(\frac {w_{n}}{z_{n}}\right)h_{n}.
$$
Also
$$
0\leq z_{n}\phi\left(\frac {w_{n}}{z_{n}}\right)
\leq z_{n}\left(\phi_{0}+\phi'(0)\frac {w_{n}}{z_{n}}\right)
\leq C(z_{n}+w_{n}),
$$
for some $C>0$. Therefore $z_{n}\phi\left(w_{n}/z_{n}\right)$ 
converges in $L^1(\Gw)$ as $n\to\infty$. Since for any $\xi\in 
C^{1,1}_{c}(\overline\Gw)$, $\xi\geq 0$, there holds
\begin {eqnarray}\label {pst8}
-\int_{\Gw}z_{n}\phi\left(\frac {w_{n}}{z_{n}}\right)\Gd\xi dx
\geq \int_{\Gw}\phi'\left(\frac {w_{n}}{z_{n}}\right)h_{n}\xi dx,
\end{eqnarray}
we derive (\ref {pst7}) by passing to the limit with Lebesgue and 
Fatou's theorems.\qeda 
\medskip

\noindent {\it Proof of \rth{KVth2}. } First, we prove the result 
when $L=-\Gd$. Since $\gs>0$, we can assume $\gs=1$ and apply \rlemma 
{BClem} with $w=u$, the solution of (\ref {pst2}),  $z=G_{-\Gd}^\Gw(\gl)$ and 
$$\phi(s)=\left\{\BA{l}(1-s^{1-q})/(q-1),\quad \mbox { if }\,s\geq 
1,\\[2mm]
s-1,\quad\qquad \qquad \qquad  \mbox {if }\,s\leq 1.
\EA\right.$$
Because $u\geq G_{-\Gd}^\Gw(\gl)$, 
\begin {eqnarray}\label {pst9}
-\Gd\left(G_{-\Gd}^\Gw(\gl)\phi\left(\frac {u}
{G_{-\Gd}^\Gw(\gl)}\right)\right)\geq \phi'\left(\frac {u}
{G_{-\Gd}^\Gw(\gl)}\right) u^q=\left(G_{-\Gd}^\Gw(\gl)\right)^q,
\end{eqnarray}
holds weakly. By the maximum principle,
\begin {eqnarray}\label {pst10}
\frac {1}{q-1}G_{-\Gd}^\Gw(\gl)-\frac {1}{q-1}u^{1-q}\left(G_{-\Gd}^\Gw(\gl)\right)^q
\geq G_{-\Gd}^\Gw\left(\left(G_{-\Gd}^\Gw(\gl)\right)^q\right),
\end{eqnarray}
which is the expected inequality in the case $L=-\Gd$. We turn now to 
the general case. By \rth {equivth}, the Green functions of $L$ and 
$-\Gd$ are equivalent in the sense that
\begin {eqnarray*}
C^{-1}G_{-\Gd}^\Gw(x,y)
\leq G_{L}^\Gw(x,y)\leq CG_{-\Gd}^\Gw(x,y),
\forevery (x,y)\in\Gw\times\Gw\setminus D_{\Gw},
\end {eqnarray*}
for some $C>0$. Thus (\ref {pst6}) follows.\qeda \medskip 

\noindent\Remark In \cite {KV}, inequality (\ref {pst6}) is proven for a 
very general class of positive kernels, not only for a Green kernel.\medskip 

The next result, proven in \cite {BVY}, exhibits a large class of measures for 
which Problem (\ref {pst2}) will be solvable by applying \rth {KVth}.
\bth {bvyth} Let $q>0$, $\ga\in [0,1]$ and 
$\gl\in\mathfrak M_{+}(\Gw;\gr^\ga_{_{\prt\Gw}})$ with 
$\norm{\gl}_{\mathfrak M_{+}(\Gw;\gr^\ga_{_{\prt\Gw}})}=1$. If
\begin {equation}\label {pst11}
q<\frac {n+\ga}{n+\ga-2},
\end {equation}
then $\mathbb G_{L}^\Gw(\gl)\in 
L^1(\Gw;\gr^\ga_{_{\prt\Gw}}dx)$, and there exists a 
positive constant $C=C(n,q,\ga,\gl,\Gw)$ such that 
\begin {equation}\label {pst12}
\mathbb G_{L}^\Gw\left(\left(G_{L}^\Gw(\gl)\right)^q\right)
\leq C G_{L}^\Gw(\gl)\quad \mbox {a.e. in }\; \Gw.
\end {equation}
\es 
\Proof As in the proof of \rth {KVth2}, it is sufficient to consider 
the case $L=-\Gd$ and then use the equivalence of Green kernels. 
\smallskip

\noindent {\it Step 1  } The case $\gl=\gd_{y}$ for $y\in \Gw$, $n\geq 3$. Since 
$G_{-\Gd}^\Gw(x,y)\leq C(n)\abs {x-y}^{2-n}$ we put $d=$diam$(\Gw)$ and 
\begin {equation}\label {pst13}h(x)=\left\{\BA{l}
\abs {x-y}^{2-(n-2)q}\quad\qquad \mbox {if }\;q>2/(n-2),\\[2mm]
d-\abs {x-y}^{2-(n-2)q}\quad\;\mbox {if }\;q<2/(n-2),\\[2mm]
\ln(d/\abs {x-y})\qquad\qquad\mbox {if }\;q=2/(n-2).
\EA\right.
\end {equation}
Hence
$$-\Gd h(.)=C_{1}\abs {.-y}^{(2-n)q}\quad\mbox {in }\;\CD'(\Gw),
$$
and consequently  
$$\mathbb G_{-\Gd}^\Gw\left(\left(G_{-\Gd}^\Gw(.,y)\right)^q\right)(x)
\leq C_{2}h(x)\leq C_{3}\abs {x-y}^{2-n},
$$
with $C_{i}=C_{i}(n,q,d)>0$. Let $r>0$ be such that 
$\overline B_{r}(y)\subset\Gw$. Clearly 
$$\mathbb G_{-\Gd}^\Gw\left(\left(G_{-\Gd}^\Gw(.,y)\right)^q\right)(x)
\leq C'_{y}\gr_{_{\prt\Gw}}(x)\leq C''_{y}G_{-\Gd}^\Gw(x,y),
$$
on $\overline B_{r}(y)\setminus \{y\}$. On $\Gw\setminus B_{r}(y)$ the 
function $\mathbb 
G_{-\Gd}^\Gw\left(\left(G_{-\Gd}^\Gw(.,y)\right)^q\right)$ is 
$C^{1}$. We get a similar inequality by Hopf boundary lemma. Finally 
there exists $C_{y}>0$ such that
\begin {equation}\label {pst14}
\mathbb 
G_{-\Gd}^\Gw\left(\left(G_{-\Gd}^\Gw(.,y)\right)^q\right)(x)\leq 
C_{y}G_{-\Gd}^\Gw(x,y),\forevery x\in\Gw\setminus \{y\}.
\end {equation}
As we shall see it in next step, $C_{y}$ is bounded independently of 
$y$.
\smallskip

\noindent {\it Step 2  } The general case. By 
\rth{estth}, $\mathbb G_{-\Gd}^\Gw(\gl)\in 
L^q(\Gw;\gr^\ga_{_{\prt\Gw}}dx)$ since (\ref {pst11}) holds. First 
assume $q\geq 1$, then
$$\mathbb G_{-\Gd}^\Gw(\gl)(x)=\int_{\Gw}G_{-\Gd}^\Gw(x,y)d\gl (y)
=\int_{\Gw}\frac {G_{-\Gd}^\Gw(x,y)}{\gr^\ga_{_{\prt\Gw}}(y)}
\gr^\ga_{_{\prt\Gw}}(y)d\gl (y).
$$
By Jensen's inequality, 
$$\left(\mathbb G_{-\Gd}^\Gw(\gl)(x)\right)^q\leq
\int_{\Gw}\left(\frac 
{G_{-\Gd}^\Gw(x,y)}{\gr^\ga_{_{\prt\Gw}}(y)}\right)^q
\gr^\ga_{_{\prt\Gw}}(y)d\gl (y),
$$
$$\mathbb G_{-\Gd}^\Gw\left(\left(\mathbb 
G_{-\Gd}^\Gw(\gl)\right)^q\right)(x)\leq
\int_{\Gw}\mathbb G_{-\Gd}^\Gw
\left(G_{-\Gd}^\Gw(.,y)\right)(x)\gr^{\ga(1-q)}_{_{\prt\Gw}}(y)d\gl (y).
$$
Now
$$\mathbb G_{-\Gd}^\Gw\left(G_{-\Gd}^\Gw(.,y)\right)(x)\gr^{\ga(1-q)}_{_{\prt\Gw}}(y)
=\int_{\Gw}G_{-\Gd}^\Gw(x,z)G_{-\Gd}^\Gw(y,z)
\left(\frac {G_{-\Gd}^\Gw(y,z)}{\gr_{_{\prt\Gw}}(y)}\right)^{q-1}dz.
$$
Because
\begin {equation}\label {pst15}
G_{-\Gd}^\Gw(y,z)\leq C\min\{\abs {y-z}^{2-n},\gr_{_{\prt\Gw}}(y)\abs {y-z}^{1-n}\},
\end {equation}
it follows
$$
G_{-\Gd}^\Gw(y,z)\leq C\gr^\ga_{_{\prt\Gw}}(y)\abs {y-z}^{2-n-\ga}.
$$
At that point of the proof we recall the following relation called the $3$-$G$ inequality (see \cite 
{DaL} for example),
\begin {equation}\label {pst16}
\frac {G_{-\Gd}^\Gw(x,z)G_{-\Gd}^\Gw(y,z)}{G_{-\Gd}^\Gw(x,y)}
\leq C\left(\abs {x-z}^{2-n}+\abs {y-z}^{2-n}\right),
\end {equation}
where $C=C(\Gw)$. It implies
$$\mathbb G_{-\Gd}^\Gw\left(G_{-\Gd}^\Gw(.,y)\right)(x)\gr^{\ga(1-q)}_{_{\prt\Gw}}(y)
\leq G_{-\Gd}^\Gw(x,y)I(x,y),
$$
for some $C=C(q,\Gw,\ga)$, and
$$I(x,y)=\int_{\Gw}\abs {y-z}^{(2-n-\ga)(q-1)}\left(\abs {x-z}^{2-n}+\abs {y-z}^{2-n}\right)dz.
$$
Since
$$I(x,y)\leq C\int_{\Gw}
\left(\abs {x-z}^{2-n+(2-n-\ga)(q-1)}+\abs 
{y-z}^{2-n+(2-n-\ga)(q-1)}\right)dz,
$$
this last quantity is clearly bounded independently of $x$ and $y$ by 
some constant depending on the various parameters and data. Notice 
that we have used 
$$q<(n+\ga)/(n+\ga-2)\leq n/(n-2).$$
Thus
\begin {equation}\label {pst17}
\mathbb G_{-\Gd}^\Gw\left(\left(\mathbb 
G_{-\Gd}^\Gw(\gl)\right)^q\right)(x)\leq 
C\int_{Gw}G_{-\Gd}^\Gw(x,y)d\gl(y)=C\mathbb G_{-\Gd}^\Gw(x).
\end {equation}
Obviously, $C=C(\Gw)$ when $q=1$.\smallskip

\noindent Next we assume $0\leq q<1$. Then
$$
\mathbb G_{-\Gd}^\Gw\left(\left(\mathbb 
G_{-\Gd}^\Gw(\gl)\right)^q\right)\leq \mathbb G_{-\Gd}^\Gw(1)+
\mathbb G_{-\Gd}^\Gw\left(\left(\mathbb 
G_{-\Gd}^\Gw(\gl)\right)\right).
$$
By Hopf boundary lemma $\mathbb G_{-\Gd}^\Gw(1)(x)\leq 
C\gr_{_{\prt\Gw}}(x)$. Let $K$ be a compact subset contained in 
the support of $\gl$ and denote by $\gl\vline_{ K}$ the restriction 
of $\gl$ to $K$. By the regularity results, $\mathbb G_{-\Gd}^\Gw(\gl\vline_{K})
\in C^{1}(\overline\Gw\setminus K)$. Then 
$\mathbb G_{-\Gd}^\Gw(\gl)\geq \mathbb 
G_{-\Gd}^\Gw(\gl\vline_{K})\geq C\gr_{_{\prt\Gw}}$ in $\overline\Gw\setminus K$.
In turn it implies $\mathbb G_{-\Gd}^\Gw(\gl)\geq C\gr_{_{\prt\Gw}}$ 
for another constant $C>0$ and (\ref {pst12}) follows.\qeda \medskip

Condition (\ref {pst11}) on $q$ is called $\ga$-subcriticality. 
However, 
as we have seen it in previous sections, there exists measures for 
which (\ref {pst2}) is solvable even if $q$ is not $\ga$-subcritical.

\bdef {}{\rm A measure $\gl\in \mathfrak M_{+}(\Gw;\gr^\ga_{_{\prt\Gw}})$ 
is called $q$-admissible if there exists some $\gs_{0}\geq 0$ such 
that Problem (\ref {pst2}) admits a solution $u\in L^{1}(\Gw)\cap 
L^q(\Gw;\gr_{_{\prt\Gw}}dx)$ whenever $0<\gs\leq\gs_{0}$.
}\es

The following theorem summarizes the results of Baras and Pierre \cite 
{BP2}, Adams and Pierre \cite {AP} and Kalton and 
Verbitsky \cite {KV} in the super-critical range of exponents.
\bth {apkvth} Let $q>1$, $\ga\in [0,1]$ and 
$\gl\in\mathfrak M_{+}(\Gw;\gr^\ga_{_{\prt\Gw}})$. Then the 
following conditions are equivalent :\smallskip

\noindent (i) $\gl$ is $q$-admissible.\smallskip

\noindent (ii) There exists some $C_{0}>0$ such that
\begin {equation}\label {(i)}
\mathbb G_{L}^\Gw\left(\left(\mathbb G_{L}^\Gw(\gl)\right)^q\right)
\leq C_{0}\mathbb G_{L}^\Gw(\gl).
\end {equation}
(iii) $\left(\mathbb G_{L}^\Gw(\gl)\right)^q$ is $q$-admissible.\smallskip

\noindent (iv) There exists $C>0$ such that 
\begin {equation}\label {(iv)}
\int_{\Gw}\mathbb G_{L}^\Gw(\gl)dx\leq C\int_{\Gw}\frac 
{g^{q'}}{\left(\mathbb G_{L}^\Gw(g)\right)^{q'-1}}dx, \forevery 
g\in L_{c}^{\ity}(\Gw),\,g\geq 0.
\end {equation}
\noindent (iv) There exists $c>0$ such that 
\begin {equation}\label {(v)}
\int_{A}d\gl\leq cC_{2,q',\ga}(A),\forevery A\subset\Gw,\,A\mbox { 
Borel},
\end {equation}
where $C_{2,q',\ga}$ is the weighted capacity defined by
\begin {equation}\label {(v')}
C_{2,q',\ga}(A)=\inf\left\{\int_{\Gw}\eta^{q'}dx:\eta\in 
L^{q'}(\Gw),\,\eta\geq 0,\, \mathbb G_{L^{*}}^\Gw(\gl)\geq 
\gr^\ga_{_{\prt\Gw}}\;\mbox {on }A\right\}.
\end {equation}
\es 
\subsection {Isolated singularities}
If one looks for radial positive solutions of 
\begin {equation}\label {singspr}
-\Gd u=\abs u^{q-1}u,
\end {equation} 
with $q>1$, in $\mathbb R^n\setminus \{0\}$ under the form $x\mapsto a\abs x^b$, one immediately finds
\begin {equation}\label {singspr2}
u(x)=u_{s}(x)=\gg_{q,n}\abs x^{-2/(q-1)},
\end {equation}
where
\begin {equation}\label {singspr3}
\gg_{q,n}=\left(\left(\frac {2}{q-1}\right)\left(n-\frac {2q}{q-1}\right)\right)^{1/(q-1)}.
\end {equation}
However such a solution exists if and only if $q>n/(n-2)$. Moreover, 
if $q\geq n/(n-2)$, it follows by \rth {remov1} that, if $\Gw$ is an open subset of 
$\mathbb R^n$ containing $0$, $\Gw^{*}=\Gw\setminus \{0\}$, and if $u\in 
L^{q}_{loc}(\Gw^{*})$ is nonnegative and satisfies
\begin {equation}\label {singspr4}
-\Gd u=u^q\quad \mbox {in }\,\CD'(\Gw^{*}),
\end {equation} 
then $u\in L^{q}_{loc}(\Gw)$, and that Equation (\ref{singspr4}) holds in $\CD'(\Gw)$. 
In this way, 
the singularity of $u$ at $0$ exists, but is not visible in the sense 
of distributions. In the subcritical range, $1<q<n/(n-2)$ it is proven 
by Brezis and Lions \cite {BL} that any positive solution of (\ref {singspr4}) 
satisfies actually
\begin {equation}\label {singspr5}
-\Gd u=u^q+C_{n}\gg\gd_{0}\quad \mbox {in }\,\CD(\Gw),
\end {equation}
for some $\gg\geq 0$ (see Step 4 in the proof of \rth {VV2th}). Furthermore $u$ admits an expansion near $0$;
\begin {equation}\label {singspr6}
u(x)=\gg\abs x^{2-n}(1+\circ(1))+C,\quad\mbox {as }\;x\to 0,
\end {equation}
if $n\geq 3$, with the usual modification if $n=2$. Finally, although 
this was noticed before by Lions \cite {Lio}, \rth {bvyth} implies that 
the Dirac mass
$\gd_{0}$ is $q$-admissible. The classification of isolated 
singularities of positive solutions of (\ref {singspr}) has been 
performed by Lions \cite {Lio} in the case $1<q<n/(n-2)$, Aviles \cite 
{Av} in the case $q=n/(n-2)$, Gidas and Spruck \cite {GS} when 
$n/(n-2)<q<(n+2)/(n-2)$ and Caffarelli, Gidas and Spruck \cite {CGS} in the case 
$q=(n+2)/(n-2)$. The case $q>(n+2)/(n-2)$ remains essentially open, 
except if the solutions are supposed to be radial.
\bth{isosing} Let $\Gw$ be an open subset of $\mathbb R^n$  
containing $0$, $\Gw^{*}=\Gw\setminus\{0\}$, $q>0$ and $u\in C^{2}(\Gw^{*})$ 
be a positive solution of (\ref {singspr}) in $\Gw^{*}$.\smallskip

\noindent (i) If $q<n/(n-2)$ : either $u\in C^{\ity}(\Gw)$, or there 
exists $\gg>0$ such that (\ref {singspr6}) and (\ref {singspr5}) 
hold.\smallskip

\noindent (ii) If $q=n/(n-2)$ : either $u\in C^{\ity}(\Gw)$, or 
\begin {equation}\label {Ave}
\lim_{x\to 0}{\abs x}^{n-2}\left(\ln (1/\abs x)\right)^{(2-n)/2}u(x)=\left(\frac {n-2}{\sqrt 
2}\right)^{n-2}.
\end {equation}

\noindent (iii) If $n/(n-2)<q<(n+2)/(n-2)$ : either $u\in C^{\ity}(\Gw)$, or 
\begin {equation}\label {GSe}
\lim_{x\to 0}{\abs x}^{2/(q-1)}u(x)=\gg_{q,n}.
\end {equation}

\noindent (iv) If $q=(n+2)/(n-2)$ : either $u\in C^{\ity}(\Gw)$, or 
\begin {equation}\label {CGSe}
\lim_{x\to 0}{\abs x}^{(n-2)/2}(u(x)-v(\abs x)=0,
\end {equation}
where $r\mapsto v(r)$ is a radial solution of (\ref {singspr}). 
\es 

Notice that in the so-called conformal case $q=(n+2)/(n-2)$, all the 
radial solutions $v$ of (\ref {singspr}) are classified by their reduced 
energy : if $v(r)=r^{(2-n)/2}w(t)$ and $t=\ln(1/r)$, then $w$ verifies
\begin {equation}\label {CGSe2}
w''-\frac {(n-2)^{2}}{4}w+\abs w^{(4)/(n-2)}w=0.
\end {equation}
Therefore the reduced energy-function
$$\CE(w)={w'}^{2}+\frac {n+2}{n}\abs {w}^{2n/(n+2)}-\frac {(n-2)^{2}}{4}w^{2}
$$
is constant. The proofs of these different results relies on 
regularity estimates and bootstrap arguments in case (i), the Lyapounov analysis as 
for \rth {isolth2} in cases (ii) and (iii), and the asymptotic symmetry method in 
the case (iv). However, there are two difficulties in case (iii) ((ii) 
being much simpler) : the first one is to prove the {\it a priori} estimate
\begin {equation}\label {GSe2}
u(x)\leq C\abs x^{2/(q-1)}\quad \mbox {near }\; 0.
\end {equation}
The second one is to identify the limit set at the end of the Lyapounov 
analysis, in which situation, it is to be proven that the only positive 
solutions to
\begin {equation}\label {GSe3}
-\Gd_{_{S^{n-1}}}\gw+\gg_{q,n}^{q-1}\gw -w^q=0
\end {equation}
on $S^{n-1}$ are the constant solutions $0$ and $\gg_{q,n}$.\medskip

\noindent \Remark Part of the results can be extended to equation
\begin {equation}\label {GSe4}
Lu=u^q,
\end {equation}
where $L$ is a general elliptic operator, satisfying condition (H). This extension is easy 
for (i), a little more complicated in case (iii) (and (ii) in the same 
way), in  particular to get (\ref {GSe2}). 
It is still completely open in case (iv).

\mysection {Boundary singularities and boundary trace}

In this chapter we shall study generalized boundary value problems for 
equation
\begin {equation}\label {Bsing}
Lu+g(x,u)=0\quad\mbox {in }\;\Gw,
\end {equation}
where $\Gw$ is an open domain in $\mathbb R^n$, $n\geq 2$, with a 
$C^{2}$ boundary, $L$ is an elliptic operator defined in $\Gw$ by (\ref {lin1}) and $g$ a 
continuous function of absorption type.
\subsection{Measures boundary data}
\subsubsection {General solvability}
Let $\gm$ be a Radon measure on $\prt\Gw$ and $g\in C(\Gw\times\mathbb 
R$. The semilinear Dirichlet problem with measure data is written 
under the form
\begin {equation}\label {Bmeas}\left.\BA{ll}
Lu+g(x,u)=0&\mbox{in }\;\Gw,\\[2mm]
\phantom {Lu+g(x,)}
u=\gm&\mbox{on }\;\prt\Gw.
\EA\right.
\end {equation}

\bdef {Bmeasdef} Let $\gm\in \mathfrak M(\prt\Gw)$. A 
function $u$ is a solution of (\ref {Bmeas}), if $u\in L^1(\Gw)$, 
$g(.,u)\in L^1(\Gw;\gr_{_{\prt\Gw}}dx)$, and
if for any $\gz\in C_{c}^{1,L}(\overline \Omega)$, there holds
\begin {eqnarray}\label {Bmeas1}
\int_{\Gw}\left( uL^*\gz +g(x,u)\gz\right)dx=-\int_{\Gw}\frac {\prt 
\gz}{\prt{{\bf n}_{L^{*}}}} d\gm.
\end {eqnarray}
\es 
\bdef {BGV} {\rm A real valued function $g\in C(\Gw\times \mathbb R)$ 
holds the {\it boundary-weak-singularity assumption}, 
if there exists $r_{0}\geq 0$ such that 
\begin {eqnarray}\label {passive''}
  rg(x,r)\geq 0, \forevery (x,r)\in\Gw\times 
(-\ity,-r_{0}]\cup[r_{0},\ity), 
\end{eqnarray}
and a nondecreasing function $\tilde g\in 
C([0,\ity))$ such that $\tilde g\geq 0$, 
\begin {eqnarray}\label {bwsa}
\int_{0}^1\tilde g(r^{1-n})r^ndr<\ity,
\end{eqnarray}
and 
\begin {eqnarray}\label {bwsa1}
\abs{g(x,r)} \leq \tilde g(\abs r),\forevery (x,r)\in 
\Gw\times \mathbb R.
\end{eqnarray}
}\es 
The following result was proven first, but under a weaker form, by Gmira and 
V\'eron \cite {GmV}.
\bth {Bmeasth}Let $\Omega$ be a $C^2$ bounded domain in $\mathbb R^n$, $n\geq 
2$, $L$ the elliptic operator defined by (\ref {lin1}) and 
$g\in C(\Gw\times \mathbb R)$ a real valued function. If $L$ 
satisfies assumptions (H) and $g$ the boundary-weak-singularity 
assumption, for any $\gm\in \mathfrak M(\prt\Gw)$ there 
exists a solution $u$ to Problem (\ref {Bmeas}).
\es 
\Proof The general idea follows the proof of \rth {genMeas}, with some 
significant changes.\smallskip

\noindent {\it Step 1 } Approximate solutions. Let $\gm_{n}$ be a 
sequence of $C^{2}(\Gw)$ functions converging to $\gm$ in the weak 
sense of measures and $m_{n}=\mathbb P_{L}^\Gw(\gm_{n})$. The 
function $g^n$ defined by
$$g^n(x,r)=g(x,r-m_{n}(x)),\forevery (x,r)\in \Gw\times\mathbb R, 
$$
is continuous in $\Gw\times\mathbb R$ and satisfies (\ref 
{passive''}) with $r_{0}$ replaced by $r_{0}+{\norm {m_{n}}}_{L^\ity}$.
By \rth{genMeas} there exists a solution to 
\begin {equation}\label {Bmeas-n}\left.\BA{ll}
Lv_{n}+g^n(x,v_{n})=0&\mbox{in }\;\Gw,\\[2mm]
\phantom {Lv_{n}+g^n(x,)}
v_{n}=0&\mbox{on }\;\prt\Gw.
\EA\right.
\end {equation}
Thus the function $u_{n}=v_{n}+m_{n}$ is a solution of
\begin {equation}\label {Bmeas-n'}\left.\BA{ll}
Lu_{n}+g(x,u_{n})=0&\mbox{in }\;\Gw,\\[2mm]
\phantom {Lu_{n}+g(x,)}
u_{n}=\gm_{n}&\mbox{on }\;\prt\Gw.
\EA\right.
\end {equation}
From the proof of  \rth{genMeas},  Steps 2-3, $u_{n}$ is bounded in 
$\Gw$ and (\ref {Bmeas1}) holds with $u_{n}$ and $m_{n}$. By \rth 
{L1-th}, for any $\gz\in C_{c}^{1,L}(\overline\Gw)$, $\gz\geq 0$,
\begin {eqnarray}\label {Bmeas3}
\int_{\Gw}\left( \abs {u_{n}}L^*\gz +{\rm 
sign(u_{n})}g(x,u_{n})\gz\right)dx\leq-\int_{\Gw}\frac {\prt\gz}{\prt 
{\bf n}_{L^{*}}}\abs{\gm_{n}} dx,
\end {eqnarray} 
which implies
\begin {eqnarray}\label {Bmeas4}
\norm {u_{n}}_{L^1(\Gw)}+\norm {\gr_{_{\prt\Gw}}g(.,u_{n})}_{L^1(\Gw)}
\leq \Gth\int_{\Gw}\gr_{_{\prt\Gw}} dx+C_{1}\norm 
{\gr_{_{\prt\Gw}}\gm_{n}}_{L^1(\prt\Gw)}.
\end {eqnarray}
Consequently, using also (\ref {marsest3}) in \rth {estth},
\begin {eqnarray}\label {Bmeas5}
{\norm{u_{n}}}_{M^{(n+\ga)/(n+\ga-2)}(\Gw;\gr^\ga_{_{\prt\Gw}})}
\leq C_{2}\norm {\gl_{n}-g(.,u_{n})}_{\mathfrak 
M(\Gw;\rho^\ga_{_{\prt\Gw}})}
\leq C_{3}\left(\Gth+\norm 
{\gr_{_{\prt\Gw}}\gm_{n}}_{L^1(\prt\Gw)}\right),
\end {eqnarray}
for $\ga=0$, $1$. \smallskip

\noindent {\it Step 2  }Convergence. By \rcor {L^1regloc} and (\ref 
{Bmeas5}), there exists a subsequence of $\{u_{n}\}$, still denoted by 
$\{u_{n}\}$ for simplicity, which converges 
to some $u$ in $ L^{1}(\Gw)$ and a.e. in $\Gw$. In order to prove 
that $g(.,u_{n)}$ converges in $L^{1}(\Gw;\gr_{_{\prt\Gw}}dx)$, we 
use Vitali's theorem and we procede as in the proof of \rth 
{genMeas}- Step 3 with $\ga=1$. \qeda \medskip

The following stability result follows from the uniform integrability 
argument.

\bcor {Bmeascor} Let  $g$ 
satisfy the boundary-weak-singularity assumption and $r\mapsto g(x,r)$
is nondecreasing, for any $x\in\Gw$. Then the solution $u$ is unique. 
If we assume that $\{\gm_{k}\}$ is a sequence of measures in  
$\mathfrak M(\Gw)$ which converges weakly to $\gm$, then the 
corresponding solutions $u_{\gm_{k}}$ of problem 
 \begin {equation}\label {Bmeask}\left.\BA{ll}
Lu_{\gm_{k}}+g(x,u_{\gm_{k}})=0&\mbox{in }\;\Gw,\\[2mm]
\phantom {Lu_{\gm_{k}}+g(x,)}
u_{\gm_{k}}=\gm_{k}&\mbox{on }\;\prt\Gw,
\EA\right.
\end {equation}
converge in $L^{1}(\Gw)$ to the solution $u$ of (\ref {Bmeas}), when 
$k\to\ity$.
\es 
\Remark If $g(x,r)=\abs r^{q-1}r$,  the boundary-weak-singularity 
assumption is satisfied if and only if
\begin {eqnarray}\label {Bpower q}
0<q<\frac {n+1}{n-1}.
\end{eqnarray}
\subsubsection {Admissible boundary measures and the $\Delta_{2}$-condition}

\bdef {Bgaddef}{\rm Let $\tilde g$ be a continuous real valued 
nondecreasing function defined in $\mathbb R_{+}$, $\tilde g\geq 0$. A 
Radon measure $\gm$ 
in $\prt\Gw$ is called $(\tilde g,k)$-boundary-admissible if
\begin {eqnarray}\label {Bg-adm}
\int_{\Gw}\tilde g(\mathbb P_{L}^\Gw(\abs\gm)+k)\rho_{_{\prt\Gw}}dx<\infty,
\end{eqnarray}
where $P_{L}^\Gw(\abs\gm)$ is the Poisson potential of $\gm$ and $k\geq 0$.}
\es

The proof of the following theorem is similar to the one of \rth{g-adm-th}.
\bth {Bg-adm-th} Let $\Gw$ be a $C^2$ bounded domain in $\mathbb 
R^n$, $n\geq 2$, $L$ an elliptic operator defined by (\ref 
{lin1}) verifying condition (H), and $g\in C(\Gw\times\mathbb R)$ 
satisfying (\ref {passive''}) for some $r_{0} \geq 0$ and (\ref 
{bwsa1}) for some function $\tilde g$ as in 
\rdef {admdef}. Then for any $(\tilde g,r_{0})$-boundary-admissible Radon measure 
$\gm\in\mathfrak M(\prt\Gw)$, Problem (\ref {Bmeas}) admits a 
solution.
\es

The proof of the next 
result, is a boundary adaptation of the one of  \rth{g-adm-th2}.

\bth {Bg-adm-th2} Let $\Gw$ and $L$ be as in \rth {Bg-adm-th}. Assume
$g\in C(\Gw\times\mathbb R)$ satisfies the $\Delta_{2}$-condition 
(\ref {delta-2}), 
$r\mapsto g(x,r)$ is nondecreasing for any $x\in\Gw$ and
(\ref {bwsa1}) holds for some nonnegative, nondecreasing function $\tilde g$. 
For any Radon measure 
$\gl\in\mathfrak M(\prt\Gw)$, 
with $\gl=\tilde\gl+\gl^*$, where $\tilde\gl\in 
L^1(\prt\Gw)$ and $\gl^*$  is 
$(\tilde g,0)$-boundary-admissible and singular with respect to the 
$(n-1)$-dimensional Hausdorff measure, problem 
(\ref {Bmeas}) admits a unique solution.
\es
\subsubsection {Sharp solvability}
The existence of a solution, necessarily unique, to
 \begin {equation}\label {Bq-power}\left.\BA{ll}
Lu+\abs u^{q-1}u=0&\mbox{in }\;\Gw,\\[2mm]
\phantom {Lu+\abs u^{q-1}}
u=\gm&\mbox{on }\;\prt\Gw,
\EA\right.
\end {equation}
where $\gm$ is a boundary measure follows unconditionaly from \rth 
{Bmeasth} in the subcritical range $0<q<(n+1)/(n-1)$. The super-critical case $q\geq 
(n+1)/(n-1)$ is treated separately according the value of $q$ with 
respect to $2$ by Le Gall \cite {LG3}, Dynkin and Kuznestov \cite 
{DK3}, \cite {DK4} and 
Marcus and V\'eron \cite {MV5}. The synthetic presentation in all the 
super-critical cases is found in \cite {MV6}.
\bth{Bq-powerth} Let $\Gw$ be a bounded domain in $\mathbb R^n$ with a 
$C^{2}$ boundary, $L$ the elliptic operator defined by (\ref {lin1}) 
satisfying condition (H), $q\geq (n+1)/(n-1)$ and $\gm\in \mathfrak 
M(\prt\Gw)$. Then Problem (\ref {Bq-power}) admits a solution 
$u=u_{\gm}$ if and 
only if $\gm$ does not charge boundary sets with 
$C_{2/q,q'}$-capacity zero. Moreover, the mapping $\gm\mapsto u_{\gm}$ 
is increasing.  
\es 

Following \rdef {Bgaddef}, a Radon measure $\gm$ on $\prt\Gw$ is called 
{\it boundary-q-admissible} for the operator $L$ if 
\begin {eqnarray}\label {Bq-adm}
\int_{\Gw}\left(\mathbb 
P_{L}^\Gw(\abs\gm)\right)^q\gr_{_{\prt\Gw}}dx<\infty.
\end {eqnarray}
However, under assumption (H), under which the Green and Poison 
kernels are constructed, this property is independent of $L$, since 
all the kernels are equivalent (see \rth {equivth}).
The proof is based upon a deep result concerning representation of 
boundary Bessel classes in terms of integrability properties of 
Poisson potentials.
\bprop{strongA}
Let the assumtions of \rth {Bq-powerth}, on $\Gw$ and the operator $L$,  be 
satisfied, $q\geq (n+1)/(n-1)$ and $\gm\in\mathfrak M(\prt\Gw)$. 
Then : \smallskip 

\noindent (i) If $\mu$ is boundary-$q$-admissible, then $\mu\in 
W^{-2/q,q}(\prt\Gw)$. \smallskip 

\noindent (ii) If $\mu\in \mathfrak M_{+}(\prt\Gw)\cap 
W^{-2/q,q}(\prt\Gw)$, then $\mu$ is boundary-$q$-admissible. Moreover 
there exists a constant $C=C(q,\Gw,L)$ such that,
\begin{equation}\label{norm-eqL}
 C^{-1} 
 \|\mu\|_{W^{-2/q,q}(\prt\Gw)}\leq\|\BBP_{L}^\Gw(\mu)\|_{L^q(\Gw;\gr_{_{\prt\Gw}} dx)}
 \leq C \|\mu\|_{W^{-2/q,q}(\prt\Gw)}.
\end{equation}
\es
\Proof The proof we present here is settled upon the interpolation 
theory between a Banach space and the domain of an analytic semigroup of 
operators.\smallskip

\noindent {\it Step 1 } The case where $\Gw$ is the unit ball $B$. 
We shall assume $n\geq 3$, 
the $2$-dimensional case requiring some easy technical modifications. 
Let $(r,\gs)$ be the spherical coordinates in $\mathbb R^n$ , $ t=-\ln 
r$. If $\mu\in W^{-2/q,q}(S^{n-1})$, we set $u=\mathbb 
P_{-\Gd}^\Gw(\gm)$, and $\tilde u(t,\gs)=u(r,\gs)$. Then relation 
(\ref {norm-eqL}) turns into
\begin{equation}\label{norm-eqL2}
 C^{-1}\norm{\mu}_{W^{-2/q,q}(S^{n-1})}\leq \int^\infty_0\int_{S^{n-1}}
 \abs{\tl u}^q(1-e^{-t})e^{-nt}d\gs\,dt\leq C\norm{\mu}_{W^{-2/q,q}(S^{n-1})}.
\end{equation}
By density it can be assumed that $\gm$ is a regular function, and let 
$f$ be the solution of
$$\mu=\frac{(n-2)^2}{4}f-\Gd_{_{S^{n-1}}} f \quad\txt{in }\, S^{n-1}.$$
By elliptic equations regularity theory, there exists $c>0$ such that
 \begin{equation}\label{norm-eqL3}
 c^{-1}\norm{\mu}_{W^{-2/q,q}(S^{n-1})}\leq {\norm f}_{W^{2-2/q,q}(S^{n-1})}
 \leq c\norm{\mu}_{W^{-2/q,q}(S^{n-1})}.
\end{equation}
Let $v=\BBP_{-\Gd}^\Gw(f)$ in $B$ and $\tl v(t,\gs)=v(r,\gs)$. Then
 \begin {equation}\label {norm-eqL4}\left.\BA{ll}
\tl L{\tl v}:= {\tl v}_{tt}-(N-2){\tl v}_t+\Gd_{_{S^{n-1}}}{\tl 
v}=0&\mbox{in }
\BBR_+\ti S^{n-1},\\[2mm]
\phantom {\tl L{\tl v}:= {\tl v}_{tt}-(N-2){\tl v}_t=,\Gd_{_{}}}
{\tl v}|_{t=0}=f&\text{on }S^{n-1}.
\EA\right.
\end {equation}
This implies 
\begin{eqnarray}\label{norm-eqL5}
  \tl L(\Gd_{_{S^{n-1}}}s\tl v)=0\quad \text{in }\,\BBR_+\ti S^{n-1},\; 
  \text{ and }\,
  \Gd_{_{S^{n-1}}}\tl v\vline_{t=0}= \Gd_{_{S^{n-1}}} f\quad\text{on 
  }S^{n-1}.
\end{eqnarray}
This problem has a unique solution which is bounded near $t=\infty$,
therefore
\begin{equation}\label{W-adm6}
\BBP_{-\Gd}^\Gw(\Gd_{_{S^{n-1}}}f)=\Gd_{_{S^{n-1}}} \tl v,
\end{equation}
and equivalently
\begin{equation}\label{norm-eqL7}
\tl 
u=\BBP_{-\Gd}^\Gw(\mu)=\BBP_{-\Gd}^\Gw\left(\myfrac{(n-2)^2}{4}f-\Gd_{_{S^{n-1}}} f\right)
=\myfrac{(n-2)^2}{4}\tl v-\Gd_{_{S^{n-1}}} \tl v.
\end{equation}
Put $v^*:=e^{-t(N-2)/2}\tl v$, then
 \begin {equation}\label {norm-eqL8}\left.\BA{ll}
v^*_{tt}-\myfrac{(n-2)^2}{4}v^* +\Gd_{_{S^{n-1}}} v^*=0&\mbox{in }
\BBR_+\ti S^{n-1},\\[2mm]
\phantom {v^*_{tt}-\frac{(n-2)^2}{4}v^* +++,}
v^*(0,\cdot)=f&\text{on } S^{n-1}.
\EA\right.
\end {equation}
One way to represent $v^{*}$ is to introduce semigroups of linear 
operators and to express the above relations in terms of interpolation 
spaces between Banach spaces. Put
$$v^*=e^{tA}(f) \q \text{where} \q A=-\paran{\frac{(n-2)^2}{4}I-\Gd_{_{S^{n-1}}}}^{1/2}.$$
It is wellknown that the square root of a densily defined closed operator $A$ defines an analytic 
semi-group in $L^q(S^{n-1})$ (see \cite {Yo} for example). The domain of $A^2$ is precisely $W^{2,q}(S^{n-1})$. 
Therefore (see 
\cite[p. 96]{Tr}),
\begin{equation}\label{norm-eqL9}\BA{rcl}
\norm{f}^q_{W^{2-2/q,q}(S^{n-1})}&\approx& \norm{f}_{L^q(S^{n-1})}^q +
\myint{0}{\infty}\paran{t^{2/q}\norm{A^2v^*}_{L^q(S^{n-1})}}^q\dfrac{dt}{t}\\
&\approx& \norm{f}_{L^q(S^{n-1})}^q +
\myint{0}{1}\paran{t^{2/q}\norm{A^2v^*}_{L^q(S^{n-1})}}^q\dfrac{dt}{t}\\
&=& \norm{f}_{L^q(S^{n-1})}^q +
\myint{0}{1}\paran{t^{2/q}e^{-t(N-2)/2}\norm{A^2\tl 
v}_{L^q(S^{n-1})}}^q\dfrac{dt}{t},
\EA\end{equation}
where the symbol $\approx$ denotes equivalence of norms. Notice that for 
$q>1$ the exponent $2-2/q$ is an integer only if $q=2$, in which 
case the Besov and Sobolev spaces coincide. Thus, by 
\R{norm-eqL3},
\begin{equation}\label{norm-eqL10}\BA{rcl}
\norm{f}^q_{W^{2-2/q,q}(S^{n-1})}&\geq& C \norm{f}_{L^q(S^{n-1})}^q +
C\myint{0}{1}\paran{t^{2/q}e^{-t(n-2)/2}\norm{\tl u}_{L^q(S^{n-1})}}^q\dfrac{dt}{t}\\
&\geq& C \norm{f}_{L^q(S^{n-1})}^q +
C\myint{0}{1}\norm{\tl u}_{L^q(S^{n-1})}^q e^{-nt}tdt.
\EA
\end{equation}
Since $u$ is an harmonic function, 
$$r\mapsto r^{1-n}\int_{\prt B_r}|u|^q dS$$
is nonincreasing on $(0,1]$. Equivalently
$$t\mapsto \int_{S^{n-1}}|\tilde u(t,.)|^q d\gs
$$
is nonincreasing on $[0,\infty)$. Furthermore
\begin{equation}\label{norm-eqL11}\BA{rcl}
\myint{0}{\ity}\norm{\tl u}_{L^q(S^{n-1})}^q(1-e^{-t}) e^{-nt}dt &\leq& C
\myint{0}{1}\norm{\tl u}_{L^q(S^{n-1})}^q (1-e^{-t})e^{-nt}dt\\
&\leq& C\myint{0}{1}\norm{\tl u}_{L^q(S^{n-1})}^q e^{-nt}tdt.
\EA
\end{equation}
This inequality
implies that
$$\int_{|x|<1}|u|^q (1-r) \,dx\leq c(\gg) \int_{\gg<|x|<1}|u|^q (1-r)\,dx,$$
 for every $\gg\in(0,1)$. Because of (\ref{norm-eqL3}),
\begin{equation}\label{norm-eqL12}
\norm{\mu}^q_{W^{-2/q,q}(S^{n-1})} \approx \norm{f}^q_{W^{2-2/q,q}(S^{n-1})}.
\end{equation}
Therefore, the right-hand side inequality in \R{norm-eqL} follows from 
\R{norm-eqL2},
\R{norm-eqL10} and \R{norm-eqL11}.\\[1mm]
\indent Next assume that $\mu$ is a distribution on $S^{n-1}$ and
$\BBP(\mu)\in L^q(B;(1-r)\,dx)$. In order to prove that $\mu\in W^{-2/q,q}(S^{n-1})$
and that the left-hand side inequality in \R{norm-eqL} holds, we can assume that $\mu\in \GTM(S^{n-1})$.
By \R{norm-eqL3}, if $f\in L^q(S^{n-1})$ then $\mu\in W^{-2/q,q}(S^{n-1})$.
Therefore, if it is proven
\begin{equation}\label{norm-eqL13}
 \norm{f}_{L^q(S^{n-1})}\leq C \norm{u}_{L^q(B; (1-r)\,dx)},
\end{equation}
the left-hand side inequality in \R{norm-eqL} follows. Equation \R{norm-eqL7} implies that
\begin{equation}\label{est-v}
  \norm{v(r,\cdot)}_{W^{2,q}(S^{n-1})}\leq C 
  \norm{u(r,\cdot)}_{L^q(S^{n-1})},\forevery r\in(0,1).
\end{equation}
for some $C=C(n)>0$. Hence
\begin{equation}\label{proofB2}
  \norm{v}_{L^q(B; (1-r)\,dx)} + \norm{\Gd_{_{S^{n-1}}} v}_{L^q(B; (1-r)\,dx)}\leq C \norm{u}_{L^q(B; (1-r)\,dx)}.
\end{equation}
We write \R{norm-eqL4} under the form
\begin{equation}\label{proofB3}
  \begin{cases}
   {\tl v}_{tt}-(N-2){\tl v}_t = \tl h := -\Gd_{_{S^{n-1}}}{\tl v} &\text{in }\BBR_+\ti S^{n-1},\\
  {\tl v}|_{t=0}=f, & \text{in }S^{n-1}.
\end{cases}
\end{equation}
Since  $u\in L^q(B; (1-r)\,dx)$,  \R{est-v} implies that $h\in L^q(B; (1-r)\,dx)$
 (where $h(x)=\tl h(t,\gs)$). Let $\gs$ be a fixed but
arbitrary  point on $S^{n-1}$. Since Equation \R{proofB3} is a first order o.d.e.
in $\tl v_t(\cdot,\gs)$ with a forcing term $\tl h(.,\gs)$, we fix some initial time
$t_0\in (0,\infty)$
 and compute the value of the solution in $(0,t_0)$. 
Integrating twice one derives
 \begin {equation}\label {proofB4}\left.\BA{l}
\tl v(t,\gs)= \myint{t_0}{t} e^{(N-2)s}\myint{t_0}{s} e^{-(N-2)\tau}\tl h(\tau,\gs)\,d\tau\,ds\\[2mm]
\phantom {\tl L{\tl v}:= {\tl v}_{tt}-(N-2){\tl v}_t=,\Gd_{_{}}}
+\myfrac{1}{N-2}(e^{(N-2)(t-t_0)}-1)\tl v_t(t_0,\gs) + \tl v(t_0,\gs).
\EA\right.
\end {equation}
Therefore
\begin{align} 
\abs{v(0,\gs)}=\abs{f(\gs)}\leq{}& C\Big( \int_0^{t_0}\int_s^{t_0}\abs{\tl h(\tau,\gs)}\,d\tau\,ds +
\abs{\tl v_t(t_0,\gs)} + \abs{\tl v(t_0,\gs)}\Big )\notag\\
={}& C\Big ( \int_0^{t_0} s\abs{\tl h(s,\gs)}\,ds + \abs{\tl v_t(t_0,\gs)} + \abs{\tl v(t_0,\gs)}\Big )\label{proofB5}\\
\leq{}&  C\Big ( \int_{e^{-t_0}}^1(1-r)\abs{h(r,\gs)}r^{N-1}dr + \abs{\tl v_t(t_0,\gs)} + \abs{\tl v(t_0,\gs)}\Big ),\notag
\end{align}
where $C$ is a constant independent of $t_0$, for $t_0\leq\ln 2$.
Taking the $q$-power and integrating over $S^{n-1}$ yields to
\begin{align}\notag
\int_{S^{n-1}}\abs{f}^q\,d\gs \leq C\Big (& 
\int_{r_0<\abs{x}<1}\abs{h}^q(x)(1-\abs x)\,dx\\
& + \int_{S^{n-1}}\abs{v_r}^q(r_0,\gs)\,d\gs + \int_{S^{n-1}}\abs{v}^q(r_0,\gs)\,d\gs \Big ), \notag
\end{align}
where $C$ is independent of $r_0$, for $r_0\geq 1/2$. We multiply the inequality by
$r_0^{N-1}$ and
integrate  with respect to $r_0$ in $(5/8,6/8)$. It follows that
 \begin {equation}\label {proofB7}\left.\BA{l}
\myint{S^{n-1}}{}\abs{f}^q\,d\gs \leq C\Big (\myint{1/2<\abs{x}<1}{}
\abs{h}^q(x)(1-\abs x)\,dx\\[2mm]
\phantom {\tl L{\tl v}:= {\tl v}_{tt}-(N-2){\tl v}_t=,\Gd_{_{}}}
+ \myint{5/8<\abs{x}<6/8}{}\abs{v_r}^q\,dx + \myint{5/8<\abs{x}<6/8}{}\abs{v}^q\,dx \Big ).
\EA\right.
\end {equation}
By interior elliptic estimates,
\begin{equation}\label{proofB8}
\int_{5/8<\abs{x}<6/8}\abs{v_r}^q\,dx \leq \int_{1/2<\abs{x}<7/8}\abs{v}^q\,dx.
\end{equation}
Finally, by \R{proofB7}, \R{proofB8} and \R{proofB2} we obtain \R{norm-eqL13}.
\smallskip

\noindent {\it Step 2 } The case of a general operator $L$ in $B$. 
Because of the equivalence property 
of \rth {equivth}  already mentioned, if $\mu\geq 0$, there exists a constant $C$ such that, 
for every measure $\mu\in\GTM_{+}(S^{N-1})$,
\begin{equation}\label{PL,PD}
 C^{-1} \BBP^\Gw_{-\Gd} (\mu)\leq  \BBP^\Gw_L (\mu)\leq  C 
 \BBP^\Gw_{-\Gd} (\mu).
\end{equation}
Therefore, if \R{norm-eqL} holds
 with respect to $\BBP^\Gw_{-\Gd} $, it holds for $\BBP^\Gw_{L} $, for every measure 
 $\mu\in W^{-2/q,q}(S^{n-1})\cap \mathfrak M_{+}(S^{n-1})$. 
If  $\mu$ is a boundary-$q$-admissible measure for $L$,
 not necessarily positive, then $\mu_+$ and $\mu_-$ are boundary-$q$-admissible.
Therefore $\mu_+,\,\mu_-\in W^{-2/q,q}(S^{N-1})$, and the same holds 
with $\gm$. Furthermore
\begin{equation}\label{PL,PD*}
C^{-1} 
 \|\mu_{\pm}\|_{W^{-2/q,q}(\prt\Gw)}\leq\|\BBP_{L}^\Gw(\mu_{\pm})\|_{L^q(\Gw;\gr_{_{\prt\Gw}} dx)}
 \leq C \|\mu_{\pm}\|_{W^{-2/q,q}(\prt\Gw)}.
\end{equation}
\noindent {\it Step 3 } The case of a general operator $L$ in a 
general bounded $C^{2}$ domain $\Gw$. There exists a finite set of 
bounded open subdomains $U_{i}$ ($1\leq i\leq k)$ of $\mathbb R^n$ such 
that 
$$\prt\Gw\subset \bigcup_{i=1}^kU_{i},$$ 
and for each $i$ there exists 
a $C^{2}$ diffeomorphism $\Phi_{i}$ from $U_{i}$ of into some open 
subset $V_{i}$ such that $\Phi_{i}(U_{i}\cap\Gw)=B$, and 
$\Phi_{i}(U_{i}\cap\prt\Gw)=\Gg_{i}\subset \prt B\approx S^{n-1}$. This 
diffeomorphism induces an isomorphism, say $\Phi^*_{i}$, between 
$\mathfrak M(U_{i}\cap\prt\Gw)$ and $\mathfrak M(\Gg_{i})$,
$W^{-2/q,q}(U_{i}\cap\prt\Gw)$ and $W^{-2/q,q}(\Gg_{i})$, 
and it preserves positivity. Moreover, by the change of variables 
$x\in U_{i}\mapsto y=\Phi_{i}(x)\in V_{i}$, the operator $L$ is 
transformed into an elliptic operator $L_{i}^{*}$ on $B$, which still 
satisfies the maximum principle, not necessarily the condition (\ref 
{uniq}), 
but this is not crucial for the equivalence property in small domains. If  $\mu\in\GTM(\prt\Gw)$ has 
its support in $U_{i}\cap\prt\Gw$, the function $u=\mathbb 
P_{\Gw}^{L}(\gm)$ satisfies
\begin {eqnarray}
Lu&=&0\quad\qquad\mbox{ in }U_{i}\cap\Gw,\notag\\
u&=&\gm\quad\qquad\mbox{ on }U_{i}\cap\prt\Gw,\label {gene1}\\
u&=&u_{c}\quad\qquad\mbox{on }\prt U_{i}\cap\Gw,\notag
\end{eqnarray}
where $u_{c}$, the restriction of $u$ to $\prt U_{i}\cap\Gw$, is $C^{1}$. 
Thus the function $v_{i}=u\circ \Phi_{i}^{-1}$ satisfies
\begin {eqnarray}
L_{i}^{*}v_{i}&=&0\quad\quad\qquad\mbox{in }\; B,\notag\\
v_{i}&=&\Phi^{*}_{i}(\gm)\quad\quad\mbox{on }\;\Gg_{i},\label {gene2}\\
v_{i}&=&u_{c}\circ \Phi_{i}^{-1}\quad\mbox{on }\;\prt B\setminus \Gg_{i}.\notag
\end{eqnarray}
Therefore, if $\gm$ is nonnegative and $\mathbb P^{\Gw}_{L}(\gm)\in 
L^q(\Gw;\rho_{_{\prt\Gw}}dx)$, $v_{i}\in L^q(B;(1-\abs 
y)dy)$, 
which leads to $\Phi^{*}_{i}(\gm)\in W^{-2/q,q}(\Gg_{i})$ and 
$\gm\in W^{-2/q,q}(\prt\Gw)$. Moreover
 \begin {equation}\label {gene3}\left.\BA{l}
\norm{\gm}_{W^{-2/q,q}(\prt\Gw)}\approx 
\norm{\Phi^{*}_{i}(\gm)}_{W^{-2/q,q}(S^{n-1})}\leq C
\norm{v_{i}}_{L^q(B;(1-\abs x)dx)}\\[2mm]
\phantom{\norm{\gm}_{W^{-2/q,q}(\prt\Gw)}\approx 
\norm{\Phi^{*}_{i}(\gm)}_{W^{-2/q,q}(S^{n-1})}}
\approx \norm{u}_{L^q(\Gw\cap U_{i};\rho_{_{\prt(\Gw\cap U_{i})}}dx)}.
\EA\right.
 \end {equation}
Since $\rho_{_{\prt(\Gw\cap U_{i})}}\leq \rho_{_{\prt\Gw}}$ in 
$U_{i}\cap\Gw$, the integral term on the right in (\ref {gene3}) is 
dominated by the norm of $\mathbb P^{\Gw}_{L}(\gm)$ in $
L^q(\Gw;\gr_{_{\prt\Gw}}dx)$.
By using a partition of unity, any measure $\gm$ on $\prt\Gw$ can be 
decomposed in the sum of measures $\gm_{i}$ with compact support in 
$\Gg_{i}$. Hence the following estimate holds when $P_{\Gw}^{L}(\gm)\in 
L^q(\Gw;\rho_{_{\prt\Gw}}dx)$ :
\begin {eqnarray}\label {gene4}
\norm{\gm}_{W^{-2/q,q}(\prt\Gw)}\leq 
C\norm {\mathbb P^{\Gw}_{L}(\gm)}_{L^q(\Gw;\gr_{_{\prt\Gw}}dx)}.
\end{eqnarray}
Conversely, if we assume that $\gm\in\mathfrak M_{+}(\prt\Gw)\cap 
W^{-2/q,q}(\prt\Gw)$ with support in some fixed 
compact $K_{i}\subset\prt\Gw\cap U_{i}$, then 
$\Phi^{*}_{i}(\gm)\in\mathfrak M_{+}(S^{n-1})\cap 
W^{-2/q,q}(S^{n-1})$ with support in $\Gg_{i}$ and equivalence of norms. 
Then $\mathbb 
P^B_{L^*}(\Phi^{*}_{i}(\gm))\in L^q(B;(1-\abs x)dx)$, with
\begin {eqnarray}\label {gene5}
\norm{\mathbb P^B_{L^*}(\Phi^{*}_{i}(\gm))}_{L^q(B;(1-\abs x)dx)} 
\leq C\norm{\Phi^{*}_{i}(\gm)}_{W^{-2/q,q}(S^{n-1})}\approx
\norm{\gm}_{W^{-2/q,q}(\prt\Gw)}.
\end{eqnarray}
But the left-hand side term in (\ref {gene5}) is comparable to 
$\norm{\mathbb P^{U_{i}\cap\Gw}_{L}(\gm)}
_{L^q(\Gw\cap U_{i};\rho_{_{\prt(\Gw\cap U_{i})}}dx)}$, and 
\begin {eqnarray}\label {gene6}
\norm{\mathbb P^B_{L^*}(u_{c}\circ \Phi_{i}^{-1}\chi_{_{\prt B\setminus \Gg_{i}}})}_{L^q(B)} 
\approx \norm{\mathbb P^{U_{i}\cap\Gw}_{L}(u_{c})}_{L^q(U_{i}\cap\Gw)}.
\end{eqnarray}
Because $u$ is an harmonic function,
\begin {eqnarray}\label {gene7}
{\norm {u_c}}_{L^{\infty}(\prt U_{i}\cap\Gw)}\leq C\norm 
{\gm}_{W^{-2/q,q}(\prt\Gw)}.
\end{eqnarray}
Finally
\begin {equation}\label {gene8}
u=\mathbb P^{\Gw}_{L}(\gm)=\left\{\begin {array}{l}\mathbb P^{\Gw\setminus U_{i}}_{L}(\gm)+
\mathbb P^{\Gw\setminus U_{i}}_{L}(u_{c})\mbox { in }\Gw\setminus 
U_{i},_{_{_{\,}}}\\\,\\
\mathbb P^{\Gw\cap U_{i}}_{L}(u_{c})\qquad \qquad\quad\mbox { in }\Gw\cap 
U_{i}.\end {array}\right.
\end{equation}
Moreover
\begin {eqnarray*}
 &(i)&\quad \norm{\mathbb P^{U_{i}\cap\Gw}_{L}(u_{c})}_{L^q(U_{i}\cap\Gw)}
 \leq C\norm {\gm}_{W^{-2/q,q}(\prt\Gw)},
 \quad \mbox {}\\
 &(ii)&\quad 
  \norm{\mathbb P^{\Gw\setminus U_{i}}_{L}(u_{c})}_{L^q(\Gw\setminus U_{i})}
 \leq C\norm {\gm}_{W^{-2/q,q}(\prt\Gw)}.\notag 
\end{eqnarray*}
Combining these inequalities with (\ref {gene5}), (\ref{gene6}) 
yields to 
\begin {eqnarray}\label {gene9}
\norm{\mathbb P^{\Gw}_{L}(\gm)}
_{L^q(\Gw;\rho_{_{\prt\Gw}}dx)}=
\norm{u}_{L^q(\Gw;\rho_{_{\prt\Gw}}dx)}\leq 
C\norm {\gm}_{W^{-2/q,q}(\prt\Gw)},
\end{eqnarray}
and we finish the proof with the help of a partition of unity. The 
proof of (\ref {norm-eqL}) is the same as in Step 2.\qeda \medskip

\noindent \Remark By using sharp estimates on the Green kernel of a 
general elliptic operator in a general smooth domain
it can be checked directely that (\ref {norm-eqL}) is valid for any signed boundary 
$q$ admissible measure. However, it is not known if the implication
\begin {equation}\label {false}
\gm\in\mathfrak M(\prt\Gw)\cap W^{-2/q,q}(\prt\Gw)\Longrightarrow \gm
\mbox { is boundary 
{\it q}-admissible,}
\end {equation}
holds.\medskip

It is proven in \cite {MV9} that \rprop {strongA} admits an extension in the 
framework of Besov spaces $B^{-s,q}$ (see e.g. \cite {Tr}). When $s$ 
is not an integer or $q=2$, the Besov space $B^{-s,q}$ coincides with 
the Sobolev space $W^{-s,q}$.
\bprop {main}Let $s>0$, $q>1$ and $\gm$ be a distribution on $S^{n-1}$. Then 
$$\mu\in B^{-s,q}(S^{n-1})\Longleftrightarrow 
\BBP^{B}_{-\Gd} (\mu)\in L^q(B;(1-\abs x)^{sq-1}dx).
$$
Moreover there exists a constant $C>0$ such that for any 
$\gm\in B^{-s,q}(S^{n-1})$,
\begin {equation}\label {equivnorm} C^{-1}{\norm \mu}_{B^{-s,q}(S^{n-1})}
\leq \left(\int_{B}{\abs {\BBP^{B}_{-\Gd} (\mu)}}^q
(1-\abs x)^{sq-1}dx\right)^{1/q}
\leq C{\norm \mu}_{B^{-s,q}(S^{n-1})}.
\end {equation}
\es\

The dual form of \rprop {strongA} is the following,
\bprop{strongB} Let $q\geq (n+1)/(n-1)$ and the assumptions on $L$ 
and $\Gw$  be satisfied as in \rprop {strongA}. Then
 $$\varphi\in L^{q'}(\Gw;\gr_{_{\prt\Gw}}^{-q'/q}dx)\Longleftrightarrow
 \myfrac {\partial}{\partial{\bf n}_{L^{*}}}\BBG_{L^{*}}^\Gw(\varphi)\in W^{2/q,q'}(\prt\Gw).
 $$
 Moreover there exists a constant $C>0$ such that, for any 
 $\varphi\in L^{q'}(\Gw;\gr_{_{\prt\Gw}}^{-q'/q})dx)$,
 \begin{equation}\label{dual5}
     C^{-1}{\norm 
     {\varphi}}_{L^{q'}(\Gw;\gr_{_{\prt\Gw}}^{-q'/q}dx)}\leq 
{\norm {\myfrac {\partial}{\partial{\bf n}_{L^{*}}}\BBG_{L^{*}}^\Gw(\varphi)}}_{W^{2/q,q'}(\prt\Gw)}
     \leq C{\norm {\varphi}}_{L^{q'}(\Gw;\gr_{_{\prt\Gw}}^{-q'/q}dx)}.
  \end{equation}    
 \es
\Proof Let $\gm\in\mathfrak M (\prt\Gw)$. By duality between $L^q(\Gw;\gr_{_{\prt\Gw}}dx)$ and 
$L^{q'}(\Gw;\gr_{_{\prt\Gw}}dx)$, we write
\begin{equation}\label{dual1}
\myint{\Gw}{}\BBP_{L}^\Gw(\mu)\psi \gr_{_{\prt\Gw}}dx=
\myint{\Gw}{}\BBP_{L}^\Gw(\mu) L^{*}\zeta dx=-\myint{\prt\Gw}{}\myfrac{\partial 
\zeta}{\partial {\bf n}_{L^{*}}}d\gm,
\end{equation}
where $\zeta=\BBG_{L^{*}}^\Gw(\psi \gr_{_{\prt\Gw}})$. 
Then the adjoint operator $\left[\BBP_{L}^\Gw\right]^{*}$ of $\BBP_{L}^\Gw$ is defined by
\begin{equation}\label{dual3}
   \left[\BBP_{L}^\Gw\right]^{*}(\psi)=-\myfrac 
   {\partial}{\partial{\bf n}_{L^{*}}}\BBG_{L^{*}}^\Gw(\psi \gr_{_{\prt\Gw}}).
 \end{equation}   
Consequently, \rprop {strongA} implies that 
 there exists a constant $C>0$ such that
\begin{equation}\label{dual4}
    C^{-1}\norm {\psi}_{L^{q'}(\Gw;\gr_{_{\prt\Gw}}dx)}
    \leq \norm {\myfrac 
   {\partial}{\partial{\bf n}_{L^{*}}}\BBG_{L^{*}}^\Gw(\gr_{_{\prt\Gw}}\psi)}_{W^{2/q,q'}(S^{n-1})}
    \leq  C\norm {\psi}_{L^{q'}(\Gw;\gr_{_{\prt\Gw}}dx)}.
 \end{equation} 
 But 
 $$\psi\in L^{q'}(\Gw;\gr_{_{\prt\Gw}}dx)\Longleftrightarrow \gr_{_{\prt\Gw}}\psi\in 
 L^{q'}(\Gw;\gr_{_{\prt\Gw}}^{(1-q')}dx).
 $$
 Putting $\varphi=\gr_{_{\prt\Gw}}\psi$, implies 
 (\ref {dual5}).\qeda 
\medskip


\noindent {\it Proof of \rth {Bq-powerth}.} (i) Assume that $u$ is a solution 
 of (\ref {Bq-power}). Then $u\in 
 L^{q}(\Gw;\gr_{_{\prt\Gw}}dx)$, and for any $\gz\in 
 C_{c}^{1,L}(\overline\Gw)$, there holds
\begin {eqnarray}
 \abs {\int_{\prt\Gw}\frac {\prt\gz}{\prt{\bf n}_{L^{*}}}d\gm}&=&\abs 
 {\int_{\Gw}\left(uL^{*}\gz+\gz\abs u^{q-1}u \right)dx},\notag\\
 &\leq& {\norm {u}_{L^q (\Gw;\rho_{_{\prt\Gw}}dx)}}{\norm {L^{*}\gz}}_{L^{q'} 
 (\Gw;\rho^{-q'/q}_{_{\prt\Gw}}dx)}
+\int_{\Gw}\abs u^q\abs\gz dx,\label {daul}\\
  &\leq& {\norm {u}_{L^q (\Gw;\rho_{_{\prt\Gw}}dx)}} 
{\norm {\myfrac {\partial\gz}{\partial{\bf n}_{L^{*}}}}}_{W^{2/q,q'}(\prt\Gw)}
  +\int_{\Gw}\abs u^q\abs\gz dx,\notag
\end {eqnarray}
since $\BBG_{L^{*}}^\Gw(L^*\gz)=\gz$. Let $\eta\in W^{2/q,q'}\prt\Gw$, 
and, for $\gd>0$, put 
$\gz=\gd^{-2}\gr_{_{\prt\Gw}}(\gd-\gr_{_{\prt\Gw}})_{+}^2\mathbb 
P_{L}^\Gw(\eta)$,
\begin {eqnarray}\label {dau2}
 \abs {\int_{\prt\Gw}\eta d\gm}\leq 
 {\norm {u}_{L^q (\Gw;\rho_{_{\prt\Gw}}dx)}}
{\norm {\eta}}_{W^{2/q,q'}(\prt\Gw)}
  + \gd^{-2}\int_{\Gw}
 \gr_{_{\prt\Gw}}(\gd-\gr_{_{\prt\Gw}})_{+}^2\abs{\mathbb P_{L}^\Gw(\eta)}
 \abs u^q dx.
\end {eqnarray}
Let $K\subset\prt\Gw$ be a compact subset such that $C_{2/q,q'}(K)=0$. 
Then there exists a sequence $\{\eta_{n}\}\subset W^{2/q,q'}(\prt\Gw)$ 
with the property that $0\leq\eta_{n}\leq 1$, $\eta_{n}\equiv 1$ in a 
neigborhood of $K$ and $\eta_{n}\to 0$ in $W^{2/q,q'}(\prt\Gw)$ as 
$n\to\infty$. We take $\eta=\eta_{n}$ in (\ref {dau2}). Since 
$u\in L^q(\Gw;\rho_{_{\prt\Gw}}dx)$ and $K$ has measure zero,  the two 
terms in the right-hand side of (\ref {dau2}) converge to $0$ when 
$n\to\infty$.
Thus $\gm$ does not charge Borel subsets with $C^{2/q,q'}$-capacity zero.
It follows that $\gm$ is the sum of an integrable function and a mesure 
in $W^{-2/q,q}(\prt\Gw)$, by \rcor{decomposition}.
 \medskip 
 
 \noindent (ii) Conversely, let $\gm$ be a boundary measure which does 
 not charge Borel subsets with $C_{2/q,q'}$-capacity zero. 
 Assuming first that $\gm$ is positive, by \rprop {BP1lem2} there exists an increasing 
 sequence $\{\gm_{j}\}$ of elements of $W^{-2/q,q}(\prt\Gw)\cap\mathfrak 
 M_{+}(\prt\Gw)$ which converges to $\gm$. By \rprop{strongA}, the 
 $\gm_{j}$ are boundary-$q$-admissible and the sequence $\{u_{j}\}$ of solutions 
 of  
 \begin {equation}\label {B-qpower-j}\left.\BA{ll}
Lu_{j}+\abs {u_{j}}^{q-1}u_{j}=0&\mbox{in }\;\Gw,\\[2mm]
\phantom {Lu_{j}+\abs {u_{j}}^{q-1}}
\,u_{j}=\gm_{j}&\mbox{on }\;\prt\Gw,
\EA\right.
\end {equation}
 is increasing. Moreover $u_{j}\geq 0$. If $u=\lim_{j\to\ity}u_{j}$, 
 then $u\geq 0$. Since
  \begin {eqnarray}\label {B-qpower-k}
\int_{\Gw}\left(u_{j}L^{*}\gz+u_{j}^q\gz\right)dx=-\int_{\prt\Gw}\frac 
{\prt\gz}{\prt{\bf n}_{L^{*}}}d\gm_{j},
 \end {eqnarray}
 for any $\gz\in C^{1,L}_{c}(\overline\Gw)$. Taking $\gz=\eta_{1}$, the 
 solution of 
  \begin {eqnarray*}
 L^{*}\eta_{1}&=&1\quad \mbox {in }\Gw,\\
\eta_{1}&=&0\quad \mbox {on }\prt\Gw,\notag
 \end {eqnarray*}
 we deduce $u\in L^{1}(\Gw)\cap L^q(\Gw;\gr_{_{\prt\Gw}}dx)$ by the 
 monotone convergence Theorem. Therefore (\ref {B-qpower-k}) implies 
 that $u$ is the solution of
  \begin {equation}\label {B-qpower-l}\left.\BA{ll}
Lu+\abs {u}^{q-1}u=0&\mbox{in }\;\Gw,\\[2mm]
\phantom {Lu+\abs {u}^{q-1}}
\,u=\gm&\mbox{on }\;\prt\Gw.
\EA\right.
\end {equation}
 If $\gm$ is a signed measure, we procede as in the proof of \rth 
 {BP1th}, 
 by truncating the nonlinearity and inroducing the solutions of (\ref {Bq-power}) associated to $\gm_{+}$ 
 and $-\gm_{-}$ on the boundary.\qeda 
\subsection {Boundary singularities}
\subsubsection {Isolated singularities}
The study of boundary singularities of solutions of semilinear 
elliptic equations started with the work of Gmira and V\'eron \cite {GmV}. 
As in the case of equations with internal singularities, the starting 
idea is to study the model case where $\Gw=\mathbb R^n_{+}$, 
$\prt\Gw=\prt\mathbb R^n_{+} \approx \mathbb R^{n-1}$ and the singularity is located at $x=0$. In 
spherical coordinates $x=(r,\gs)$ where $r>0$, $\gs\in S^{n-1}$, the 
existence of a solution $u$ to
\begin {equation}\label {model1}
-\Gd u+\abs u^{q-1}u=0,
\end {equation}
in $\mathbb R^n_{+}$ ($q>1$) which vanishes on $\prt\mathbb 
R^n_{+}\setminus \{0\}$ is enlighted if we look for it under the 
separable form $u(r,\gs)=r^{\ga}\gw(\gs)$. Then $\ga=-2/(q-1)$ and 
$\gw$ is a solution of 
\begin {equation}\label {model2}
-\Gd _{_{S^{n-1}}}\gw-
\left(\frac {2}{q-1}\right)\left(\frac {2q}{q-1}-n\right)\gw+\abs\gw^{q-1}\gw=0\quad\mbox {on 
}\;S^{n-1}_{+}=S^{n-1}\cap\mathbb R^n_{+},
\end {equation}
which vanishes on the equator $\prt S^{n-1}_{+}\approx S^{n-2}$. 
Since the first nonzero eigenvalue of the Laplace-Beltrami operator 
in $W^{1,2}_{0}(S^{n-1}_{+})$ is $n-1$, it is clear, by multiplying 
(\ref {model2}) by $\gw$ and integrating over $S^{n-1}_{+}$, that no nontrivial solution 
of (\ref{model2}) exists whenever $(2/(q-1))(2q/(q-1)-n)\leq n-1$. 
Equivalently $q\geq (n+1)/(n-1)$. Conversely, if $(2/(q-1))(2q/(q-1)-n)< 
n-1$ solutions to (\ref {model2}) exist. The stable solutions are 
obtained by minimizing the functional
\begin {equation}\label {model2'}
\eta\mapsto J(\eta)=\int_{S^{n-1}_{+}}\left(\frac {1}{2}\abs{\nabla_{_{S^{n-1}}} \eta}^{2}
-\left(\frac {1}{q-1}\right)\left(\frac 
{2q}{q-1}-n\right)\eta^{2}+\frac {1}{q+1}\abs\eta^{q+1}\right) d\gs,
\end {equation}
over the space $W^{1,2}_{0}(S^{n-1}_{+})$, where $\nabla_{_{S^{n-1}}}$ 
denotes the covariant derivative identified with the tangential 
gradient 
thanks to the isometrical imbedding of $S^{n-1}$ into $\mathbb R^n$. 
Put
$$\Gl_{q,n}=\left(\left(\frac {2}{q-1}\right)\left(\frac 
{2q}{q-1}-n\right)\right)=\ell^{q-1}_{q,n},
$$
and let $\mathcal S_{+}$ be the set of solutions of (\ref{model2}) in 
$S^{n-1}_{+}$ which vanishes on $\prt S^{n-1}$. As we have already 
seen it, if $q\geq (n+1)/(n-1)$ this set is reduced to $\{0\}$. 
Conversely, if  
$1<q< (n+1)/(n-1)\Longleftrightarrow \Gl_{q,n}>n-1$, there exist 
minimizing solutions to (\ref{model2}). Besides this fact, the 
positive solutions are unique. Moreover, if $\Gl_{q,n}\leq 2n$, 
which is the second eigenvalue of $\Gd_{_{S^{n-1}}}$
in $W^{1,2}_{0}(S^{n-1}_{+})$, all the solutions of (\ref{model2}) 
vanishing on the equator have constant sign. Finally, if  $\Gl_{q,n}> 
2n$ there exist changing sign solutions. \medskip

Let $\Gw$ be an open subset 
of $\mathbb R^n$ with a boundary of class $C^{2,\gth}$ for some 
$\gth\in (0,1)$, and $0\in\prt\Gw$. It can be performed an orthogonal 
change of coordinates in $\mathbb R^n$ in order the axis 
$\{x:\,x_{j}=0,\;\forall j=1,...,n-1\} $ be the normal direction to $\prt\Gw$, 
${\bf e}_{n}$ be the unit outward normal vector at $0$ and 
$\prt\mathbb R^n_{+}\approx \mathbb R^{n-1}$ the tangent plane to 
$\prt\Gw$ at $0$. Let $u$ be any solution to (\ref {model1}) in $\Gw$ 
which is continuous in $\overline\Gw\setminus\{0\}$ and coincides on 
$\prt\Gw\setminus\{0\}$ with a function $g\in C(\prt\Gw)$. For 
$R>0$ small enough and $m_{+}=\max\{g(x):\,x\in\prt\Gw\cap B_{R}\}$, 
the function 
$$x\mapsto \tilde u(x)=\left\{\BA {l} (u(x)-m_{+})_{+}\quad \mbox {if }x\in 
\Gw\cap B_{R},\\[2mm]
0\quad\qquad\qquad \quad \mbox {if }x\in B_{R}\setminus\overline\Gw,
\EA\right.$$ 
is a subsolution of (\ref {model1}) in $B_{R}\setminus\{0\}$. But the 
Keller-Osserman estimate implies
$$u(x)\leq m_{+}+C\abs x^{-2/q-1},\forall x\in \overline\Gw\cap 
B_{R}\setminus\{0\},
$$
for some $C=C(n,q,R)>0$. In the same way, $u$ is bounded from below 
in the same set by $m_{-}-C\abs x^{-2/q-1}$, where 
$m_{-}=\min\{g(x):\,x\in\prt\Gw\cap B_{R}\}$. Hence the function 
$x\mapsto \abs x^{2/q-1}u(x)$ is uniformly bounded in 
$\overline\Gw\cap B_{R}\setminus\{0\}$. We perform a change of coordinates
$y=\phi (x)$ which transforms $\Gw\cap B_{R}$ into
$\mathbb R_{+}^n\cap B_{R}$ and $\prt\Gw\cap B_{R}$ into 
$R^{n-1}\cap B_{R}$. We define $v$ by
$$y\mapsto v(y)=v(r,\gs)=\abs{\phi^{-1}(y)}^{2/(q-1)}u(\phi^{-1}(y)), 
\quad (r,\gs)\in (0,R)\times S_{+}^{n-1},
$$
and put $w(t,\gs)=v(r,\gs)$ with $t=\ln r$. Then $w$ satisfies an equation of the type
  \begin {equation}\label {model3}\left.\BA{ll}
0=(1+\ge_{1}(t))w_{tt}+\left(n-2\myfrac {q+1}{q-1}+\ge_{2}(t)\right)w_{t}
+\left(\Gl_{q,n}+\ge_{3}(t)\right)w+\Gd_{_{S^{n-1}}}w\\[2mm]
+\langle\nabla_{_{S^{n-1}}}w.\ge_{4}(t)\rangle+
+\langle\nabla_{_{S^{n-1}}}w_{t}.\ge_{5}(t)\rangle
+\langle\nabla_{_{S^{n-1}}}\langle\nabla_{_{S^{n-1}}}w.{\bf e}_{n}\rangle.\ge_{6}(t)\rangle
+\abs w^{q-1}w\EA\right.
\end {equation}
in $(-\ity,\ln R]\times S_{+}^{n-1}$, where the $\ge_{j}(t)$ depend 
on the change of coordinates and verify
\begin {equation}\label {model4}
\abs {\ge_{j}(t,\gs)}\leq C_{j}e^{t},\forevery 
(t,\gs)\in (-\ity,\ln R]\times S_{+}^{n-1},\quad j=1,...,6.
\end {equation}
Since $\abs {w(t,\gs)}\leq Ce^{2qt/(q-1)}$, we can use the elliptic 
equations regularity theory and a Lyapounov style analysis at 
$-\infty$. The following result is due to Gmira and V\'eron \cite 
{GmV}.
\bth {GmVth1}Suppose $1<q<(n+1)/(n-1)$. Then, with the previous notations, 
there exists a compact connected subset $\CE_{+}$ of the set of the solutions of (\ref{model2}) in 
$S^{n-1}_{+}$ which vanish on $\prt S_{+}^{n-1}$, such that
\begin {equation}\label {model5}
\lim_{t\to-\infty}\rm{dist}_{C^2(S_{+}^{n-1})}(w(t,.),\CE_{+})=0,
\end {equation}
where dist$_{C^2(S_{+}^{n-1})}$ denotes the distance associated with 
the $C^2(S_{+}^{n-1})$-norm. Moreover, the set $\CE_{+}$ is reduced 
to a singleton in the following cases :\smallskip

\noindent (i) $u$ is nonnegative,\smallskip

\noindent (ii) $(n+2)/2n\leq q<(n+1)/(n-1)$,\smallskip

\noindent (iii) $n=2$.
\es

When $\CE_{+}=\{0\}$ it is possible to make more precise the way the 
function $w(t,.)$ converges to $0$ as $t\to-\infty$. By adapting the 
method developed in \cite {CMV}, it is proven in \cite {GmV} that the 
following result holds,
\bth {GmVth2} Suppose $1<q<(n+1)/(n-1)$ and let $w$ be the solution of 
(\ref{model3}) associated to $u$, solution of (\ref {model1}). Assume
$$\lim_{t\to-\infty}\norm {w(t,.)}_{C^2(S_{+}^{n-1})}=0.$$
Then, if one of the following conditions holds :\smallskip

\noindent (a) $u$ is nonnegative,\smallskip

\noindent (b) $n=2$ and $\prt\Gw$ is locally a straight line near $0$,\smallskip

\noindent (c) $2/(q-1)$ is not an integer,\smallskip

\noindent (i) either $u$ can be extended to $\overline\Gw$ as continuous 
function solution of the Dirichlet problem 
  \begin {equation}\label {model6}\left.\BA{ll}
-\Gd u+\abs u^{q-1}u=0&\mbox {in }\Gw,\\[2mm]
\phantom {-\Gd u+\abs u^{q-1}}
\,u=g&\mbox {on }\prt\Gw,
\EA\right.
\end {equation}
(ii) or there exists an integer $k\in[n-1,2/(q-1))$ and a nonzero 
solution $\psi$ of 
  \begin {equation}\label {model7}\left.\BA{ll}
\Gd_{_{S^{n-1}}}\psi +k(n+k-2)\psi =0&\mbox {in }S^{n-1}_{+},\\[2mm]
\phantom {\Gd_{_{S^{n-1}}}\psi +k(n+k-2)}
\,\psi=0&\mbox {on }\prt S^{n-1}_{+},
\EA\right.
\end {equation}
such that
\begin {eqnarray}\label {model8}
\lim_{t\to-\infty}e^{(k-2/(q-1))t}w(t,.)=\psi,
\end {eqnarray}
in the $C^2(S_{+}^{n-1})$-topology.
\es 

The meaning of this result is the following : either $u$ has a strong 
boundary singularity which is described thanks to the set $\CS_{+}$ 
of solutions of (\ref {model2}) vanishing on the equator, 
either there exists a spherical harmonic of 
degree $k$ such that 
\begin {eqnarray}\label {model9}
\lim_{\tiny{\begin {array} {l}\;\;x\to 0\\ x/\abs x\to \gs\end {array}}}\abs x^ku(x)=\psi (\gs),
\quad\mbox {uniformly for }\gs\in S_{+}^{n-1},
\end {eqnarray}
or $u$ is regular function. \medskip

When $-\Gd$ is replaced by an elliptic operator $L$ with variable 
Lipschitz continuous coefficients, most of the above results extend in the same way as for 
the isolated internal singularities (see the section on isolated 
singularities).

\subsubsection {Removable singularities}
The first result on removability (see \cite {GmV}) is the following.
\bth {GmVth}Let $\Gw$ be a $C^{2}$ domain in $\mathbb R^n$, $x_{0}$ a boundary 
point, and $g$ a continuous real valued function defined on $\Gw\times 
\mathbb R$, such 
that
\begin {equation}\label {GVcond}
\liminf_{r\to\infty}\frac {g(x,r)}{r^{(n+1)/(n-1)}}>0\quad\mbox {and 
}\;\limsup_{r\to\infty}\frac {g(x,r)}{\abs r^{(n+1)/(n-1)}}<0,\forevery 
x\in\Gw,
\end {equation}
uniformly with respect to $x\in\Gw$. If $u\in C^{2}(\Gw)\cap C(\overline\Gw\setminus\{x_{0}\})$ is a solution 
of 
\begin {equation}\label {Bsing1}
-\Gd u+g(x,u)=0\quad\mbox {in }\;\Gw,
\end {equation}
which coincides on $\prt\Gw\setminus\{x_{0}\}$ 
with some $\phi\in C(\prt\Gw)$, then $u$ can be extended as a $C(\overline\Gw)$ 
function, which verifies
  \begin {equation}\label {Bsingphi}\left.\BA{ll}
-\Gd u+g(x,u)=0&\mbox {in }\;\Gw,\\[2mm]
\phantom {-\Gd u+g(x,)}
u=\phi&\mbox {on }\;\prt\Gw.
\EA\right.
\end {equation}
\es 
Actually, their proof could have been adapted, without any deep 
modification, to Equation (\ref {Bsing}). A much more general result 
will be given later on. 

\bdef {Rem}\rm {Let $\Gw$ be a $C^{2}$ domain in $\mathbb R^n$ and 
$q\geq (n+1)/(n-1)$. \smallskip

\noindent (i) A Borel subset $K$ of $\prt\Gw$ is said  
$q$-removable if any nonnegative function $u\in C^{2}(\Gw)\cap 
C({\overline\Gw\setminus K})$  solution of (\ref {model1})
which vanishes on $\prt\Gw$ is identically zero.\smallskip

\noindent (ii) A Borel subset $K$ of $\prt\Gw$ is said  
conditionally $q$-removable if any nonnegative function $u\in C^{2}(\Gw)\cap 
C({\overline\Gw\setminus K})$  solution of (\ref {model1}) belongs to 
$L^q_{loc}(\overline\Gw;\rho_{_{\prt\Gw}}dx)$.}
\es

The condition $q\geq (n+1)/(n-1)$ is necessary, since, below this 
value, only the empty set is removable by \rth {Bg-adm-th}. The main 
removability result is the folllowing,

\bth {Bremovth}Let $\Gw$ be a $C^{2}$ bounded domain in $\mathbb R^n$, 
$q\geq (n+1)/(n-1)$ and $K\subset\Gw$ be compact. Then the following assertions are 
equivalent.\smallskip

\noindent (i) $K$ is $q$-removable.\smallskip

\noindent (ii)  $K$ is conditionally $q$-removable.\smallskip

\noindent (iii) $C_{2/q,q'}(K)=0$. 
\es

This result was first proven by Le Gall \cite {LG3} in the case 
$q=2$, by probabilistic methods, then by Dynkin and Kuznetsov \cite 
{DK1} in the 
case $q\leq 2$, by a combination of analytic and probabilistic methods 
and by Marcus and V\'eron \cite {MV5} when $q>2$ with purely analytic 
tools. All the proof are based upon the construction of suitable 
lifting operators which transform functions defined on the boundary into 
functions defined in $\Gw$.  In \cite {MV6} the first unified proof, 
valid in all the cases $q\geq (n+1)/(n-1)$ is given. We shall present 
a sketch of it below.
\bdef {lift}\rm { A linear map $R: C^{2}(\Gw)\mapsto C^{2}(\overline\Gw)$ 
is called a positive lifting if 
\begin {equation}
R(\eta)\vline_{\prt\Gw}=\eta\quad\mbox {and }\;\eta\geq 
0\Longrightarrow R(\eta)\geq 0.
\end {equation}
}\es
\blemma{liftlem} Let $\phi$ be the first eigenfunction of $-\Gd$ 
in $W^{1,2}_{0}(\Gw)$ and $q\geq (n+1)/(n-1)$. There exists a positive lifting operator 
$R:\eta\mapsto R(\eta)=R_{\eta}$ 
with the additional property 
$$\norm{R_{\eta}}_{L^\ity(\Gw)}\leq 
\norm{\eta}_{L^\ity(\prt\Gw)},
$$
and 
\begin {equation}\label{Leta}
\norm{\abs {\phi^{1/q'}\Gd R_{\eta}} +2\abs 
{\phi^{-1/q}\langle\nabla R_{\eta}.\nabla\phi\rangle}}_{L^{q'}(\Gw)}
\leq C{\norm 
{\eta}}_{W^{2/q,q'}(\prt\Gw)},\forevery\eta\in W^{2/q,q'}(\prt\Gw).
\end {equation}
Furthermore 
  \begin {equation}\label {Meta}\left.\BA{ll}
\norm{\abs {\phi^{1/q'} R_{\eta}\Gd R_{\eta}} +2\abs 
{\phi^{-1/q}R_{\eta}\langle\nabla R_{\eta}.\nabla\phi\rangle}
+\phi^{1/q'}\abs {\nabla{R_{\eta}}}^{2}}_{L^{q'}(\Gw)}\\[2mm]
\phantom {\abs {\phi^{1/q'} R_{\eta}\Gd R_{\eta}} +2\abs 
{\phi^{-1/q}R_{\eta}\langle\nabla R_{\eta}.\nabla\phi\rangle}}
\leq C(1+{\norm 
{\eta}}_{W^{2/q,q'}(\prt\Gw)}),\forevery\eta\in \CT^{*},
\EA\right.
\end {equation}
where $\CT^{*}=\{\eta\in W^{2/q,q'}(\prt\Gw): 0\leq \eta\leq 1\}$.
\es
\Proof In Section 2.4 we have already introduced the foliation of 
$\prt\Gw$ by the $\Gs_{\gb}$
$$
\Gs_{\gb}:=\{x\in\Gw:\gr_{_{\Gw}}(x)=\gb\},\quad 0<\gb\leq\gb_{0},
$$
for $\gb_{0}$ depending on the curvature of $\prt\Gw$, with $\Gs_{0}=\Gs=\prt\Gw$, 
$\Gw_{\gb}=\{x\in\Gw:\rho_{_{\prt\Gw}}(x)>\gb\}$ and 
$G_{\gb}=\Gw\setminus\overline\Gw_{\gb}$.
For every $0<\gb\leq\gb_{0}$ and $x\in G_{\gb}$ 
there exists a unique $\gs(x)\in\Gs$ such that $\abs 
{x-\gs(x)}=\gr_{_{\prt\Gw}}(x)$, and the correspondence 
$x\longleftrightarrow (\gr_{_{\prt\Gw}}(x),\gs(x))$ defines a smooth change of 
coordinates near the boundary called the flow coordinates. In terms of 
flow ccordinates, the Laplacian has the following form
$$\Gd=\frac {\prt^{2}}{\prt\gr^{2}}+b_{0}\frac 
{\prt}{\prt\gr}+\Gl_{\Gs},
$$
where $\gr$ stands for $\gr_{_{\prt\Gw}}$, $b_{0}$ depends on $x$ 
and $\Gl_{\Gs}$ is a second order elliptic operator on $\Gs$ with 
coefficients depending also on $x$. Moreover
$$\Gl_{\Gs}\to\Gd_{\Gs}\quad\mbox {and }\;b_{0}\to\gk\quad\mbox {as 
}\;\gr_{_{\prt\Gw}}(x)\to 0,
$$
where $\Gd_{\Gs}$ is the Laplace-Beltrami operator on $\Gs$, and $\gk$ 
the mean curvature of $\Gs$ (see \cite {BK}). If $\eta\in C(\Gs)$, 
let $H=H_{\eta}$ be the solution of the initial value problem
\begin{equation}\label{defH}\BA{ll}
\dfrac{\prt H}{\prt\gt}=\Gd_{\Gs}H & \text{in } \BBR_+\ti\Gs,\\[2mm]
H(0,\cdot)=\eta(\cdot) & \text{in } \Gs.
\EA
\end{equation}
We can express $H$ in terms of the two coordinates $(\gt,\gs)$. Let 
$h\in C^{\ity}(\mathbb R_{+})$ be a truncation function with value in 
$[0,1]$, $h\equiv 1$, on $[0,\gb_{0}/2]$ and $h\equiv 0$, on 
$[\gb_{0},\ity)$. The lifting $R=R_\eta$ of $\eta$ is defined by
\begin{equation}\label{goodlift}
R_\eta(x)=\begin{cases}
H_\eta(\phi^2(x),\gs(x))h(\gr_{_{\prt\Gw}}(x)), &\forevery x\in G_{\gb_{0}}\\[2mm]
0, &\forevery x\in \Gw_{\gb_0}. \end{cases}
\end{equation}
Clearly the positivity and contraction principle in uniform norms hold. 
The proof of (\ref {Leta}) and (\ref {Meta}) is much more elaborated and 
settled upon analytic semigroups theory and delicate interpolation results (see \cite {MV6} for a 
detailled proof).\qeda \medskip

\noindent {\it Proof of \rth {Bremovth}. }(iii)$\Longrightarrow$ (ii) 
Let 
$$\CT_{K}=\{\eta\in C^{2}(\prt\Gw):0\leq \eta\leq 1,\,\eta\equiv 
0\mbox{ in an open relative neighborhood of } K\}.$$
Put $\gz_{\eta}:=\phi R_{\eta}^{2q'}$. Then $0\leq \gz\leq \phi$, and 
$\gz_{\eta}(x)=O\left((\gr_{_{\prt\Gw}}(x))^{1+2q'}\right)$ in a 
neighborhood $V_{\eta}$
of $K$. Since in the case of 
Equation (\ref {model1}), the Keller-Osserman {\it a priori} bound implies
\begin{equation}\label{KOdist}
\abs{u(x)}\leq C(N,q)(\rho_{_{\prt\Gw}}(x))^{-2/(q-1)},\forevery 
x\in\Gw,
\end{equation}
and $u(x)=O\left(\gr_{_{\prt\Gw}}(x)\right)$ if 
$\gr_{_{\prt\Gw}}(x)\to 0$, outside $V_{\eta}$,
 we derive
\begin{equation}\label{Odist}
u^q(x)\gz_{\eta}(x)=O\left(\gr_{_{\prt\Gw}}(x)\right)\quad \mbox 
{in }\,\Gw.
\end{equation}
Moreover, if $\gl_{1}$ is the eigenvalue corresponding to $\phi$, 
  \begin {equation}\label {Odist2}\left.\BA{ll}
\Gd \gz_{\eta}&=-\gl_{1}\phi R_{\eta}^{2q'}+\phi\Gd 
R_{\eta}^{2q'}+2\langle\nabla \phi.\nabla R_{\eta}^{2q'}\rangle\\[2mm]
\phantom {} 
&=-\gl_{1}\gz_{\eta}+2q'\phi R_{\eta}^{2q'-1}\Gd 
R_{\eta}+2q'(2q'-1)R_{\eta}^{2q'-2}\abs{\nabla {R_{\eta}}}^{2}
+2q'R_{\eta}^{2q'-1}\langle\nabla \phi.\nabla R_{\eta}\rangle.
\EA\right.
\end {equation}
Therefore
$$u\abs{\Gd\gz_{\eta}}\leq C(\eta)uR_{\eta}^{2q'-2}.
$$
Because $\eta\in\CT_{K}$, $u\Gd\gz_{\eta}$ remains bounded in $\Gw$. 
For $0<\gb\leq\gb_{0}$, 
\begin{equation}\label{int-gb}
  \int_{\Gw\setminus G_{\gb}}\gz_\eta\Gd u\,dx= \int_{\Gw\setminus G_{\gb}}u\Gd\gz_\eta\,dx +
\int_{\Gs_\gb}\Big(\gz_\eta\frac{\prt u}{\prt{\bf n}} - u \frac{\prt \gz_\eta}{\prt{\bf n}}\Big)\,dS,
\end{equation}
and combining (\ref {KOdist}) with Schauder estimates, 
$$\frac{\prt u}{\prt{\bf n}}
\vline_{\Gs_{\gb}}=O(\gb^{-(q+1)/(q-1)}),
$$
hence
$$\Big(\gz_\eta\frac{\prt u}{\prt{\bf n}} - u \frac{\prt \gz_\eta}{\prt{\bf n}}\Big)\vline_{\Gs_{\gb}}=O(\gb).
$$
Letting $\gb\to 0$ in (\ref{int-gb}) implies 
\begin{equation}\label{int-0}
   \int_{\Gw}\Big(-u\Gd\gz_\eta +u^q \gz_\eta\Big)dx =0.
\end{equation}
By H\"older's inequality,
  \begin {equation}\label {est-t2}\left.\BA{ll}
\myint{\Gw}{}u \abs{\Gd\gz_\eta}\, dx \leq 
\left(\myint{\Gw}{} u^q \gz_\eta\,dx\right)^{1/q}
\left(\myint{\Gw}{}  \gz_\eta^{-q'/q}\abs{\Gd\gz_\eta}^{q'}\,dx\right)^{1/q'}\\[2mm]
\phantom {\myint{\Gw}{}u \abs{\Gd\gz_\eta}\, dx } 
\leq c\left(\myint{\Gw}{} u^q \gz_\eta\,dx\right)^{1/q}
\left(\myint{\Gw}{} (\gz_\eta +  M(\eta)^{q'})\,dx\right)^{1/q'}, 
\EA\right.
\end {equation}
where 
$$M(\eta)=\abs {\phi^{1/q'} R_{\eta}\Gd R_{\eta}} +2\abs 
{\phi^{-1/q}R_{\eta}\langle\nabla R_{\eta}.\nabla\phi\rangle}
+\phi^{1/q'}\abs {\nabla{R_{\eta}}}^{2}.
$$
Since by \rlemma {liftlem},
$$\norm{M(\eta)}_{L^{q'}(\Gw)}\leq 
C_{1}(1+\norm\eta_{W^{2/q,q'}(\prt\Gw)}),
$$
it follows from (\ref {int-0}) and (\ref {est-t2}),
\begin{equation}\label{int- u^q}
   \int_{\Gw}u^q \gz_\eta dx \leq 
   C_{2}(1+\norm\eta_{W^{2/q,q'}(\prt\Gw)})^{q'}.
\end{equation}
If we put $\eta^{*}=1-\eta$, then 
$\norm{\eta}^{q'}_{W^{2/q,q'}(\prt\Gw)}\leq 
C'+\norm{\eta^{*}}^{q'}_{W^{2/q,q'}(\prt\Gw)}$. If $K$ has 
$C_{2/q,q'}$-capacity zero, there exists a sequence 
$\{\eta^{*}_{n}\}\subset C^{2}(\prt\Gw)$ such that 
$0\leq\eta_{n}^{*}\leq 1$, $\eta_{n}^{*}\equiv 1$ in a relatively open 
neighborhood of $K$ and 
$$\norm{\eta_{n}^{*}}_{W^{2/q,q'}(\prt\Gw)}\to 0\quad\mbox {as 
}\,n\to\infty.$$
Since a boundary set with $C_{2/q,q'}$-capacity zero  has zero $(n-1)$
-Hausdorff measure, $\eta_{n}^{*}\to 0$ as $n\to\infty$. Thus 
$\gz_{\eta_{n}^{*}}\to \phi$. If we let $n\to\ity$ in (\ref{int- 
u^q}) we finally obtain
\begin{equation}\label{int- u^q+1}
   \int_{\Gw}u^q \phi dx \leq 
   C_{2},
 \end{equation}  
 with $C_{2}=C_{2}(K)$. Thus $K$ is conditionally 
 $q$-removable.\medskip
 
 \noindent (ii) $\Longrightarrow$ (i) Since $u^q\in 
 L^1(\Gw;\rho_{_{\prt\Gw}}dx)$, $u\geq 0$ and
 $$-\Gd u=-u^q,
 $$
the function $v=u+\mathbb G_{{-\Gd}}^\Gw(u^q)$ is positive and 
harmonic in $\Gw$, thus it admits a boundary trace $\gm\in \mathfrak 
M_{+}(\prt\Gw)$. Since the boundary trace of 
$\mathbb G_{{-\Gd}}^\Gw(u^q)$ is the zero measure, it is infered that $u$
admits the same boundary trace $\gm$,  the support of which
is included into the set $K$. Moreover
$$0\leq u=\mathbb P_{{-\Gd}}^\Gw(\gm)-
\mathbb G_{{-\Gd}}^\Gw(u^q)
\leq \mathbb P_{{-\Gd}}^\Gw(\gm).
$$
Therefore $u=u_{\gm}$, solution of Problem (\ref {Bq-power}) 
with $L=-\Gd$. Consequently $\gm$ does not charge boundary sets with 
$C_{2/q,q'}$-capacity zero and the same property is shared by 
$k\gm$, for any $k\in\mathbb N_{*}$. Put $u_{k}=u_{k\gm}$. If $\gm$ 
is not zero, the sequence of solutions $\{u_{k}\}$ is increasing and 
converges to some $u_{\ity}$ when $k\to\ity$. Because $u_{k}$ 
vanishes on $\prt\Gw\setminus K$, it follows from the Keller-Osserman 
construction that $u_{\ity}$ inherits the same property.  Furthermore
\begin {equation}\label {condition}
\int_{\Gw}\left(-u_{k}\Gd\gz_{\eta^*}+\gz_{\eta^*}u_{k}^q\right)dx
 =-k\int_{\prt\Gw}\frac{\prt\gz_{\eta^*}}{\prt{\bf n}}d\gm,
\end {equation}
 where $\eta\in\CT$, $\eta^*=1-\eta$ and 
 $\gz_{\eta^*}=\phi R^{2q'}_{\eta^*}$.
Because $\gm$ is not zero, the right-hand side of (\ref {condition}) 
 tends to infinity with $k$. Since $K$ is conditionally q-removable
 $u_{\ity}\in  L^1(\Gw;\rho_{_{\prt\Gw}}dx)$. Moreover, as we have 
 seen it before,
 $$\abs {\int_{\Gw}u_{k}\Gd\gz_{\eta^*} dx}
 \leq C\left(\int_{\Gw}u_{k}^q\phi dx\right)^{1/q}
\left(1+{\norm {\eta^*}}_{W^{2/q,q'}(\prt\Gw)}\right).
 $$
 Hence, the right-hand side of (\ref {condition}) is bounded 
  independently of $k$, which is a contradiction.\medskip
 
 \noindent (i) $\Longrightarrow$ (iii). If we assume 
 $C_{2/q,q'}(K)>0$, there exists a measure $\gm_{K}
\in\mathfrak M_{+}(\prt\Gw)\cap W^{-2/q,q}(\prt\Gw)$, satisfying 
$\gm_{K}(\prt\Gw\setminus K)=0$
 and $C_{2/q,q'}(K)=\gm_{K}(K)$. This measure is an extremal for the 
dual definition of the capacity of $K$ (already introduced in 
(\ref{Bessel4'} with Bessel potentials) :
$$C_{2/q,q'}(K)=\sup_{\tiny{\BA {l}\gm\in\mathfrak M_{+}(\prt\Gw)\\
\gm(\prt\Gw\setminus K)=0\EA}}\left(\frac{\gm (K)}{\norm{\mathbb 
P^\Gw_{-\Gd}(\gm)}_{L^q(\Gw;\gr_{_{\prt\Gw}}dx)}}\right)^{q'},
$$
see \cite [Th. 2.2.7]{AH}. 
Hence Problem (\ref {Bq-power}) with $L=-\Gd$ is solvable with 
$\gm=\gm_{K}$, thus $K$ is not conditionally q-removable.\qeda 

\subsection {The boundary trace problem}
One of the most striking aspects in the study on positive solutions of (\ref {Bq-power})
in a domain $\Gw$ relies on the possibility of defining a boundary 
trace which is no longer a Radon measure, but a generalized Borel 
measure, that is a measure which can take infinite values on compact 
boundary subsets. The second important task of the theory of boundary 
trace is to analyse the
connection between the set of all the boundary traces and the set of solutions. 
These notions were first studied by Le Gall \cite {LG1}, \cite {LG2} in 
the case $L=-\Gd$, 
$q=n=2$, and then extended by Marcus and V\'eron \cite {MV2}, \cite 
{MV3}, \cite {MV4}. For simplicity we shall consider first the model 
case 
\begin {equation}\label {tr1}
-\Gd u+\abs u^{q-1}u=0\quad\mbox {in }\,\Gw.
\end{equation}
We adopt the notations of Section 2.4
\bth {trth1}Let $\Gw\subset\mathbb R^n$ be a smooth domain and 
$q>1$. Let $u$ be a positive solution of (\ref {tr1}). Then for any 
$a\in\prt\Gw$ the following dichotomy holds :\\
(i) either for every relatively open subset $\CO\subset\Gw$ 
containing $a$,
  \begin {equation}\label {tr2}
\lim_{t\to 0}\int_{\CO_{t}}u(x)dS_{t}=\ity,
\end {equation}
(ii) or there exist a relatively open subset $\CO\subset\Gw$ 
containing $a$ and a positive linear functional $\ell$ on 
$C_{c}^\ity(\CO)$ such that for every $\theta\in C^\ity_{c}(\CO)$,
  \begin {equation}\label {tr3}
\lim_{t\to 0}\int_{\CO_{t}}u(x)\theta(x)dS=\ell (\theta).
\end {equation}
\es
\Proof The proof of this result is settled upon the following 
alternative which holds for every boundary point $a$ :
\smallskip 

\noindent (I) either there exists an open ball $B_{r_{0}}(a)$ such that
 \begin {equation}\label 
 {t5}\int_{B_{r_{0}}(a)\cap\Gw}u^q\rho_{_{\prt\Gw}} dx<\infty,
 \end {equation}
 \noindent (II) or for any $r>0$, 
  \begin {equation}\label {t6}
  \int_{B_{r}(a)\cap\Gw}u^q\rho_{_{\prt\Gw}}  dx=\infty.
 \end {equation}
 If (I) holds, let $\ge>0$ and $\CU_{\ge}$ be a smooth open 
 subdomain of 
 $\Gw\cap B_{r_{0}}(a)$ containing $\overline B_{r-\ge}(a)\Gw$
 and  such that 
 $$\overline B_{r-\ge}(a)\cap\prt\Gw\subset\overline \CU_{\ge}
 \cap\prt\Gw\subset B_{r}(a)\cap\prt\Gw.
 $$
 The function $\tilde u=u\vline_{\CU_{\ge}}$ is a nonnegative solution
 of (\ref {tr1}) in $\CU_{\ge}$ with 
 $\tilde u^q\in L^1(\CU_{\ge};\gr_{_{\prt\CU_{\ge}}}dx)$. Thus it 
 admits a boundary trace on $\prt\CU_{\ge}$ which belongs to 
 $\mathfrak M_{+}(\prt\CU_{\ge})$. Therefore,
 for any $\theta\in C^\ity_{c}(\prt\CU_{\ge})$, there holds
  \begin {equation}\label {tr3-ge}
\lim_{t\to 0}\int_{\prt\CU_{\ge\,t}}u(x)\theta(x)dS=\ell_{\ge} (\theta).
\end {equation}
Since $\ge$ is arbitrary and $\ell_{\ge}$ is uniquely determined on 
$\prt\CU_{\ge}$, assertion (ii) follows.\medskip

\noindent If (II) holds, let $\eta\in C^\ity_{c}(\prt\Gw\cap B_{r}(a))$ 
such that $0\leq \eta\leq 1$, $\eta\equiv 1$ on $\prt\Gw\cap 
B_{r/2}(a)$. For $t\in (0,\gb_{0}/2)$ small enough, we define 
$\gz_{\eta,t}$ in the set $\Gw_{t}\setminus \Gw_{\gb_{0}}$ by
$$\gz_{\eta,t}(x)=\gz_{\eta,t}(\gr_{_{\prt\Gw}}(x)-t,\gs(x))
=(\phi R_{\eta}^{2q'})(\gr_{_{\prt\Gw}}(x)-t,\gs(x)).
$$
Then
  \begin {equation}\label {tr4-ge}
\int_{\Gw_{t}\setminus 
\Gw_{\gb_{0}}}\left(-u\Gd\gz_{\eta,t}+u^q\gz_{\eta,t}\right)dx=
\int_{\Gs_{t}}\eta^{2q'}udS-\int_{\Gs_{\gb_{0}}}
\frac {\prt \gz_{\eta,t}}{\prt{\bf n}}(\gb_{0}-t,\gs)dS.
\end {equation}
As we have already seen it
$$\int_{\Gw_{t}\setminus 
\Gw_{\gb_{0}}}\abs {u\Gd\gz_{\eta,t}}dx\leq C{\norm\eta}_{W^{2/q,q'}}
\left(\int_{\Gw_{t}\setminus \Gw_{\gb_{0}}}u^q\gz_{\eta,t}dx\right)^{1/q}.
$$
Because the surface integral term in (\ref {tr4-ge}) on $\Gs_{\gb_{0}}$ 
is bounded independently of $t$, it follows
  \begin {equation}\label {tr5-ge}
\int_{\Gs_{t}}\eta^{2q'}udS\geq \int_{\Gw_{t}\setminus\Gw_{\gb_{0}}} 
u^q\gz_{\eta,t}dx-C_{1}{\norm\eta}_{W^{2/q,q'}}
\left(\int_{\Gw_{t}\setminus 
\Gw_{\gb_{0}}}u^q\gz_{\eta,t}dx\right)^{1/q}-C_{2}.
\end {equation}
Moreover, as $\eta\equiv 1$ on $\prt\Gw\cap 
B_{r/2}(a)$, there exists $\gd>0$ such that 
$\phi R_{\eta}^{2q'}\geq \gd $ on $\Gw\cap B_{r/2}(a)$. Hence, by (\ref 
{t6}) and the 
Beppo-Levi Theorem, 
$$\lim_{t\to 0}\int_{\Gw_{t}\setminus\Gw_{\gb_{0}}} 
u^q\gz_{\eta,t}dx=\ity,
$$
which implies
  \begin {equation}\label {tr6-ge}
\lim_{t\to 0}\int_{\Gs_{t}}\eta^{2q'}udS=\ity,
\end {equation}
and assertion (i) follows.\qeda\medskip

We write $\prt\Gw=\mathcal S(u)\cup\mathcal R(u)$ where $\mathcal 
S(u)$ is the closed subset of boundary points where (i) occurs, and 
$\mathcal R(u)=\prt\Gw\setminus\mathcal S(u)$. By using a partition 
of unity, there exists a unique positive Radon measure $\gm$ on $\mathcal R(u)$
such that 
 \begin {equation}\label {t3}
 \lim_{t\downarrow 0}\int_{\mathcal R(u)}u(\gs,t)\zeta_{t} (\gs,t)dS_{t}=\int_{\mathcal R(u)}\zeta 
 (\gs)d\gm,
 \end {equation}
 for every $\zeta\in C_{c}(\mathcal R(u))$.
 Thus we define the boundary trace by the following identification
  \begin {equation}\label {t4}
Tr_{\prt\Gw}(u)=(\mathcal S(u),\gm).
 \end {equation}
The set $\mathcal S(u)$ is called the {\it singular part} of the boundary trace of 
$u$, while $\gm\in\GTM_{+}(\mathcal R(u))$ is the {\it regular part}. 
The couple $(\mathcal S(u),\gm)$ defines in a unique way an outer 
regular positive Borel measure $\gn$ (an element of 
$\GTB_{+}^{reg}(\prt\Gw))$, with singular part $\mathcal S(u)$ 
and regular part $\gm$.\medskip

In the subcritical case, an important pointwise characterization of the singular part
 is the following minoration,
 \bprop {minprop} Let $\Gw$ be a bounded domain in $\mathbb R^n$ 
 with a $C^{2}$ boundary $\prt\Gw$, $1<q<(n+1)/(n-1)$ and $u$ be a positive solution 
of (\ref {tr1}) in $\Gw$ with boundary trace $(\CS(u),\gm)$. If $a\in 
\CS(u)$, then
  \begin {equation}\label {min1}
u(x)\geq u_{\ity a}(x),\forevery x\in \Gw,
 \end {equation}
 where $u_{\ity a}=\lim_{k\to\ity}u_{k\gd_{a}}$, and $u_{k\gd_{a}}$ 
 is the solution of
  \begin {equation}\label {min2}\left.\BA{ll}
-\Gd u_{k\gd_{a}}+\abs {u_{k\gd_{a}}}^{q-1}u_{k\gd_{a}}=0&\mbox 
{in } \Gw,\\[2mm]
\phantom {-\Gd u_{k\gd_{a}}+\abs {u_{k\gd_{a}}}^{q-1}} 
\,u_{k\gd_{a}}=k\gd_{a}&\mbox {on }\, \prt\Gw.
\EA\right.
\end {equation}
\es 
\Proof Since for any $r>0$, there holds
$$\lim_{t\to 0}\int_{B_{r}(a)\cap \Gs_{t}}u(x)dS_{t}=\infty,
$$
for any $k>0$ and $t=t_{k}=1/k$, there exists $r_{k,t}>0$ such that
$$\int_{B_{r_{k,t}}(a)\cap \Gs_{t}}u(x)dS_{t}\geq k.
$$
Let $m_{k}$ be such that
$$\int_{B_{r_{k,t}}(a)\cap \Gs_{t}}\min\{m_{k},u(x)\}dS_{t}= k, 
$$
and denote by $v_{k}$ the solution of 
  \begin {equation}\label {min3}\left.\BA{ll}
-\Gd v_{k}+\abs {v_{k}}^{q-1}v_{k}=0&\mbox 
{in } \Gw_{t},\\[2mm]
\phantom {-\Gd v_{k}+\abs {v_{k}}^{q-1}} 
\,v_{k}=\chi_{_{B_{r_{k,t}}(a)\cap \Gs_{t}}}&\mbox {on }\Gs_{t}.
\EA\right.
\end {equation}
 By the maximum principle, $v_{k}\leq u$ in $\Gw_{t}$ and by the 
 stability result of \rcor {Bmeascor}, $v_{k}$ converges to $u_{k\gd_{a}}$ locally 
 uniformly in $\Gw$ (actually the proof is given for a fixed 
 domain $\Gw$, but the adaptation to a sequence of expanding smooth 
 domains is straightforward). Thus $u_{k\gd_{a}}\leq u$ in $\Gw$. 
 Since $k$ is arbitrary, (\ref {min1}) follows.\qeda \medskip
 
 \noindent \Remark Notice that the boundary behaviour of $u_{\ity a}$ 
 is given by \rth {GmVth1} : with an apropriate rotation in the space, 
 it is
   \begin {eqnarray}\label {min4}
\lim_{\tiny {\BA {c}x\to a\\
(x-a)/\abs {x-a}\to\gs\EA}}\abs {x-a}^{2/(q-1)}u_{\ity a}(x)=\gw(\gs),\quad \mbox 
{uniformly on }S^{n-1}_{+},
 \end {eqnarray}
 where $\gw$ is the unique solution of (\ref {model2}) on 
 $S^{n-1}_{+}$  which vanishes on the equator $\prt 
 S^{n-1}_{+}$.\medskip
 
The most general boundary value problem concerning 
 positive solutions of (\ref {tr1}) is to solve the Dirichlet boundary 
 value problem with a given outer regular Borel measure as boundary 
 trace. If $\gn\in \mathfrak B_{+}^{reg}(\prt\Gw)$, we put
 $$\CS=\CS_{\gn}=\{\gs\in\prt\Gw:\,\gn (U)=\ity\mbox { for every 
 relatively open neighborhood $U$ of $\gs$}\}.
 $$
 Clearly $\CS_{\gn}$ is closed and the restriction $\gm$ of $\gn$ to 
 $\CR_{\gn}=\prt\Gw\setminus \CS_{\gn}$ is a Radon measure. This 
 establishes a one to one correspondence between 
 $\mathfrak B_{+}^{reg}(\prt\Gw)$ and the set of couples $(\CS,\gm)$, 
 where $\CS$ is a closed subset of $\prt\Gw$ and $\gm$ a positive 
 Radon measure on $\CR=\prt\Gw\setminus \CS$. 
 The following result is proven in \cite {MV4}.

\bth {trth2}Let $\Gw\subset\mathbb R^n$ be a smooth domain and 
$1<q<(n+1)/(n-1)$. Then for any 
$\gn\in\mathfrak B^{reg}_{+}(\prt\Gw)$ with $\gn\approx (\mathcal S,\gm)$, 
where $\mathcal S$ is a closed subset of $\prt\Gw$ and $\gm$ a 
positive Radon measure on $\prt\Gw\setminus \mathcal S$,
there exists a unique solution of 
  \begin {equation}\label {min5}\left.\BA{ll}
-\Gd u+\abs u^{q-1}u=0&\mbox 
{in } \Gw,\\[2mm]
\phantom {-\Gd v_{k}i} 
\,Tr_{\prt\Gw}(u)=\gn.&
\EA\right.
\end {equation}
\es
\Proof The proof is long and technical, and we shall just indicate 
the main steps : \smallskip

\noindent (1) By approximation, a minimal solution 
$\underline u_{\mathcal S,\gm}$ and a maximal solution 
$\overline u_{\mathcal S,\gm}$ of Problem (\ref {min5}) are constructed, 
so any other solution $u$ satisfies
 \begin {eqnarray}\label {min6}
\underline u_{\mathcal S,\gm}\leq u\leq \overline u_{\mathcal S,\gm}.
\end{eqnarray}

\smallskip

\noindent (2) Using convexity and the approximations of the minimal 
and the maximal solutions, it is proven that 
 \begin {eqnarray}\label {min7}
\overline u_{\mathcal S,\gm}-\underline u_{\mathcal S,\gm}
\leq \overline u_{\mathcal S,0}-\underline u_{\mathcal S,0}.
\end{eqnarray}
\smallskip

\noindent (3) Using (\ref {KOdist}), (\ref {min1}), (\ref {min4}) and Hopf 
boundary lemma, there 
exists $K=K(q,\Gw)>1$ such that 
 \begin {eqnarray}\label {min7'}
\overline u_{\mathcal S,0}\leq K\underline u_{\mathcal S,0}.
\end{eqnarray}
\smallskip

\noindent (4) Assuming that $\underline u_{\mathcal S,0}\neq\overline u_{\mathcal S,\gm}$
 (and the strict inequality follows by the strong maximum principle), 
a convexity argument implies that the function
 $$w=\underline u_{\mathcal S,0}-\frac {1}{2K}
 (\overline u_{\mathcal S,0}-\underline u_{\mathcal S,0}),
 $$
 is a supersolution of (\ref {min5}) with $\gn\approx (\CS,0)$. Since 
 for $0<\ga<1/(2K)$
 $\ga\underline u_{\mathcal S,0}$ is a subsolution of the same 
 problem with the same boundary trace, and 
 $$\ga\underline u_{\mathcal S,0}\leq w,
 $$
 it follows by (\rth {subsup}) that there exists a solution $u$ of 
 (\ref {tr1}) in $\Gw$ and 
  \begin {eqnarray}\label {min8}
 \ga\underline u_{\mathcal S,0}\leq u\leq w<\underline u_{\mathcal S,0}.
 \end{eqnarray}
 Because both $\ga\underline u_{\mathcal S,0}$ and $w$ have the same 
 boundary trace $(\mathcal S,0)$ in the sense of \rth {trth1}, $u$ is 
 a solution of Problem (\ref {min5}) with $\gn\approx (\mathcal S,0)$.
 This fact contradicts the 
 minimality of $\underline u_{\mathcal S,0}$, thus 
 $\overline u_{\mathcal S,0}=\underline u_{\mathcal S,0}$, which, in 
 turn,  implies
 $\overline u_{\mathcal S,\gm}=\underline u_{\mathcal S,\gm}$.
 \qeda \medskip
 
When $q\geq (n+1)/(n-1)$ neither any positive 
Radon measure on $\prt\Gw$ is eligible for 
being the regular part of the boundary trace of a positive solution 
of  (\ref {tr1}), nor any closed boundary subset for being the 
singular part : these facts follow from \rth {Bq-powerth} and \rth 
{Bremovth}. 
\bdef {singbdry}{\rm (i) Let $\CA$ be a relatively open subset of 
$\prt\Gw$ and $\gm\in\mathfrak M_{+}(\CA)$. Then {\it the singular 
boundary of $\CA$ relative to $\gm$} is defined by
  \begin {eqnarray}\label {min9}
\prt_{\gm}\CA=\{\gs\in\overline\CA:\gm(U\cap\CA)=\ity,\quad\mbox {for 
every neighborhood }U\mbox { of }\gs\}.
 \end{eqnarray}
(ii) Let $\CA$ be a Borel subset of $\prt\Gw$. A boundary point $\gs$ 
is 
{\it q-accumulation point of $\CA$} if, for every relatively open neighborhood
$U$ of $\gs$, $C_{2/q,q'}(\CA\cap U)>0$. The set of $q$-accumulation 
points of $\CA$ will be denoted by $\CA^*_{q}$.}\es 

The following result, announced (under a slighly different 
form) in \cite {MV3}, is proven in \cite {MV5} (see also \cite {DK3}, \cite 
{DK4}). 
\bth {trth3}Let $\Gw\subset\mathbb R^n$ be a smooth domain, 
$q\geq (n+1)/(n-1)$ and $\gn\approx (\mathcal S,\gm)$
 an element of $\mathfrak B^{reg}_{+}(\prt\Gw)$. Then Problem (\ref 
 {min5}) admits a solution if and only if  the following condition is 
 fulfilled :
 \begin {equation}\label {CNS}
 \left.\BA {l}
 (i) \mbox { For every Borel subset }\CA\subset \CR= \prt\Gw\setminus\mathcal 
 S,\,C_{2/q,q'}(\CA)=0\Longrightarrow \gm(\CA)=0,\\[2mm]
 (ii) \,\CS=\CS_{q}^*\cup \prt_{\gn}(\CR).\EA 
 \right.\end {equation}
\es

One of the most striking aspect of the super-critical case is the 
loss of uniqueness. It has been proven by Le Gall \cite {LG4} in the case $q=2$ 
and extended by Marcus and V\'eron \cite {MV5} that there exist 
infinitely many solutions of Problem (\ref {min5}) whenever the 
singular set $\CS$ has a non-empty relative interior. Actually there 
exists a maximal solution, but no minimal solution. This fact has led 
Dynkin and Kuznetsov in \cite {DK5} to introduce a thiner notion of 
boundary trace called the {\it fine trace}. However their definition 
is only working when $q\leq 2$. When $q=2$ and with a 
fundamental use of probability techniques (the Brownian snake), Mselati proved in \cite {Ms} 
the one to one correpondence between positive solutions of (\ref 
{tr1}) and the fine trace. The extension of this result in the general 
case remains open.

\subsection {General nonlinearities}
\subsubsection {The exponential}
There are many extensions of the nonlinear boundary value problems 
when the nonlinearity in no longer of a power type. In \cite {GrV} 
the boundary trace of the prescribed Gaussian curvature equation is studied 
\begin {equation}\label {gauss}
-\Gd u=K(x)e^{2u},
\end {equation}
in a $2$-dimensional bounded domain $\Gw$. In this equation, $K$ is 
a given function ; the question is to find out a new metric conformal to  the 
standard metric of a subdomain on the 
hyperbolic plane $\mathbb H^2$ so that $K$ is the Gaussian curvature of 
this metric (see \cite {RRV} for example). The existence of boundary trace in the set 
of outer regular Borel measures on $\prt\Gw$ is proven. In the case of 
a Radon measure the following existence 
result is obtained :
\bth {GrVth}Suppose $\gb\leq K(x)\leq \ga<0$ is a continuous function in a smooth bounded 
domain $\Gw$ of the plane and $\gm\in\mathfrak M(\prt\Gw)$ with 
Lebesgue decomposition
$$\gm=\gm_{R}dH_{1}+\gm_{s},
$$
where $\gm_{R}\in L^1(\prt\Gw)$ and $\gm_{s}\perp\gm_{R}$. If there 
exists some $p\in (1,\ity]$ such that 
 \begin {equation}\label {gauss2}
 \left.\BA {l}
 (i)\, \exp \left(2P_{-\Gd}^\Gw(\gm_{s})\right)\in 
 L^{p'}(\Gw;\rho_{_{\prt\Gw}}dx),\\[2mm]
 (ii) \,\exp (2\gm_{R})\in L^{p-1}(\prt\Gw)
 ,\EA 
 \right.\end {equation}
 then there exists a unique $u\in L^1(\Gw)$ with 
 $e^{2u}\in L^1(\Gw;\rho_{_{\prt\Gw}}dx)$ solution of
 \begin {equation}\label {gauss3}\left.\BA{ll}
-\Gd u-K(x)e^{2u}=0&\mbox {in }\,\Gw,\\[2mm]
\phantom {-\Gd u unn +gf} 
\,u=\gm.&
\EA\right.
\end {equation}
\es  

As for the power case, sufficient conditions for solving
\begin {equation}\label {gauss4}\left.\BA{ll}
-\Gd u-K(x)e^{2u}=0&\mbox {in }\,\Gw,\\[2mm]
\phantom {-\Gd unn} 
\,Tr_{\prt\Gw}(u)=\gn.&
\EA\right.
\end {equation}
where $\gn\in\mathbb B_{+}^{reg}(\prt\Gw)$ are given. They are 
expressed in terms of a boundary logarithmic capacity.
\subsubsection {The case of a general nonlinearity}
For general semilinear equations of the form
 \begin {equation}\label {gen1}
 -\Gd u+g(x,u)=0\quad\mbox {in }\,\Gw,
\end {equation}
where $\Gw$ is a smooth domain in $\mathbb R^n$, not necessarily 
bounded, and $g$ a continuous function defined on $\Gw\times \mathbb R$, 
a new approach of the boundary trace problem is provided by 
Marcus and V\'eron in \cite {MV7}. As it has already been 
observed in the implication [(i) $\Longrightarrow$ (ii)] 
in  the proof of \rth {Bremovth}, 
if $u$ is a positive solution of (\ref {gen1}) with $g(x,u)\geq 0$, and 
if for some $a\in\prt\Gw$ there exists $r>0$ such that
 \begin {equation}\label {gen1'}
\int_{B_{r}(a)\cap\Gw}g(x,u)\gr_{_{\prt\Gw}}dx<\ity,
\end {equation}
then $u\in L^1(B_{r'}(a)\cap\Gw)$ for any $0<r'<r$ and there exists
a positive linear functional $\ell$ on $C^{\ity}_{c}(\Gs\cap 
B_{r}(a))$ such that, for any $\gth$ in this space,   
 \begin {equation}\label {gen1''}
\lim_{t\to 0}\int_{B_{r}(a)\cap\Gs_{t}}u(x)\gth (x)dS_{t}=\ell (\gth).
\end {equation}

This result leads to the notion of regular and singular points if it 
is assumed for example that $g$ satisfies
 \begin {equation}\label {gen13}
g(x,r)\geq 0,\forevery (x,r)\in\Gw\times \mathbb R_{+}.
\end {equation}
\bdef {RP}{\rm Let $u$ be a continuous nonnegative solution of (\ref 
{gen1}). A point $a\in\prt\Gw$ is called a {\it regular point} of $u$ if 
there exists an open neighborhood $U$ of $a$ such that (\ref{gen1'}) 
holds. The set of regular points is denoted by $\CR (u)$. It is a
relatively open subset of $\prt\Gw$. Its complement, $\CS 
(u)=\prt\Gw\setminus \CR (u)$ is the {\it singular set} of $u$.
}\es

Using  a partition of unity, it exists a positive Radon measure $\gm$ on $\CR 
(u)$ such that 
   \begin {equation}\label {gen2}
 \lim_{t\downarrow 0}\int_{\CR (u)_{t}}u(\gs,t)\zeta_{t} (\gs,t)dS_{t}
 =\int_{\CR (u)}\zeta (\gs)d\gm,
 \end {equation}
 for every $\zeta\in C_{c}(\CR (u))$. 
 
 \bdef {CL}{\rm A function $g$ is a {\it coercive nonlinearity} in $\Gw$ 
 if, for every compact subset $K\subset\Gw$, the set of positive 
 solutions of (\ref {gen1}) is uniformly bounded on $K$}.
 \es
 
 An example of coercive nonlinearity is the following :
    \begin {equation}\label {cnl}
g(x,r)\geq h(x)g(r),\forevery (x,r)\in\Gw\times\mathbb R_{+},
 \end {equation}
 where $h\in C(\Gw)$ is continuous and positive, and $f\in C(\mathbb R_{+})$ 
 is nondecreasing and satisfies the Keller-Osserman assumption :
\begin {equation}\label {KO2}
 \myint{\theta}\infty\left(\myint{0}tf(s)ds\right)^{-1/2}dt<\infty, 
  \forevery\gth >0.
 \end {equation}
 The verification of this property is based upon the maximum 
 principle and the construction of local super solutions by the 
 Keller-Osserman method.
\bdefÊ{}{\rm A function $g$ possesses the {\it strong 
  barrier property} at $a\in\prt\Gw$ 
 if there exists $r_{0}>0$ such that, for any $0<r\leq r_{0}$, there 
 exists a positive super solution $v=v_{a,r}$ of (\ref {gen1}) in $B_{r}(a)\cap\Gw$ 
 such that $v\in C(B_{r}(a)\cap\overline{\Gw})$ and}
\begin {equation}\label {sbp}
\lim_{\tiny {\begin {array}{l} y\to 
x\\y\in\Gw\end{array}}}v(y)=\infty,\forevery x\in\Gw\times \prt 
B_{r}(a).
 \end {equation}
 \es
 
 If $g(x,r)= f(r)$ where $f$ satisfies the Keller-Osserman 
 assumption, then it possesses the strong barrier property at any 
 boundary point. If 
 $$ g(x,r)=\left(\gr_{_{\prt\Gw}}(x)\right)^\ga r^q,\forevery (x,r)\in 
 \Gw\times\mathbb R_{+}
 $$
 for some $\ga>-2$ and $q>1$, it possesses also the strong barrier 
 property, but the proof, due to Du and Guo \cite {DG}, is 
 difficult in the case $\ga>0$ (the nonlinearity is degenerate at the 
 boundary).

 \bprop{toinfty}Let $u\in C(\Gw)$ be a positive solution of
(\ref {gen1}) and suppose that $a\in \CS(u)$. Suppose that at least one
of the following sets of conditions holds :\\[1mm]
 (i)  There exists an open neighborhood  $U'$ of $a$ such that 
 $u\in L^1(U'\cap\Gw)$.\\[1mm]
(ii) (a) $g(x,\cdot)$ is non-decreasing in $\BBR_{+}$, for every $x\in\Gw$;\\
\phantom{(i). }(b) $\exists U_a$, an open neighborhood of $a$, such that $g$
is coercive in $U_a\cap\Gw$;\\
\phantom{(i). }(c) $g$ possesses the strong barrier property at $a$.\\[1mm]
 Then, for every open neighborhood $U$ of $a$,
\begin{equation}\label{loc-inf}
\lim_{t\to 0}\int_{U\cap\Gs_t}u(x)\,dS_{t}=\infty.
\end{equation}
\es
 
 This result, jointly with (\ref {gen1''}),  yields to the 
 following trace theorem.
 \bth {THGEN} Let $g$ be a coercive nonlinearity which has 
 the strong barrier property at any boundary point. Assume also 
 that $r\mapsto g(x,r)$ is nondecreasing on $\mathbb R_{+}$ for every 
 $x\in\Gw$. Then any continuous nonnegative solution $u$ of (\ref {gen1}) 
 possesses a boundary trace $\gn$ in $\GTB^{reg}_{+}(\prt\Gw)$ with 
\begin {equation}\label {gt3}
\gn=Tr_{\prt\Gw}(u)\approx (\CS(u),\gm),\mbox { where 
}\gm\in\GTM_{+}(\CR(u)).
 \end {equation}
 \es
 
 This result applies in the particular case where 
 $g(x,r)=\gr_{_{\prt\Gw}}(x)^\ga 
 r^q$. Moreover a complete extension of \rth {trth2} in 
 the subcritical range
 $$1<q<\frac {n+1+\ga}{n-1},\quad \ga>-2,
 $$
is valid. The super critical case is still completely open. 

 


\end {document}